\documentclass[11pt]{article}
\usepackage[a4paper,margin=25mm]{geometry}
\usepackage[T1]{fontenc}
\usepackage{lmodern}
\usepackage{amsmath,amssymb,amsthm,mathtools}
\usepackage{microtype,array,longtable,booktabs,mathrsfs}

\usepackage[unicode,hidelinks]{hyperref}
\usepackage{enumitem}
\setlist{nosep,leftmargin=*}
\allowdisplaybreaks[1]
\numberwithin{equation}{section}
\renewcommand{\[}{\begin{equation}}
\renewcommand{\]}{\end{equation}}
\newtheorem{lemma}{Lemma}[section]
\newtheorem{theorem}[lemma]{Theorem}
\newtheorem{proposition}[lemma]{Proposition}

\theoremstyle{definition}

\newtheorem{example}[lemma]{Example}
\theoremstyle{remark}
\newtheorem{remark}[lemma]{Remark}
\newcommand{\Hh}{\mathbb H}
\newcommand{\R}{\mathbb R}

\newcommand{\Z}{\mathbb Z}
\newcommand{\C}{\mathbb C}
\newcommand{\dd}{\,\mathrm d}
\newcommand{\supp}{\operatorname{supp}}
\newcommand{\sech}{\operatorname{sech}}
\newcommand{\dist}{\operatorname{dist}}
\newcommand{\Ree}{\operatorname{Re}}
\newcommand{\one}{\mathbf 1}
\newcommand{\norm}[1]{\lVert #1\rVert}
\newcommand{\ip}[2]{\langle #1,#2\rangle}
\newcommand{\G}{\mathcal G}

\newcommand{\J}{\mathcal J}
\newcommand{\K}{\mathcal K}

\newcommand{\cp}{c_+}
\newcommand{\lam}{\lambda_+}
\newcommand{\E}{\mathbb E}
\newcommand{\Prb}{\mathbb P}
\newcommand{\phiZ}{\varphi_0}
\hypersetup{pdftitle={Periodic KPP equations in hyperbolic space},pdfauthor={Yuanyang Hu},pdfsubject={Extinction, persistence and spreading for periodic Fisher--KPP reactions in hyperbolic space}}
\title{Periodic KPP equations in hyperbolic space}
\author{Yuanyang Hu\\[0.4em]
\small School of Mathematics and Statistics, Henan University\\
\small 475004 Kaifeng, China\\
\small \href{mailto:yuanyhu@henu.edu.cn}{\texttt{yuanyhu@henu.edu.cn}}}
\date{}
\begin{document}
\maketitle
\begin{abstract}
We study Fisher-KPP equations in hyperbolic space with reactions periodic in a
Busemann coordinate, in the distance from a point, or in the distance from a
complete geodesic. For the Busemann and radial reactions, solutions with nonzero
compactly supported initial data become uniformly extinct precisely when the
principal eigenvalue of the periodic linearisation with the curvature correction
is nonpositive. When this eigenvalue is positive, we determine the logarithmic
delay of spreading up to a bounded error. For compactly supported data we bound the Busemann superlevel sets between
horoballs in fixed transverse windows and between balls in sectors away from
the endpoint of the defining ray. For data between a positive multiple of one horoball indicator
and the indicator of a containing horoball, the horoball inclusions hold
throughout hyperbolic space. For reactions
periodic in the distance from a geodesic, extinction is governed by the top of the
transverse spectrum. In the hyperbolic plane we construct a smooth periodic
reaction for which compactly supported populations die out in the Busemann and
radial cases but persist near the geodesic in the geodesic-distance case. For the
geodesic-distance reaction in this example, the smallest positive stationary
solution tends to zero with distance from the geodesic, and a variational formula
gives the spreading speed along it. For every fixed level below the value of this solution on the geodesic,
we determine the logarithmic delay up to a bounded error in both directions.
\end{abstract}


\section{Introduction}\label{sec:intro}

We study reactions periodic in a Busemann coordinate, distance from
a point, or distance from a complete geodesic in hyperbolic space.
How does the choice of coordinate affect the long-time behaviour
of solutions from compactly supported initial data and the positions
of their superlevel sets?
We determine the spectral thresholds and propagation speeds for these
three coordinates. In the hyperbolic plane, we construct a smooth
periodic reaction with a positive stationary solution that tends to zero
as the distance from a geodesic tends to infinity. For the Busemann and radial
reactions, positivity of the periodic principal eigenvalue governing extinction
yields a logarithmic delay for nonzero compactly supported data, even when the
linear growth rate varies periodically in space.

Let $(\Hh^d,g)$ be the complete simply connected Riemannian manifold
of dimension $d\ge2$ and constant sectional curvature $-1$.
We denote its volume measure by $\dd V$, set $\nabla=\nabla_g$, and write
$\Delta=\Delta_{\Hh^d}=\operatorname{div}_g\nabla_g$ for its
Laplace--Beltrami operator, where $\nabla_g$ and
$\operatorname{div}_g$ are the Riemannian gradient and divergence;
see~\cite[Section~2]{GGQ}.
Throughout, $D>0$ is a fixed diffusion coefficient and
$u_0\in L^\infty(\Hh^d)$ satisfies $0\le u_0\le1$ almost everywhere.
For each reaction $f:\Hh^d\times[0,1]\to\R$ specified below,
we consider the Cauchy problem
\[\label{eq:KPP}
 \begin{cases}
u_t=D\Delta_{\Hh^d}u+f(x,u),&\text{for }x\in\Hh^d,\ t>0,\\
 u(x,0)=u_0(x),&\text{for almost every }x\in\Hh^d.
 \end{cases}
\]
Equalities and inequalities between initial data are understood
almost everywhere with respect to $\dd V$.

For the reactions considered here, \eqref{eq:KPP} has a unique
bounded mild solution taking values in $[0,1]$, as established in
Lemma~\ref{lem:wellposed} and the discussion following its proof.
Here a mild solution means a function satisfying the Duhamel integral
identity \eqref{eq:duhamel} with initial datum $u_0$.
For positive times, $u$ is a classical solution of~\eqref{eq:KPP}:
$u\in C^{2,1}(\Hh^d\times(0,\infty))$ and satisfies the equation pointwise.
In local coordinates, its spatial derivatives through order two and its
first time derivative are continuous. It attains its initial datum in the sense that
$u(\cdot,t)\to u_0$ in $L^1_{\mathrm{loc}}(\Hh^d)$ as
$t\downarrow0$.
We write $u(x,t;u_0,f)$ when the initial datum or the reaction
needs to be specified. The level-position results describe the sets
$\{x\in\Hh^d:u(x,t)\ge\theta\}$ for $0<\theta<1$ and $t>0$.

For $1\le p<\infty$, we use
$\|v\|_{L^p(\Hh^d)}=(\int_{\Hh^d}|v|^p\dd V)^{1/p}$.
For functions on $\Hh^d$, we write $\|v\|_\infty=\|v\|_{L^\infty(\Hh^d)}$ and
$\|v\|_2=\|v\|_{L^2(\Hh^d)}$. The $L^\infty$ norm is the
essential supremum and agrees with the supremum for continuous functions. Compact support means vanishing
almost everywhere outside a compact set. For a measurable set $E$,
$\one_E$ denotes its indicator, equal to one on $E$ and zero on its
complement. A stationary solution with values in $[0,1]$ is a function
$q\in C^2(\Hh^d)$ satisfying $0\le q\le1$ and
$D\Delta q+f(x,q)=0$ on $\Hh^d$. A smallest positive stationary
solution $p_*$ is a stationary solution with $0<p_*\le1$ and
$p_*(x)\le q(x)$ for every stationary solution $q$ with $0<q\le1$
and every $x\in\Hh^d$.

One class of reactions has the form $f(x,v)=F(z(x),v)$, where $z:\Hh^d\to\R$ is a
Busemann function, distance from a fixed point, or distance from a fixed
complete geodesic. For a fixed $L>0$, the function
$F:\R\times[0,1]\to\R$ satisfies $F(\zeta+L,v)=F(\zeta,v)$ for
$(\zeta,v)\in\R\times[0,1]$. The choices $z=b,r,s$ give the reactions in~\eqref{eq:patterns}.
We assume \eqref{eq:F1} and, for the two distance reactions, the axis
condition $F_z(0,v)=0$ for $0\le v\le1$ in~\eqref{eq:axis-regularity}. The geometric conventions follow
\cite[Chapters~3--5]{Ratcliffe}, \cite[Chapters~5--6]{Lee18}
and \cite{Helgason}. Lemma~\ref{lem:coordinates} gives the coordinate expressions for
the diffusion operators and the volume elements entering the Rayleigh
quotients \eqref{eq:sigma} and \eqref{eq:rayleigh-full}.

Fisher \cite{Fisher} and Kolmogorov, Petrovskii and Piskunov
\cite{KPP} introduced the classical equation in 1937. Aronson and
Weinberger \cite{AW75,AW78} established spreading for scalar
reaction--diffusion equations. Periodic habitats and monotone
evolution systems were subsequently treated in~\cite{Weinberger82,Weinberger02,LiangYiZhao06,LiangZhao,LiangZhao10,
FangZhao14,FangYuZhao17}. For periodic reaction--diffusion equations,
Berestycki and Hamel \cite{BH02} proved pulsating-front existence,
and Berestycki, Hamel and Nadirashvili \cite{BHN05} characterised
minimal KPP speeds by periodic principal eigenvalues. Hamel \cite{Hamel08} studied exponential decay of pulsating fronts,
Hamel and Roques \cite{HamelRoques11} treated uniqueness and stability,
and Giletti \cite{Giletti14} proved convergence to minimal-speed fronts
in periodic KPP problems. For compactly supported data whose solutions
converge locally to $1$, Rossi \cite{RossiFG17} related the spreading
speed along a ray to the speeds of planar fronts by the
Freidlin--G\"artner formula.
The principal-eigenvalue and maximum-principle results of~\cite{BNV94,BR15} provide the elliptic background.

For the homogeneous hyperbolic equation with unit diffusion, Matano, Punzo and Tesei
\cite[Theorems~3.2 and~3.6]{MPT15} identified the threshold
$(d-1)^2/4$. For continuous initial data and KPP reactions satisfying
$f(v)\le f'(0)v$, compact support and $f'(0)<(d-1)^2/4$ imply
$\|u(\cdot,t)\|_{L^\infty(\Hh^d)}\to0$ as $t\to\infty$. If
$f'(0)>(d-1)^2/4$, the solution from every nonzero continuous
datum converges locally uniformly to $1$, and compactly supported data spread with speed
$2\sqrt{f'(0)}-(d-1)$. Their Corollary~3.8 bounds the width of a
fixed level set in radial distance, and they also construct
horospheric waves. Punzo \cite{Punzo16} studied propagation and
extinction on Riemannian manifolds. For heat diffusion on $\Hh^d$, kernel estimates and formulas are
given in~\cite{DaviesMandouvalos,GrigoryanNoguchi}, and large-time
behaviour is studied in~\cite{Vazquez22}.

Gonz\'alez, Gonz\'alvez and Quir\'os \cite{GGQ} treat critical
extinction and precise wave behaviour for the homogeneous hyperbolic
KPP equation with unit diffusion. Their Theorem~1.3 considers initial
data depending only on distance from a point, a Busemann function,
or signed distance from a totally geodesic hyperplane.
Writing $\rho$ for the chosen coordinate, the initial datum has the
form $u_0(x)=\widetilde u_0(\rho(x))$, where $\widetilde u_0\not\equiv0$
and $\widetilde u_0(\rho)=0$ for every $\rho\in I$ with $\rho\ge X$, for some
$X\in I$. Here $I=[0,\infty)$ for distance from a point and $I=\R$
for either of the other two coordinates \cite[(1.23)--(1.24)]{GGQ}. For wave
convergence they use the stronger bound $f'(v)\le f'(0)$ for $v\in[0,1]$ and
$2\sqrt{f'(0)}-(d-1)>0$. Under these hypotheses,
\cite[Theorem~1.3]{GGQ} gives convergence to a fixed translate of a
minimal Euclidean wave in the moving coordinate, uniformly on
$\{x\in\Hh^d:\rho(x)\ge0\}$.
The moving coordinate in~\cite[Theorem~1.3]{GGQ} uses speed
$2\sqrt{f'(0)}-(d-1)$ and logarithmic coefficient
$3/(2\sqrt{f'(0)})$.

Berestycki and Chapuisat \cite{BC13} studied a Euclidean model in which
mortality increases without bound away from a fixed axis. They characterised when a positive stationary profile exists and studied
travelling fronts and spreading along the axis in that regime.

In the hyperbolic plane we ask whether a reaction that varies periodically
with distance from a geodesic can sustain a population near that geodesic
and support spreading along it. Bouin, Henderson and Ryzhik~\cite{BHR17}
established a logarithmic delay in a model where the spatial diffusion rate depends on a bounded trait variable.

Logarithmic corrections have also been studied for Euclidean KPP equations.
Uchiyama \cite{Uchiyama78} studied large-time wave behaviour, and
Bramson \cite{Bramson78,Bramson83} determined the logarithmic delay
and selected phase. For $v_t=v_{xx}+v(1-v)$ on $\R$
with initial datum $\one_{(-\infty,0]}$, each fixed level in $(0,1)$
lies at $2t-\tfrac32\log t+O(1)$. Hamel, Nolen, Roquejoffre and
Ryzhik \cite{HNRR13} obtained the position $2t-\tfrac32\log t+O(1)$ by PDE methods.
Convergence to a fixed critical-wave translate and finer corrections
were developed in~\cite{NRRsingle,NRRrefined19,Graham19}; weighted
critical-front stability was studied in~\cite{FayeHolzer19}.
The effects of initial tails and front tails are treated in~\cite{AlfaroGilettiXiao25,Giletti22}. Multidimensional front locations,
direction-dependent phases and level-set geometry are studied in~\cite{Gartner82,Ducrot15,RRR19,Rossi17,RRM18}.

For one-dimensional periodic media, Hamel, Nolen, Roquejoffre and
Ryzhik \cite[Theorems~1.1--1.2]{HNRR16} prove logarithmic delay
from nonzero compactly supported data and convergence on the right
half-line to a minimal pulsating front with an additional time shift $\xi(t)=O(1)$
after the logarithmic delay. In their equation $f(x,\cdot)$ is concave,
$g(x)=f_u(x,0)$ may change sign, and the principal periodic eigenvalue of
$-\partial_{xx}-g$ is negative. For the reaction on $v\ge0$ in~\cite[Section~1.1]{HNRR16}, there are $m,V>0$ such that $f(x,v)\le g(x)v-mv^2$ for every
$x\in\R$ and $v\ge V$. Under these assumptions the stationary equation
$p''+f(x,p)=0$ has a positive $1$-periodic solution
\cite[Section~1.1]{HNRR16}. Shabani
\cite[Theorem~1.3]{Shabani19} obtains multidimensional logarithmic
positions for a positive periodic spatial factor multiplying a
scalar KPP nonlinearity. In Euclidean dimension $d$, the delay along a unit vector $e$ is
$\frac{d/2+1}{\lambda_*(e')\,e\cdot e'}\log t+O(1)$.
Here $e'$ minimises $c_*(\nu)/(e\cdot\nu)$ over unit vectors $\nu$ with
$e\cdot\nu>0$, $c_*(\nu)$ is the minimal front speed normal to $\nu$,
and $\lambda_*(e')$ is the critical exponential decay rate in direction $e'$
\cite[Theorem~1.3]{Shabani19}. Periodic branching processes
give related descriptions of maxima and pulsating waves
\cite{LTZ22,RenSongYangExist,RenSongYangUnique}.

Related wave and spreading results for homogeneous, periodic and
almost-periodic lattice equations appear in~\cite{ChenGuo03,ChenFuGuo06,GuoWu08,MaZhao08,GuoHamel06,GuoWu10,
GuoWu09,ChengLiLin15,WuChenHsu21,LiangZhou20,LiangWangZhouZhou24}.
Besse, Faye, Roquejoffre and Zhang \cite{BFRZ} prove a logarithmic
delay for the Fisher--KPP equation on $\Z$. For propagation on trees, see~\cite{HoffmanHolzer,AlvarezTree26}. Nonlocal diffusion,
nonlocal competition and transport with a nonlocal boundary reaction
lead to further logarithmic corrections
\cite{Graham22,BouinHendersonRyzhik20,RoquejoffreZhang26,
FayeRoquejoffreZhang26}.

We apply the identity in Lemma~\ref{lem:ground-transform} to positive
solutions and supersolutions of the linearised equations, and construct
compactly supported tests far from the fixed point or geodesic.
Positive Dirichlet eigenfunctions give the smallest positive stationary
solution in Lemma~\ref{lem:minimal-state}; the integral identity in
Lemma~\ref{lem:localized-uniqueness} proves the uniqueness assertion of
Theorem~\ref{thm:localized}(ii). Periodic exponential modes determine the front
speeds. For logarithmic positions, scalar reactions with the same
derivative at zero bound the spatially dependent reaction; shifted
one-dimensional comparisons absorb the exponentially decaying
geometric drift error. Proposition~\ref{prop:longitudinal-martingale} identifies the
asymptotic longitudinal variance of a diffusion whose transverse
coordinate is stationary, using the
martingale method for additive functionals in~\cite{KV86}.

\subsection{The equation and spectral thresholds}

For $x,y\in\Hh^d$, let $\dist(x,y)$ denote their hyperbolic distance.
For $x\in\Hh^d$ and a nonempty set $A\subset\Hh^d$, put
$\dist(x,A)=\inf_{y\in A}\dist(x,y)$.

Fix $o\in\Hh^d$, a unit-speed geodesic ray $\gamma_+:[0,\infty)\to\Hh^d$
with $\gamma_+(0)=o$, and a complete unit-speed geodesic
$\gamma:\R\to\Hh^d$. Put $\G=\gamma(\R)$ and
\[\label{eq:coordinates}
 \begin{aligned}
 b(x)&=\lim_{R\to\infty}\bigl(\dist(x,\gamma_+(R))-R\bigr),\\
 r(x)&=\dist(o,x),\qquad s(x)=\dist(x,\G) \quad {\text{for }x\in\Hh^d}.
 \end{aligned}
\]
The nearest point of $x$ on $\G$ is unique; the distance formula in
Lemma~\ref{lem:coordinates} proves this assertion. Define $\ell(x)$ by
$\dist(x,\gamma(\ell(x)))=\min_{t\in\R}\dist(x,\gamma(t))$. The parametrisation $\gamma$ fixes the sign
of $\ell$. For $R\in\R$, we use
\[\label{eq:sets}
 B(o,R)=\{x\in\Hh^d:r(x)\le R\},\quad
 \mathcal H_R=\{x\in\Hh^d:b(x)\le R\},\quad
 \mathcal T_R=\{x\in\Hh^d:s(x)\le R\}.
\]
For $R<0$, $B(o,R)=\mathcal T_R=\varnothing$. The balls in~\eqref{eq:sets} are closed; when a ball is used as a domain for a
differential equation or a Dirichlet eigenvalue problem, its interior
is understood. In statements on $\{x\in\Hh^d:|\ell(x)|\le M\}$, the parameter
$M\ge0$ is the half-width in the $\ell$ coordinate. A tube has infinite volume when its radius is positive.

Let $L>0$. Write $F_u(z,v)=\partial_vF(z,v)$ and
$F_z(z,v)=\partial_zF(z,v)$; mixed derivatives use the same convention.
The three periodic reactions in~\eqref{eq:patterns} share the following
hypotheses on $F$. For the general reaction $f(x,v)$ in
Theorem~\ref{thm:common}, the hypotheses are stated separately in that theorem.
\[\label{eq:F1}
 \begin{gathered}
 F\in C^2(\R\times[0,1]),\qquad F(z+L,v)=F(z,v)\quad\text{for }z\in\R,\ v\in[0,1],\\
 F(z,0)=F(z,1)=0,\qquad a(z):=F_u(z,0)>0,\qquad F_u(z,1)<0\quad\text{for }z\in\R,\\
 0<F(z,v)\le a(z)v\quad\text{for }z\in\R,\ 0<v<1.
 \end{gathered}
\]
Fix a constant $C_2\ge0$ such that
\[\label{eq:uniform-taylor}
 a(z)v-C_2v^2\le F(z,v)\le a(z)v
 \quad\text{for }z\in\R,\ 0\le v\le1.
\]
For example, $C_2=\tfrac12\sup_{[0,L]\times[0,1]}|F_{uu}|$ suffices.
We write $R_F(z,v)=a(z)v-F(z,v)$ for the quadratic remainder; thus
$0\le R_F(z,v)\le C_2v^2$ for $z\in\R$ and $0\le v\le1$.
Define
\[\label{eq:patterns}
 \begin{gathered}
 f_b(x,v)=F(b(x),v),\qquad f_o(x,v)=F(r(x),v),\\
 f_{\G}(x,v)=F(s(x),v)
 \quad {\text{for }(x,v)\in\Hh^d\times[0,1]}.
 \end{gathered}
\]
For $f=f_o$ and $f=f_\G$, we also assume
\[\label{eq:axis-regularity}
 F_z(0,v)=0\quad\text{for }0\le v\le1.
\]
Condition \eqref{eq:axis-regularity} makes $x\mapsto F(r(x),v)$ and $x\mapsto F(s(x),v)$
$C^2$ across $o$ and $\G$, respectively, as verified in~\eqref{eq:distance-composition-hessian}. Periodicity and \eqref{eq:F1} give uniform bounds for $F$
and its derivatives of order at most two on $\R\times[0,1]$. Set
\[\label{eq:growth-extrema}
 a_{\min}=\min_{z\in[0,L]}a(z)>0,\qquad
 a_{\max}=\max_{z\in[0,L]}a(z)<\infty.
\]
Write $\xi=\gamma_+(\infty)$ for the ideal endpoint of the ray $\gamma_+$.

For a closed linear operator $\mathscr L$ on a specified
Hilbert space, let $\sigma(\mathscr L)$ denote its spectrum.
If $\mathscr L$ is self-adjoint and bounded above, we call
$\sup\sigma(\mathscr L)$ its upper spectral endpoint.

Set $\beta=d-1$. For $\lambda\in\R$, let $k(\lambda)$ be the principal
periodic eigenvalue of
\[\label{eq:cell}
 \K_\lambda\phi
 =D\phi''+D(\beta-2\lambda)\phi'
       +\bigl(D\lambda^2-D\beta\lambda+a(z)\bigr)\phi,
 \qquad \phi(z+L)=\phi(z)\quad\text{for }z\in\R.
\]
The operator acts on $L^2(0,L;\dd z)$ with domain
$H^2_{\mathrm{per}}(0,L)=\{\phi\in H^2(0,L):
\phi(0)=\phi(L),\ \phi'(0)=\phi'(L)\}$.
Its principal eigenvalue is the real eigenvalue admitting a strictly
positive $L$-periodic eigenfunction. We denote this eigenfunction by
$\phi_\lambda$, so that
$\K_\lambda\phi_\lambda=k(\lambda)\phi_\lambda$, and normalise it by
$\int_0^L\phi_\lambda^2\dd z=1$. In particular, define
\[\label{eq:parameters}
 \lambda_0=\frac\beta2,\quad
 \mu=k(\lambda_0),\quad \phiZ=\phi_{\lambda_0},\quad
 \cp=\min_{\lambda>0}\frac{k(\lambda)}\lambda,\quad
 \lam=\operatorname*{argmin}_{\lambda>0}\frac{k(\lambda)}\lambda.
\]
Lemma~\ref{lem:cell} proves existence and uniqueness of the
minimiser and gives the identities
\[\label{eq:symmetry-intro}
 \begin{gathered}
 k(\lambda)=k(\beta-\lambda)\quad\text{for }\lambda\in\R,\\
 \mu=\min_{\lambda\in\R}k(\lambda),\qquad
 \operatorname{sgn}\cp=\operatorname{sgn}\mu.
 \end{gathered}
\]

Recall from \eqref{eq:coordinates} that $s(x)=\dist(x,\G)$ for
$x\in\Hh^d$. Define
\[\label{eq:tube-weight}
 m_d(s)=(\sinh s)^{d-2}\cosh s,\qquad
 B_d(s)=\frac{m_d'(s)}{m_d(s)}=(d-2)\coth s+\tanh s
 \quad {\text{for }s>0}.
\]
For $v\in C^2((0,\infty))$, the function $v\circ s$ is defined on
$\Hh^d\setminus\G$. By \eqref{eq:fermi-metric},
$\Delta_{\Hh^d}(v\circ s)(x)
=v''(s(x))+B_d(s(x))v'(s(x))$ for $x\in\Hh^d\setminus\G$.
The factor $m_d(s)$ in~\eqref{eq:fermi-metric} gives the measure
$m_d(s)\dd s$ in the Rayleigh quotient below.
For $d=2$, $m_2(s)=\cosh s$ and $B_2(s)=\tanh s$ extend smoothly
to $s=0$, with $m_2(0)=1$.

Let $\mathcal D_d$ be the set of functions $v:[0,\infty)\to\R$
such that $y\mapsto v(|y|)$ belongs to
$C_c^\infty(\R^{d-1};\R)$. We define
\[\label{eq:sigma}
 \sigma_d=\sup_{0\ne v\in\mathcal D_d}
 \frac{\displaystyle\int_0^\infty
       \bigl(a(s)v(s)^2-Dv'(s)^2\bigr)m_d(s)\dd s}
      {\displaystyle\int_0^\infty v(s)^2m_d(s)\dd s}.
\]
The form domain is the completion of $\mathcal D_d+i\mathcal D_d$
under the norm
$\bigl(\int_0^\infty(|v(s)|^2+|v'(s)|^2)m_d(s)\dd s\bigr)^{1/2}$.
On this domain, the lower-semibounded closed quadratic form
$v\mapsto\int_0^\infty(D|v'|^2-a(s)|v|^2)m_d(s)\dd s$
defines a self-adjoint operator $-\mathscr L_d$ on
$L^2((0,\infty),m_d(s)\dd s)$;
see~\cite[Chapter~VI, Theorems~2.1 and~2.6]{Kato}.
For $s>0$, the differential expression for $\mathscr L_d$ is
$\mathscr L_dv=Dm_d^{-1}(m_dv')'+a(s)v$.
By \eqref{eq:sigma}, the upper spectral endpoint of $\mathscr L_d$
is $\sigma_d$.

For $d=2$, the form domain is
$\{v\in H^1(0,\infty):
\int_0^\infty(|v|^2+|v'|^2)\cosh s\dd s<\infty\}$.
Its elements may have any trace value $v(0)$.
The operator domain of $\mathscr L_2$ imposes the Neumann condition
$v'(0)=0$; in particular, every eigenfunction of $\mathscr L_2$
satisfies this condition. For $d\ge3$, the form closure specifies
$\mathscr L_d$ at the singular endpoint $s=0$, and every function
$v$ in its operator domain satisfies
$\lim_{s\downarrow0}m_d(s)v'(s)=0$.

For complex-valued functions on a measure space $(X,\nu)$,
we use the Hilbert-space inner product
$\langle v,w\rangle=\int_X\overline{v(x)}w(x)\dd\nu(x)$.
The symbols $\Ree z$ and $\operatorname{Im}z$ denote the real and
imaginary parts of a complex number. Each use of a closed form or a
spectral endpoint below specifies its measure and form domain.
The notation $O(1)$ as $t\to\infty$ denotes a quantity bounded
uniformly for all sufficiently large real $t$. Constants may depend on the fixed
diffusion coefficient, reaction, period, dimension, initial datum, and level. Uniformity in any
additional parameter is stated where it is used.

\begin{theorem}\label{thm:threshold}
Assume \eqref{eq:F1}, and choose $f=f_b$ or $f=f_o$.
In the case $f=f_o$, assume also \eqref{eq:axis-regularity}.
Let $u$ solve \eqref{eq:KPP}.
\begin{enumerate}[label=(\roman*),beginpenalty=10000]
\item \emph{Extinction.} If $\mu\le0$ and $u_0\in L^2(\Hh^d)$, then
\[\label{eq:threshold-extinction}
 \lim_{t\to\infty}\norm{u(\cdot,t;u_0,f)}_{L^\infty(\Hh^d)}=0.
\]
If $\mu<0$ and $u_0\in L^2(\Hh^d)$, then for every $\tau>0$,
\[\label{eq:threshold-rate}
 \norm{u(\cdot,t;u_0,f)}_{L^{\infty}(\Hh^d)}
 \le C_\tau e^{\mu(t-\tau)}\norm{u_0}_{L^{2}(\Hh^d)}
 \quad\text{for }t\ge\tau.
\]
\item \emph{Persistence.} If $\mu>0$, then every nonzero datum in $[0,1]$
satisfies
\[\label{eq:threshold-persistence}
 \lim_{t\to\infty}\sup_{x\in K}|u(x,t;u_0,f)-1|=0
 \quad\text{for every compact }K\subset\Hh^d.
\]
\end{enumerate}
\end{theorem}

Let $\mathscr L_{\G}=D\Delta_{\Hh^d}+a(s(x))$ on
$L^2(\Hh^d)$, with the self-adjoint realisation specified in~\eqref{eq:rayleigh-full}. Then
$\sigma_d=\sup\sigma(\mathscr L_{\G})\ge\mu$,
where $\sigma_d$ is defined by \eqref{eq:sigma}; see
Lemma~\ref{lem:tube-spectrum}.
\begin{theorem}\label{thm:dimension}
Assume \eqref{eq:F1} and \eqref{eq:axis-regularity}, and take $f=f_{\G}$.
\begin{enumerate}[label=(\roman*),beginpenalty=10000]
\item \emph{Extinction.} If $\sigma_d\le0$, every
$u_0\in L^2(\Hh^d)$ with $0\le u_0\le1$ satisfies
\eqref{eq:threshold-extinction}. If $\sigma_d<0$ and
$u_0\in L^2(\Hh^d)$ with $0\le u_0\le1$, then
\eqref{eq:threshold-rate} holds with $\mu$ replaced by $\sigma_d$,
for every $\tau>0$.
\item \emph{Persistence.} If $\sigma_d>0$, there is a smallest
positive stationary solution $p_*\le1$, and every nonzero datum satisfies
\[\label{eq:minimal-persistence}
 \liminf_{t\to\infty}u(x,t;u_0,f_{\G})\ge p_*(x)>0
 \quad\text{for }x\in\Hh^d.
\]
The function $p_*$ depends only on $s(x)$. If, in addition, $\mu>0$,
then $p_*\equiv1$ and \eqref{eq:threshold-persistence} holds.
\item  If $d\ge3$ and $\mu\le0$, then
\[\label{eq:no-localization}
 \sigma_d=\mu\le0.
\]
\end{enumerate}
\end{theorem}

The dimension restriction follows from the drift inequality
$B_d(s)>\beta$ for $s>0$ when $d\ge3$; see the proof of
Theorem~\ref{thm:dimension}.

\subsection{Stationary solutions near a geodesic}

For the rest of this subsection $d=2$, so $\beta=1$. Write
$\sigma=\sigma_2$ and
\[\label{eq:h-weight}
 h(s)=(\cosh s)^{-1/2}\phiZ(s),\quad\text{for }s\ge0.
\]

For uniqueness and convergence, we shall also assume

\[\label{eq:quotient}
 \frac{F(z,v)}v>\frac{F(z,w)}w
 \quad\text{for }z\in\R,\ 0<v<w\le1.
\]

\begin{theorem}\label{thm:localized}
Assume \eqref{eq:F1} and \eqref{eq:axis-regularity}, with $d=2$,
$f=f_{\G}$ and $\mu\le0<\sigma$.
\begin{enumerate}[label=(\roman*),beginpenalty=10000]
\item \emph{Profile and confinement.} Let $p_*$ be the smallest
positive stationary solution in Theorem~\ref{thm:dimension}(ii).
There exist $p\in C^2([0,\infty))$ and $C>0$ such that
$p_*(x)=p(s(x))$ for every $x\in\Hh^2$ and
\[\label{eq:localized-bounds}
 \left\{\begin{aligned}
 &Dp''(s)+D\tanh s\,p'(s)+F(s,p(s))=0
       &&\text{for }s>0,\\
 &0<p(s)<1,\quad p(s)\le C h(s)
       &&\text{for }s\ge0,\\
 &p'(s)<0&&\text{for }s>0,\\
 &p'(0)=0.
 \end{aligned}\right.
\]
For every compactly supported $u_0:\Hh^2\to[0,1]$, there is a
constant $C_{u_0}>0$ such that
\[\label{eq:confinement}
 u(x,t;u_0,f_{\G})\le\min\{1,C_{u_0}h(s(x))\}
 \quad\text{for }x\in\Hh^2,\ t\ge0.
\]
\item \emph{Uniqueness and convergence.} Assume, in addition,
\eqref{eq:quotient}. Let $q\in C^2(\Hh^2)$ satisfy
$D\Delta q(x)+F(s(x),q(x))=0$ and $0<q(x)\le1$ for every
$x\in\Hh^2$. If there is a constant $C_q>0$ such that
$q(x)\le C_qh(s(x))$ for every $x\in\Hh^2$, then $q=p_*$.
For every measurable datum $u_0\not\equiv0$ satisfying
$0\le u_0(x)\le\min\{1,C_0h(s(x))\}$ for almost every
$x\in\Hh^2$ and some constant $C_0>0$,
\[\label{eq:strip-convergence}
 \lim_{t\to\infty}\sup_{\{x\in\Hh^2:|\ell(x)|\le M\}}
 |u(x,t;u_0,f_{\G})-p_*(x)|=0\quad\text{for }0\le M<\infty.
\]
\item \emph{Critical tail.} Suppose $\mu=0$, and let $p$ be the function
in part~(i). Then there is $A>0$ such that
\[\label{eq:critical-tail}
 p(s)=A e^{-s/2}\phiZ(s)+O(e^{-s})\quad\text{as }s\to\infty.
\]
\end{enumerate}
\end{theorem}

The profile $p(s)$ describes decay away from $\G$.
Along $\G$ its value is $p(0)$, and Theorem~\ref{thm:longitudinal}(iv)
locates each level $0<\theta<p(0)$ on the two half-geodesics.
The next theorem shows why the bound $q\le C_q h$ in
Theorem~\ref{thm:localized}(ii) cannot simply be omitted from the uniqueness assertion. There are other
positive stationary solutions whose ratio to $h$ is unbounded.

\begin{theorem}\label{thm:other-states}
Assume \eqref{eq:F1}, \eqref{eq:axis-regularity} and \eqref{eq:quotient},
with $d=2$, $f=f_{\G}$ and $\mu\le0<\sigma$.
Let $p\in C^2([0,\infty))$ be the function in
Theorem~\ref{thm:localized}(i), so that $p_*(x)=p(s(x))$ for
$x\in\Hh^2$. For each $\alpha\in(p(0),1)$, there exists
$p_\alpha\in C^2([0,\infty))$ such that
$x\mapsto p_\alpha(s(x))$ is a positive stationary solution
with reaction $f_{\G}$ and
\[\label{eq:other-states}
 \begin{gathered}
 p_\alpha(0)=\alpha,\quad p_\alpha'(0)=0,\quad
 p_\alpha'(s)<0\quad\text{for }s>0,\quad
 \lim_{s\to\infty}p_\alpha(s)=0,\\
 p_\alpha(s)>p(s)\quad\text{for }s\ge0,\qquad
 \sup_{s\ge0}\frac{p_\alpha(s)}{h(s)}=\infty.
 \end{gathered}
\]
\end{theorem}

\begin{example}\label{ex:periodic-reaction}
There is a $C^\infty$, $8$-periodic function $a:\R\to[D/100,D]$, constant
near zero, such that the logistic reaction $F(z,v)=a(z)v(1-v)$ satisfies
\[\label{eq:example-endpoints}
 \mu<0<\sigma_2\qquad\text{in dimension }d=2.
\]

One can take $a=D$ on $[0,1]$, $a=D/100$ on
$[26/25,199/25]$, and monotone smooth transitions on $[1,26/25]$ and $[199/25,8]$. For this choice, $\sigma_2>D/200$.
\end{example}

\subsection{Fronts and level positions}

For $u(x,t)=w(b(x),t)$, equation \eqref{eq:KPP} with $f=f_b$ reduces to
\[\label{eq:line}
 w_t=Dw_{zz}+D\beta w_z+F(z,w)\quad\text{for }z\in\R,\ t>0.
\]
An entire classical solution of~\eqref{eq:line} satisfies its
differential equation for every $(z,t)\in\R^2$.
Following \cite{BH02,BHN05}, for $c\in\R\setminus\{0\}$, a
pulsating front of speed $c$ for \eqref{eq:line} is an entire classical
solution $W:\R^2\to(0,1)$ admitting a function
$\Phi:\R^2\to(0,1)$ such that
\[\label{eq:front-definition}
 \begin{gathered}
 W(z,t)=\Phi(z,z-ct)\quad\text{for }z,t\in\R,\\
 \Phi(z+L,\xi)=\Phi(z,\xi)\quad\text{for }z,\xi\in\R,\\
 \lim_{\xi\to-\infty}\sup_{z\in\R}|\Phi(z,\xi)-1|=0,\qquad
 \lim_{\xi\to\infty}\sup_{z\in\R}\Phi(z,\xi)=0.
 \end{gathered}
\]
By the periodicity of $\Phi$ in its first variable, the two suprema
in~\eqref{eq:front-definition} can equivalently be taken over $[0,L]$.

We call $P\in C^2(\R)$ a stationary front of speed zero for
\eqref{eq:line} if
\[\label{eq:stationary-front-definition}
 \left\{\begin{aligned}
 &DP''(z)+D\beta P'(z)+F(z,P(z))=0,\quad
 0<P(z)<1\quad\text{for }z\in\R,\\
 &\lim_{z\to-\infty}P(z)=1,\qquad
 \lim_{z\to\infty}P(z)=0.
 \end{aligned}\right.
\]
The function $W:\R^2\to(0,1)$ defined by $W(z,t)=P(z)$ for
$z,t\in\R$ is then an entire classical solution of~\eqref{eq:line}.

\begin{theorem}\label{thm:fronts}
Assume \eqref{eq:F1}.
\begin{enumerate}[label=(\roman*),beginpenalty=10000]
\item \emph{Nonzero speeds.} For each
$c\in[\cp,\infty)\setminus\{0\}$, equation \eqref{eq:line} has a
pulsating front $W(z,t)=\Phi(z,z-ct)$ of speed $c$ satisfying
\eqref{eq:front-definition}, with $\Phi_\xi(z,\xi)<0$ for
$z,\xi\in\R$. Conversely, every pulsating front of~\eqref{eq:line}
with nonzero speed $c$ satisfies $c\ge\cp$.
\item \emph{Speed zero.} Equation \eqref{eq:line} has a stationary
front $P$ satisfying \eqref{eq:stationary-front-definition} if and
only if $\mu\le0$. Every such $P$ satisfies $P'(z)<0$ for $z\in\R$.
\end{enumerate}
Thus the set of front speeds for \eqref{eq:line} is $[\cp,\infty)$.
For a pulsating front in part~(i), the function
$u(x,t)=\Phi(b(x),b(x)-ct)$ solves \eqref{eq:KPP} with $f=f_b$
and initial datum $u_0(x)=\Phi(b(x),b(x))$.
For a stationary front $P$ in part~(ii), the function
$u(x,t)=P(b(x))$ solves \eqref{eq:KPP} with $f=f_b$ and
initial datum $u_0(x)=P(b(x))$.
\end{theorem}

By \eqref{eq:front-definition}, a pulsating front of speed $c\ne0$
satisfies $W(z+L,t+L/c)=W(z,t)$ for $z,t\in\R$.

Stability uses the derivative bound
\[\label{eq:derivative}
 F_u(z,v)\le a(z)\quad\text{for }z\in\R,\ 0\le v\le1.
\]
Under \eqref{eq:F1}, condition \eqref{eq:quotient} gives
$\partial_v(F(z,v)/v)\le0$ for $z\in\R$ and $0<v<1$, and hence
$vF_u(z,v)\le F(z,v)\le a(z)v$ on this range.
Continuity then gives \eqref{eq:derivative} also at $v=0$ and $v=1$.
The uniqueness statement in Theorem~\ref{thm:localized}(ii) uses
\eqref{eq:quotient}.
For a pulsating front $W(z,t)=\Phi(z,z-ct)$ of~\eqref{eq:line}
with speed $c\ne0$ and for $\eta\in\R$, define
$W^\eta(z,t)=\Phi(z,z-ct+\eta)$.
Since $W^\eta(z,t)=W(z,t-\eta/c)$, each $W^\eta$ solves
\eqref{eq:line}. For front stability, we require the following initial
inequalities for some real $\eta_1\le\eta_2$:
\[\label{eq:front-bracket-explicit}
 W^{\eta_2}(b(x),0)\le u_0(x)\le W^{\eta_1}(b(x),0)
 \quad\text{for almost every }x\in\Hh^d.
\]
By the $L$-periodicity of $F$, if $P$ satisfies
\eqref{eq:stationary-front-definition}, then $z\mapsto P(z+jL)$
is a stationary solution of~\eqref{eq:line} for every $j\in\Z$.

\begin{theorem}\label{thm:stability}
Assume \eqref{eq:F1} and \eqref{eq:derivative}, and fix
$c\in[\cp,\infty)$. For $c\ne0$, let
$W(z,t)=\Phi(z,z-ct)$ be a pulsating front of speed $c$ for
\eqref{eq:line}, as defined in~\eqref{eq:front-definition}.
For $c=0$, let $P\in C^2(\R)$ satisfy
\eqref{eq:stationary-front-definition} and set $W(z,t)=P(z)$ for
$z,t\in\R$.
\begin{enumerate}[label=(\roman*),beginpenalty=10000]
\item \emph{Weighted stability on $\Hh^d$.}
Let $u$ solve \eqref{eq:KPP} with $f=f_b$.
For $c\ne0$, assume \eqref{eq:front-bracket-explicit} for some
real $\eta_1\le\eta_2$. For $c=0$, assume that, for some
integers $j_1\le j_2$,
$P(b(x)+j_2L)\le u_0(x)\le P(b(x)+j_1L)$ for almost every
$x\in\Hh^d$.
In either case, suppose that, for some $1\le p<\infty$,
\[\label{eq:weighted-data}
 B_p=\norm{
  \frac{e^{\lam b(\cdot)}}{\phi_{\lam}(b(\cdot))}
  \bigl(u_0(\cdot)-W(b(\cdot),0)\bigr)
 }_{L^p(\Hh^d)}<\infty.
\]
For example, \eqref{eq:weighted-data} holds if $u_0(x)=W(b(x),0)$
outside a compact subset of $\Hh^d$. The integral in~\eqref{eq:weighted-data} includes all horospherical directions.
Under the hypotheses of part~(i),
\[\label{eq:weighted-stability}
 \sup_{x\in\Hh^d}
 \frac{e^{\lam(b(x)-ct)}}{\phi_{\lam}(b(x))}
 |u(x,t)-W(b(x),t)|
 \le C_p t^{-1/(2p)}e^{-\lam(c-\cp)t}B_p\quad\text{for }t\ge1.
\]
\item \emph{Uniform stability on $\Hh^d$.}
Under the assumptions of part~(i),
\[\label{eq:uniform-stability}
 \lim_{t\to\infty}
 \norm{u(\cdot,t)-W(b(\cdot),t)}_{L^\infty(\Hh^d)}=0.
\]
\item \emph{Stability on $\R$.}
Let $w$ solve \eqref{eq:line} with measurable initial datum
$w_0:\R\to[0,1]$. For $c\ne0$, assume that, for some
real $\eta_1\le\eta_2$,
$W^{\eta_2}(z,0)\le w_0(z)\le W^{\eta_1}(z,0)$ for almost
every $z\in\R$. For $c=0$, assume that, for some integers
$j_1\le j_2$, $P(z+j_2L)\le w_0(z)\le P(z+j_1L)$ for almost
every $z\in\R$. Suppose that, for some $1\le p<\infty$,
\[\label{eq:line-weighted-data}
 B_p^{(1)}=
 \left\|\frac{e^{\lam z}}{\phi_{\lam}(z)}
       \bigl(w_0(z)-W(z,0)\bigr)\right\|_{L^p(\R;\dd z)}<\infty.
\]
Then
\[\label{eq:line-stability}
 \begin{gathered}
 \sup_{z\in\R}
 \frac{e^{\lam(z-ct)}}{\phi_{\lam}(z)}|w(z,t)-W(z,t)|
 \le C_p t^{-1/(2p)}e^{-\lam(c-\cp)t}B_p^{(1)}
 \quad\text{for }t\ge1,\\
 \lim_{t\to\infty}\norm{w(\cdot,t)-W(\cdot,t)}_{L^\infty(\R)}=0.
 \end{gathered}
\]
\end{enumerate}
\end{theorem}

\begin{theorem}\label{thm:periodic-speed}
Assume \eqref{eq:F1}. Then we have
\begin{enumerate}[label=(\roman*),beginpenalty=10000]
\item \emph{Horoballs.} For $f=f_b$ and data satisfying
$\one_{\mathcal H_A}\le u_0\le\one_{\mathcal H_B}$ for real $A\le B$,
\[\label{eq:periodic-horo-speed}
 \begin{aligned}
 \lim_{t\to\infty}\sup_{\{b(x)\ge ct\}}u(x,t)&=0&&\text{for }c>\cp,\\
 \lim_{t\to\infty}\sup_{\{b(x)\le ct\}}|u(x,t)-1|&=0&&\text{for }c<\cp.
 \end{aligned}
\]
The speed $\cp$ may have either sign or be zero.
\item \emph{Balls.} Suppose $\mu>0$. For $f=f_o$ with
\eqref{eq:axis-regularity} and nonzero compactly supported data, the
solution $u$ of~\eqref{eq:KPP} satisfies
\[\label{eq:periodic-ball-speed}
 \begin{aligned}
 \lim_{t\to\infty}\sup_{\{r(x)\ge ct\}}u(x,t)&=0&&\text{for }c>\cp,\\
 \lim_{t\to\infty}\sup_{\{r(x)\le ct\}}|u(x,t)-1|&=0&&\text{for }0<c<\cp.
 \end{aligned}
\]
\item \emph{Geodesic tubes.} Suppose $\mu>0$. For $f=f_{\G}$ with
\eqref{eq:axis-regularity} and
$\eta\one_{\mathcal T_{R_1}}\le u_0\le\one_{\mathcal T_{R_2}}$ with
$0<\eta\le1$ and $0<R_1\le R_2<\infty$, the two limits in~\eqref{eq:periodic-ball-speed} hold with $r$ replaced by $s$.
\end{enumerate}
\end{theorem}

For the periodic reactions, define
\[\label{eq:periodic-position}
 m_*(t)=\cp t-\frac{3}{2\lam}\log(1+t)
 \quad\text{for }t\ge0.
\]
For the Busemann reaction, $0<\theta<1$ and $t>0$, write $E_\theta(t)=\{x\in\Hh^d:u(x,t)\ge\theta\}$.
Use upper half-space coordinates $(y,q)$ with $o=(0,1)$ and
$\gamma_+(t)=(0,e^t)$, so $z=b(x)=-\log q$. Write $X(z,y)=(y,e^{-z})$.
For $A\in\R$ and $R>0$ put
\[
 \Omega_{A,R}=\{X(z,y):z\ge A,\ y\in\R^{d-1},\ |y|\le R\},
\]
where $|y|$ is the Euclidean norm. Let $T_o\Hh^d$ be the tangent space at $o$
and $g_o$ its hyperbolic inner product. The unit vector $e=\dot\gamma_+(0)$
points towards the ideal endpoint $\xi=\gamma_+(\infty)$.
For a unit vector $\omega\in T_o\Hh^d$, $\exp_o(r\omega)$ is the point at
distance $r$ along the geodesic starting at $o$ with initial velocity $\omega$.
For $0<\delta<\pi$, define the closed sector
\[
 \mathcal C_\delta=\{o\}\cup\{\exp_o(r\omega):r>0,\ \omega\in T_o\Hh^d,\ 
       g_o(\omega,\omega)=1,\ g_o(\omega,e)\le\cos\delta\}.
\]
The last inequality means that the angle between $\omega$ and $e$ is at least
$\delta$ radians. The vertex $o$ is included separately.

\begin{theorem}\label{thm:horo-compact-periodic}
Assume \eqref{eq:F1} and $\mu>0$, and take $f=f_b$ in \eqref{eq:KPP}.
Fix $0<\theta<1$.
\begin{enumerate}[label=(\roman*)]
\item Suppose $0\le u_0\le1$ is nonzero and compactly supported.
The set $E_\theta(t)$ is compact for $t>0$. For each $A\in\R$ and $R>0$
there are $\tilde{C}_1,\tilde{T}_1>0$ such that, for $t\ge \tilde{T}_1$,
\[\label{eq:busemann-window}
 \Omega_{A,R}\cap\mathcal H_{m_*(t)-\tilde{C}_1}
 \subset E_\theta(t)\cap\Omega_{A,R}
 \subset \Omega_{A,R}\cap\mathcal H_{m_*(t)+\tilde{C}_1}.
\]
For every $0<\delta<\pi$ there are $\tilde{C}_2,\tilde{T}_2>0$ such that $m_*(t)\ge \tilde{C}_2$ and
\[\label{eq:busemann-sector}
 B(o,m_*(t)-\tilde{C}_2)\cap\mathcal C_\delta
 \subset E_\theta(t)\cap\mathcal C_\delta
 \subset B(o,m_*(t)+\tilde{C}_2)\cap\mathcal C_\delta
 \qquad(t\ge \tilde{T}_2).
\]
The constants $\tilde{C}_1,\tilde{T}_1$ may depend on
$u_0,\theta,A,R$, and $\tilde{C}_2,\tilde{T}_2$ on $u_0,\theta,\delta$.
There are $\tilde{C}_3,\tilde{T}_3>0$ such that $E_\theta(t)$ is nonempty for $t\ge \tilde{T}_3$ and
\begin{equation*}
 \bigl|H_\theta^b(t)-m_*(t)\bigr|\le \tilde{C}_3
 \quad\text{for }t\ge \tilde{T}_3,\qquad
 H_\theta^b(t):=\max_{x\in E_\theta(t)}b(x).
\end{equation*}
\item Fix $0<\eta\le1$ and $A_-\le A_+$ in $\R$. If
$\eta\one_{\mathcal H_{A_-}}\le u_0\le\one_{\mathcal H_{A_+}}$
almost everywhere, then there are $\tilde{C}_4,\tilde{T}_4>0$ such that
\[\label{eq:busemann-global}
 \mathcal H_{m_*(t)-\tilde{C}_4}\subset E_\theta(t)\subset\mathcal H_{m_*(t)+\tilde{C}_4}
 \qquad(t\ge \tilde{T}_4).
\]
The constants $\tilde{C}_4,\tilde{T}_4$ can be chosen in terms of
$\theta,\eta,A_-,A_+$ uniformly over these data.
\end{enumerate}
\end{theorem}

\begin{theorem}\label{thm:periodic-bramson-distance}
Assume \eqref{eq:F1}, \eqref{eq:axis-regularity}, and $\mu>0$.
Let $u$ solve \eqref{eq:KPP} with $0\le u_0\le1$.
Use $m_*$ from \eqref{eq:periodic-position}.
\begin{enumerate}[label=(\roman*),beginpenalty=10000]
\item Suppose $f=f_o$, and let $u_0$ be nonzero and compactly
supported. For every $0<\theta<1$, there are
$\tilde{C}_{\theta,1},\tilde{T}_{\theta,1}>0$ such that
\[\label{eq:periodic-bramson-balls}
 B(o,m_*(t)-\tilde{C}_{\theta,1})
 \subset\{x\in\Hh^d:u(x,t)\ge\theta\}
 \subset B(o,m_*(t)+\tilde{C}_{\theta,1})
 \quad\text{for }t\ge \tilde{T}_{\theta,1}.
\]
\item Suppose $f=f_\G$, and assume that
\[
 \eta\one_{\mathcal T_{R_1}}\le u_0\le\one_{\mathcal T_{R_2}}
 \quad\text{on }\Hh^d,
 \qquad 0<\eta\le1,\quad 0<R_1\le R_2<\infty.
\]
For every $0<\theta<1$, there are $\tilde{C}_{\theta,2},\tilde{T}_{\theta,2}>0$ such that
\[\label{eq:periodic-bramson-tubes}
 \mathcal T_{m_*(t)-\tilde{C}_{\theta,2}}
 \subset\{x\in\Hh^d:u(x,t)\ge\theta\}
 \subset\mathcal T_{m_*(t)+\tilde{C}_{\theta,2}}
 \quad\text{for }t\ge \tilde{T}_{\theta,2}.
\]
\end{enumerate}
\end{theorem}

Write $f_u(x,v)=\partial_v f(x,v)$. Theorem~\ref{thm:common} also applies to spatially nonperiodic reactions
with $f_u(x,0)=a_0>0$ for every $x\in\Hh^d$. Under the hypotheses below
it gives an explicit decay rate for square-integrable data when $a_0\le D(d-1)^2/4$.
Assume that $\sup_{x\in\Hh^d}\|f(x,\cdot)\|_{C^2([0,1])}<\infty$.
The uniform spatial H\"older hypothesis means that there are
$\alpha\in(0,1)$ and $C_{\rm H}<\infty$ such that
\[\label{eq:common-spatial-regularity}
 |f(x,v)-f(y,v)|\le C_{\rm H}\dist(x,y)^\alpha
 \quad {\text{for }x,y\in\Hh^d,\quad\dist(x,y)\le1,\quad0\le v\le1}.
\]
The remaining reaction assumptions are
\[\label{eq:common-assumptions}
 \begin{gathered}
 f(x,0)=f(x,1)=0,\qquad f_u(x,0)=a_0>0
 \quad {\text{for }x\in\Hh^d},\\
 0<f(x,v)\le a_0v
 \quad {\text{for }x\in\Hh^d,\quad0<v<1},\\
 \sup_{x\in\Hh^d}f_u(x,1)<0,\qquad
 \inf_{x\in\Hh^d,\ v\in[\delta,1-\delta]}f(x,v)>0
 \quad {\text{for }0<\delta<1/2}.
 \end{gathered}
\]
In every subsequent reference to \eqref{eq:common-assumptions}, the
uniform $C^2$ bound in the population variable and the spatial
H\"older condition \eqref{eq:common-spatial-regularity} are included.
The reactions in~\eqref{eq:patterns} satisfy these assumptions
when $a(z)\equiv a_0$.

Define
\[\label{eq:common-position}
 \begin{gathered}
 \rho_d=\frac{(d-1)^2}{4},\quad
 \lambda_{a_0}=\sqrt{\frac{a_0}{D}},\quad
 c_{a_0}=2\sqrt{Da_0}-D(d-1),\\
 m_{a_0}(t)=c_{a_0}t-\frac{3}{2\lambda_{a_0}}\log(1+t)
 \quad\text{for }t\ge0.
 \end{gathered}
\]
When $a(z)\equiv a_0$, equation \eqref{eq:cell} gives
$k(\lambda)=D\lambda^2-D\beta\lambda+a_0$ for $\lambda\in\R$.
Hence $\mu=a_0-D\rho_d$, $\lam=\lambda_{a_0}$ and
$\cp=c_{a_0}$ by \eqref{eq:parameters}.
\begin{theorem}\label{thm:common}
Assume \eqref{eq:common-assumptions},
\eqref{eq:common-spatial-regularity}, and
$\sup_{x\in\Hh^d}\|f(x,\cdot)\|_{C^2([0,1])}<\infty$.
\begin{enumerate}[label=(\roman*),beginpenalty=10000]
\item \emph{Extinction.} If $a_0\le D\rho_d$ and $u_0\in L^2(\Hh^d)$, then
the solution of~\eqref{eq:KPP} satisfies \eqref{eq:threshold-extinction}.
More precisely,
\[\label{eq:common-decay}
 \norm{u(\cdot,t)}_{L^{\infty}(\Hh^d)}
 \le C t^{-3/4}e^{(a_0-D\rho_d)t}\norm{u_0}_{L^{2}(\Hh^d)}\quad\text{for }t\ge1.
\]
\item \emph{Balls.} If $a_0>D\rho_d$ and $u_0$ is nonzero and compactly
supported, then for each $0<\theta<1$ there are
$\hat{C}_{\theta,1},\hat{T}_{\theta,1}<\infty$ such that
\[\label{eq:common-balls}
 B(o,m_{a_0}(t)-\hat{C}_{\theta,1})
 \subset\{x\in\Hh^d:u(x,t)\ge\theta\}
 \subset B(o,m_{a_0}(t)+\hat{C}_{\theta,1})\quad\text{for }t\ge \hat{T}_{\theta,1}.
\]
\item \emph{Geodesic tubes.} If $a_0>D\rho_d$ and
$\eta\one_{\mathcal T_{R_1}}\le u_0\le\one_{\mathcal T_{R_2}}$ for
$0<\eta\le1$ and $0<R_1\le R_2<\infty$, then for each $0<\theta<1$
there are $\hat{C}_{\theta,2},\hat{T}_{\theta,2}<\infty$ such that
\[\label{eq:common-tubes}
 \mathcal T_{m_{a_0}(t)-\hat{C}_{\theta,2}}
 \subset\{x\in\Hh^d:u(x,t)\ge\theta\}
 \subset\mathcal T_{m_{a_0}(t)+\hat{C}_{\theta,2}}\quad\text{for }t\ge \hat{T}_{\theta,2}.
\]
\item \emph{Horoballs.} If
$\one_{\mathcal H_A}\le u_0\le\one_{\mathcal H_B}$ for $A\le B$, then
for every $a_0>0$ and each $0<\theta<1$ there are
$\hat{C}_{\theta,3},\hat{T}_{\theta,3}<\infty$ such that
\[\label{eq:common-horoballs}
 \mathcal H_{m_{a_0}(t)-\hat{C}_{\theta,3}}
 \subset\{x\in\Hh^d:u(x,t)\ge\theta\}
 \subset\mathcal H_{m_{a_0}(t)+\hat{C}_{\theta,3}}\quad\text{for }t\ge \hat{T}_{\theta,3}.
\]
\end{enumerate}
\end{theorem}

For fixed $D$ and $a_0$, the logarithmic coefficient
$3/(2\lambda_{a_0})=\tfrac32\sqrt{D/a_0}$ is independent of
dimension, while $c_{a_0}=2\sqrt{Da_0}-D(d-1)$.

\subsection{Spreading along a geodesic}

In this subsection $d=2$, $\beta=1$, and $\sigma=\sigma_2$.
Choose one side of $\G$ as positive and write
$\rho=s$ on that side and $\rho=-s$ on the other. The coordinate $\ell(x)$ is determined by the nearest point
$\gamma(\ell(x))\in\G$ specified after \eqref{eq:coordinates}.
Multiplication by $\sqrt{\cosh s}$ maps
$L^2((0,\infty),\cosh s\dd s)$ unitarily onto $L^2(0,\infty;\dd s)$.
Conjugating $v\mapsto Dv''+D\tanh s\,v'+a(s)v$, with $v'(0)=0$,
by this multiplication gives the half-line operator
\[\label{eq:J-intro}
 \J=D\frac{\mathrm d^2}{\mathrm ds^2}+a(s)-\frac D4-\frac D4\sech^2s
 \quad {\text{on }(0,\infty)},\qquad v'(0)=0.
\]
The operator in~\eqref{eq:J-intro} is realised on $L^2(0,\infty;\dd s)$
with domain $\{v\in H^2(0,\infty):v'(0)=0\}$. Its upper spectral
endpoint is $\sigma$. In the regime $\mu\le0<\sigma$ put
\[\label{eq:longitudinal-intro}
 \J_\parallel(\lambda)=\J+D\lambda^2\sech^2s,\qquad
 \Lambda(\lambda)=\sup\sigma(\J_\parallel(\lambda))
 \quad {\text{for }\lambda\in\R},\qquad
 c_{\parallel}=\min_{\lambda>0}\frac{\Lambda(\lambda)}\lambda.
\]
For every $\lambda\in\R$, the term $D\lambda^2\sech^2s$ acts
by pointwise multiplication, and $\J_\parallel(\lambda)$ has the
same operator domain as $\J$. In the next theorem, $p_*$ is the
smallest positive stationary solution supplied by
Theorem~\ref{thm:dimension}(ii); its representation
$p_*(x)=p(s(x))$ and transverse decay are given in
Theorem~\ref{thm:localized}(i).

On $\G$ the profile equals $p(0)\in(0,1)$ by~\eqref{eq:localized-bounds}.
For a solution $u$ and $0<\theta<p(0)$, write
$E_\theta^\pm(t)=\{r\ge0:u(\gamma(\pm r),t)\ge\theta\}$.
With $\lambda_\parallel$ denoting the minimiser in~\eqref{eq:longitudinal-intro},
characterised in part~\textup{(i)} below, put
$m(t)=c_\parallel t-\frac{3}{2\lambda_\parallel}\log(1+t)$ for $t\ge0$.

\begin{theorem}\label{thm:longitudinal}
Assume \eqref{eq:F1} and \eqref{eq:axis-regularity}, with $d=2$
and $\mu\le0<\sigma$.
Let $u$ solve
\eqref{eq:KPP} with $f=f_\G$ and $0\le u_0\le1$.
\begin{enumerate}[label=(\roman*),beginpenalty=10000]
\item \emph{Spectral characterisation.} For every real $\lambda$,
$\Lambda(\lambda)$ is a simple isolated eigenvalue of
$\J_\parallel(\lambda)$, with a positive eigenfunction in $L^2(0,\infty;\dd s)$. The function $\lambda\mapsto\Lambda(\lambda)$ is even, real
analytic, and strictly convex.
The minimum defining $c_{\parallel}>0$ has a unique minimiser
$\lambda_{\parallel}>0$, characterised by
\[\label{eq:longitudinal-critical}
 \lambda_{\parallel}\Lambda'(\lambda_{\parallel})
 =\Lambda(\lambda_{\parallel}).
\]
\item \emph{Outer bound.} If $u_0$ is compactly supported, then
\[\label{eq:longitudinal-outer}
 \lim_{t\to\infty}\sup_{\{x\in\Hh^2:|\ell(x)|\ge ct\}}
 u(x,t;u_0,f_{\G})=0\quad\text{for }c>c_{\parallel}.
\]
\item \emph{Inner bound.} If $u_0$ is nonzero
and compactly supported, then
\[\label{eq:longitudinal-inner}
 \liminf_{t\to\infty}\inf_{\{x\in\Hh^2:|\ell(x)|\le ct\}}
       [u(x,t;u_0,f_{\G})-p_*(x)]\ge0
 \quad {\text{for }0<c<c_\parallel}.
\]
\item \emph{Level positions on $\G$.} Suppose $u_0$ is nonzero and compactly supported.
For every $0<\theta<p(0)$ there are $C_{\theta,u_0},T_{\theta,u_0}>0$ such that for both signs
and every $t\ge T_{\theta,u_0}$,
\[\label{eq:longitudinal-level-speed}
 [0,m(t)-C_{\theta,u_0}]\subset E_\theta^\pm(t)
 \subset [0,m(t)+C_{\theta,u_0}].
\]
\item \emph{Profile convergence.} If \eqref{eq:quotient} also holds and
$u_0$ is nonzero and compactly supported, then
\[\label{eq:longitudinal-profile}
 \lim_{t\to\infty}\sup_{\{x\in\Hh^2:|\ell(x)|\le ct\}}
    |u(x,t;u_0,f_{\G})-p_*(x)|=0\quad\text{for }0<c<c_{\parallel}.
\]
\end{enumerate}
\end{theorem}

\begin{remark}\label{rem:outermost-axial-level}
Under the assumptions of Theorem~\ref{thm:longitudinal}(iv), let
$0<\theta<p(0)$ and $t\ge T_{\theta,u_0}$. The sets $E_\theta^\pm(t)$
are closed by continuity of $u(\cdot,t)$ and are nonempty and bounded
by~\eqref{eq:longitudinal-level-speed}. Taking their maxima
$L_\theta^\pm(t)=\max E_\theta^\pm(t)$ in the two inclusions gives,
for each $\varepsilon\in\{+,-\}$,
\begin{equation*}
 \left|L_\theta^\varepsilon(t)-c_\parallel t+
 \frac{3}{2\lambda_\parallel}\log(1+t)\right|\le C_{\theta,u_0}.
\end{equation*}
\end{remark}

Section~\ref{sec:basic} derives the coordinate operators and spectral
thresholds. Section~\ref{sec:fronts} proves the existence and stability of fronts
and determines spreading speeds and logarithmic level positions. Section~\ref{sec:localized} treats the planar stationary states and proves
Theorem~\ref{thm:longitudinal}(i)--(iii) and (v).
Section~\ref{sec:axial-delay} proves part~(iv), the logarithmic delay
on each half-geodesic.

\section{Geometry, extinction and persistence}\label{sec:basic}

\subsection{Coordinates and geometric formulas}\label{sec:coordinates}

Let $\nabla^g$ be the Levi-Civita connection of $g$. We use
$\nabla^2v(Y,Z)=g(\nabla^g_Y\nabla v,Z)$ and
$\Delta v=\operatorname{tr}_g\nabla^2v$. In local coordinates,
write $(g^{ij})=(g_{ij})^{-1}$ and $|g|=\det(g_{ij})$. Then
\[\label{eq:geometric-operators}
 \dd V=\sqrt{|g|}\,\mathrm dx^1\cdots\mathrm dx^d,\qquad
 \Delta v=\frac1{\sqrt{|g|}}
  \sum_{i,j=1}^d\partial_i\bigl(\sqrt{|g|}\,g^{ij}\partial_jv\bigr).
\]
For a tangent vector $X$, write $|X|_g=\sqrt{g(X,X)}$ and $\langle X,Y\rangle_g=g(X,Y)$ for tangent vectors $X,Y$. Then
$\int_{\Hh^d}v\Delta v\dd V
=-\int_{\Hh^d}|\nabla v|_g^2\dd V$
for $v\in C_c^\infty(\Hh^d;\R)$; see~\cite[Chapters~2 and~5]{Lee18}. We derive the coordinate
operators and the nearest-point formula defining $\ell$ in the
standard models \cite[Chapters~3--5]{Ratcliffe}; see also
\cite{Helgason}.

The ideal boundary $\partial_\infty\Hh^d$ is the set of equivalence
classes of unit-speed geodesic rays. Two such rays
$\gamma_1,\gamma_2:[0,\infty)\to\Hh^d$ are equivalent when
$\sup_{R\ge0}\dist(\gamma_1(R),\gamma_2(R))<\infty$.
Each class has a unique representative starting at any prescribed
point $o$. We write $\xi=\gamma_+(\infty)$ for the class determined
by the ray in~\eqref{eq:coordinates}.

In Lorentz space let
$\langle X,Y\rangle_{\rm Lor}=-X_0Y_0+\sum_{j=1}^dX_jY_j$.
The hyperboloid model is
\[\label{eq:hyperboloid}
 \Hh^d=\{X\in\R^{d+1}:-X_0^2+X_1^2+\cdots+X_d^2=-1,\ X_0>0\},
\]
with the positive definite metric obtained by restricting $\langle\cdot,\cdot\rangle_{\rm Lor}$ to
the tangent spaces of $\Hh^d$. The hyperbolic distance satisfies
$\cosh\dist(X,Y)=-\langle X,Y\rangle_{\rm Lor}$.
We denote the standard round metric on $\mathbb S^m$ by
$g_{\mathbb S^m}$ and its area measure by $\dd\omega$. For $m\ge1$, $\Delta_{\mathbb S^m}$ denotes the Laplace--Beltrami operator of this round metric. For $m=0$,
$\mathbb S^0=\{-1,1\}$ has counting measure and no differential angular
term. For an integrable even function $v:\R\to\R$,
$\int_\R v(\rho)\dd\rho=2\int_0^\infty v(s)\dd s$. The factor $2$ cancels between the numerator and denominator of a
Rayleigh quotient; integration with respect to $\dd V$ retains it.
We use the real branches of $\operatorname{arsinh}$ and
$\operatorname{artanh}$.

\begin{lemma}\label{lem:coordinates}
The Busemann function in~\eqref{eq:coordinates} satisfies
$|\nabla b|_g=1$ and $\Delta b=\beta$. There are global coordinates
$(z,y)\in\R\times\R^{d-1}$, with $z=b(x)$, in which
\[\label{eq:horo-metric}
 \begin{gathered}
 g=\dd z^2+e^{2z}|\dd y|^2,\qquad
 \dd V=e^{\beta z}\dd z\dd y,\\
 \Delta=\partial_{zz}+\beta\partial_z+e^{-2z}\Delta_y
 \quad {\text{for }(z,y)\in\R\times\R^{d-1}}.
 \end{gathered}
\]
Here $\Delta_y=\sum_{j=1}^{d-1}\partial^2/\partial y_j^2$ is the Euclidean Laplacian in $y$. In polar coordinates about $o$,
\[\label{eq:polar-metric}
 \begin{gathered}
 g=\dd r^2+\sinh^2r\,g_{\mathbb S^{d-1}},\qquad
 \dd V=\sinh^{d-1}r\dd r\dd\omega,\\
 \Delta=\partial_{rr}+\beta\coth r\,\partial_r
                +\sinh^{-2}r\,\Delta_{\mathbb S^{d-1}}
 \quad {\text{for }r>0,\quad\omega\in\mathbb S^{d-1}}.
 \end{gathered}
\]
About a complete geodesic, away from the geodesic itself,
\[\label{eq:fermi-metric}
 \begin{gathered}
 g=\dd s^2+\cosh^2s\,\dd\ell^2+
                         \sinh^2s\,g_{\mathbb S^{d-2}},\qquad
 \dd V=m_d(s)\dd s\dd\ell\dd\omega,\\
 \Delta=\partial_{ss}+B_d(s)\partial_s+
          \sech^2s\,\partial_{\ell\ell}
          +\sinh^{-2}s\,\Delta_{\mathbb S^{d-2}} \quad {\text{for }s>0,\quad\ell\in\R,\quad\omega\in\mathbb S^{d-2}}.
 \end{gathered}
\]
When $d=2$, let $\rho$ denote signed distance from $\G$. In the
global coordinates $(\rho,\ell)\in\R^2$, the metric and Laplacian are
\[\label{eq:plane-fermi}
 g=\dd\rho^2+\cosh^2\rho\,\dd\ell^2,\qquad
 \Delta=\partial_{\rho\rho}+\tanh\rho\,\partial_\rho
                     +\sech^2\rho\,\partial_{\ell\ell} \quad {\text{for }(\rho,\ell)\in\R^2}.
\]
Here $s=|\rho|$.
\end{lemma}

\begin{proof}
By the isometry invariance of the hyperbolic metric, we may take
$o=(0,1)$ and $\gamma_+(R)=(0,e^R)$ in the upper half-space model
$(y,q)\in\R^{d-1}\times(0,\infty)$.
The metric is $q^{-2}(|\dd y|^2+\dd q^2)$, and the upper half-space distance formula gives
\[\label{eq:upper-halfspace-distance}
 \cosh\dist\bigl((y,q),(0,e^R)\bigr)
 =1+\frac{|y|^2+(q-e^R)^2}{2qe^R}
 \quad {\text{for }y\in\R^{d-1},\quad q>0,\quad R\in\R}.
\]
Using $\operatorname{arcosh}A=\log(2A)+o(1)$ in~\eqref{eq:upper-halfspace-distance} gives
$\dist((y,q),(0,e^R))-R\to-\log q$, hence $b(y,q)=-\log q$. Substitution of $z=-\log q$ therefore yields
the metric and volume element in~\eqref{eq:horo-metric}, while
\eqref{eq:geometric-operators} gives the Laplacian and hence
$|\nabla b|_g=1$ and $\Delta b=\beta$. For the upper half-space model see~\cite[Section~4.6]{Ratcliffe}; its Chapter~3 treats the Lorentz model.

For polar coordinates, take $o=(1,0,\ldots,0)$ and
$X=(\cosh r,\sinh r\,\omega)$, where $\omega\in\mathbb S^{d-1}$.
The Lorentz distance formula gives $\dist(o,X)=r$; taking inner
products of the coordinate derivatives and applying
\eqref{eq:geometric-operators} then proves \eqref{eq:polar-metric}.

By \cite[Chapter~5]{Ratcliffe}, an isometry sends the complete geodesic
to $\gamma(\ell)=(\cosh\ell,\sinh\ell,0,\ldots,0)$. For
$X\in\Hh^d$, put $X_\perp=(X_2,\ldots,X_d)$. Since
\eqref{eq:hyperboloid} gives
$X_0^2-X_1^2=1+|X_\perp|^2>0$, the coordinates
\[\label{eq:nearest-geodesic-coordinates}
 s=\operatorname{arsinh}|X_\perp|,\qquad
 \ell=\operatorname{artanh}(X_1/X_0)
 \quad {\text{for }X\in\Hh^d}
\]
are well defined. For every $t\in\R$, the Lorentz distance formula gives
\[\label{eq:nearest-geodesic-distance}
 \cosh\dist(X,\gamma(t))=\cosh s\cosh(t-\ell)
 \quad {\text{for }X\in\Hh^d,\quad t\in\R}.
\]
By \eqref{eq:nearest-geodesic-distance}, the function
$t\mapsto\dist(X,\gamma(t))$ has its unique minimum at $t=\ell$, proving that $s=\dist(X,\G)$ and that
$\gamma(\ell)$ is the unique nearest point. For $s>0$, setting $\omega=X_\perp/|X_\perp|$ yields
\[\label{eq:fermi-embedding}
 X=(\cosh s\cosh\ell,\cosh s\sinh\ell,\sinh s\,\omega)
 \quad {\text{for }s\ge0,\quad\ell\in\R,\quad\omega\in\mathbb S^{d-2}}.
\]
Taking inner products of the coordinate derivatives in~\eqref{eq:fermi-embedding} gives the metric factors
$1$, $\cosh^2s$ and $\sinh^2s$, and hence the volume density
$m_d(s)=\cosh s\,\sinh^{d-2}s$. Formula \eqref{eq:geometric-operators} then
gives $B_d=m_d'/m_d$ and \eqref{eq:fermi-metric}.
For fixed $\ell$ and $\omega$, the curves $s\mapsto X(s,\ell,\omega)$
in~\eqref{eq:fermi-embedding} are unit-speed geodesics orthogonal to $\G$, giving the Fermi coordinates of~\cite[Chapter~5, pp.~133--136]{Lee18}; their global coverage away
from $\G$ follows from \eqref{eq:nearest-geodesic-distance}.
Finally, when $d=2$, the choice $\rho=\operatorname{arsinh}X_2$
combines the two signs in $\mathbb S^0$ and gives
\eqref{eq:plane-fermi}, with $s=|\rho|$.
\end{proof}

In the horospherical coordinates, the maps
\[\label{eq:horo-shift-isometry}
 T_\tau(z,y)=(z+\tau,e^{-\tau}y),\quad\text{for }(z,y)\in\R\times\R^{d-1},\ \tau\in\R,
\]
preserve the metric in~\eqref{eq:horo-metric} and satisfy
$b(T_\tau x)=b(x)+\tau$. In particular, $F(b(T_\tau x),v)=F(b(x)+\tau,v)$ for
$\tau\in\R$ and $v\in[0,1]$.
The formula \eqref{eq:nearest-geodesic-coordinates} shows that $\ell$ is
smooth even on $\G$, and the inverse metric gives
\[\label{eq:longitudinal-geometric-identities}
 |\nabla\ell|_g^2=\sech^2s,\qquad \Delta\ell=0
 \quad {\text{for }x\in\Hh^d}.
\]
The functions $r$ and $s$ are one-Lipschitz by the triangle inequality;
the same holds for $b$ by its defining limit. Integrating $|\nabla\ell|_g=\sech s\le1$ along a minimising
geodesic gives $|\ell(x)-\ell(y)|\le\dist(x,y)$, so $\ell$
is also one-Lipschitz. Closed bounded sets are compact: in the hyperboloid model,
the ball of radius $R$ centred at $(1,0,\ldots,0)$ lies in
$1\le X_0\le\cosh R$, and an isometry takes any other centre to $(1,0,\ldots,0)$
\cite[Chapter~6]{Lee18}.

The metrics in~\eqref{eq:horo-metric}--\eqref{eq:fermi-metric}
give the Hessians used in Lemma~\ref{lem:coordinate-limit}. On $s>0$ set $\vartheta_\ell=\cosh s\,\dd\ell$;
this is a unit covector in the longitudinal direction. Then
\[\label{eq:coordinate-hessians}
 \begin{aligned}
 \nabla^2 b&=g-\dd b\otimes\dd b,&&\text{for }x\in\Hh^d,\\
 \nabla^2 r&=\coth r\,(g-\dd r\otimes\dd r),&&\text{for }r>0,\\
 \nabla^2 s&=\tanh s\,\vartheta_\ell\otimes\vartheta_\ell
 +\coth s\,(g-\dd s\otimes\dd s
                   -\vartheta_\ell\otimes\vartheta_\ell),&&\text{for }s>0.
 \end{aligned}
\]
In local coordinates $(z,\theta^1,\ldots,\theta^{d-1})$, consider a metric
$\dd z^2+\sum_{A,B=1}^{d-1}g_{AB}(z,\theta)\dd\theta^A\dd\theta^B$.
Write $\Gamma^k_{ij}$ for the connection coefficients defined by
$\nabla^g_{\partial_i}\partial_j={\sum_{k=1}^d}\Gamma^k_{ij}\partial_k$.
The formula for the Levi-Civita connection gives
$\Gamma^z_{AB}=-\tfrac12\partial_zg_{AB}$ and
$(\nabla^2 z)_{AB}=\tfrac12\partial_zg_{AB}$. The components involving the normal coordinate satisfy
$(\nabla^2z)_{zz}=(\nabla^2z)_{zA}=0$. Substitution of $e^{2z}$, $\sinh^2r$, $\cosh^2s$ and
$\sinh^2s$ in this connection formula proves
\eqref{eq:coordinate-hessians}; see also the treatment of distance
Hessians in~\cite[Chapters~6, 10 and~11]{Lee18}.
For the hypersurface $\{z=c\}$ with unit normal $n=\nabla z$,
the second fundamental form on tangent vectors $Y,Z$ is
$g(\nabla^g_Yn,Z)$. Thus the normal
drift $\Delta z$ is the sum of the principal curvatures, and
$\Delta z/(d-1)$ is their average \cite[Chapter~8]{Lee18}.

If $x\mapsto U(r(x))$ is $C^2$ at $o$, then $U_r(0)=0$.
Likewise, if $x\mapsto V(s(x))$ is $C^2$ across $\G$, invariance
under orthogonal transformations of the normal space gives $V_s(0)=0$.
To check the coefficient regularity, use geodesic normal coordinates
$\xi\in\R^d$ at $o$, or normal coordinates $\xi\in\R^{d-1}$ to $\G$.
Then $r=|\xi|$ or $s=|\xi|$, respectively. Let $i,j$ range over
$1,\ldots,d$ in the first case and $1,\ldots,d-1$ in the second, and
write $\delta_{ij}$ for the Kronecker delta. For $q=|\xi|>0$,
\[\label{eq:distance-composition-hessian}
 \partial_{\xi_i\xi_j}F(q,v)
 =F_{zz}(q,v)\frac{\xi_i\xi_j}{q^2}
  +\frac{F_z(q,v)}q\left(\delta_{ij}-\frac{\xi_i\xi_j}{q^2}\right) \quad {\text{for }q=|\xi|>0,\quad v\in[0,1]}.
\]
Assumption~\eqref{eq:axis-regularity} and continuity of $F_{zz}$ make
the right-hand side of~\eqref{eq:distance-composition-hessian}
tend to $F_{zz}(0,v)\delta_{ij}$ as $q\downarrow0$,
uniformly in $v\in[0,1]$. The mixed derivatives
$\partial_{\xi_i}\partial_v F(|\xi|,v)
=F_{zu}(q,v)\xi_i/q$ tend to $0$ as $q\downarrow0$,
since differentiating \eqref{eq:axis-regularity} in $v$ gives
$F_{zu}(0,v)=0$.
Thus $x\mapsto F(r(x),v)$ and $x\mapsto F(s(x),v)$
extend with their derivatives of order at most two across $o$ and
$\G$, respectively.
For initial data depending only on $b$, $r$ or $s$, respectively,
take $f=f_b$, $f=f_o$ or $f=f_{\G}$ as defined in~\eqref{eq:patterns}.
The solutions of~\eqref{eq:KPP} with these respective reactions
and initial data have the forms
$u(x,t)=w(b(x),t)$, $u(x,t)=U(r(x),t)$ or $u(x,t)=V(s(x),t)$.
Take measurable functions $w_0:\R\to[0,1]$ and
$U_0,V_0:[0,\infty)\to[0,1]$. Initial data $w_0(b(x))$, $U_0(r(x))$
or $V_0(s(x))$ give, respectively, $w(z,0)=w_0(z)$,
$U(r,0)=U_0(r)$ or $V(s,0)=V_0(s)$. For $t>0$, the equations and boundary conditions are
\[\label{eq:three-reductions}
 \begin{aligned}
 w_t&=Dw_{zz}+D\beta w_z+F(z,w)
       &&\text{for }z\in\R,\ t>0,\\
 U_t&=DU_{rr}+D\beta\coth r\,U_r+F(r,U)
       &&\text{for }r>0,\ t>0,\\
 V_t&=DV_{ss}+DB_d(s)V_s+F(s,V)
       &&\text{for }s>0,\ t>0,\\
 U_r(0,t)&=V_s(0,t)=0&&\text{for }t>0.
 \end{aligned}
\]
For the geodesic-distance equation in~\eqref{eq:three-reductions},
regularity at $s=0$ means that $\xi\mapsto V(|\xi|,t)$
is $C^2$ on the $(d-1)$-dimensional normal space to $\G$ for $t>0$. In particular,
$B_d(s)V_s(s,t)\to(d-2)V_{ss}(0,t)$ as $s\downarrow0$ for $d\ge3$.
The Neumann conditions in~\eqref{eq:three-reductions} and this
limit of $B_dV_s$ hold for $t>0$; the initial datum may
be merely measurable, as in Section~\ref{sec:intro}.

\begin{lemma}\label{lem:drifts}
For $s>0$,
\[\label{eq:drift-comparison}
 \beta\tanh s\le B_d(s)\le\beta\coth s
 \quad {\text{for }s>0,\quad d\ge2}.
\]
If $d\ge3$, then
\[\label{eq:drift-excess}
 B_d(s)-\beta
 =2\left(\frac{d-2}{e^{2s}-1}-\frac1{e^{2s}+1}\right)>0
 \quad {\text{for }s>0,\quad d\ge3}.
\]
For $d=2$, $B_2(s)-1=\tanh s-1<0$. In every dimension,
$B_d(s)-\beta=O(e^{-2s})$ as $s\to\infty$.
\end{lemma}

\begin{proof}
By \eqref{eq:tube-weight}, subtracting $\beta\tanh s$ from $B_d(s)$
and $B_d(s)$ from $\beta\coth s$ gives, respectively,
$(d-2)(\coth s-\tanh s)$ and $\coth s-\tanh s$, which proves
\eqref{eq:drift-comparison}. Expressing $\coth s-1$ and
$1-\tanh s$ in terms of $e^{2s}$ gives
\eqref{eq:drift-excess}; after the positive common denominator is
removed, its sign is that of $(d-3)e^{2s}+d-1$, which is positive for
$d\ge3$. The same expressions give $B_2(s)-1=\tanh s-1<0$ and
$|B_d(s)-\beta|\le C_d e^{-2s}$ for $s\ge1$ in every dimension $d\ge2$.
\end{proof}

Lemma~\ref{lem:tube-sandwich} uses signed distance from a
totally geodesic hyperplane to construct an upper comparison;
see~\cite[Section~2]{GGQ}.
A totally geodesic hyperplane $\mathcal P$ is an isometric copy of
$\Hh^{d-1}$ whose intrinsic geodesics are ambient geodesics. Put it in
the form $X_d=0$ in the hyperboloid. Every point has a unique
representation
$X=(\cosh\varrho\,Y,\sinh\varrho)$ with
$Y\in\Hh^{d-1}$ and $\varrho\in\R$. The distance formula shows that
$ (Y,0)$ is the nearest point on $\mathcal P$ and that $|\varrho|$ is the
distance to $\mathcal P$; the sign specifies its side. Differentiation
gives
\[\label{eq:hyperplane-metric-full}
 \begin{gathered}
 g=\dd\varrho^2+\cosh^2\varrho\,g_{\Hh^{d-1}},\qquad
 \dd V=\cosh^\beta\varrho\,\dd\varrho\dd V_{\Hh^{d-1}},\\
 \Delta=\partial_{\varrho\varrho}+\beta\tanh\varrho\,\partial_\varrho
                       +\sech^2\varrho\,\Delta_{\Hh^{d-1}} \quad {\text{for }\varrho\in\R}.
 \end{gathered}
\]
Here $g_{\Hh^{d-1}}$, $\dd V_{\Hh^{d-1}}$ and $\Delta_{\Hh^{d-1}}$ are the metric, volume element and Laplace--Beltrami operator on the hyperplane $\mathcal P$ with its induced metric. See also \cite[Section~2.3]{BirindelliMazzeo}.
The metric induced on $\{\varrho=c\}$ is
$\cosh^2c\,g_{\Hh^{d-1}}$, a constant rescaling of the curvature
$-1$ metric. For $d\ge3$, its sectional curvature is therefore
$-\sech^2c$; at $c=0$ the induced metric is $g_{\Hh^{d-1}}$.

The coordinate $\varrho$ ranges over $\R$. For an even $C^2$ function $W(\varrho)$, $W'(0)=0$; its
restriction to $[0,\infty)$ therefore satisfies the Neumann condition. For $d=2$, $\mathcal P$ is a geodesic and
$W(\varrho)=W(s)$ with $s=|\varrho|$. The drift inequalities required for comparison are
given by Lemma~\ref{lem:drifts}. An isometric action has cohomogeneity
one when its orbit space is one-dimensional. For example, let
$1\le j\le d-2$ and let $\mathcal P_j$ be a totally geodesic copy of
$\Hh^j$. In the hyperboloid, the parametrisation
$X=(\cosh s\,Y,\sinh s\,\omega)$, with $Y\in\Hh^j$,
$\omega\in\mathbb S^{d-j-1}$ and $s>0$, gives
\[\label{eq:subspace-tube-metric}
 g=\dd s^2+\cosh^2s\,g_{\Hh^j}
                   +\sinh^2s\,g_{\mathbb S^{d-j-1}}
 \quad {\text{for }s>0,\quad1\le j\le d-2}.
\]
The orthogonal coordinate derivatives give the displayed metric,
and the Lorentz distance formula identifies $(Y,0)$ as the nearest
point, so $s=\dist(X,\mathcal P_j)$. Isometries of $\mathcal P_j$ and rotations of the
normal sphere act transitively on each positive-distance level. For
$j=1$, the levels $\{X\in\Hh^d:\dist(X,\mathcal P_1)=s\}$ are geodesic
tubes. Their principal curvatures include both $\tanh s$ and $\coth s$
when $d\ge3$. For the isometric actions on the levels
$\{X\in\Hh^d:\dist(X,\mathcal P_j)=s\}$, see~\cite[Sections~1 and~4]{BerndtTamaru}.

\subsection{Solutions and comparison}

For $t\ge0$, let $P_t=e^{tD\Delta}$ denote the heat semigroup
generated by $D\Delta$. For $t>0$, its minimal positive heat kernel
$p_t(x,y)$ is taken with respect to $\dd V(y)$, so that
$(P_t\psi)(x)=\int_{\Hh^d}p_t(x,y)\psi(y)\dd V(y)$ for bounded
measurable $\psi:\Hh^d\to\R$ and $x\in\Hh^d$. By isometry invariance, the diagonal value $p_t(x,x)$ is independent of $x$; we denote it by $p_t(0)$ for $t>0$. The local stochastic equations construct Brownian motion with generator
$D\Delta$ up to its lifetime; see~\cite[Chapters~1 and~3]{Hsu}.
The following estimate shows that this lifetime is infinite.
The smooth function $\mathscr V(x)=1+r(x)^2$ satisfies
\[\label{eq:radial-Lyapunov}
 \Delta\mathscr V=2+2\beta r\coth r
       \le (2+3\beta)\mathscr V\quad\text{for }x\in\Hh^d,
\]
where the value at $r=0$ is understood by continuity. The inequality
uses $r\coth r\le r+1$ and $2r\le1+r^2$. More generally, for a smooth
vector field $Y$ with $|Y|_g\le K$, the generator
$\mathscr B=D\Delta+\langle Y,\nabla\,\cdot\,\rangle_g$ obeys
$\mathscr B\mathscr V\le C\mathscr V$, with $C=D(2+3\beta)+K$.
For a diffusion $X_t$ with generator $\mathscr B$ and $X_0=x$, let
$\tau_R=\inf\{t\ge0:r(X_t)\ge R\}$, where $R>r(x)$.
It\^o's formula applied up to $\tau_R$ to
$e^{-Ct}\mathscr V(X_t)$ gives
\[\label{eq:nonexplosion-exit-bound}
 \mathbb P_x(\tau_R\le T)
 \le \frac{e^{CT}(1+r(x)^2)}{1+R^2}
       \quad\text{for }T>0.
\]
The stopped stochastic integral has zero expectation because its
coefficients are bounded on $\overline{B(o,R)}\times[0,T]$.
Letting $R\to\infty$ excludes a finite lifetime. For $Y=0$, this proves
stochastic completeness, that is,
$\int_{\Hh^d}p_t(x,y)\dd V(y)=1$.

Isometry invariance makes $p_t(x,y)$ a function of
$R=\dist(x,y)$; we also denote it by $p_t(R)$. If $p^{(1)}_\tau$
is the heat kernel of $e^{\tau\Delta}$, then $p_t=p^{(1)}_{Dt}$.
We apply the unit-diffusion bounds of~\cite{DaviesMandouvalos,Davies}, and the explicit formulas of~\cite{GrigoryanNoguchi}, through this time change. For large-time heat behaviour, including the
distinction between radial and general integrable data, see~\cite{Vazquez22}.

We use the following weak inequality for the piecewise smooth comparison
functions. On an open region $\mathcal O\subset\Hh^d\times(0,\infty)$,
a function $v:\mathcal O\to[0,1]$ with locally square-integrable
spatial weak gradient is a weak subsolution of~\eqref{eq:KPP} if, for
every nonnegative $\eta\in C_c^\infty(\mathcal O)$,
\[\label{eq:weak-subsolution-definition}
 \int_{\mathcal O}
 \bigl[-v\eta_t+D\langle\nabla v,\nabla\eta\rangle_g
                    -f(x,v)\eta\bigr]\dd V\dd t\le0
 \quad\text{for }0\le\eta\in C_c^\infty(\mathcal O).
\]
A weak supersolution satisfies the reverse inequality. For a linear equation $v_t=D\Delta v+A(x,t)v$, the weak inequality
\eqref{eq:weak-subsolution-definition} uses $A(x,t)v$ in place of
$f(x,v)$ and allows nonnegative $v$ of arbitrary size. For comparison on a space--time region, the subsolution is required
to lie below the supersolution on its initial and lateral boundary. For measurable
initial data, the initial trace is taken in
$L^1_{\mathrm{loc}}(\Hh^d)$ as $t\downarrow0$; positive-time lateral
values of the comparison functions are continuous.
The comparison functions below are continuous at positive times and
locally Lipschitz across the hypersurfaces where their defining formulas change. Across a moving joining hypersurface, continuity of the function
prevents a surface measure in its first time derivative.
A jump in the spatial normal derivative can contribute a surface
measure to the Laplacian; its sign is checked at each join.
Whenever a comparison function is replaced by its minimum with $1$
or its maximum with $0$, we check
\eqref{eq:weak-subsolution-definition}, or its reverse, for that function.

\begin{lemma}\label{lem:wellposed}
Assume \eqref{eq:F1} and choose $f\in\{f_b,f_o,f_\G\}$,
with \eqref{eq:axis-regularity} for $f=f_o$ or $f=f_\G$.
Then \eqref{eq:KPP} has a unique bounded mild solution. It takes values in
$[0,1]$ and satisfies
\[\label{eq:duhamel}
 \begin{aligned}
 u(x,t)={}&\int_{\Hh^d}p_t(x,y)u_0(y)\dd V(y)\\
 &+\int_0^t\int_{\Hh^d}p_{t-\tau}(x,y)
            f(y,u(y,\tau))\dd V(y)\dd\tau
 \quad {\text{for }x\in\Hh^d,\quad t>0}.
 \end{aligned}
\]
It is a classical solution on $\Hh^d\times(0,\infty)$. If $u_0$ is nonzero, $u(x,t)>0$ for $t>0$.
If $u_0$ is not identically one, $u(x,t)<1$ for $t>0$.
If a bounded subsolution lies below a bounded supersolution at
$t=0$, that order holds for all $t>0$.
The same comparison holds for continuous, piecewise smooth weak barriers
whose differential inequalities hold in the distributional sense.
\end{lemma}

\begin{proof}
Extend $f(x,\cdot)$ to a globally Lipschitz function agreeing with
$f$ on $[0,1]$, and let $L_f$ be a Lipschitz constant independent of $x$.
Since $\|P_t g\|_\infty\le\|g\|_\infty$, the right-hand side of~\eqref{eq:duhamel} with fixed datum $u_0$ has Lipschitz constant at most
$TL_f$ on $L^\infty(\Hh^d\times(0,T))$. Take
$0<T\le1/[2(1+L_f)]$, so $TL_f<1/2$.
The fixed point of this contraction satisfies \eqref{eq:duhamel}
and is therefore a bounded mild solution of~\eqref{eq:KPP}. By the interior parabolic estimates in~\cite{Friedman}, this solution has a spatial weak gradient and a
classical representative on every compact cylinder separated from
time zero. Approximation of $u_0$ on compact neighbourhoods, together
with the Gaussian heat-kernel tail, gives
$P_tu_0\to u_0$ in $L^1_{\mathrm{loc}}$ as $t\downarrow0$.
Since the integral term in~\eqref{eq:duhamel} is $O(t)$ in the
supremum norm, the mild solution has the same initial trace.
The comparison argument below bounds the mild solution between zero
and one, allowing it to be continued for all time.

Let $v_1$ be a bounded subsolution and $v_2$ a bounded supersolution
with $v_1(\cdot,0)\le v_2(\cdot,0)$. Put $z=v_1-v_2$ and $w=z_+$.
By the weak inequality \eqref{eq:weak-subsolution-definition},
$z_t-D\Delta z\le c_z(x,t)z$ in distributions, where $c_z$ is $[f(x,v_1)-f(x,v_2)]/(v_1-v_2)$ where $v_1\ne v_2$,
and zero where $v_1=v_2$. Then $|c_z|\le K$ for a common
Lipschitz constant $K\ge0$. The positive part $w$ has zero
$L^1_{\rm loc}$ initial trace. Choose $\alpha>\beta$ and set
\[\label{eq:comparison-integrable-weight}
 \omega(x)=e^{-\alpha\sqrt{1+r(x)^2}},\qquad
 \int_{\Hh^d}\omega\dd V<\infty,\qquad
 |\nabla\omega|_g\le\alpha\omega \quad {\text{for }x\in\Hh^d,\quad\alpha>\beta}.
\]
The smoothness of $r^2$ at $o$ makes $\omega$ smooth there, while
the polar volume element in~\eqref{eq:polar-metric} gives its
integrability. Take a smooth cutoff $\chi_R$ equal to one on
$B(o,R)$ and zero outside $B(o,2R)$, with $0\le\chi_R\le1$ and
$|\nabla\chi_R|_g\le C/R$. Testing the difference inequality against
$w\omega\chi_R^2$ gives, in the time-integrated sense,
\[\label{eq:weak-comparison-energy}
 \begin{split}
 \frac12 E_R'(t)+\frac D2\int_{\Hh^d}
                  \omega\chi_R^2|\nabla w|_g^2\dd V
 &\le (K+D\alpha^2)E_R(t)
       +4D\int_{\Hh^d}\omega w^2|\nabla\chi_R|_g^2\dd V,\\
 E_R(t)&=\int_{\Hh^d}\omega\chi_R^2w(x,t)^2\dd V.
 \end{split}
\]
Indeed, \eqref{eq:comparison-integrable-weight} bounds the term
containing $\nabla\omega$ by
$\tfrac D4{\int_{\Hh^d}}\omega\chi_R^2|\nabla w|_g^2\dd V+D\alpha^2E_R$. The term containing $\nabla\chi_R$ is bounded by
$\tfrac D4{\int_{\Hh^d}}\omega\chi_R^2|\nabla w|_g^2\dd V
+4D{\int_{\Hh^d}}\omega w^2|\nabla\chi_R|_g^2\dd V$.

For weak barriers, we justify testing with $z_+\omega\chi_R^2$ by time averaging.
Set $z_h(x,t)=h^{-1}\int_t^{t+h}z(x,\tau)\dd\tau$ for
$0<t<T-h$, whose time derivative is
$(z(x,t+h)-z(x,t))/h$, locally in $L^2$. Approximation by smooth spatial functions in $H^1$ makes
$(z_h)_+\omega\chi_R^2$ an admissible test in~\eqref{eq:weak-subsolution-definition}.
On each compact cylinder $\Omega\times(t_0,t_1)$ in the comparison
region, continuity of time translation gives $z_h\to z$ in
$L^2((t_0,t_1);H^1(\Omega))$ and $(c_zz)_h\to c_zz$ in
$L^2(\Omega\times(t_0,t_1))$. Hence $(c_zz)_h(z_h)_+$ converges in $L^1$ on the cylinder
to $c_zw^2\le Kw^2$, and the Sobolev chain rule and lower
semicontinuity of the gradient integral yield
\eqref{eq:weak-comparison-energy} as $h\downarrow0$.
The inequality holds first for almost every pair of time endpoints;
continuity of the solutions and barriers gives it for all positive
endpoints. The distributional derivatives include the contributions
at the joining sets.

For fixed $T$, boundedness of $w$ and
\eqref{eq:comparison-integrable-weight} bound the cutoff error in~\eqref{eq:weak-comparison-energy} by $C_T/R^2$. Letting
$R\to\infty$, using dominated convergence for $E_R$ and discarding
the nonnegative gradient term, gives
$E(t)\le E(s)+2(K+D\alpha^2)\int_s^tE(\tau)\dd\tau$ for
$0<s<t<T$, where $E(t)={\int_{\Hh^d}}\omega w(x,t)^2\dd V$.
The zero local trace, uniform boundedness of $w$ and integrability
of $\omega$ imply $E(s)\to0$ as $s\downarrow0$.
Gronwall's inequality gives $E(t)=0$, hence $w=0$ and
$v_1\le v_2$ on $\Hh^d\times(0,T)$. Applying comparison in both directions proves uniqueness.

Applying this comparison to the constant solutions zero and one
keeps the mild solution in $[0,1]$ and completes its continuation.
For comparison in singular coordinates, see also \cite{PPT11}. By \eqref{eq:F1}, $u\ge P_tu_0$, so strict positivity
of the heat kernel gives $u(x,t)>0$ for nonzero $u_0$ and $t>0$.
Likewise, $f(x,1)=0$ and the Lipschitz bound give
$f(x,u)\le K(1-u)$, whence
$1-u(x,t)\ge e^{-Kt}P_t(1-u_0)(x)>0$ when $u_0\not\equiv1$.
The same heat-kernel approximation gives
$u(x,t)\to u_0(x)$ as $t\downarrow0$ at every continuity point $x$ of $u_0$.
\end{proof}

Existence and comparison require a uniform Lipschitz bound
in the population variable, local spatial H\"older bounds and
$f(x,0)=f(x,1)=0$. Conditions
\eqref{eq:common-assumptions} and \eqref{eq:common-spatial-regularity},
together with $\sup_{x\in\Hh^d}\|f(x,\cdot)\|_{C^2([0,1])}<\infty$,
give these bounds and endpoint values. Thus the solutions in
Theorem~\ref{thm:common} have the mild and positive-time classical
interpretation of Lemma~\ref{lem:wellposed}.

Let $u$ solve \eqref{eq:KPP} with initial datum $u_0$ and
reaction $f$. For an isometry $T$, put $u^T(x,t)=u(Tx,t)$. Invariance of the
Laplacian gives
\[\label{eq:isometry-transformed-equation}
 \left\{\begin{aligned}
 (u^T)_t&=D\Delta u^T+f(Tx,u^T)
       &&\text{for }x\in\Hh^d,\ t>0,\\
 u^T(x,0)&=u_0(Tx)&&\text{for almost every }x\in\Hh^d.
 \end{aligned}\right.
\]
Writing $f_T(x,v)=f(Tx,v)$ and $u_{0,T}(x)=u_0(Tx)$, uniqueness gives
\[\label{eq:isometry-covariance}
 u(Tx,t;u_0,f)=u(x,t;u_{0,T},f_T)
 \quad {\text{for }x\in\Hh^d,\quad t>0}.
\]
Thus $u^T=u$ whenever $f_T=f$ and $u_{0,T}=u_0$.
For $f=f_b$ and $T_\tau$ in~\eqref{eq:horo-shift-isometry},
$f_{T_\tau}(x,v)=F(b(T_\tau x),v)=F(b(x)+\tau,v)$.
For $\tau\in L\Z$, one has $f_{T_\tau}=f_b$, and hence
$u^{T_\tau}=u$ whenever $u_{0,T_\tau}=u_0$.
Identity \eqref{eq:isometry-covariance} itself holds for every
isometry $T$. See~\cite[Chapters~2 and~5]{Lee18} for metric covariance and
\cite[Corollary~3.3]{GGQ} for its use with homogeneous reactions.
The isometries preserving $b$, $r$ or $s$ act transitively on each
level set $\{b=z\}$, $\{r=R\}$ or $\{s=S\}$, respectively,
for $z\in\R$ and $R,S\ge0$. In the coordinates of
Lemma~\ref{lem:coordinates}, one may use, respectively, Euclidean
motions in $y$, rotations about $o$, and longitudinal translations
together with normal orthogonal transformations, including reflections.
Thus initial data depending only on $b$, $r$ or $s$
yield, respectively, the first, second or third evolution equation
in~\eqref{eq:three-reductions}. For general initial data, the coordinate solutions of~\eqref{eq:three-reductions}, composed with $b$, $r$ or $s$, provide
upper and lower comparisons for \eqref{eq:KPP}. Support and integrability refer to the functions
$x\mapsto w(b(x),t)$, $x\mapsto U(r(x),t)$ or
$x\mapsto V(s(x),t)$ on $(\Hh^d,\dd V)$.

\begin{lemma}\label{lem:compactness}
Fix $0<\alpha<1$ and an open time interval $I$. Let
$u_n:\Hh^d\times I\to[0,1]$ solve
$(u_n)_t=D\Delta u_n+f_n(x,u_n)$. Assume that, for each compact spatial
set, the functions $f_n$ have uniformly bounded $C^\alpha$ norms in $x$,
uniformly for $0\le v\le1$, and a uniform Lipschitz constant in $v$.
After passage to a subsequence, there are $f_\infty$ and $u_\infty$ such
that $f_n\to f_\infty$ locally uniformly on $\Hh^d\times[0,1]$ and
\[\label{eq:compactness-topology}
 u_n\to u_\infty
 \quad\hbox{in }C^{2+\alpha',\,1+\alpha'/2}_{\rm loc}
                  (\Hh^d\times I)\text{ as }n\to\infty,
 \quad\text{for }0<\alpha'<\alpha.
\]
Here the convergence is on compact cylinders contained in the open
space--time domain. The limit $u_\infty$ solves
$(u_\infty)_t=D\Delta u_\infty+f_\infty(x,u_\infty)$
on $\Hh^d\times I$.
For isometries $T_n$, the assertion also applies to $u_n(T_nx,t)$
and the reactions $(x,v)\mapsto f_n(T_nx,v)$, provided that on every compact spatial set their $C^\alpha$ norms in $x$
and Lipschitz constants in $v\in[0,1]$ are bounded independently of $n$.
For a fixed reaction $f$, assume that $f(\cdot,v)$ is locally
$C^\alpha$ uniformly in $v\in[0,1]$ and $f(x,\cdot)$ is Lipschitz
uniformly for $x$ in compact sets. A solution $0\le u\le1$ of
$u_t=D\Delta u+f(x,u)$ on $\Hh^d\times(0,\infty)$ has a locally uniform stationary
limit as $t\to\infty$ whenever it is monotone in time.
\end{lemma}

\begin{proof}
The hypotheses give a common bound for $f_n$ on each compact spatial
set times $[0,1]$. Interior parabolic $W^{2,1}_p$ estimates on nested
compact cylinders therefore bound $u_n$ for every finite $p$.
Choose $p>d+2$ with $1-(d+2)/p>\alpha$; by parabolic Sobolev
embedding, the resulting $C^{\alpha,\alpha/2}$ bound for $u_n$
and the spatial H\"older bounds and the Lipschitz bounds in the second
variable of $f_n$ give uniform $C^{\alpha,\alpha/2}$ bounds
for $f_n(x,u_n(x,t))$. Applying interior Schauder estimates on a
smaller cylinder now bounds $u_n$ in
$C^{2+\alpha,1+\alpha/2}$ \cite{Friedman}. Compact embedding
and a diagonal extraction yield \eqref{eq:compactness-topology},
while Arzel\`a--Ascoli on compact subsets of
$\Hh^d\times[0,1]$ gives $f_n\to f_\infty$ locally uniformly.
Passing to the limit gives the equation for $u_\infty$.

For $f(x,v)=F(z(x),v)$, choose $z\in\{b,r,s\}$ and an
isometry $T$, and put $K_z=\|F_z\|_\infty$ and
$K_u=\|F_u\|_\infty$, with both suprema over
$[0,L]\times[0,1]$. The one-Lipschitz property of $z$ gives
\[\label{eq:pullback-coefficient-control}
 |F(z(Tx),v)-F(z(Ty),w)|
 \le K_z\dist(x,y)+K_u|v-w|
 \quad {\text{for }x,y\in\Hh^d,\;0\le v,w\le1}.
\]
Thus $(x,v)\mapsto F(z(Tx),v)$ has H\"older and Lipschitz
bounds independent of the isometry $T$. The parabolic estimates on isometric balls of fixed radius
consequently have uniform constants.

Finally, let $u:\Hh^d\times(0,\infty)\to[0,1]$ solve
$u_t=D\Delta u+f(x,u)$ and be monotone in $t$. Set
$q(x)=\lim_{t\to\infty}u(x,t)$. By
\eqref{eq:compactness-topology}, every sequence of time translates
tending to infinity has a subsequence converging on compact
space--time cylinders. For every fixed $\tau\in\R$, monotonicity gives
$u(x,t_n+\tau)\to q(x)$ as $t_n\to\infty$. Thus every
locally convergent sequence of time translates has the same limit
$q$, which is independent of $\tau$ and satisfies $D\Delta q+f(x,q)=0$.
Local equicontinuity makes the convergence uniform on each compact
spatial set. Otherwise, a sequence of time translates violating uniform
convergence would have a subsequence converging locally uniformly to
$q$, a contradiction.
\end{proof}

\subsection{Linear estimates and stationary solutions}

Fix a real potential $A\in L^\infty(\Hh^d)$. The lower-semibounded
closed form $\int_{\Hh^d}(D|\nabla v|_g^2-A|v|^2)\dd V$ on
$H^1(\Hh^d)$ defines $-\mathscr L_A$, where
$\mathscr L_A=D\Delta+A$ is self-adjoint on $L^2(\Hh^d)$.
Here $H^1(\Hh^d)$ consists of locally weakly differentiable
functions $v$ with $\int_{\Hh^d}(|v|^2+|\nabla v|_g^2)\dd V<\infty$.
For $g\in L^2(\Hh^d)$, $v(x,t)=(e^{t\mathscr L_A}g)(x)$ solves
$v_t=D\Delta v+A(x)v$ on $\Hh^d\times(0,\infty)$ with
$\|e^{t\mathscr L_A}g-g\|_2\to0$ as $t\downarrow0$.
The space $C_c^\infty(\Hh^d)$ is dense in $H^1(\Hh^d)$
for the form norm: multiply by the ball cutoffs used in
Lemma~\ref{lem:wellposed}, then smooth in finitely many coordinate charts.
We use this realisation for each potential in this subsection; see~\cite{Davies,ReedSimon}. Write $s_A=\sup\sigma(\mathscr L_A)$.
The Rayleigh principle gives
\[\label{eq:rayleigh-full}
 s_A=\sup_{0\ne\zeta\in C_c^\infty(\Hh^d)}
 \frac{{\int_{\Hh^d}}(A\zeta^2-D|\nabla\zeta|_g^2)\dd V}{{\int_{\Hh^d}}\zeta^2\dd V}.
\]
We use real tests in Rayleigh quotients and identities containing
$\zeta^2$; the complex identities follow by summing over real and
imaginary parts. When applying the smoothing estimate \eqref{eq:smoothing} to the nonlinear problem,
we assume $u_0\in L^2(\Hh^d)$.

\begin{lemma}\label{lem:smoothing}
Let $A\in L^\infty(\Hh^d)$ be real with $A(x)\le A_+$,
and let $g\in L^2(\Hh^d)$. For $\tau>0$,
\[\label{eq:smoothing}
 \norm{e^{\tau\mathscr L_A}g}_{L^{\infty}(\Hh^d)}
 \le e^{A_+\tau}p_{2\tau}(0)^{1/2}\norm{g}_{L^{2}(\Hh^d)}
 \quad {\text{for }\tau>0,\quad g\in L^2(\Hh^d)}.
\]
If $s_A\le0$ and $\ker\mathscr L_A=\{0\}$, then
\[\label{eq:strong-extinction}
 \lim_{t\to\infty}\norm{e^{t\mathscr L_A}g}_{L^{\infty}(\Hh^d)}=0
 \quad\text{for }g\in L^2(\Hh^d).
\]
If $f(x,v)\le A(x)v$ for $x\in\Hh^d$ and $0\le v\le1$,
the solution $u$ of~\eqref{eq:KPP} with $0\le u_0\le1$ and
$u_0\in L^2(\Hh^d)$ satisfies \eqref{eq:smoothing} with
$e^{\tau\mathscr L_A}g$ replaced by $u(\cdot,\tau)$ and $g$ by $u_0$.
If also $s_A\le0$ and $\ker\mathscr L_A=\{0\}$, then
\eqref{eq:strong-extinction} holds with $e^{t\mathscr L_A}g$
replaced by $u(\cdot,t)$.
\end{lemma}

\begin{proof}
For bounded compactly supported $g$, linear comparison, or equivalently
the bounded-potential Feynman--Kac formula \cite{Pinsky}, gives
$|e^{\tau\mathscr L_A}g|\le e^{A_+\tau}P_\tau|g|$.
By symmetry and the semigroup property of the heat kernel,
\[\label{eq:kernel-L2}
 \int_{\Hh^d}p_\tau(x,y)^2\dd V(y)=p_{2\tau}(x,x)=p_{2\tau}(0)
 \quad {\text{for }x\in\Hh^d,\quad\tau>0}.
\]
Cauchy--Schwarz applied to \eqref{eq:kernel-L2} proves
\eqref{eq:smoothing} for such data. For arbitrary
$g\in L^2(\Hh^d)$, approximate $g$ in $L^2$ by bounded compactly
supported functions. For fixed $\tau>0$, their images under $e^{\tau\mathscr L_A}$
form a Cauchy sequence in $L^\infty(\Hh^d)$ by \eqref{eq:smoothing}, and boundedness of $e^{\tau\mathscr L_A}$
on $L^2$ identifies the limit with $e^{\tau\mathscr L_A}g$.
Thus \eqref{eq:smoothing} holds throughout $L^2$.

By the spectral theorem,
$\norm{e^{t\mathscr L_A}g}_{L^{2}(\Hh^d)}\le e^{s_A t}\norm g_{L^{2}(\Hh^d)}$.
Combining this bound with \eqref{eq:smoothing} over the final
$\tau$ units of time yields
\[\label{eq:smoothing-spectral-bound}
 \norm{e^{t\mathscr L_A}g}_{L^{\infty}(\Hh^d)}
 \le e^{A_+\tau}p_{2\tau}(0)^{1/2}
             e^{s_A(t-\tau)}\norm{g}_{L^{2}(\Hh^d)}\quad\text{for }t\ge\tau>0.
\]
For $s_A<0$, decay follows directly from
\eqref{eq:smoothing-spectral-bound}. When $s_A=0$ and
$\ker\mathscr L_A=\{0\}$, the same spectral theorem expresses
$\norm{e^{t\mathscr L_A}g}_{L^{2}(\Hh^d)}^2$ as the integral of $e^{2t\lambda}$
against the finite spectral measure associated with $g$, supported
on $(-\infty,0]$ and assigning zero mass to $\{0\}$.
Dominated convergence therefore gives decay in $L^2$, and
\eqref{eq:smoothing}, applied over the final unit of time, proves
\eqref{eq:strong-extinction}. The nonlinear assertions follow by
linear comparison, since
$0\le u(\cdot,t)\le e^{t\mathscr L_A}u_0$ when $f(x,v)\le A(x)v$ for $x\in\Hh^d$ and $0\le v\le1$.
\end{proof}

For the unit-diffusion kernel, the estimate of~\cite{DaviesMandouvalos,Davies} is
\[
 p^{(1)}_\tau(R)\asymp
 \tau^{-d/2}(1+R)(1+R+\tau)^{(d-3)/2}
 \exp\left(-\rho_d\tau-\frac\beta2R-\frac{R^2}{4\tau}\right)
 \quad\text{for }R\ge0,\ \tau>0.
\]
Substituting $\tau=Dt$ gives the heat-kernel estimate used below:
\[\label{eq:heat-estimate}
 p_t(R)\asymp
 (Dt)^{-d/2}(1+R)(1+R+Dt)^{(d-3)/2}
 \exp\left(-D\rho_d t-\frac\beta2 R-\frac{R^2}{4Dt}\right),
 \quad\text{for }R\ge0,\ t>0,
\]
Here $\asymp$ means that the ratio of the left-hand side to the
right-hand side lies in $[C^{-1},C]$, with $C\ge1$ depending only on $d$.
In particular, for $t\ge1$,
\[\label{eq:heat-diagonal}
 p_t(0)\le C_{d,D}t^{-3/2}e^{-D\rho_d t}
 \quad {\text{for }t\ge1}.
\]
For $0<t\le1$, \eqref{eq:heat-estimate} and
\eqref{eq:kernel-L2} give
\[\label{eq:small-time-L2-smoothing}
 \|P_t\|_{L^2(\Hh^d)\to L^\infty(\Hh^d)}
       =p_{2t}(0)^{1/2}\asymp t^{-d/4}
       \quad\text{for }0<t\le1.
\]
The constants in~\eqref{eq:small-time-L2-smoothing} may depend on
$d$ and the fixed $D>0$.
For fixed $0<t\le1$ and $x_0\in\Hh^d$, \eqref{eq:kernel-L2} gives
$\|p_t(x_0,\cdot)/p_{2t}(0)^{1/2}\|_2=1$ and
$(P_t[p_t(x_0,\cdot)/p_{2t}(0)^{1/2}])(x_0)=p_{2t}(0)^{1/2}$.
Hence the datum attains the operator norm in~\eqref{eq:small-time-L2-smoothing}. Continuity makes this point value a lower bound for the
essential supremum. 

A positive Dirichlet eigenfunction and the uniform quadratic lower
bound on the reaction give the smallest positive stationary solution in Lemma~\ref{lem:minimal-state}. On a smooth
bounded domain, the principal eigenfunction is positive in the
interior and has strictly negative outward normal derivative; see~\cite[Chapters~3, 6 and~8]{GilbargTrudinger} and \cite{BNV94}.

\begin{lemma}\label{lem:minimal-state}
Let $f:\Hh^d\times[0,1]\to\R$ satisfy the spatial H\"older
condition \eqref{eq:common-spatial-regularity}, with
$f(x,\cdot)\in C^2([0,1])$ for every $x\in\Hh^d$ and
$\sup_{x\in\Hh^d}\|f(x,\cdot)\|_{C^2([0,1])}<\infty$.
Assume
\begin{equation*}
\begin{gathered}
 f(x,0)=f(x,1)=0,\qquad
 0<f(x,v)\le f_u(x,0)v
 \quad\text{for }x\in\Hh^d,\quad0<v<1,\\
 \inf_{x\in\Hh^d}f_u(x,0)>0,\qquad
 \sup_{x\in\Hh^d}f_u(x,1)<0,\\
 \inf_{x\in\Hh^d,\,v\in[\delta,1-\delta]}f(x,v)>0
 \quad\text{for }0<\delta<1/2.
\end{gathered}
\end{equation*}
If the upper spectral endpoint of $D\Delta+f_u(x,0)$ on
$L^2(\Hh^d)$ is positive, the following conclusions hold. There is a pointwise
smallest positive stationary solution $p$ with $0<p\le1$. It is the limit, as $t\to\infty$, of solutions that are
nondecreasing in time and start from arbitrarily small positive
multiples of one fixed, nonzero, compactly supported stationary
weak subsolution. Every nonzero solution of~\eqref{eq:KPP} satisfies
\[\label{eq:liminf-local}
 \liminf_{t\to\infty}\inf_{x\in K}(u(x,t)-p(x))\ge0
 \quad\text{for every compact }K\subset\Hh^d.
\]
Every isometry preserving $f$ also preserves $p$.
\end{lemma}

For the reactions in~\eqref{eq:patterns}, the regularity, comparison
and uniform Taylor bound used below follow from \eqref{eq:F1},
\eqref{eq:axis-regularity} where required, and
\eqref{eq:pullback-coefficient-control}.

\begin{proof}
By \eqref{eq:rayleigh-full} and positivity of the upper spectral
endpoint, a compactly supported test has positive Rayleigh quotient.
Place its support in a smooth bounded connected domain $\Omega$;
the principal Dirichlet eigenvalue $\nu$ of $D\Delta+f_u(x,0)$ on
$\Omega$ is then positive. Let $\varphi>0$ be its principal
eigenfunction, normalised by $\max_{\overline\Omega}\varphi=1$. The Taylor bound $f(x,v)\ge f_u(x,0)v-C_2v^2$ for
$x\in\Hh^d$ and $0\le v\le1$, with $C_2$ independent of $x$,
permits the choice
\[\label{eq:seed-size}
 \varepsilon_0=\min\{1/2,\nu/(2(C_2+1))\}.
\]
For $0<\varepsilon\le\varepsilon_0$, the expression
$D\Delta(\varepsilon\varphi)+f(x,\varepsilon\varphi)$ on $\Omega$ is at least
$\varepsilon(\nu-C_2\varepsilon)\varphi\ge0$.
Define $\widetilde\varphi=\varphi$ on $\Omega$ and
$\widetilde\varphi=0$ on $\Hh^d\setminus\Omega$. Then
Green's formula gives
\[\label{eq:zero-extended-seed}
 \Delta\widetilde\varphi
 =\one_\Omega\Delta\varphi
       -(\partial_{\mathbf n}\varphi)\delta_{\partial\Omega}
 \quad\hbox{in distributions on }\Hh^d.
\]
Here $\langle\delta_{\partial\Omega},\eta\rangle
=\int_{\partial\Omega}\eta\dd S_g$, where $\dd S_g$ is the
induced area measure and $\mathbf n$ is the outward unit normal.
By the Hopf lemma, $-\partial_{\mathbf n}\varphi\ge0$, so
\eqref{eq:zero-extended-seed} and
$D\Delta(\varepsilon\varphi)+f(x,\varepsilon\varphi)\ge0$ in $\Omega$ prove that
$\varepsilon\widetilde\varphi$ is a stationary weak subsolution.
Let $u_\varepsilon$ solve \eqref{eq:KPP} with reaction $f$ and
$u_\varepsilon(x,0)=\varepsilon\widetilde\varphi(x)$. By comparison, $u_\varepsilon$ is nondecreasing in
time and remains bounded by $1$; Lemma~\ref{lem:compactness}
therefore gives the positive stationary limit
$p_\varepsilon(x)=\lim_{t\to\infty}u_\varepsilon(x,t)$.

To show that $p_\varepsilon$ is independent of $\varepsilon$, fix
$0<\varepsilon_1<\varepsilon_2\le\varepsilon_0$ and put
$\alpha(t)=\varepsilon_1e^{\nu t/2}$ until
$T_{12}=2\nu^{-1}\log(\varepsilon_2/\varepsilon_1)$.
By \eqref{eq:seed-size},
$(\partial_t-D\Delta)(\alpha(t)\varphi)-f(x,\alpha(t)\varphi)
\le\alpha(-\nu/2+C_2\alpha)\varphi\le0$
on $\Omega\times[0,T_{12}]$.
On $\partial\Omega\times[0,T_{12}]$,
$\alpha(t)\varphi=0\le u_{\varepsilon_1}$, so comparison yields
$u_{\varepsilon_1}(x,T_{12})\ge\varepsilon_2\varphi(x)$ in
$\Omega$, and $u_{\varepsilon_1}(x,T_{12})\ge0
=\varepsilon_2\widetilde\varphi(x)$ on $\Hh^d\setminus\Omega$. Applying comparison from $T_{12}$
onward and then passing to the stationary limits gives
$p_{\varepsilon_1}\ge p_{\varepsilon_2}$. The reverse inequality $p_{\varepsilon_1}\le p_{\varepsilon_2}$
follows by comparison from
$\varepsilon_1\widetilde\varphi\le\varepsilon_2\widetilde\varphi$, proving
$p_{\varepsilon_1}=p_{\varepsilon_2}=:p$.

Every positive stationary solution $q\le1$ has a positive minimum
on $\overline\Omega$, and hence dominates
$\varepsilon\widetilde\varphi$ for a sufficiently small
$\varepsilon>0$. Comparison with $u_\varepsilon$ gives $q\ge p$,
which proves minimality. Likewise, by Lemma~\ref{lem:wellposed},
a solution $u$ of~\eqref{eq:KPP} from a nonzero datum is strictly
positive at time one and satisfies
$u(x,1)\ge\varepsilon\widetilde\varphi(x)$ for a sufficiently
small $\varepsilon$. The local uniform convergence of
$u_\varepsilon$ to $p$ then proves \eqref{eq:liminf-local}.
Finally, if an isometry $T$ preserves $f$, then $p\circ T$
is a positive stationary solution, so minimality gives $p\le p\circ T$.
Applying the same argument to $T^{-1}$ gives $p\circ T\le p$,
and hence $p\circ T=p$.
\end{proof}

\begin{lemma}\label{lem:ground-transform}
Let $A\in L^\infty(\Hh^d)$ be real, let
$H\in C^2(\Hh^d)$ be positive, and let
$\zeta\in C_c^\infty(\Hh^d;\R)$. Then
\[\label{eq:ground-identity}
 {\int_{\Hh^d}}(A\zeta^2-D|\nabla\zeta|_g^2)\dd V
 ={\int_{\Hh^d}}\frac{(D\Delta+A)H}{H}\zeta^2\dd V
       -D{\int_{\Hh^d}} H^2\left|\nabla\frac\zeta H\right|_g^2\dd V\quad {\text{for }\zeta\in C_c^\infty(\Hh^d)}.
\]
In particular a positive weak supersolution
$(D\Delta+A)H\le\nu H$ implies $s_A\le\nu$. For the positive, locally Lipschitz functions
$H(x)=H_0(\dist(x,\Sigma))$ used below, where
$H_0:[0,\infty)\to(0,\infty)$, $\Sigma$ is a point or a smooth
submanifold of codimension at least two, and $H$ is $C^2$ away
from $\Sigma$, the identity holds as the limit of the identities on
$\{x\in\Hh^d:\dist(x,\Sigma)>\varepsilon\}$ as $\varepsilon\downarrow0$.
\end{lemma}

\begin{proof}
Set $v=\zeta/H$. Then
$|\nabla(Hv)|_g^2=H^2|\nabla v|_g^2+
\nabla H\cdot\nabla(Hv^2)$. Multiplication by $D$ and integration
by parts prove
\eqref{eq:ground-identity}. If $(D\Delta+A)H\le\nu H$ weakly,
testing this inequality with $\zeta^2/H$, with compactly supported approximation where required,
gives $s_A\le\nu$ by \eqref{eq:rayleigh-full}.

For $H(x)=H_0(\dist(x,\Sigma))$ as in the statement, work
on a compact tubular neighbourhood on which $H$, $|\nabla H|_g$ and
$H^{-1}$ are bounded. By the coordinate metrics in
Lemma~\ref{lem:coordinates}, the boundary at distance $\varepsilon$
from a codimension-$m$ submanifold has area
$O(\varepsilon^{m-1})$ over each fixed compact portion.
The boundary term in the preceding integration by parts therefore tends
to zero when $m\ge2$. The integrals of $A\zeta^2$, $|\nabla\zeta|_g^2$ and
$H^2|\nabla(\zeta/H)|_g^2$ converge by the local bounds for $H$, $H^{-1}$
and $|\nabla H|_g$. Hence the integral containing $(\Delta H)/H$ has a finite
excision limit, which proves the stated extension of~\eqref{eq:ground-identity}. When $(\Delta H)/H\in L^1_{\mathrm{loc}}(\Hh^d)$,
this limit is the Lebesgue integral.
In codimension one, choose a local signed normal coordinate
$\rho$ and write $H(x)=H_0(|\rho(x)|)$ for a $C^2$ half-line
function $H_0$. In distributions on $\R$,
$\frac{\mathrm d^2}{\mathrm d\rho^2}H_0(|\rho|)
=H_0''(|\rho|)+2H_0'(0+)\delta_0$,
where $\delta_0$ is unit mass at $\rho=0$. The weak inequality is
determined by the regular part and this surface contribution,
which vanishes when $H_0'(0)=0$.
\end{proof}

\begin{lemma}\label{lem:endpoint-cutoff}
Let $A\in L^\infty(\Hh^d;\R)$ and $\nu\in\R$.
Suppose $H:\Hh^d\to(0,\infty)$ is locally Lipschitz and is $C^2$ on
$\Hh^d\setminus\Sigma$, where $\Sigma$ is a point or a smooth submanifold
of codimension at least two, as in Lemma~\ref{lem:ground-transform}.
Assume $(D\Delta+A)H=(\nu-\mathcal E)H$ weakly, with
$0\le\mathcal E\in L^1_{\mathrm{loc}}(\Hh^d)$.
If $v\in L^2(\Hh^d)$ is an eigenfunction of $\mathscr L_A$ at $\nu$,
then, for every smooth compactly supported cutoff $\chi$,
\[\label{eq:endpoint-cutoff}
 D\int_{\Hh^d}H^2\left|\nabla\frac{\chi v}{H}\right|_g^2\dd V
 +\int_{\Hh^d}\mathcal E\chi^2|v|^2\dd V
 =D\int_{\Hh^d}|\nabla\chi|_g^2|v|^2\dd V
 \quad\text{for }\chi\in C_c^\infty(\Hh^d;\R).
\]
Consequently $v/H$ is constant, and $\mathcal E|v|^2=0$ almost everywhere.
In particular, a nonzero eigenfunction at $\nu$ is excluded if
$\mathcal E>0$ almost everywhere or if $H\notin L^2$.
\end{lemma}

\begin{proof}
For real $v$, the operator-domain assumption permits testing the
weak eigenvalue equation with $\chi^2v$, giving
\[
 \int_{\Hh^d}\bigl(A|\chi v|^2-D|\nabla(\chi v)|_g^2\bigr)\dd V
 =\nu\int_{\Hh^d}|\chi v|^2\dd V
       -D\int_{\Hh^d}|\nabla\chi|_g^2|v|^2\dd V.
\]
Subtracting \eqref{eq:ground-identity}, applied to $\chi v$,
from this equality yields \eqref{eq:endpoint-cutoff}.
To justify the application, iterate local elliptic estimates from
$v\in L^2$: boundedness of $A$ gives
$v\in W^{2,p}_{\mathrm{loc}}$ for every finite $p$, and hence
local boundedness \cite[Chapters~8--9]{GilbargTrudinger}.
Compact approximation therefore applies, while $\mathcal E\in L^1_{\mathrm{loc}}$
ensures local integrability of $\mathcal E|v|^2$. The excision argument in Lemma~\ref{lem:ground-transform}
justifies passage across $\Sigma$.
Applying this real calculation to the real and imaginary parts
proves \eqref{eq:endpoint-cutoff} for complex eigenfunctions as well.

Choose $\chi_R$ equal to one on $B(o,R)$ and zero outside
$B(o,2R)$, with $0\le\chi_R\le1$ and
$|\nabla\chi_R|_g\le C/R$. By \eqref{eq:endpoint-cutoff}, the
sum of its two nonnegative left-hand integrals is at most
$CR^{-2}\norm v_{L^{2}(\Hh^d)}^2$. On every fixed compact set,
$\chi_R=1$ for all sufficiently large $R$, so passage to the limit
gives $\nabla(v/H)=0$ and $\mathcal E|v|^2=0$ almost everywhere there.
Connectedness then yields $v=cH$ on $\Hh^d$. If $\mathcal E>0$ almost
everywhere, the second identity forces $v=0$; if
$H\notin L^2$, the identity $v=cH\in L^2$ forces $c=0$.
\end{proof}

Applying Lemmas~\ref{lem:ground-transform} and
\ref{lem:endpoint-cutoff} with $A=0$, and dividing the quadratic
identity by $D$, gives the geometric spectral normalisation
of the unperturbed hyperbolic Laplacian. Set
$H_b(x)=e^{-\beta b(x)/2}$. Since $|\nabla b|_g=1$ and
$\Delta b=\beta$, we have $\Delta H_b=-\rho_dH_b$, where
$\rho_d=\beta^2/4$. For every real $\zeta\in C_c^\infty(\Hh^d)$,
\[\label{eq:hyperbolic-spectral-normalization}
 \int_{\Hh^d}\bigl(|\nabla\zeta|_g^2-\rho_d\zeta^2\bigr)\dd V
 =\int_{\Hh^d}H_b^2
             \left|\nabla\left(\frac\zeta{H_b}\right)\right|_g^2\dd V\quad {\text{for }\zeta\in C_c^\infty(\Hh^d)}.
\]
The constant $\rho_d$ in~\eqref{eq:hyperbolic-spectral-normalization} is sharp. Indeed, in the horospherical chart choose nonzero
$\chi\in C_c^\infty(0,1)$ and
$\eta\in C_c^\infty(\R^{d-1})$, and put
$\zeta_N=H_b\chi((z-R_N)/N)\eta(y)$, with $R_N\to\infty$.
The quotient of the right side of~\eqref{eq:hyperbolic-spectral-normalization} by $\|\zeta_N\|_{L^2(\Hh^d)}^2$
is at most $C(N^{-2}+e^{-2R_N})$.
Since $H_b^2\dd V=\dd z\dd y$ has infinite mass,
Lemma~\ref{lem:endpoint-cutoff}, applied at $\nu=-D\rho_d$,
excludes an $L^2$ eigenfunction of $D\Delta$ at $-D\rho_d$.
Thus $\sup\sigma(\Delta)=-\rho_d$, with no eigenfunction
at that endpoint, and $\sup\sigma(D\Delta)=-D\rho_d$.
The geometric normalisation agrees with that in~\cite{DaviesMandouvalos,Helgason}.

\paragraph{Feynman--Kac formulas.}

Let $X_t$ be Brownian motion on $\Hh^d$ with generator
$D\Delta$. This is the time change by a factor $2D$ of the convention
with generator $\frac12\Delta$ in~\cite{Hsu}. For bounded real measurable
$A$ and bounded measurable $g$,
\[\label{eq:FK}
 e^{t(D\Delta+A)}g(x)
 =\mathbb E_x\left[g(X_t)
             \exp\left(\int_0^tA(X_\tau)\dd\tau\right)\right]\quad {\text{for }x\in\Hh^d,\quad t>0}.
\]
For smooth $A$ and $g$, apply It\^o's formula to
$w(X_\tau,t-\tau)\exp(\int_0^\tau A(X_r)\dd r)$, where
$w(x,t)=(e^{t(D\Delta+A)}g)(x)$ solves
$w_t=D\Delta w+A(x)w$ with $w(x,0)=g(x)$. Stop when $X$ exits $B(o,R)$ and let $R\to\infty$. Boundedness of $A$ and
$g$, and nonexplosion, justify this limit. For bounded measurable data, the identity follows from the
heat-kernel integral equation or its absolutely convergent iterated series. Nonnegative unbounded data are admitted
only by monotone truncation, allowing an infinite value unless a norm
estimate gives finiteness. For signed $L^2$ data, apply the formula
to the positive and negative parts and use
Lemma~\ref{lem:smoothing} at positive time. The Brownian motion throughout has generator $D\Delta$;
see~\cite{Hsu,Pinsky} for these semigroup and stochastic formulas.

For the solution $u(x,t;u_0,f)$ of~\eqref{eq:KPP}, put
$A_f(x)=f_u(x,0)$ for $x\in\Hh^d$ and define $N(x,v)$ by
\[\label{eq:nonlinear-killing}
 \left\{\begin{aligned}
 N(x,v)&=A_f(x)-\frac{f(x,v)}v
       &&\text{for }x\in\Hh^d,\ 0<v\le1,\\
 N(x,0)&=0&&\text{for }x\in\Hh^d,
 \end{aligned}\right.
\]
For $f(x,v)=F(z(x),v)$, the definition gives
$A_f(x)=a(z(x))$ for $x\in\Hh^d$. Taylor's formula and the KPP upper inequality give
$0\le N(x,v)\le C_2v$. For $f(x,v)=F(z(x),v)$, we have $vN(x,v)=R_F(z(x),v)$.
The nonlinear Feynman--Kac formula is
\[\label{eq:nonlinear-FK}
 u(x,t)=\mathbb E_x\left[u_0(X_t)
  \exp\left(\int_0^t
   [A_f(X_\tau)-N(X_\tau,u(X_\tau,t-\tau))]\dd\tau\right)\right] \quad {\text{for }x\in\Hh^d,\quad t>0}.
\]
This identity also applies to measurable initial data. Once $u$ is
fixed, its coefficient $c(x,t)=A_f(x)-N(x,u(x,t))$ is bounded and
measurable. Iterating \eqref{eq:duhamel} with
$f(x,u(x,t))=c(x,t)u(x,t)$ gives integrals over
$\{(\tau_1,\ldots,\tau_n)\in\R^n:0<\tau_1<\cdots<\tau_n<t\}$.
The $n$th term has absolute value at most
$\|u_0\|_\infty\|c\|_\infty^nt^n/n!$, where $\|c\|_\infty$ is the
essential supremum on $\Hh^d\times(0,t)$.
The Markov property identifies their sum with the exponential in~\eqref{eq:nonlinear-FK}. The argument $u(X_\tau,t-\tau)$ in~\eqref{eq:nonlinear-FK}
comes from evaluating the coefficient along the backward space--time path.
The next lemma gives the stopped Feynman--Kac formula for a
moving domain, including its initial and lateral boundary terms.

\begin{lemma}\label{lem:moving-domain-FK}
Let $T>0$, let $b_0\in C^1([0,T])$, and put
$\Omega_t=\{x\in\Hh^d:b(x)>b_0(t)\}$. Suppose a bounded classical function $v$
solves $v_t=D\Delta v+A(x,t)v$ for $x\in \Omega_t$, $0<t\le T$, with
bounded continuous data $v_0$ and $g_0$ satisfying
$v(x,0)=v_0(x)$ for $x\in\overline{\Omega_0}$ and
$v(x,t)=g_0(x,t)$ for $x\in\partial \Omega_t$, $0\le t\le T$. Assume $g_0(x,0)=v_0(x)$ on $\{x\in\Hh^d:b(x)=b_0(0)\}$ and
that $v$ is continuous
up to the parabolic boundary; let $A$ be bounded and continuous.
For $x\in \Omega_t$, let $X$ have generator $D\Delta$, start at $x$, and set
\[\label{eq:backward-exit}
 \tau=\inf\{r\in[0,t]:b(X_r)\le b_0(t-r)\}\wedge t,
 \quad\text{for }0<t\le T,\ X_0=x\in \Omega_t,
 \qquad \inf\varnothing=+\infty.
\]
Then
\[\label{eq:moving-FK-boundary}
 \begin{aligned}
 v(x,t)={}&\mathbb E_x\bigg[
  e^{\int_0^\tau A(X_r,t-r)\dd r}
  \Big(\one_{\{\tau=t\}}v_0(X_t)\\
 &\qquad\qquad+\one_{\{\tau<t\}}g_0(X_\tau,t-\tau)\Big)\bigg]
 \quad {\text{for }x\in \Omega_t,\quad0<t\le T}.
 \end{aligned}
\]
The same formula holds for $\Omega_t=\{x\in\Hh^d:\ell(x)>b_0(t)\}$,
with $b$ replaced by $\ell$ in~\eqref{eq:backward-exit}. For zero lateral data only the term with
$\tau=t$ remains.
\end{lemma}

\begin{proof}
Let $h(r)=b(X_r)-b_0(t-r)$ and, for sufficiently large $n$, define
$\tau_n=(t-1/n)\wedge\inf\{r\in[0,t]:h(r)\le1/n
\text{ or }\dist(o,X_r)\ge n\}$.
By continuity of the paths and nonexplosion,
$\tau_n\uparrow\tau$ almost surely, with $\tau$ given by
\eqref{eq:backward-exit}. Before $\tau_n$, the space--time path
remains in a compact subset of the classical equation's domain,
so It\^o's formula applies to
\[
 e^{\int_0^r A(X_q,t-q)\dd q}v(X_r,t-r)
 \quad\text{for }0\le r\le\tau_n.
\]
The equation $D\Delta v-v_t+Av=0$ cancels its drift, and the
stopped stochastic integral has expectation zero. Letting
$n\to\infty$, continuity of the paths and of $v$ up to the
parabolic boundary gives $v_0(X_t)$ on $\{\tau=t\}$ and
$g_0(X_\tau,t-\tau)$ on $\{\tau<t\}$, as in~\eqref{eq:moving-FK-boundary}. This passage is justified by
dominated convergence: the exponential is bounded by
$e^{T\norm A_{L^{\infty}(\Hh^d)}}$, where
$\|A\|_\infty=\sup_{0<s\le T,\,x\in \Omega_s}|A(x,s)|$, and the solution
and its boundary data are bounded. If exit occurs exactly at time $t$, compatibility identifies
the initial and lateral values. Replacing $b$ by the smooth
coordinate $\ell$ gives the same stopping argument and proves the
formula for $\Omega_t=\{x\in\Hh^d:\ell(x)>b_0(t)\}$; see~\cite{Hsu,Pinsky,FitzsimmonsPitman}
for the stopped Feynman--Kac formula.
\end{proof}

Let $p^{\mathrm{Dir}}(t,x,y)$ denote the kernel, relative to $\dd y$, of the
heat equation $w_t=Dw_{xx}$ on $(0,\infty)$ with $w(0,t)=0$;
the superscript indicates the Dirichlet boundary condition. Reflection gives
\[\label{eq:reflection-kernel}
 p^{\mathrm{Dir}}(t,x,y)=\frac{e^{-(x-y)^2/(4Dt)}}{\sqrt{4\pi Dt}}
                  \left(1-e^{-xy/(Dt)}\right),\quad\text{for }t,x,y>0.
\]
The difference of the Gaussian kernels centred at $y$ and $-y$
solves the heat equation on $x>0$, has zero boundary value and
initial trace $\delta_y$, proving \eqref{eq:reflection-kernel}
by uniqueness. Since
$(1-e^{-1})\min\{1,r\}\le1-e^{-r}\le\min\{1,r\}$ for $r\ge0$,
\[\label{eq:reflection-endpoints}
 p^{\mathrm{Dir}}(t,x,y)\asymp (Dt)^{-1/2}e^{-(x-y)^2/(4Dt)}
                                  \min\{1,xy/(Dt)\}
 \quad {\text{for }t>0,\quad x,y>0}.
\]
For fixed positive $x,y$, \eqref{eq:reflection-endpoints} gives
$p^{\mathrm{Dir}}(t,x,y)\asymp xy\,(Dt)^{-3/2}$ as $t\to\infty$.

\subsection{Spectral thresholds}

We first establish the properties of $k(\lambda)$ and then
identify the spectral thresholds by comparing quadratic forms.

\begin{lemma}\label{lem:cell}
For every $\lambda\in\R$, the operator $\K_\lambda$ in~\eqref{eq:cell} has a simple principal eigenvalue $k(\lambda)$ and a
strictly positive $L$-periodic eigenfunction $\phi_\lambda$. The functions
$k$ and $\phi_\lambda$ are real analytic in $\lambda$, with
$\phi_\lambda$ normalised as specified after \eqref{eq:cell}.
The eigenvalue $k(\lambda)$ satisfies
\[\label{eq:cell-bounds}
 D(\lambda^2-\beta\lambda)+a_{\min}
 \le k(\lambda)\le
 D(\lambda^2-\beta\lambda)+a_{\max}
 \quad\text{for }\lambda\in\R,
\]
and $k''(\lambda)>0$ for $\lambda\in\R$.
All assertions in~\eqref{eq:symmetry-intro} hold.
The unique minimiser $\lam>0$ is characterised by
\[\label{eq:critical-cell}
 \lam k'(\lam)=k(\lam).
\]
\end{lemma}

\begin{proof}
Use the fixed domain $H^2_{\rm per}(0,L)$ specified after
\eqref{eq:cell}.
For $\alpha>\max_{[0,L]}[D(\lambda^2-\beta\lambda)+a]$,
the maximum principle makes $(\alpha-\K_\lambda)^{-1}$ positive
on continuous periodic functions, and the strong maximum principle
makes its image strictly positive for every nonzero nonnegative
datum. Elliptic regularity and compact embedding make the resolvent
compact. The principal-eigenvalue theorem consequently supplies a
positive eigenfunction and algebraic simplicity.
Analytic perturbation on $H^2_{\rm per}(0,L)$ gives real analytic
dependence of $k(\lambda)$ and the normalised $\phi_\lambda$
on $\lambda$ \cite{Kato}. Evaluating
$\K_\lambda\phi_\lambda=k(\lambda)\phi_\lambda$ at a maximum
and a minimum of $\phi_\lambda$, using \eqref{eq:cell}, proves
\eqref{eq:cell-bounds}.

Periodic integration by parts in~\eqref{eq:cell} gives
$\K_\lambda^*=\K_{\beta-\lambda}$; since an operator and its
adjoint have the same principal eigenvalue,
$k(\lambda)=k(\beta-\lambda)$. To prove strict convexity, fix $\lambda\in\R$ and set
\[\label{eq:Hlambda-line}
 \begin{gathered}
 H_\lambda(z)=e^{-\lambda z}\phi_\lambda(z),\qquad
 D_\lambda(z)=D\left(\beta-2\lambda+
                  2\frac{\phi_\lambda'(z)}{\phi_\lambda(z)}\right),\\
 Q_\lambda=D\partial_{zz}+D_\lambda\partial_z
 \quad \text{for }z\in\R.
 \end{gathered}
\]
By \eqref{eq:cell}, the identity
$(D\partial_{zz}+D\beta\partial_z+a)H_\lambda=k(\lambda)H_\lambda$
is equivalent to
\[\label{eq:log-cell}
 D(\log H_\lambda)''+D((\log H_\lambda)')^2
       +D\beta(\log H_\lambda)'+a=k(\lambda)
 \quad \text{for }z\in\R.
\]
Put $v_\lambda=\partial_\lambda\log\phi_\lambda$.
Differentiating \eqref{eq:log-cell} once and twice with respect
to $\lambda$ gives
\[\label{eq:cell-derivatives}
 Q_\lambda(-z+v_\lambda)=k'(\lambda),\qquad
 Q_\lambda(\partial_\lambda v_\lambda)
        +2D(-1+v_\lambda')^2=k''(\lambda)\quad \text{for }z\in\R.
\]
The positive eigenfunction $\phi_{\beta-\lambda}$ of
$\K_\lambda^*=\K_{\beta-\lambda}$ gives the density
\[\label{eq:cell-invariant-density}
 \pi_\lambda(z)=
 \frac{\phi_\lambda(z)\phi_{\beta-\lambda}(z)}
      {\displaystyle\int_0^L
           \phi_\lambda(\zeta)\phi_{\beta-\lambda}(\zeta)\dd\zeta}\quad \text{for }z\in\R,
\]
where $\pi_\lambda(z)>0$ for $z\in\R$ and
$\int_0^L\pi_\lambda(z)\dd z=1$. Indeed, for periodic
$v\in C^2$, the identity
$Q_\lambda v=\phi_\lambda^{-1}
(\K_\lambda-k(\lambda))(\phi_\lambda v)$ and periodic
integration by parts give
$\int_0^LQ_\lambda v\,\pi_\lambda\dd z=0$.
Thus the density $\pi_\lambda$ defined by
\eqref{eq:cell-invariant-density} satisfies
$D\pi_\lambda''-(D_\lambda\pi_\lambda)'=0$ and is the invariant
probability density for $Q_\lambda$ on $\R/L\Z$; see~\cite{Pinsky}.
Since $Q_\lambda(-z+v_\lambda)=-D_\lambda+Q_\lambda v_\lambda$
and $v_\lambda$ is periodic, integration of~\eqref{eq:cell-derivatives} against $\pi_\lambda$ yields
\[\label{eq:cell-variance}
 k'(\lambda)=-\int_0^L D_\lambda\pi_\lambda\dd z,
 \qquad
 k''(\lambda)=2D\int_0^L(1-v_\lambda')^2\pi_\lambda\dd z>0\quad {\text{for }\lambda\in\R}.
\]
The integral
$\int_0^L(1-v_\lambda')^2\pi_\lambda\dd z$ in~\eqref{eq:cell-variance} is strictly positive, since its vanishing would give $v_\lambda'\equiv1$,
whereas periodicity gives $\int_0^Lv_\lambda'\dd z=0$.
At $\lambda=\lambda_0$, the normalisation of $\phiZ$ reduces
\eqref{eq:cell-invariant-density} to $\pi_{\lambda_0}=\phiZ^2$.

By the symmetry and strict convexity just proved, $k$ has its unique
minimum at $\beta/2$, and \eqref{eq:cell-bounds} gives
$k(0)\ge a_{\min}>0$. For $\lambda>0$, define
$G(\lambda)=\lambda k'(\lambda)-k(\lambda)$.
Equation \eqref{eq:cell-variance} implies
$G'(\lambda)=\lambda k''(\lambda)>0$, while
$G(0+)=-k(0)<0$. Since \eqref{eq:cell-bounds} also gives
$k(\lambda)/\lambda\to\infty$ as $\lambda\downarrow0$ and as $\lambda\to\infty$,
this quotient attains an interior minimum, where $G=0$.
Strict increase of $G$ proves uniqueness and
\eqref{eq:critical-cell}. Finally, $k$ is negative somewhere
exactly when $\mu<0$; when $\mu=0$, its unique zero is
$\beta/2>0$. Together with the attained minimum of the quotient,
the identities $\cp=\min_{\lambda>0}k(\lambda)/\lambda$ and
$\mu=k(\beta/2)$ prove the sign relation in~\eqref{eq:symmetry-intro}.
\end{proof}

Set $H_{\mathrm{sym}}(z)=H_{\lambda_0}(z)$ for $z\in\R$, with
$H_{\lambda_0}$ given by \eqref{eq:Hlambda-line}.
The following derivative bound holds when $\mu\le0$.

\begin{lemma}\label{lem:decreasing-ground}
Suppose $\mu\le0$. The function
\[\label{eq:bulk-ground}
 H_{\mathrm{sym}}(z)=e^{-\beta z/2}\phiZ(z)
 \quad {\text{for }z\in\R}
\]
satisfies $H_{\mathrm{sym}}'(z)<0$ for every $z\in\R$. There is $\eta>0$ such that
$H_{\mathrm{sym}}'/H_{\mathrm{sym}}\le-\eta$ on $\R$.
\end{lemma}

\begin{proof}
By \eqref{eq:bulk-ground} and periodicity of $\phiZ$, the logarithmic
derivative $g=H_{\mathrm{sym}}'/H_{\mathrm{sym}}$ is periodic
with mean $-\beta/2$. Substitution of the cell equation
\eqref{eq:cell} at $\lambda=\beta/2$ gives
\[\label{eq:riccati}
 g'=\frac{\mu-a}{D}-\beta g-g^2
 \quad {\text{for }z\in\R}.
\]
At a zero of $g$, \eqref{eq:riccati} yields $g'=(\mu-a)/D<0$,
so every zero would be a strict crossing from positive to negative.
Periodicity would then require a crossing in the reverse direction
within the same period, a contradiction. Since the mean is negative,
$g<0$ everywhere, and compactness of $[0,L]$ gives
$g\le-\eta$ for some $\eta>0$.
\end{proof}

The positive periodic eigenfunction $\phiZ$ of
$\K_{\lambda_0}=D\partial_{zz}+a-D\beta^2/4$, defined in~\eqref{eq:cell}--\eqref{eq:parameters}, gives the identity
\[\label{eq:periodic-ground-form}
 \int_\R\left[(a-D\beta^2/4)\zeta^2-D|\zeta'|^2\right]\dd z
 =\mu\int_\R\zeta^2\dd z
 -D\int_\R\phiZ^2\left|\left(\frac\zeta{\phiZ}\right)'\right|^2\dd z
 \quad\text{for }\zeta\in C_c^\infty(\R;\R).
\]
Consequently the upper spectral endpoint of
$D\partial_{zz}+a-D\beta^2/4$ on $L^2(\R;\dd z)$, with domain
$H^2(\R)$, is $\mu$: the upper bound is
immediate, and cutoffs of $\phiZ$ over intervals of length tending to
infinity give the lower bound.

For the changes of unknown below, $m$ is a positive $C^2$ function
on $(0,\infty)$. Multiplication by $\sqrt m$ maps
$L^2((0,\infty),m\dd s)$ unitarily onto $L^2(0,\infty;\dd s)$.
The conjugation identity \eqref{eq:liouville} is first used on compact intervals
contained in $(0,\infty)$.

\begin{lemma}\label{lem:liouville}
Let $m\in C^2((0,\infty))$ be positive, put $B=m'/m$, and set
$v=m^{-1/2}w$ for $w\in C^2((0,\infty))$. Then
$D v''+DBv'+a(s)v$ equals
\[\label{eq:liouville}
 m^{-1/2}\left[Dw''+
       \left(a-D\frac{B'}2-D\frac{B^2}4\right)w\right].
\]
For $B(r)=\beta\coth r$ and $B(s)=B_d(s)$, respectively,
the term $B'/2+B^2/4$ is
\[\label{eq:geometric-potentials}
 \begin{gathered}
 \frac{(\beta\coth r)'}2+\frac{\beta^2\coth^2r}4
   =\frac{\beta^2}4+
          \frac{\beta(\beta-2)}{4\sinh^2r}
 \quad {\text{for }r>0},\\
 \frac{B_d'}2+\frac{B_d^2}4
   =\frac{\beta^2}4+
      \frac{(d-2)(d-4)}{4\sinh^2s}+\frac1{4\cosh^2s}
 \quad {\text{for }s>0}.
 \end{gathered}
\]
The differences between the two expressions on the left of~\eqref{eq:geometric-potentials} and $\beta^2/4$ are $O(e^{-2r})$ as $r\to\infty$
in the spherical case and $O(e^{-2s})$ as $s\to\infty$ in the tube case.
\end{lemma}

\begin{proof}
Differentiating $v=m^{-1/2}w$ twice and using $m'/m=B$ cancels
the coefficient of $w'$ and leaves the zeroth-order diffusion term
$-D(B'/2+B^2/4)$, proving \eqref{eq:liouville}. Substituting the
spherical drift $\beta\coth r$ and the tube drift
\eqref{eq:tube-weight}, and using
$\coth^2s=1+\sinh^{-2}s$, $\tanh^2s=1-\cosh^{-2}s$ and
$\coth s\tanh s=1$, gives \eqref{eq:geometric-potentials}.
\end{proof}

For the transverse problem, the domain at zero is fixed by the
closed form in~\eqref{eq:sigma}; multiplication by $\sqrt{m_d}$ transports
this closed form to $L^2((0,\infty);\dd s)$, where the associated operator
has differential expression $w\mapsto Dw''+[a-D(B_d'/2+B_d^2/4)]w$ for $s>0$.
Tests supported away from zero are used in the estimates at infinity.

\begin{lemma}\label{lem:far-tests}
The operators $D\Delta+a(b(x))$ and $D\Delta+a(r(x))$ on
$L^2(\Hh^d)$ have
upper spectral endpoint at least $\mu$. The transverse endpoint
$\sigma_d$ is also at least $\mu$.
For the transverse self-adjoint operator determined by
\eqref{eq:sigma}, the upper endpoint of the essential spectrum is $\mu$.
\end{lemma}

\begin{proof}
In the horospherical coordinates of~\eqref{eq:horo-metric}, choose nonzero functions
$\chi\in C_c^\infty(0,1)$ and
$\eta\in C_c^\infty(\R^{d-1})$.
Let $(R_N)_{N\ge1}$ be a sequence of positive numbers tending to infinity. For $(z,y)\in\R\times\R^{d-1}$ define
\[\label{eq:horo-tests}
 v_N(z,y)=e^{-\beta z/2}\phiZ(z)
                    \chi\left(\frac{z-R_N}{N}\right)\eta(y).
\]
These compactly supported tests belong to the form domain and may
be approximated in $H^1$ by smooth compactly supported functions
with the same limiting Rayleigh quotients. Since the volume factor
$e^{\beta z}$ cancels the exponential square, applying
\eqref{eq:ground-identity} with $H_{\mathrm{sym}}(b)$ bounds
the difference between $\mu$ and the Rayleigh quotient of $v_N$
by $CN^{-2}+Ce^{-2R_N}$. Hence that quotient tends to $\mu$.

For the radial and transverse tube operators, apply
\eqref{eq:liouville} away from zero with
$w_N(s)=\phiZ(s)\chi((s-R_N)/N)$.
By \eqref{eq:periodic-ground-form}, the cutoff error is at most
$CN^{-2}$, and \eqref{eq:geometric-potentials} bounds the
geometric error by $Ce^{-2R_N}$. In the spherical case the functions
$x\mapsto m_{\rm rad}(r(x))^{-1/2}w_N(r(x))$, where
$m_{\rm rad}(r)=\sinh^\beta r$, are admissible on $\Hh^d$;
in the tube case $m_d^{-1/2}w_N$ is admissible in~\eqref{eq:sigma}. This proves all three
lower spectral bounds.

To determine the upper endpoint of the essential spectrum of the
transverse operator, cut the half-line at a fixed $S>0$ and impose a
Dirichlet condition there. On $(0,S)$, the ratio
$m_d(s)/s^{d-2}$ is smooth and positive up to zero, so the weighted
radial $H^1$ norm is equivalent to the radial Sobolev norm on a
Euclidean ball of dimension $d-1$. Rellich compactness gives compact
resolvent for this bounded part \cite[Chapter~7]{GilbargTrudinger}.
Removing the artificial Dirichlet condition at $S$ and restoring
the full half-line operator changes the resolvent by finite rank:
the difference of the resolvent solutions satisfies the homogeneous
second-order equation on each side and is determined by the
finitely many traces at $S$.

On $(S,\infty)$, \eqref{eq:geometric-potentials} gives the
potential $a(s)-D\beta^2/4+DE(s)$, where
$E(s)=-(d-2)(d-4)/(4\sinh^2s)-1/(4\cosh^2s)$ and
$|E(s)|\le Ce^{-2s}$ for $s\ge S$. Multiplication by $DE$ is relatively compact
with respect to the Dirichlet second derivative, by compact
$H^2$ embedding on bounded intervals and
$\sup_{s\ge A}|E(s)|\to0$ as $A\to\infty$.
The Dirichlet operator $D\partial_{ss}+a(s)-D\beta^2/4$ on
$L^2(S,\infty)$ has upper spectral bound $\mu$ by
\eqref{eq:periodic-ground-form}. For the reverse essential bound,
choose the supports of $w_N$ mutually disjoint and set
$\mathscr P=D\partial_{ss}+a-D\beta^2/4$.
With $\|\cdot\|_2$ denoting the norm on
$L^2((S,\infty);\dd s)$ in the next estimate, the equation for $\phiZ$ gives
\[\label{eq:tail-Weyl-residual}
 \begin{aligned}
 (\mathscr P-\mu)w_N
 &=\frac{2D}{N}\phiZ'\chi'\left(\frac{s-R_N}{N}\right)
  +\frac{D}{N^2}\phiZ\chi''\left(\frac{s-R_N}{N}\right)
       \quad\text{for }s>S,\\
 \frac{\|(\mathscr P-\mu)w_N\|_2}{\|w_N\|_2}&\le\frac CN
       \quad\text{for }N=1,2,\ldots.
 \end{aligned}
\]
In view of~\eqref{eq:tail-Weyl-residual}, the functions
$w_N/\|w_N\|_{L^2(S,\infty)}$ form a Weyl sequence at $\mu$;
their disjoint supports give weak convergence to zero. Relative compactness of the
perturbation $DE$ and the finite-rank change of the resolvent at $S$ therefore
preserve the upper endpoint of the essential spectrum \cite{ReedSimon}, completing
the proof.
\end{proof}

\begin{lemma}\label{lem:horo-spectrum}
The upper spectral endpoint of $D\Delta+a(b(x))$ is exactly $\mu$.
It is not an $L^2(\Hh^d)$ eigenvalue.
\end{lemma}

\begin{proof}
By \eqref{eq:horo-metric} and \eqref{eq:bulk-ground},
$(D\Delta+a(b))H_{\mathrm{sym}}(b)=\mu H_{\mathrm{sym}}(b)$.
Lemma~\ref{lem:ground-transform} therefore gives the upper bound
$\mu$, while Lemma~\ref{lem:far-tests} gives the matching lower
bound. If $v\in L^2(\Hh^d)$ were an eigenfunction at $\mu$,
Lemma~\ref{lem:endpoint-cutoff}, applied with $\mathcal E=0$, would yield
$\nabla(v/H_{\mathrm{sym}}(b))=0$, and hence
$v=cH_{\mathrm{sym}}(b)$. However,
\[\label{eq:horo-nonsquare}
 H_{\mathrm{sym}}(b)^2\dd V=\phiZ(z)^2\dd z\dd y
\]
has infinite total mass, since $\phiZ$ is positive and periodic.
Thus \eqref{eq:horo-nonsquare} and $v\in L^2$ force $c=0$,
which excludes a nonzero eigenfunction at $\mu$.
\end{proof}

\begin{lemma}\label{lem:radial-spectrum}
If $\mu\le0$, the upper spectral endpoint of $D\Delta+a(r(x))$ is $\mu$,
and it is not an $L^2$ eigenvalue.
\end{lemma}

\begin{proof}
By the radial Laplacian in~\eqref{eq:polar-metric} and
Lemma~\ref{lem:decreasing-ground}, for $r>0$ we have
\[\label{eq:radial-supersolution}
 (D\partial_{rr}+D\beta\coth r\,\partial_r+a(r))H_{\mathrm{sym}}(r)
 =\mu H_{\mathrm{sym}}(r)+D\beta(\coth r-1)H_{\mathrm{sym}}'(r)
 <\mu H_{\mathrm{sym}}(r)\quad {\text{for }r>0}.
\]
The function $H_{\mathrm{sym}}(r(x))$ is locally Lipschitz.
Excision of $B(o,\varepsilon)$ leaves a boundary term
$O(\varepsilon^{d-1})$, because $H_{\mathrm{sym}}'$ is bounded
near zero; this term vanishes for $d\ge2$. Thus
\eqref{eq:radial-supersolution} holds weakly on the full space.
Applying Lemma~\ref{lem:ground-transform} gives the spectral
upper bound $\mu$, and Lemma~\ref{lem:far-tests} gives equality.

If $(D\Delta+a(r(x)))v=\mu v$ with $v\in L^2(\Hh^d)$,
apply Lemma~\ref{lem:endpoint-cutoff} with
$H(x)=H_{\mathrm{sym}}(r(x))$. The second term on the left of~\eqref{eq:endpoint-cutoff} is the following integral with $\chi^2$
inserted in its integrand:
\[\label{eq:radial-strict-term}
 \int_{\Hh^d}
  \left[-D\beta(\coth r-1)\frac{H_{\mathrm{sym}}'(r)}{H_{\mathrm{sym}}(r)}\right]|v|^2\dd V.
\]
By Lemma~\ref{lem:decreasing-ground},
$-D\beta(\coth r-1)H_{\mathrm{sym}}'(r)/H_{\mathrm{sym}}(r)>0$
for $r>0$. Letting the cutoff radius tend to infinity
in~\eqref{eq:endpoint-cutoff} forces
\eqref{eq:radial-strict-term} to vanish, and hence $v=0$ almost
everywhere. Thus $\mu$ is not an $L^2(\Hh^d)$ eigenvalue.
\end{proof}

\begin{lemma}\label{lem:tube-spectrum}
The upper spectral endpoint of $D\Delta+a(s(x))$ on $L^2(\Hh^d)$ equals
$\sigma_d$. If $\sigma_d=0$, zero is not an $L^2$ eigenvalue. If $d\ge3$
and $\mu\le0$, then $\sigma_d=\mu$.
\end{lemma}

\begin{proof}
First take $v\in C_c^\infty(\Hh^d)$. In the Fermi coordinates
\eqref{eq:fermi-metric}, apply the transverse form inequality \eqref{eq:sigma} to the $s$ variable,
retain the longitudinal gradient term, and discard the nonpositive
angular contribution. This gives
\[\label{eq:tube-form-upper}
 \int_{\Hh^d}(a(s)|v|^2-D|\nabla v|_g^2)\dd V
 \le\sigma_d\norm v_{L^{2}(\Hh^d)}^2
       -D\int_{\Hh^d}\sech^2s\,|\partial_\ell v|^2\dd V\quad {\text{for }v\in C_c^\infty(\Hh^d)}.
\]
For $d=2$, the half-line inequality applies separately on the two
sides of $\G$, with unrestricted form trace at zero. For $d\ge3$,
for each fixed $(\ell,\omega)$, the function
$s\mapsto v(s,\ell,\omega)$ is approximated in the weighted form domain
by replacing its values on $[0,\varepsilon]$ with its value at
zero, interpolating smoothly on $[\varepsilon,2\varepsilon]$,
and leaving it unchanged for $s\ge2\varepsilon$.
The derivatives remain bounded, and
$m_d(s)\asymp s^{d-2}$ gives a squared $H^1$ error
$O(\varepsilon^{d-1})$. Smooth approximation away from zero
then permits application of~\eqref{eq:sigma} to the slice and
integration in the other variables, proving
\eqref{eq:tube-form-upper} for smooth compactly supported $v$.
Density of $C_c^\infty(\Hh^d)$ in $H^1(\Hh^d)$ extends
\eqref{eq:tube-form-upper} to every $v\in H^1(\Hh^d)$.
Conversely, choose $0\ne v\in\mathcal D_d$ from
\eqref{eq:sigma} and $0\ne\chi\in C_c^\infty(\R)$, and use
$x\mapsto v(s(x))\chi(\ell(x)/T)$ for $T>0$.
By \eqref{eq:fermi-metric}, the quotient of the longitudinal
Dirichlet integral of $x\mapsto v(s(x))\chi(\ell(x)/T)$ by its
squared $L^2(\Hh^d)$ norm is at most $CT^{-2}$, since
$\sech^2s\le1$.
Letting $T\to\infty$ and optimising over the transverse tests
proves equality of the full-space and transverse endpoints.

If $\sigma_d=0$ and an $L^2$ function satisfies
$(D\Delta+a(s))v=0$, \eqref{eq:tube-form-upper} implies
$\partial_\ell v=0$ almost everywhere. The volume element in~\eqref{eq:fermi-metric} is independent of $\ell\in\R$,
so an $L^2$ function with this property is zero. This proves the
assertion that $0$ is not a full-space eigenvalue, even when
$0$ is an eigenvalue of the transverse operator.

For $d\ge3$ and $\mu\le0$, Lemmas~\ref{lem:drifts} and
\ref{lem:decreasing-ground} give
\[\label{eq:tube-supersolution}
 (D\partial_{ss}+DB_d\partial_s+a)H_{\mathrm{sym}}
 =\mu H_{\mathrm{sym}}+D(B_d-\beta)H_{\mathrm{sym}}'\le\mu H_{\mathrm{sym}} \quad {\text{for }s>0,\quad d\ge3,\quad\mu\le0},
\]
The normal codimension of $\G$ is $d-1\ge2$; thus deleting
$s<\varepsilon$ produces a boundary term
$O(\varepsilon^{d-2})$ on each compact longitudinal interval.
By the excision argument of Lemma~\ref{lem:ground-transform},
this term vanishes, and $x\mapsto H_{\mathrm{sym}}(s(x))$
satisfies \eqref{eq:tube-supersolution} weakly on $\Hh^d$. The same lemma gives the upper spectral bound $\mu$ for $D\Delta+a(s(x))$, while
Lemma~\ref{lem:far-tests} gives the lower bound $\mu$. For $d=2$, the Neumann realisation is specified in~\eqref{eq:J-intro}.
\end{proof}

For $\mu>0$, translated positive Dirichlet eigenfunctions yield a
uniform lower bound for stationary solutions. We first quantify the
limits of the distance coordinates that enter this translation.

\begin{lemma}\label{lem:coordinate-limit}
Choose $z\in\{r,s\}$, and let $x_n\in\Hh^d$ satisfy
$z(x_n)\to\infty$ as $n\to\infty$. There are isometries $T_n$
with $T_no=x_n$ such that, after passage to a subsequence,
\[\label{eq:coordinate-limit}
 z\circ T_n-z(x_n)\to b_*
 \quad\text{in }C^2(K)\text{ as }n\to\infty
 \quad\text{for every compact }K\subset\Hh^d,
\]
where $b_*$ is a Busemann function with $b_*(o)=0$.
If $z(x_n)$ tends to $\theta$ modulo $L$, then
$a(z(T_nx))\to a(b_*(x)+\theta)$ locally uniformly for
$x\in\Hh^d$ as $n\to\infty$.
\end{lemma}

\begin{proof}
Write $d_n=r$ in the point-distance case and $d_n=s$ in the
geodesic-distance case, and put $R_n=d_n(x_n)\to\infty$.
By the transitivity of rotations on the unit tangent sphere
\cite[Chapter~5]{Ratcliffe}, we may choose $T_no=x_n$ so that
$\nabla(d_n\circ T_n)(o)$ is a fixed unit vector $e$.
Set $z_n=d_n\circ T_n-R_n$. On each fixed ball $B(o,M)$,
the one-Lipschitz property gives
$d_n(T_nx)\ge R_n-M$ and $|z_n(x)|\le M$.
Thus $T_n^{-1}\{o\}$ or $T_n^{-1}\G$, respectively, is
disjoint from $B(o,M)$ for $R_n>M$, and
\eqref{eq:coordinate-hessians} gives
\[\label{eq:distance-hessian-limit}
 \begin{gathered}
 |\nabla z_n|_g=1,\qquad
 \nabla^2z_n=g-\dd z_n\otimes\dd z_n+E_n
 \quad\text{for }x\in B(o,M),\ R_n\ge M+1,\\
 \sup_{B(o,M)}|E_n|_g\le C_Me^{-2R_n}
 \quad\text{for }R_n\ge M+1,
 \end{gathered}
\]
where $|E_n|_g$ is the pointwise norm of the covariant
two-tensor $E_n$ induced by $g$. The error bound follows
from $|\coth q-1|+|\tanh q-1|\le Ce^{-2q}$ for $q\ge1$.

Choose the Busemann function $b_*$ with $b_*(o)=0$ and
$\nabla b_*(o)=e$ by rotating the function of
Lemma~\ref{lem:coordinates}. Put $v_n=e^{z_n}$,
$v_*=e^{b_*}$, $w_n=v_n-v_*$ and $Z_n=\nabla w_n$.
Combining \eqref{eq:distance-hessian-limit} with the Busemann
Hessian in~\eqref{eq:coordinate-hessians} yields
\[\label{eq:exponential-hessian-error}
 \nabla^2w_n=w_ng+v_nE_n\quad\text{for }x\in B(o,M),\ R_n\ge M+1,
 \qquad w_n(o)=0,\quad Z_n(o)=0.
\]
Along a unit-speed geodesic $\gamma$ from $o$, of length at most
$M$, this identity becomes
\[\label{eq:geometric-error-ODE}
 \frac{\mathrm d}{\mathrm dt}w_n(\gamma(t))
      =g(Z_n,\dot\gamma),\qquad
 \nabla^g_{\dot\gamma}Z_n=w_n\dot\gamma+v_nE_n^\sharp\dot\gamma.
\]
Here $g(E_n^\sharp\dot\gamma,Y)=E_n(\dot\gamma,Y)$.
In a parallel orthonormal frame, \eqref{eq:geometric-error-ODE}
is a linear system with bounded coefficients and forcing
$O_M(e^{-2R_n})$, because $v_n\le e^M$ on the ball.
Gronwall's inequality therefore gives
$|w_n|+|Z_n|\le C_Me^{-2R_n}$ on $B(o,M)$, and
\eqref{eq:exponential-hessian-error} gives
$|\nabla^2w_n|_g\le C_Me^{-2R_n}$ there.
Since $v_n,v_*\ge e^{-M}$, the chain rule applied to their
logarithms now yields
\[\label{eq:coordinate-C2-rate}
 \|z_n-b_*\|_{C^2(B(o,M))}\le C_Me^{-2R_n},
 \quad\text{for }R_n\ge M+1.
\]
The $C^2$ norm is the supremum of the absolute value, gradient norm
and covariant Hessian norm. As $M$ is arbitrary,
\eqref{eq:coordinate-C2-rate} proves the claimed local convergence.

For $v,\theta\in\R$, define
$d_L(v,\theta)=\min_{j\in\Z}|v-\theta-jL|$.
By periodicity of $a$, its $C^1$ bound and
\eqref{eq:coordinate-C2-rate},
\[\label{eq:phase-potential-rate}
 \sup_{x\in B(o,M)}
 |a(d_n(T_nx))-a(b_*(x)+\theta)|
 \le \|a'\|_\infty\bigl(C_Me^{-2R_n}+d_L(R_n,\theta)\bigr)\quad {\text{for }R_n\ge M+1}.
\]
After passage to a subsequence, there is $\theta\in[0,L]$
such that $d_L(R_n,\theta)\to0$ as $n\to\infty$. Equation \eqref{eq:phase-potential-rate}
then gives $a(d_n(T_nx))\to a(b_*(x)+\theta)$ locally
uniformly for $x\in\Hh^d$.
\end{proof}

\begin{lemma}\label{lem:positive-boxes}
Assume $\mu>0$. For each $f\in\{f_b,f_o,f_{\G}\}$ in~\eqref{eq:patterns}, every positive stationary solution $q\le1$ of~\eqref{eq:KPP} satisfies $\inf_{x\in\Hh^d}q(x)>0$.
Consequently $q\equiv1$.
\end{lemma}

\begin{proof}
By Lemma~\ref{lem:horo-spectrum}, there is a compactly supported
test with strictly positive Rayleigh quotient for $D\Delta+a(b)$.
Fix a ball containing its support and choose $\delta>0$ smaller
than that quotient. The isometry $T_\theta$ in~\eqref{eq:horo-shift-isometry} translates $b$ by $\theta$ and
moves the origin by $|\theta|\le L$. Thus, for every
$\theta\in[0,L]$, the potential $a(b+\theta)$ has a test
with the same quotient in one common larger ball $B(o,R)$.
Its Dirichlet principal eigenvalue is consequently at least
$\delta$, uniformly in the phase.

For bounded real potentials $A_1,A_2$ on a fixed bounded domain
$\Omega$, let $\nu_\Omega(A_i)$ denote the principal Dirichlet
eigenvalue of $D\Delta+A_i$. The Rayleigh principle gives
\[\label{eq:Dirichlet-potential-Lipschitz}
 |\nu_\Omega(A_1)-\nu_\Omega(A_2)|
                 \le\|A_1-A_2\|_{L^\infty(\Omega)}
 \quad\text{for }A_1,A_2\in L^\infty(\Omega;\R),
\]
since the potential difference changes each normalised quotient
by at most its supremum norm. In~\eqref{eq:phase-potential-rate},
take $M=R$ and choose $\theta\equiv r(x)\pmod L$ or
$\theta\equiv s(x)\pmod L$, respectively. The resulting error is at most $C_Re^{-2r(x)}$
for point distance and $C_Re^{-2s(x)}$ for geodesic distance.
Together with \eqref{eq:Dirichlet-potential-Lipschitz}, this gives
a lower bound $\delta/2$ for the principal Dirichlet
eigenvalue on $B(x,R)$ whenever the corresponding distance $r(x)$ or $s(x)$ is at least
$\max\{R+1,\frac12\log(2C_R/\delta)\}$, where $C_R>0$ bounds the coefficient of $e^{-2r(x)}$ or $e^{-2s(x)}$
in~\eqref{eq:phase-potential-rate} with $M=R$ and includes the factor
$\|a'\|_\infty$.
Indeed, \eqref{eq:phase-potential-rate} then applies on $B(o,R)$ and the
potential error is at most $\delta/2$.
For the Busemann reaction, the principal Dirichlet eigenvalue on $B(x,R)$ is at least $\delta/2$ for every centre $x$: compose \eqref{eq:horo-shift-isometry} with a horizontal
translation to obtain $To=x$ and
$b(Ty)=b(y)+b(To)$, and then reduce $b(To)$ modulo $L$.

For $z=b,r$ or $s$, respectively, let $\nu_x$ be the principal
Dirichlet eigenvalue of $D\Delta+a(z(\cdot))$ on $B(x,R)$.
On each ball for which $\nu_x\ge\delta/2$, normalise its positive
eigenfunction $\varphi_x$ by $\max_{\overline{B(x,R)}}\varphi_x=1$.
The Rayleigh principle gives $\nu_x\le a_{\max}$.
The eigenfunction satisfies
$D\Delta\varphi_x(y)=[\nu_x-a(z(y))]\varphi_x(y)$
on $B(x,R)$. Pulling these functions back to $B(o,R)$, global Dirichlet
$W^{2,p}$ estimates for $p>d$ give a common $C^{1,\alpha}$
bound for some $\alpha>0$.
The uniform $C^{1,\alpha}$ bound, the normalisation
$\max\varphi_x=1$ and the zero boundary trace keep every
maximising point a fixed positive distance from $\partial B(x,R)$. A finite interior
Harnack chain then gives $\varphi_x(x)\ge\kappa>0$, uniformly
in the centre; the relevant estimates on these isometric copies
of one smooth ball are given in~\cite[Chapters~8--9]{GilbargTrudinger}.

On each ball with $\nu_x\ge\delta/2$, choose a sufficiently small
multiple $\alpha(0)\varphi_x$ below $q$ and let
$\alpha'=\delta\alpha/4$ until $\alpha$ reaches
\[\label{eq:uniform-seed-size}
 \alpha_* =\min\left\{\frac12,
             \frac{\delta}{4(C_2+1)}\right\}.
\]
By \eqref{eq:uniform-taylor} and \eqref{eq:uniform-seed-size},
for $y\in B(x,R)$,
$(\partial_t-D\Delta_{\Hh^d})[\alpha(t)\varphi_x(y)]
-f(y,\alpha(t)\varphi_x(y))$ is at most
$(\delta/4-\nu_x+C_2\alpha_*)\alpha(t)\varphi_x(y)\le0$.
On $\partial B(x,R)$, $\alpha(t)\varphi_x=0\le q$.
Comparison with the stationary solution $q$ therefore gives
$q(x)\ge\alpha_*\kappa$, independently of the initial small
multiple.

For the radial pattern, positivity on the remaining compact ball
completes the uniform lower bound. For the tube pattern, fix the
region $\mathcal T_S=\{x\in\Hh^d:s(x)\le S\}$. Since $\sigma_d\ge\mu>0$,
a compactly supported test with positive Rayleigh quotient lies in
a sufficiently large ball about $\gamma(0)$, so its principal
Dirichlet eigenvalue is positive; enlarge the radius to exceed $S+1$.
Write $\nu>0$ for the eigenvalue on the enlarged ball and normalise
its positive eigenfunction $\varphi$ by $\max\varphi=1$.
Set $m_S=\min_{\overline{B(\gamma(0),S)}}\varphi>0$.
The comparison leading to \eqref{eq:uniform-seed-size}, with
$\delta=2\nu$, gives $q\ge\min\{1/2,\nu/[2(C_2+1)]\}\,m_S$
on this smaller ball. Translations along the geodesic preserve the
operator, so this bound holds on every ball of radius $S$ centred on
the geodesic. These balls cover $\mathcal T_S$ when taken closed. We have therefore proved $\inf q>0$ for all
three patterns.

Let $m=\inf_{\Hh^d}q$. If $m<1$, choose $q(x_n)\to m$
and isometries $T_no=x_n$, and put $q_n(x)=q(T_nx)$,
$f_n(x,v)=f(T_nx,v)$. Then
$D\Delta q_n+f_n(x,q_n)=0$. By
\eqref{eq:pullback-coefficient-control} and local elliptic
estimates, a diagonal subsequence satisfies
$q_n\to q_\infty$ in local $C^2$ and $f_n\to f_\infty$
locally uniformly on $\Hh^d\times[0,1]$.
Passing to the equation gives
$D\Delta q_\infty+f_\infty(x,q_\infty)=0$, with
$q_\infty\ge m$ and $q_\infty(o)=m$.
By \eqref{eq:F1} and periodicity,
$f_n(x,m)\ge\min_{\zeta\in[0,L]}F(\zeta,m)>0$
for $x\in\Hh^d$ and $n\ge1$. Locally uniform convergence therefore gives
$D\Delta q_\infty(o)=-f_\infty(o,m)<0$.
This contradicts the nonnegative Laplacian at a minimum of a
$C^2$ function, proving $m=1$ and hence $q\equiv1$.
\end{proof}

\begin{proof}[Proof of Theorem~\ref{thm:threshold}]
For $f=f_b$, Lemma~\ref{lem:horo-spectrum} identifies $\mu$ as
the upper spectral endpoint and shows that $\mu$ is not an eigenvalue; for $f=f_o$
and $\mu\le0$, the same conclusions follow from
Lemma~\ref{lem:radial-spectrum}. Linear comparison and
Lemma~\ref{lem:smoothing} therefore yield
\eqref{eq:threshold-extinction} and \eqref{eq:threshold-rate}.

When $\mu>0$, Lemma~\ref{lem:far-tests} gives a positive
full-space spectral endpoint. By Lemma~\ref{lem:minimal-state},
the equation with the chosen reaction $f_b$ or $f_o$ has a
smallest positive stationary solution, and
\eqref{eq:liminf-local} holds for every nonzero datum.
Lemma~\ref{lem:positive-boxes} identifies the smallest positive
stationary solution with the constant function $1$;
combining \eqref{eq:liminf-local} with $u\le1$ proves
\eqref{eq:threshold-persistence}.
\end{proof}

\begin{proof}[Proof of Theorem~\ref{thm:dimension}]
Lemma~\ref{lem:tube-spectrum} identifies the full-space spectral
endpoint with $\sigma_d$ and shows that $0$ is not an eigenvalue
when $\sigma_d=0$. Applying Lemma~\ref{lem:smoothing} and linear
comparison gives \eqref{eq:threshold-extinction} when $\sigma_d\le0$. For $\sigma_d<0$, the same comparison gives \eqref{eq:threshold-rate} with $\sigma_d$ in place of $\mu$.
For $\sigma_d>0$, Lemma~\ref{lem:minimal-state} supplies the
smallest positive stationary solution $p_*$ and the lower bound
\eqref{eq:minimal-persistence}. Every isometry preserving $\G$
preserves the reaction, so minimality makes $p_*$ invariant under
this group. Since the group acts transitively on each level of
$s$, it follows that $p_*$ depends only on $s(x)$.
When $\mu>0$, Lemma~\ref{lem:positive-boxes} gives $p_*=1$
and hence \eqref{eq:threshold-persistence}. Finally,
\eqref{eq:no-localization} follows directly from the last assertion
of Lemma~\ref{lem:tube-spectrum}.
\end{proof}

\begin{remark}
For measurable $q:[0,\infty)\to\R$ with
$\int_0^\infty|q(s)|^2m_d(s)\dd s>0$, \eqref{eq:fermi-metric} gives
$\int_{\Hh^d}|q(s(x))|^2\dd V=\infty$.
Likewise, for measurable $h:\R\to\R$ with
$\int_\R|h(z)|^2e^{\beta z}\dd z>0$, \eqref{eq:horo-metric} gives
$\int_{\Hh^d}|h(b(x))|^2\dd V=\infty$.
These conclusions also hold when $q$ or $h$ has compact support.
The $L^2(\Hh^d)$ hypothesis in Lemma~\ref{lem:smoothing} concerns
the functions $x\mapsto q(s(x))$ and $x\mapsto h(b(x))$ with respect to $\dd V$.
\end{remark}

\section{Fronts and spreading}\label{sec:fronts}

Periodic exponential solutions give the fronts and spreading bounds.
The transition-kernel estimates yield stability. We then prove the
periodic logarithmic positions and Theorem~\ref{thm:common} for
reactions with a common derivative at zero.

\subsection{Construction of fronts}

For measurable $v_0:\R\to[0,1]$, define $(S_tv_0)(z)=w(z,t)$,
where $w$ solves \eqref{eq:line} with $w(z,0)=v_0(z)$. The operators $S_t$ preserve order and satisfy
$S_t[v_0(\cdot+L)](z)=(S_tv_0)(z+L)$.
Put $K_F=\max_{[0,L]\times[0,1]}|F_u|$. Linear comparison gives, for
any two such data $v_0,w_0$ and $t>0$,
\[\label{eq:line-data-continuity}
 |(S_tv_0)(z)-(S_tw_0)(z)|
 \le \frac{e^{K_Ft}}{\sqrt{4\pi Dt}}
       \int_\R e^{-(z+D\beta t-y)^2/(4Dt)}|v_0(y)-w_0(y)|\dd y\quad {\text{for }z\in\R,\quad t>0}.
\]
By \eqref{eq:line-data-continuity}, if $v_{0,n}\to v_0$ almost
everywhere and $0\le v_{0,n},v_0\le1$, then $S_tv_{0,n}\to S_tv_0$
locally uniformly on $\R\times(0,\infty)$; interior parabolic estimates give local convergence of the spatial
derivatives through order two and the first time derivative on
$\R\times(0,\infty)$.

\begin{lemma}\label{lem:front-barriers}
Fix $c>\cp$. There are positive numbers $\lambda<\lambda'$ with
$\lambda'<2\lambda$, and positive periodic functions
$\phi_\lambda,\phi_{\lambda'}$, such that
\[\label{eq:barrier-roots}
 k(\lambda)=c\lambda,\qquad
 \delta:=c\lambda'-k(\lambda')>0.
\]
For all sufficiently large $M$, the functions
\[\label{eq:front-barriers}
 \begin{gathered}
 H(z,t)=e^{-\lambda(z-ct)}\phi_\lambda(z),\qquad
 G(z,t)=e^{-\lambda'(z-ct)}\phi_{\lambda'}(z),\\
 \overline W(z,t)=\min\{1,H(z,t)\},\qquad
 \underline W(z,t)=(H(z,t)-MG(z,t))_+
 \quad {\text{for }(z,t)\in\R^2}
 \end{gathered}
\]
satisfy $0\le\underline W\le\overline W\le1$, with
$\overline W$ a weak supersolution and $\underline W$ a weak
subsolution of~\eqref{eq:line}.
\end{lemma}

\begin{proof}
By Lemma~\ref{lem:cell} and $c>\cp$, the strictly convex function
$k(\lambda)-c\lambda$ is positive at zero and at large positive
$\lambda$, and negative somewhere between. Choose its smaller zero as
$\lambda$ and then $\lambda'$ between the two zeros with
$\lambda'<2\lambda$. These choices give \eqref{eq:barrier-roots}.

Put $\mathcal L=D\partial_{zz}+D\beta\partial_z+a(z)$.
By \eqref{eq:cell} and \eqref{eq:barrier-roots}, the functions
in~\eqref{eq:front-barriers} satisfy
$(\partial_t-\mathcal L)H=0$ and
$(\partial_t-\mathcal L)G=\delta G$.
By \eqref{eq:F1}, $H$ is a supersolution where $H\le1$.
The function $\overline W=\min\{1,H\}$ remains a weak
supersolution of~\eqref{eq:line}:
the spatial derivative jumps downward at regular joins, and smooth
concave approximation gives the general case.

On $\{H>MG\}$, the definitions in~\eqref{eq:front-barriers}
give $\xi=z-ct\ge(\log M-C)/(\lambda'-\lambda)$.
Consequently $H\le1$ there for large $M$, and the positive minima and finite maxima of $\phi_\lambda$ and
$\phi_{\lambda'}$ on $[0,L]$ yield
\[\label{eq:barrier-nonlinearity}
 C_2H^2\le M\delta G
 \quad\hbox{on }\{H>MG\},
\]
Indeed, $H^2/G\le C\exp(-(2\lambda-\lambda')\xi)$ on this
set, with $2\lambda-\lambda'>0$.
Using \eqref{eq:uniform-taylor} and
\eqref{eq:barrier-nonlinearity}, we obtain
\[
 (\partial_t-D\partial_{zz}-D\beta\partial_z)(H-MG)
       -F(z,H-MG)
 \le-M\delta G+C_2H^2\le0.
\]
For $\underline W=(H-MG)_+$, the spatial derivative jumps
upward at regular boundary points of $\{H>MG\}$; the positive-part
inequality gives the weak subsolution inequality also at the other
boundary points. The definitions in~\eqref{eq:front-barriers} then
give $0\le\underline W\le\overline W\le1$.
\end{proof}

\begin{lemma}\label{lem:supercritical-front}
For every $c>\cp$ with $c\ne0$, there is a front of speed $c$ in the
sense of~\eqref{eq:front-definition}, and its profile satisfies
$\Phi_\xi<0$.
\end{lemma}

\begin{proof}
For $s<t$, let $W_s$ solve \eqref{eq:line} on
$\R\times(s,\infty)$ with
$W_s(z,s)=\overline W(z,s)$ from \eqref{eq:front-barriers}.
Thus $W_s(z,t)=S_{t-s}\overline W(\cdot,s)(z)$, and comparison with $\underline W$ and $\overline W$ from
Lemma~\ref{lem:front-barriers} gives
$\underline W(z,t)\le W_s(z,t)\le\overline W(z,t)$.
For $s_1<s_2<t$, the supersolution property and comparison give
$W_{s_1}(\cdot,s_2)\le\overline W(\cdot,s_2)$ and then
$W_{s_1}(z,t)\le W_{s_2}(z,t)$. Thus the limit
\[\label{eq:entire-front-limit}
 W(z,t)=\lim_{s\to-\infty}S_{t-s}\overline W(\cdot,s)(z)
 \quad {\text{for }(z,t)\in\R^2}
\]
exists. On each fixed space--time cylinder, the functions $W_s$
solve \eqref{eq:line} for sufficiently negative $s$ and satisfy
$0\le W_s\le1$. Local parabolic compactness identifies the pointwise limit $W$
as a classical entire solution.
Since, for each fixed $t$, $\underline W(z,t)>0$ and
$\overline W(z,t)<1$ for all sufficiently large $z$, the strong
maximum principle, applied from any earlier time, gives $0<W<1$
on $\R^2$.

By \eqref{eq:front-barriers},
$\overline W(z+L,t+L/c)=\overline W(z,t)$.
Uniqueness and spatial periodicity therefore give
$W_{s+L/c}(z+L,t+L/c)=W_s(z,t)$.
For either sign of $c$, both starting times tend to $-\infty$
as $s\to-\infty$, so \eqref{eq:entire-front-limit} yields
\[\label{eq:front-period-relation}
 W(z+L,t+L/c)=W(z,t)
 \quad {\text{for }(z,t)\in\R^2,\quad c\ne0}.
\]
It follows that $\Phi(z,\xi)=W(z,(z-\xi)/c)$ is $L$-periodic
in $z$, and $\Phi(z,\xi)\le e^{-\lambda\xi}\phi_\lambda(z)$ gives
$\sup_{z\in[0,L]}\Phi(z,\xi)\to0$ as $\xi\to\infty$.

For $\tau>0$, a change of starting time in~\eqref{eq:entire-front-limit} gives
\[
 W(z,t+\tau)=\lim_{s\to-\infty}
                      S_{t-s}\overline W(\cdot,s+\tau)(z).
\]
By \eqref{eq:front-barriers}, $\overline W(z,s+\tau)$
is nondecreasing in $\tau$ for $c>0$ and nonincreasing for $c<0$.
Order preservation thus gives $cW_t\ge0$, or
$\Phi_\xi\le0$. By monotone convergence and local parabolic
compactness, the limits $\lim_{t\to-\infty}W(\cdot,t)$ and
$\lim_{t\to\infty}W(\cdot,t)$ are stationary solutions
of $DP''+D\beta P'+F(z,P)=0$, and
\eqref{eq:front-period-relation} makes them periodic.
A periodic stationary solution in $[0,1]$ is either zero or one:
the strong maximum principle makes a nonzero solution positive,
and \eqref{eq:F1} excludes a minimum strictly below one.
For sufficiently large fixed $\xi$,
$\underline W(z,(z-\xi)/c)>0$ for every $z\in[0,L]$.
Since $\Phi$ is nonincreasing in $\xi$, its limit as
$\xi\to-\infty$ is nonzero and hence equals $1$. Dini's theorem on $[0,L]$, applied to monotone sequences
of profiles, gives uniform convergence at both ends; monotonicity
extends these limits to real $\xi$.

Differentiating \eqref{eq:line} in time gives
$(W_t)_t=D(W_t)_{zz}+D\beta(W_t)_z+F_u(z,W)W_t$.
The coefficient $F_u(z,W)$ is bounded by \eqref{eq:F1}.
Since $cW_t\ge0$, the strong maximum principle gives either
$cW_t>0$ everywhere or $W_t\equiv0$. If $W_t\equiv0$, \eqref{eq:front-period-relation} would make
$W$ periodic and stationary, contradicting $\lim_{\xi\to-\infty}\Phi(z,\xi)=1$ and
$\lim_{\xi\to\infty}\Phi(z,\xi)=0$.
Thus $\Phi_\xi<0$.
\end{proof}

\begin{lemma}\label{lem:critical-front}
If $\cp\ne0$, there is a front $W(z,t)=\Phi(z,z-\cp t)$, with
$\xi\mapsto\Phi(z,\xi)$ strictly decreasing for every $z\in\R$.
\end{lemma}

\begin{proof}
Choose $c_n>\cp$ with $c_n\to\cp$ and $c_n$ of the same
sign as $\cp$. By Lemma~\ref{lem:supercritical-front}, there
are fronts $W_n(z,t)=\Phi_n(z,z-c_nt)$, which we normalise by
time translation so that $W_n(0,0)=1/2$.
Their common equation \eqref{eq:line} and bounds $0<W_n<1$
give local parabolic compactness. Passing to a subsequence yields
an entire solution $W$ with $W(0,0)=1/2$ and $\cp W_t\ge0$.
Since $\cp\ne0$, local uniform convergence and $L/c_n\to L/\cp$
allow passage in~\eqref{eq:front-period-relation}, giving
$W(z+L,t+L/\cp)=W(z,t)$.

As in Lemma~\ref{lem:supercritical-front}, the two monotone time
limits are periodic stationary solutions, hence zero or one.
The two limits are distinct, since equality would make $W$
stationary and periodic, contrary to $W(0,0)=1/2$. Since $\cp W_t\ge0$,
$W$ is nondecreasing in $t$ when $\cp>0$ and nonincreasing when $\cp<0$; hence
$\Phi(z,\xi)=W(z,(z-\xi)/\cp)$ has the required uniform
limits. Applying the strong maximum principle to $W_t$ as in
Lemma~\ref{lem:supercritical-front} gives $\Phi_\xi<0$.
\end{proof}

\begin{lemma}\label{lem:stationary-front}
A stationary front of~\eqref{eq:line} exists exactly when $\mu\le0$.
Every stationary front is strictly decreasing.
\end{lemma}

\begin{proof}
When $\mu<0$, the sign relation in~\eqref{eq:symmetry-intro}
gives $0>\cp$, so Lemma~\ref{lem:front-barriers} applies with
$c=0$. Let $w(z,t)$ solve \eqref{eq:line} with
$w(z,0)=\overline W(z,0)$ from \eqref{eq:front-barriers}.
By comparison, $w$ decreases in time and remains above
$\underline W(z,0)$. Its limit
$P(z)=\lim_{t\to\infty}w(z,t)$ therefore satisfies
$DP''+D\beta P'+F(z,P)=0$, $0<P<1$ and $P(+\infty)=0$.
Multiplication of this equation by $D^{-1}e^{\beta z}$ gives
\[\label{eq:stationary-flux}
 (e^{\beta z}P'(z))'=-D^{-1}e^{\beta z}F(z,P(z))<0
 \quad {\text{for }z\in\R}.
\]
If $P'(z_0)\ge0$ for some $z_0\in\R$, \eqref{eq:stationary-flux}
gives $z_1<z_0$ with $P'(z_1)>0$ and
$P'(z)\ge e^{\beta(z_1-z)}P'(z_1)$ for $z\le z_1$.
Integrating this inequality from $z$ to $z_1$ contradicts boundedness
of $P$ as $z\to-\infty$.
Hence $P'<0$. Local ODE bounds give bounded $P'$, so integration of~\eqref{eq:stationary-flux} from $-\infty$ yields
$P'(z)=-D^{-1}\int_{-\infty}^z e^{-\beta(z-y)}F(y,P(y))\dd y$.
If $\lim_{z\to-\infty}P(z)$ were in $(0,1)$, \eqref{eq:F1} and
periodicity would bound $F(y,P(y))$ below by a positive constant
on a left half-line. The integral would then force $P'$ to stay
below a negative constant, contrary to boundedness. Thus $P(-\infty)=1$.

For $\mu=0$, put $F_\varepsilon=(1-\varepsilon)F$,
$0<\varepsilon<1/2$. By the periodic Rayleigh principle, the upper spectral endpoint
of $D\partial_{zz}+(1-\varepsilon)a(z)-D\beta^2/4$ on
$L^2(0,L)$ with domain $H^2_{\rm per}(0,L)$ is at most
$-\varepsilon a_{\min}<0$. The preceding construction, applied
to \eqref{eq:line} with $F_\varepsilon$, gives a decreasing
stationary front $P_\varepsilon$.
Translate it by an integer multiple of $L$ so that its unique
crossing of $1/2$ lies in $[0,L]$; periodicity preserves the
reaction under these translations. The stationary ODE and
boundedness give uniform local derivative bounds, so a subsequence
converges in $C^2_{\rm loc}(\R)$ to a nonincreasing stationary
solution $P$ satisfying $P(z_*)=1/2$ for some $z_*\in[0,L]$. By \eqref{eq:stationary-flux} and periodic positivity in~\eqref{eq:F1},
a limit $P(-\infty)$ or $P(+\infty)$ in $(0,1)$ would force $P'$ below a negative constant
on a half-line. The crossing and monotonicity therefore give
$P(-\infty)=1$ and $P(+\infty)=0$; the flux identity makes the
decrease strict.

If $\mu>0$ and $P$ were a stationary front,
$x\mapsto P(b(x))$ would be a positive stationary solution of~\eqref{eq:KPP} with $f=f_b$, different from $1$, contradicting Lemma~\ref{lem:positive-boxes}.
For every stationary front, \eqref{eq:stationary-flux} gives
$P'<0$ by the same boundedness argument, completing the proof.
\end{proof}

To prove the speed bound, we take the Laplace transform in $\xi$
of the nonnegative profile $\Phi$ and use the periodic elliptic
operator $\K_\lambda$ in~\eqref{eq:cell}.

\begin{lemma}\label{lem:tail-average}
Let $c\ne0$ and let $W(z,t)=\Phi(z,z-ct)$ be any front satisfying
\eqref{eq:front-definition}. Set
\[\label{eq:average-profile}
 \bar\Phi(\xi)=\frac1L\int_0^L\Phi(z,\xi)\dd z,\qquad \beta_c=\beta+c/D
 \quad {\text{for }\xi\in\R}.
\]
Then $\beta_c>0$, $\bar\Phi'(\xi)<0$, and there are $R,C_-,C_+,\alpha>0$ such that
\[\label{eq:average-tail-bounds}
 C_-e^{-\beta_c\xi}\le \bar\Phi(\xi)\le C_+e^{-\alpha\xi}
 \quad\text{for }\xi\ge R.
\]
For some constants $C,R_H>0$,
\[\label{eq:phase-Harnack}
 \max_{z\in[0,L]}\Phi(z,\xi)\le C\bar\Phi(\xi-R_H)
 \quad\text{for }\xi\in\R.
\]
For integers $j,k\ge0$ with $1\le j+k\le2$, the derivatives
$\partial_z^j\partial_\xi^k\Phi$ satisfy
$|\partial_z^j\partial_\xi^k\Phi(z,\xi)|\le C\bar\Phi(\xi-R_H)$,
with possibly larger $C,R_H$.
\end{lemma}

\begin{proof}
Substituting $W(z,t)=\Phi(z,z-ct)$ into \eqref{eq:line} gives
\[\label{eq:profile-PDE}
 D(\partial_z+\partial_\xi)^2\Phi+
 D\beta\partial_z\Phi+(D\beta+c)\partial_\xi\Phi+F(z,\Phi)=0
 \quad {\text{for }(z,\xi)\in\R^2}.
\]
Dividing by $D$ and integrating over one period in $z$ eliminates
the periodic derivatives and yields
\[\label{eq:average-ODE}
 \bar\Phi''+\beta_c\bar\Phi'+G=0,
 \qquad G(\xi)=\frac1{DL}\int_0^L F(z,\Phi(z,\xi))\dd z>0
 \quad {\text{for }\xi\in\R}.
\]
Since $\bar\Phi(-\infty)=1$ and $\bar\Phi(+\infty)=0$, we have $\bar\Phi'<0$
somewhere. If $\beta_c\le0$, \eqref{eq:average-ODE} implies that
$e^{\beta_c\xi}\bar\Phi'$ is strictly decreasing, forcing $\bar\Phi'$ to remain
negative with magnitude bounded below toward $+\infty$.
This contradicts $\bar\Phi\ge0$, so $\beta_c>0$.
The identity $(e^{\beta_c\xi}\bar\Phi')'=-e^{\beta_c\xi}G<0$, applied toward
$-\infty$, excludes a nonnegative value of $\bar\Phi'$, proving
$\bar\Phi'<0$ on $\R$.

By the uniform right limit in~\eqref{eq:front-definition} and
the Taylor estimate \eqref{eq:uniform-taylor},
$G(\xi)\ge(a_{\min}/(2D))\bar\Phi(\xi)$ for sufficiently large $\xi$.
Put $Q=-\bar\Phi'>0$. Equation \eqref{eq:average-ODE} gives
$Q'+\beta_c Q=G$, and integration over the preceding unit interval gives
\[
 Q(\xi)\ge\int_{\xi-1}^{\xi}e^{-\beta_c(\xi-v)}G(v)\dd v
 \ge\frac{a_{\min}(1-e^{-\beta_c})}{2D\beta_c}\bar\Phi(\xi).
\]
It follows that $\bar\Phi$ has the upper exponential bound in~\eqref{eq:average-tail-bounds}, with
$\alpha=a_{\min}(1-e^{-\beta_c})/(2D\beta_c)$.
The equation $Q'+\beta_c Q=G>0$ gives
$Q(\xi)\ge e^{-\beta_c(\xi-R)}Q(R)$; integrating from $\xi$ to
infinity proves the lower bound in~\eqref{eq:average-tail-bounds}.

To prove \eqref{eq:phase-Harnack}, write \eqref{eq:line} as
$W_t=DW_{zz}+D\beta W_z+c_W(z,t)W$, where
$c_W(z,t)=F(z,W(z,t))/W(z,t)$ is bounded by \eqref{eq:F1}.
For $z,z'\in[0,L]$, fix a time lag $\delta>0$ and set
$R_c=\operatorname{sgn}(c)(2L+|c|\delta)$.
The later time $(z'-\xi+R_c)/c$ minus $(z-\xi)/c$ lies in
$[\delta+L/|c|,\delta+3L/|c|]$. A finite parabolic Harnack chain with
spatial length at most $L$ therefore gives
$\Phi(z,\xi)\le C\Phi(z',\xi-R_c)$.
Averaging in $z'$ and using $\bar\Phi'<0$ to replace $R_c$ by
$|R_c|$ proves \eqref{eq:phase-Harnack}.

Interior H\"older estimates for
$W_t=DW_{zz}+D\beta W_z+c_W(z,t)W$ bound
the H\"older norm of $W$ on a rectangle by its supremum on a
larger fixed rectangle. Since $F(z,0)=F_z(z,0)=0$ and $F_{zu}$
is bounded, \eqref{eq:F1} gives a H\"older bound of the same
size for $F(z,W)$. Schauder estimates then control $W_z$,
$W_{zz}$ and $W_t$ by that supremum.
Time difference quotients in~\eqref{eq:line} give
$(W_t)_t=D(W_t)_{zz}+D\beta(W_t)_z+F_u(z,W)W_t$;
a second interior estimate controls $W_{zt}$ and $W_{tt}$,
which enter the derivatives of $\Phi(z,\xi)=W(z,(z-\xi)/c)$.
The coefficients have uniform local H\"older bounds by
\eqref{eq:F1}, so transforming back and applying
\eqref{eq:phase-Harnack} with a larger fixed shift gives the
bounds $|\partial_z^j\partial_\xi^k\Phi|\le C\bar\Phi(\xi-R_H)$
for $1\le j+k\le2$; see~\cite{Friedman}.
\end{proof}

\begin{lemma}\label{lem:positive-transform-singularity}
Let $p_0:[R,\infty)\to[0,\infty)$ be measurable, and suppose the abscissa of convergence
\[
 \overline\lambda=\sup\left\{\lambda\in\R:
             \int_R^\infty e^{\lambda\xi}p_0(\xi)\dd\xi<\infty\right\}
\]
is finite. The function
$\lambda\mapsto\int_R^\infty e^{\lambda\xi}p_0(\xi)\dd\xi$,
defined for $\operatorname{Re}\lambda<\overline\lambda$,
has a nonremovable singularity at $\overline\lambda$.
\end{lemma}

\begin{proof}
After translating $\xi$, take $R=0$ and suppose that the transform
extends analytically to a neighbourhood of $\overline\lambda$.
Choose $\rho>0$ so that the extension is analytic on
$|\lambda-\overline\lambda|<\rho$. At $\lambda_1=\overline\lambda-\rho/4$
its Taylor series has radius at least $3\rho/4$ and hence converges at
$\lambda_2=\overline\lambda+\rho/4$. Within the original half-plane,
the $n$th Taylor coefficient is
$\int_0^\infty\xi^ne^{\lambda_1\xi}p_0(\xi)\dd\xi/n!$,
which is nonnegative. Tonelli's theorem therefore identifies the
sum at $\lambda_2$ with
$\int_0^\infty e^{\lambda_2\xi}p_0(\xi)\dd\xi$.
Convergence of the Taylor series makes this integral finite,
contradicting the definition of $\overline\lambda$.
\end{proof}

\begin{lemma}\label{lem:speed-necessity}
Every front \eqref{eq:front-definition} with $c\ne0$ satisfies $c\ge\cp$.
\end{lemma}

\begin{proof}
Let $\bar\Phi$ be defined by \eqref{eq:average-profile}, and
let $\overline\lambda$ be the abscissa of convergence of
$\int_R^\infty e^{\lambda\xi}\bar\Phi(\xi)\dd\xi$.
By Lemma~\ref{lem:tail-average},
$0<\alpha\le\overline\lambda\le \beta_c<\infty$.
For $\lambda\in\C$, define $\K_\lambda$ by the differential expression in~\eqref{eq:cell}, acting on $L^2(0,L;\dd z)$ with domain $H^2_{\rm per}(0,L)$.
For $\operatorname{Re}\lambda<\overline\lambda$, define
\[\label{eq:V-transform}
 V_\lambda(z)=\int_R^\infty e^{\lambda\xi}\Phi(z,\xi)\dd\xi\quad {\text{for }z\in\R,\quad\lambda\in\C,\quad\Ree\lambda<\overline\lambda}.
\]
By \eqref{eq:phase-Harnack}, $\lambda\mapsto V_\lambda$
is holomorphic on $\{\operatorname{Re}\lambda<\overline\lambda\}$
with values in $C_{\rm per}([0,L])$. The derivative bounds of Lemma~\ref{lem:tail-average}
justify two integrations by parts, and the bound
$0\le R_F(z,v)\le C_2v^2$ in~\eqref{eq:uniform-taylor} gives,
from \eqref{eq:profile-PDE},
\[\label{eq:transform-equation}
 (\K_\lambda-c\lambda)V_\lambda=R_\lambda\quad {\text{for }z\in\R,\quad\lambda\in\C,\quad\Ree\lambda<\overline\lambda},
\]
where
\[\label{eq:transform-R}
 \begin{aligned}
 R_\lambda(z)={}&\int_R^\infty e^{\lambda\xi}R_F(z,\Phi(z,\xi))\dd\xi\\
 &+e^{\lambda R}
 [D\Phi_\xi(z,R)+2D\Phi_z(z,R)+(D\beta+c-D\lambda)\Phi(z,R)].
 \end{aligned}
\]
The identity, initially obtained for real $\lambda$, extends
to the half-plane by holomorphy.
The estimates \eqref{eq:average-tail-bounds} and \eqref{eq:phase-Harnack}
give $\Phi^2\le Ce^{-\alpha\xi}\Phi$ for large $\xi$.
The nonlinear integral in~\eqref{eq:transform-R} is therefore
holomorphic for
$\operatorname{Re}\lambda<\overline\lambda+\alpha$,
while the second line of~\eqref{eq:transform-R}, evaluated at
$\xi=R$, is entire in $\lambda$.

Consider the analytic elliptic family
$A(\lambda)=\K_\lambda-c\lambda:
H^2_{\rm per}(0,L)\to L^2(0,L)$, of Fredholm index zero and
with compact resolvent. To obtain an invertible member, suppose
$A(\lambda)v=0$ for real positive $\lambda$ and set
$w(z)=e^{(\beta/2-\lambda)z}v(z)$.
By \eqref{eq:cell},
\[\label{eq:large-lambda-ODE}
 \left\{\begin{aligned}
 &Dw''+[a(z)-D\beta^2/4-c\lambda]w=0,\\
 &(w,w')(L)=e^{(\beta/2-\lambda)L}(w,w')(0).
 \end{aligned}\right.
\]
For $\lambda\ge1$, the first-order system for $(w,w'/\sqrt\lambda)$
has coefficient norm at most $C\sqrt\lambda$. Its period transfer matrix
maps $(w(0),w'(0)/\sqrt\lambda)$ to $(w(L),w'(L)/\sqrt\lambda)$.
Gronwall's inequality bounds this matrix and its inverse by
$e^{C\sqrt\lambda L}$, so each of its eigenvalues has modulus at least
$e^{-C\sqrt\lambda L}$.
The boundary multiplier $e^{(\beta/2-\lambda)L}$ in~\eqref{eq:large-lambda-ODE} is smaller than $e^{-C\sqrt\lambda L}$
for sufficiently large $\lambda$, a contradiction.
Fix an invertible member $A(\lambda_{\mathrm{inv}})$.
On $L^2(0,L)$ we have
\[\label{eq:analytic-Fredholm-factorization}
 A(\lambda)A(\lambda_{\mathrm{inv}})^{-1}
 =I+[A(\lambda)-A(\lambda_{\mathrm{inv}})]
                                  A(\lambda_{\mathrm{inv}})^{-1}.
\]
Here $I$ is the identity on $L^2(0,L)$. The operator
$[A(\lambda)-A(\lambda_{\mathrm{inv}})]A(\lambda_{\mathrm{inv}})^{-1}$
is analytic and compact: $A(\lambda_{\mathrm{inv}})^{-1}$ maps
$L^2(0,L)$ into $H^2_{\rm per}(0,L)$, the difference
$A(\lambda)-A(\lambda_{\mathrm{inv}})$ has order at most one,
and $H^1(0,L)\hookrightarrow L^2(0,L)$ is compact.
The analytic Fredholm theorem \cite{Kato,ReedSimon} applied to
\eqref{eq:analytic-Fredholm-factorization} thus gives $\lambda\mapsto A(\lambda)^{-1}$ as a
meromorphic family of bounded operators from $L^2(0,L)$ to
$H^2_{\rm per}(0,L)$.

It follows from \eqref{eq:transform-equation} that $V_\lambda$
has a meromorphic continuation near $\overline\lambda$.
Since $L^{-1}\int_0^L V_\lambda(z)\dd z
=\int_R^\infty e^{\lambda\xi}\bar\Phi(\xi)\dd\xi$,
Lemma~\ref{lem:positive-transform-singularity} forces a pole,
of some order $m\ge1$ for the meromorphic function $V_\lambda$.
Let $v$ be its coefficient of $(\overline\lambda-\lambda)^{-m}$.
Then $v=\lim_{\lambda\uparrow\overline\lambda}
(\overline\lambda-\lambda)^mV_\lambda$ in
$H^2_{\rm per}(0,L)$, hence uniformly on $[0,L]$.
For real $\lambda<\overline\lambda$, $V_\lambda(z)>0$,
so convergence in $C([0,L])$ gives $v\ge0$. Since the pole has order $m$, we have $v\not\equiv0$;
comparing the highest-order terms in~\eqref{eq:transform-equation} yields
\[
 (\K_{\overline\lambda}-c\overline\lambda)v=0.
\]
Elliptic regularity and the strong maximum principle give $v>0$.
By the principal-eigenvalue characterisation in
Lemma~\ref{lem:cell},
$k(\overline\lambda)=c\overline\lambda$.
The definition \eqref{eq:parameters} therefore implies
$c=k(\overline\lambda)/\overline\lambda\ge\cp$.
\end{proof}

\begin{proof}[Proof of Theorem~\ref{thm:fronts}]
By Lemmas~\ref{lem:supercritical-front} and
\ref{lem:critical-front}, every nonzero speed at or above $\cp$
is attained, while Lemma~\ref{lem:speed-necessity} gives the
converse restriction. Lemma~\ref{lem:stationary-front} treats
speed zero, and the sign identity in~\eqref{eq:symmetry-intro}
therefore gives the set of admissible speeds $[\cp,\infty)$.
Finally, by \eqref{eq:horo-metric}, every entire classical solution
$W$ of~\eqref{eq:line} gives a solution $u(x,t)=W(b(x),t)$ of~\eqref{eq:KPP} with $f=f_b$ and initial datum $u_0(x)=W(b(x),0)$.
Applying this identity to $W(z,t)=\Phi(z,z-ct)$ and to
$W(z,t)=P(z)$ proves the two assertions on $\Hh^d$.
\end{proof}

\subsection{Spreading with a periodic reaction}

For spreading, we construct compactly supported moving subsolutions
whose powers remain subsolutions under the small drift perturbations
in Lemma~\ref{lem:drifts}.

\begin{lemma}\label{lem:complex-root}
Let $K$ be a real analytic strictly convex function on $\R$ with
$K''>0$ and $K'(\R)=\R$. Let $I$ be a compact interval of real numbers
such that
\[\label{eq:positive-dispersion-gap}
 m(v)=\min_{\lambda\in\R}[K(\lambda)-v\lambda]>0
 \quad\text{for }v\in I.
\]
Write $\lambda_v$ for the unique solution of $K'(\lambda_v)=v$.
For every sufficiently small $\eta>0$, there are real numbers
$\alpha(v,\eta)=\lambda_v+O(\eta^2)$, uniformly for $v\in I$, such that
\[\label{eq:complex-dispersion}
 K(\alpha(v,\eta)+i\eta)-v(\alpha(v,\eta)+i\eta)
       =E(v,\eta)\in\R,
 \qquad \inf_{v\in I}E(v,\eta)>0.
\]
For the complex-variable argument, it suffices that $K$ extend
holomorphically to a neighbourhood of $\{\lambda_v:v\in I\}$.
\end{lemma}

\begin{proof}
Extend $K$ holomorphically to a neighbourhood of the compact real
interval traced by $\lambda_v$, and for real $\alpha,\eta$ put
\[
 A(\alpha,\eta)=
 \begin{cases}
 \operatorname{Im}K(\alpha+i\eta)/\eta,&\eta\ne0,\\
 K'(\alpha),&\eta=0.
 \end{cases}
\]
Taylor expansion gives $A(\alpha,\eta)=K'(\alpha)+O(\eta^2)$,
even in $\eta$. Since
$\partial_\alpha A(\lambda_v,0)=K''(\lambda_v)>0$, the implicit
function theorem gives $A(\alpha(v,\eta),\eta)=v$ and
$\alpha(v,\eta)=\lambda_v+O(\eta^2)$, uniformly on $I$.
The imaginary part in~\eqref{eq:complex-dispersion} is then zero,
and its real part $m(v)+O(\eta^2)$ is uniformly positive by
\eqref{eq:positive-dispersion-gap} and compactness.
\end{proof}

\begin{lemma}\label{lem:periodic-moving-bumps}
Let $I$ be a compact interval of velocities such that
$\min_{\lambda\in\R}[k(\lambda)-v\lambda]>0$ for $v\in I$.
There exist $A,M,\kappa,E_*>0$, positive numbers $E_v$
for $v\in I$, and nonnegative functions
$\varphi_v(z,t;R)$, for $v\in I$ and every centre $R\in\R$, with
\[\label{eq:bump-support}
 \begin{gathered}
 \supp\varphi_v(\cdot,t;R)\subset[R+vt-A,R+vt+A]
       \quad\text{for }t,R\in\R,\ v\in I,\\
 0\le\varphi_v(z,t;R)\le M
       \quad\text{for }z,t,R\in\R,\ v\in I,\\
 \varphi_v(z,t;R)\ge\kappa
       \quad\text{for }z,t,R\in\R,\ v\in I
       \text{ with }|z-R-vt|\le L,
 \end{gathered}
\]
and $E_v\ge E_*$ for every $v\in I$. For $z,t,R\in\R$ and $v\in I$
such that $\varphi_v(z,t;R)>0$, we have
\[\label{eq:bump-linear}
 (\partial_t-D\partial_{zz}-D\beta\partial_z-a(z))\varphi_v=-E_v\varphi_v.
\]
Extending $\varphi_v$ by zero outside its positive set gives
a weak subsolution of
$(\partial_t-D\partial_{zz}-D\beta\partial_z-a(z))\varphi_v
\le-E_v\varphi_v$.
There exist $p>1$, $\varepsilon_*>0$, and $\delta_*>0$ such that
$\delta\varphi_v^p$ is a weak subsolution of
\[\label{eq:perturbed-line}
 w_t=Dw_{zz}+DB(z)w_z+F(z,w)\quad\text{for }(z,t)\in\R^2
\]
whenever $0<\delta\le\delta_*$ and $|B(z)-\beta|\le\varepsilon_*$ on
its support. The constants $A,M,\kappa,E_*,p,\varepsilon_*,\delta_*$
can be chosen independently of $v\in I$, $R\in\R$ and $t\in\R$.
\end{lemma}

\begin{proof}
By the quadratic bounds and strict convexity in
Lemma~\ref{lem:cell}, $k'(\R)=\R$, so
Lemma~\ref{lem:complex-root} applies with $K=k$.
Continue the simple eigenpair to
$\zeta_v=\alpha(v,\eta)+i\eta$, using the holomorphic
normalisation $\int_0^L\phi_\zeta(z)^2\dd z=1$ and the
branch agreeing with the positive real eigenfunction.
The local eigenpairs agree on overlaps along
$\{\lambda_v:v\in I\}$ and hence extend to a common complex
neighbourhood of this interval. Positivity of the real eigenfunctions and analytic
perturbation give
$\phi_{\zeta_v}/\phi_{\lambda_v}=1+O(\eta)$ in
$C^2([0,L])$, uniformly in $v$; the $C^2$ estimate follows
from the $H^2_{\rm per}$ dependence and the eigenvalue ODE.
By the uniform estimate, choose $\eta_0>0$ such that
\eqref{eq:complex-dispersion} holds and
$\sup_{v\in I,\,z\in[0,L]}|\phi_{\zeta_v}(z)/\phi_{\lambda_v}(z)-1|\le1/2$
whenever $0<\eta\le\eta_0$. Choose its continuous complex argument to equal $0$ at $\eta=0$.
With this continuous choice,
$|\arg(\phi_{\zeta_v}(z)/\phi_{\lambda_v}(z))|\le C\eta$
uniformly for $z\in[0,L]$ and $v\in I$.
For
\[
 Z_v(z,t;R)=e^{-\zeta_v(z-vt-R)}\phi_{\zeta_v}(z)
\]
the cell equation and \eqref{eq:complex-dispersion} give
$(\partial_t-D\partial_{zz}-D\beta\partial_z-a)Z_v=-E_vZ_v$,
where $E_v=E(v,\eta)$ is real and satisfies $E_v\ge E_*>0$,
with $E_*$ independent of $v$.

Set $\varphi_v(z,t;R)=(\operatorname{Re}Z_v(z,t;R))_+$
for $|z-vt-R|<\pi/\eta$, and $\varphi_v(z,t;R)=0$ otherwise. Near $z-vt-R=\pm\pi/\eta$, the complex argument of $Z_v$
is $\mp\pi+O(\eta)$, respectively, so $\operatorname{Re}Z_v<0$ and the positive part already vanishes in a
neighbourhood of the cutoffs. The positive-part inequality therefore
gives a weak linear subsolution, with equality
\eqref{eq:bump-linear} on its positive set.
Fix $0<\eta\le\min\{\eta_0,\pi/[3(L+C)]\}$.
On $|z-vt-R|\le L$ the phase of $Z_v$ has absolute value at most
$(L+C)\eta$, so $\operatorname{Re}Z_v\ge |Z_v|/2$.
Positivity of $\phi_{\lambda_v}$ and the uniform perturbation estimate
give $\varphi_v\ge\kappa>0$ there, uniformly in $v\in I$. On $|z-vt-R|\le\pi/\eta$, compactness of $I$ and the bounds
for $\phi_{\zeta_v}$ give the upper bound in~\eqref{eq:bump-support}.

On $\{\varphi_v>0\}$, the chain rule and
\eqref{eq:uniform-taylor} bound
$w_t-Dw_{zz}-DB(z)w_z-F(z,w)$, for $w=\delta\varphi_v^p$, by
\[\label{eq:powered-bump-residual}
 \delta\{[(p-1)a-pE_v]\varphi_v^p
   -Dp(p-1)\varphi_v^{p-2}|\partial_z\varphi_v|^2
   -Dp(B-\beta)\varphi_v^{p-1}\partial_z\varphi_v
            +C_2\delta\varphi_v^{2p}\}.
\]
Young's inequality bounds the absolute value of the drift term by
\[
 \frac{Dp(p-1)}2\varphi_v^{p-2}|\partial_z\varphi_v|^2
   +\frac{Dp|B-\beta|^2}{2(p-1)}\varphi_v^p.
\]
Choose $p>1$ close enough to one that
$(p-1)a_{\max}\le E_*/4$, then choose $\varepsilon_*$ with
$Dp\varepsilon_*^2/[2(p-1)]\le E_*/4$, and finally choose
$\delta_*$ with $C_2\delta_*M^p\le E_*/4$ and
$\delta_*M^p\le1$.
With these choices, the expression displayed in~\eqref{eq:powered-bump-residual} is nonpositive. At a zero of $\varphi_v$, the first derivative of
$\varphi_v^p$ vanishes and its second derivative is locally
integrable, since $p>1$. Integration over the connected components of $\{\varphi_v>0\}$,
or smooth convex approximation, extends the inequality across
the boundary of this set. 
\end{proof}

\begin{lemma}\label{lem:plateau-amplification}
There is a function $g\in C^2([0,1])$ with $g(0)=g(1)=0$,
$g>0$ on $(0,1)$, and $g(v)\le F(z,v)$ for $z\in\R$ and $v\in[0,1]$.
For every $\eta\in(0,1)$ and $\varepsilon>0$, there are finite $T,A>0$
with the following property: if $u$ solves \eqref{eq:KPP} with $f$ one of the reactions in~\eqref{eq:patterns} and, for some $t_0\ge0$, satisfies $u(y,t_0)\ge\eta$
for $y\in B(x,A)$, then
\[\label{eq:plateau-amplification}
 u(x,t_0+T)\ge1-\varepsilon.
\]
The assertion also holds for the solution of~\eqref{eq:line},
with $B(x,A)$ replaced by an interval of radius $A$ in $\R$.
\end{lemma}

\begin{proof}
By \eqref{eq:F1} and periodicity,
$\min_{(z,v)\in[0,L]\times[\vartheta,1-\vartheta]}F(z,v)>0$
for $0<\vartheta<1/2$. Together with the uniform Taylor estimates at
zero and one, this gives $c_0,C>0$ and $0<\delta<1/2$ such that
\[\label{eq:plateau-envelope-bounds}
 \begin{aligned}
 \frac{F(z,v)}{v(1-v)}&\ge c_0
       &&\text{for }z\in\R,\ 0<v<1,\\
 \frac{F(z,v)}{v(1-v)}&\ge a_{\min}-Cv
       &&\text{for }z\in\R,\ 0<v\le\delta.
 \end{aligned}
\]
Decrease $\delta$ so that $C\delta/a_{\min}\le1/2$, and
choose $M\ge2C/a_{\min}$ with
$a_{\min}e^{-M\delta}\le c_0$.
Using $\log(1-x)\ge-2x$ for $0\le x\le1/2$ in~\eqref{eq:plateau-envelope-bounds} shows that
$g(v)=a_{\min}v(1-v)e^{-Mv}\le F(z,v)$ both for
$v\le\delta$ and $v\ge\delta$.
The function $g$ belongs to $C^2([0,1])$, satisfies
$g(0)=g(1)=0$ and is positive for $0<v<1$.

Let $y'=g(y)$, $y(0)=\eta$. Since $y(t)\to1$ as
$t\to\infty$, choose $T$ with $y(T)>1-\varepsilon/2$.
Let $U_A$ solve \eqref{eq:KPP} with $f(x,v)=g(v)$ and
$U_A(x,0)=\eta\one_{B(o,A)}(x)$.
With $K_g=\max_{[0,1]}|g'|$ and
$p_T(o,z)=p_T(\dist(o,z))$ from \eqref{eq:heat-estimate},
linear comparison gives
\[\label{eq:plateau-heat-tail}
 0\le y(T)-U_A(o,T)
 \le \eta e^{K_gT}
       \int_{\Hh^d\setminus B(o,A)}p_T(o,z)\dd V(z)
 \to0\quad\text{as }A\to\infty.
\]
Since $\int_{\Hh^d}p_T(o,z)\dd V(z)=1$,
\eqref{eq:plateau-heat-tail} gives
$U_A(o,T)\ge1-\varepsilon$ for large $A$. Isometry invariance and
comparison with $F$ prove \eqref{eq:plateau-amplification} at every
centre. The conservative, translation-invariant semigroup generated by
$D\partial_{zz}+D\beta\partial_z$ gives the interval version
by the same heat comparison, with reaction $g$.
\end{proof}

\begin{proof}[Proof of Theorem~\ref{thm:periodic-speed}]
For the horospherical pattern, choose
$C\ge e^{\lam B}/\min_{[0,L]}\phi_{\lam}$.
By \eqref{eq:Hlambda-line} and $k(\lam)=\cp\lam$,
$\min\{1,Ce^{-\lam(z-\cp t)}\phi_{\lam}(z)\}$ is a weak
supersolution of~\eqref{eq:line} whose value at $t=0$ dominates
$\one_{(-\infty,B]}$.
Composing it with $b$ and applying comparison gives
\[\label{eq:periodic-speed-exponential-upper}
 \sup_{\{x\in\Hh^d:b(x)\ge ct\}}u(x,t)
 \le C\max_{[0,L]}\phi_{\lam}\,
                  e^{-\lam(c-\cp)t}\quad\text{for }c>\cp,\ t\ge0.
\]
This proves the outer limit in~\eqref{eq:periodic-horo-speed}
whether $\cp$ is positive, zero or negative.

For the inner limit, fix $c<\cp$. Strict convexity and
$k(0)>0$ from Lemma~\ref{lem:cell} give $k'(0)<\cp$.
Choose $\max\{c,k'(0)\}<v<\cp$.
For $\lambda<0$, the tangent bound
$k(\lambda)\ge k(0)+k'(0)\lambda$ gives
$k(\lambda)-v\lambda>0$; for $\lambda>0$, the same inequality
follows from the definition of $\cp$, and at zero it is
$k(0)>0$. By \eqref{eq:cell-bounds}, the minimum over $\lambda$
is therefore strictly positive.
Lemma~\ref{lem:periodic-moving-bumps} supplies
$\delta\varphi_v(z,t;R)^p$, uniformly positive when
$|z-R-vt|\le L$. Choose $R$ so negative that its initial support
lies in $(-\infty,A]$, and decrease $\delta>0$ so its initial
value is at most $1$. The functions
$\delta\varphi_v(z-jL,t;R)^p$, $j\in\Z$, $j\le0$,
are weak subsolutions of~\eqref{eq:line} initially below
$\one_{(-\infty,A]}$, and the intervals
$[R+vt+jL-L,R+vt+jL+L]$, $j\le0$,
cover $(-\infty,R+vt+L]$.
Comparison therefore gives a fixed $\eta>0$ such that the
solution $w_A$ of~\eqref{eq:line} with
$w_A(z,0)=\one_{(-\infty,A]}(z)$ is at least $\eta$ for
$z\le vt-C$, uniformly in $t\ge0$.

For a centre $z\le ct$, the interval $[z-A_0,z+A_0]$, where $A_0$ is the radius in
Lemma~\ref{lem:plateau-amplification} for this $\eta$ and a fixed
error $\varepsilon>0$, lies in
$(-\infty,v(t-T)-C]$ once $t$ is sufficiently large, since
$c<v$. Lemma~\ref{lem:plateau-amplification} therefore gives uniform
convergence to one for $z\le ct$.
Let $w_B$ solve \eqref{eq:line} with
$w_B(z,0)=\one_{(-\infty,B]}(z)$. Comparison gives
$w_A(b(x),t)\le u(x,t)\le w_B(b(x),t)$, so the preceding
limits prove \eqref{eq:periodic-horo-speed} for $u_0$ between the
horoball indicators in Theorem~\ref{thm:periodic-speed}(i).

Suppose now $\mu>0$, so $\cp>0$, and fix
$0<c<c_1<c_2<\cp$.
For $v\in[0,c_2]$ and $\lambda\le0$,
$k(\lambda)-v\lambda\ge\mu>0$; for $\lambda>0$,
$k(\lambda)-v\lambda\ge(\cp-c_2)\lambda>0$.
These inequalities, continuity and the quadratic bounds in~\eqref{eq:cell-bounds} give
$\inf_{v\in[0,c_2]}\min_{\lambda\in\R}
[k(\lambda)-v\lambda]>0$.
Take $\delta\varphi_v^p$ from
Lemma~\ref{lem:periodic-moving-bumps} for $v\in[0,c_2]$.
By Lemma~\ref{lem:drifts}, their common initial centre $R$ can
be chosen so far out that the spherical or tube geometric drift coefficient differs
from $\beta$ by at most $\varepsilon_*$ on every support.
Their nonnegative velocities preserve this drift bound for all $t\ge0$.

For a nonzero compactly supported datum in the spherical case,
Lemma~\ref{lem:wellposed} gives a positive minimum of
$u(x,1;u_0,f_o)$ on the compact annulus containing all initial
supports. Decrease the common $\delta>0$ so that every
$\delta\varphi_v(\cdot,0;R)^p$, composed with $r(x)$, lies
below $u(x,1;u_0,f_o)$. Comparison then gives a uniform
positive lower bound at all radii in $[R,R+c_2(t-1)]$.
By Theorem~\ref{thm:threshold}, the solution converges to one
uniformly on the ball of radius $R$, so for all sufficiently large
$t$ it has a fixed positive lower bound throughout $r(x)\le c_1t$.

For the data in Theorem~\ref{thm:periodic-speed}(iii),
let $V$ solve the geodesic-distance problem in~\eqref{eq:three-reductions} with
$V(s,0)=\eta\one_{[0,R_1]}(s)$ and $V_s(0,t)=0$.
Comparison gives $V(s(x),t)\le u(x,t;u_0,f_{\G})$.
For every compact interval $J\subset[0,\infty)$,
$\min_{s\in J}V(s,1)>0$. Thus one common $\delta>0$ places
all $\delta\varphi_v(s,0;R)^p$, $v\in[0,c_2]$, below $V(s,1)$.
For $0\le s\le S$, representatives with $\ell=0$ and a fixed normal
direction in~\eqref{eq:fermi-embedding} form a compact curve.
Coordinate invariance and Theorem~\ref{thm:dimension} give
$\sup_{0\le s\le S}|V(s,t)-1|\to0$ as $t\to\infty$
for every fixed $S<\infty$.
Together with the moving comparisons, this gives a fixed positive
lower bound throughout $s(x)\le c_1t$ for large $t$.

In either geometry, the distance coordinate is one-Lipschitz.
Hence the ball of fixed radius supplied by
Lemma~\ref{lem:plateau-amplification}, centred at any point with
$r(x)\le ct$ or $s(x)\le ct$, lies in $\{r\le c_1(t-T)\}$ or $\{s\le c_1(t-T)\}$,
respectively, for large $t$. Lemma~\ref{lem:plateau-amplification}
therefore gives uniform convergence to $1$ on $\{r\le ct\}$
or $\{s\le ct\}$.

For the outer limit, fix $\cp<c_0<c'$ and put
$H(s)=e^{-\lam s}\phi_{\lam}(s)$, whose logarithmic derivative
is bounded. With $B(s)=\beta\coth s$ for point distance and
$B(s)=B_d(s)$ for geodesic distance, choose $R$ so large that
\[
 \lam c_0-k(\lam)-D|B(s)-\beta|\left|\frac{H'(s)}{H(s)}\right|>0
 \quad\text{for }s\ge R.
\]
Such a choice follows from Lemma~\ref{lem:drifts} and
$k(\lam)=\cp\lam$. The cell equation then makes
$Ce^{\lam c_0t}H(s)$ a supersolution of
$v_t=Dv_{ss}+DB(s)v_s+a(s)v$ on $s>R$.
By \eqref{eq:three-reductions}, its composition with the
corresponding coordinate $z(x)=r(x)$ or $z(x)=s(x)$ is a
supersolution of $v_t=D\Delta v+a(z(x))v$ on $\{z>R\}$,
and \eqref{eq:F1} makes the original solution a subsolution.
Take $R$ so that $u_0=0$ on $\{z>R\}$, and choose $C$ so
that $Ce^{\lam c_0t}H(R)\ge1$ for every $t\ge0$;
this is possible because $c_0>0$.
Exterior comparison yields
$u(x,t)\le C'e^{-\lam(c'-c_0)t}$ whenever $z(x)\ge c't$,
proving the outer limit for each $c'>\cp$.
\end{proof}

\subsection{Kernel bounds and front stability}

Fix a real $\lambda$ for this subsection. On the full space put
$H_\lambda(x)=e^{-\lambda b(x)}\phi_\lambda(b(x))$.
Let $K_a(t,x,y)$ be the kernel of
$e^{t(D\Delta+a(b(\cdot)))}$, relative to $\dd V(y)$.
For bounded initial datum $v_0$, the function
$v(\cdot,t)=e^{t(D\Delta+a(b(\cdot)))}v_0$
solves $v_t=D\Delta v+a(b(x))v$ on $\Hh^d$. Define
\[\label{eq:Doob-kernel}
 q_\lambda(t,x,y)=e^{-k(\lambda)t}
                    \frac{H_\lambda(y)}{H_\lambda(x)}K_a(t,x,y)\quad {\text{for }t>0,\quad x,y\in\Hh^d,\quad\lambda\in\R}.
\]
All densities on the full space in this formula are relative to
hyperbolic volume. The transition operators
$v_0\mapsto\int_{\Hh^d}q_\lambda(t,\cdot,y)v_0(y)\dd V(y)$
have generator
\[\label{eq:Doob-full}
 \mathcal Q_\lambda=D\Delta+2D\nabla\log H_\lambda\cdot\nabla
 =D\partial_{zz}+D_\lambda(z)\partial_z+De^{-2z}\Delta_y.
\]
For functions depending only on $z=b(x)$, $\mathcal Q_\lambda$
reduces to the operator $Q_\lambda$ in~\eqref{eq:Hlambda-line}.

\begin{lemma}\label{lem:periodic-kernel}
The transition density $q_\lambda^{(1)}(t,z,z')$ of $Q_\lambda$ on
$\R$, relative to $\dd z'$, satisfies
\[\label{eq:line-kernel-bound}
 \sup_{z,z'\in\R}q_\lambda^{(1)}(t,z,z')\le C_\lambda t^{-1/2}
 \quad\text{for }t>0.
\]
Consequently its transition operator maps $L^p(\R)$ to $L^\infty(\R)$
with norm at most $C_{p,\lambda}t^{-1/(2p)}$ for $1\le p<\infty$.
\end{lemma}

\begin{proof}
By \eqref{eq:Hlambda-line}, $D_\lambda$ is bounded, $C^1$ and
periodic. The positive $L$-periodic density $\pi_\lambda$ in~\eqref{eq:cell-invariant-density} is bounded above and away from zero on $\R$.
Write $\dd m_\lambda(z)=\pi_\lambda(z)\dd z$.
The equation $D\pi_\lambda''-(D_\lambda\pi_\lambda)'=0$
integrates to
$D\pi_\lambda'-D_\lambda\pi_\lambda=J_\lambda$ for a
constant $J_\lambda$. Consequently, integration by parts for
real compactly supported $v$ gives
\[\label{eq:line-energy}
 \int_\R v Q_\lambda v\,\pi_\lambda\dd z
       =-D\int_\R|v'|^2\pi_\lambda\dd z
 \quad {\text{for }v\in C_c^\infty(\R)}.
\]
The term involving $J_\lambda$ is a constant multiple of
$\int_\R(v^2)'\dd z=0$. The adjoint in $L^2(m_\lambda)$ is
$D\partial_{zz}+(2D\pi_\lambda'/\pi_\lambda-D_\lambda)\partial_z$;
its drift is bounded and periodic, and its quadratic form
satisfies \eqref{eq:line-energy}.

The diffusions generated by $Q_\lambda$ and its
$L^2(m_\lambda)$ adjoint are conservative because both drifts
are bounded.
To establish invariance of $m_\lambda$ on $\R$, let
$v_t=Q_\lambda v$ with $v(\cdot,0)\in C_c^\infty(\R)$,
and integrate against $\pi_\lambda$ and a compact cutoff.
The finite-time Gaussian tails for bounded drift, together with
interior derivative estimates at positive time, make the cutoff
errors vanish. The stationary equation for $\pi_\lambda$
therefore gives
$\int_\R v(z,t)\dd m_\lambda(z)
=\int_\R v(z,0)\dd m_\lambda(z)$ for $t>0$.
Approximation and monotone convergence extend this identity to
nonnegative data and hence to the transition semigroup.
Applying the cutoff integration to the adjoint equation proves
invariance of $m_\lambda$ for the adjoint transition semigroup. In what follows,
$L^p(\pi_\lambda)$ means $L^p(\R,m_\lambda)$.

The one-dimensional Nash inequality takes the form
\[\label{eq:Nash}
 \norm v_{L^2(\pi_\lambda)}^6
 \le C\norm{v'}_{L^2(\pi_\lambda)}^2
           \norm v_{L^1(\pi_\lambda)}^4.
\]
For Lebesgue measure on $\R$, the following Fourier argument applies.
Plancherel's identity and the splitting of
$\int_\R|\widehat v(\xi)|^2\dd\xi$ into $|\xi|\le R$
and $|\xi|>R$ bound these two integrals by $CR\norm v_{L^{1}(\R;\dd z)}^2$ and
$R^{-2}\norm{v'}_{L^{2}(\R;\dd z)}^2$. Optimisation in $R$ gives the usual
Nash inequality, and the upper and lower bounds for
$\pi_\lambda$ transfer it to \eqref{eq:Nash}; see~\cite{Nash,Davies}.

For $v_0\in C_c^\infty(\R)$, let $v(z,t)$ solve
$v_t=Q_\lambda v$ with $v(z,0)=v_0(z)$.
Positivity and invariance give $L^1(m_\lambda)$ contraction.
Combining \eqref{eq:line-energy} with \eqref{eq:Nash} yields
$Y'\le-cY^3/\norm{v_0}_1^4$, where
$Y(t)=\norm{v(\cdot,t)}_{L^2(\pi_\lambda\dd z)}^2$
and all norms in the differential inequality use $m_\lambda$.
For nonzero $v_0$, integration gives
$Y(t)\le[Y(0)^{-2}+2ct/\norm{v_0}_1^4]^{-1/2}$;
taking square roots proves the $L^1\to L^2$ bound $Ct^{-1/4}$.
The bound also holds for $v_0=0$. Applying the energy and Nash
inequalities to the adjoint and then using duality gives an
$L^2(m_\lambda)\to L^\infty$ norm bound $Ct^{-1/4}$. Composing the two estimates at time $t/2$
gives an $L^1(m_\lambda)\to L^\infty$ norm bound $Ct^{-1/2}$.

The kernel relative to $m_\lambda$ is
$q_\lambda^{(1)}(t,z,z')/\pi_\lambda(z')$.
Comparability of the measures consequently proves
\eqref{eq:line-kernel-bound} for the Lebesgue density.
Since $\int_\R q_\lambda^{(1)}(t,z,z')\dd z'=1$,
for $p>1$ and $p'=p/(p-1)$,
$\|q_\lambda^{(1)}(t,z,\cdot)\|_{L^{p'}(\R;\dd z')}
\le\bigl(\sup_{\zeta,\zeta'\in\R}q_\lambda^{(1)}(t,\zeta,\zeta')\bigr)^{1/p}$.
Thus integration against $q_\lambda^{(1)}(t,z,z')\,\dd z'$ maps
$L^p(\R,\dd z)$ to $L^\infty(\R)$ with norm at most
$C_{p,\lambda}t^{-1/(2p)}$ by H\"older's inequality. For $p=1$
the same estimate follows from \eqref{eq:line-kernel-bound}.
Density of $C_c^\infty(\R)$ gives the extension for every $1\le p<\infty$.
\end{proof}

\begin{lemma}\label{lem:full-kernel}
The kernel $q_\lambda$ defined by \eqref{eq:Doob-kernel}
satisfies $\int_{\Hh^d}q_\lambda(t,x,y)\dd V(y)=1$
for $x\in\Hh^d$ and $t>0$, and
\[\label{eq:full-kernel-bound}
 \sup_{x,y\in\Hh^d}q_\lambda(t,x,y)\le C_\lambda t^{-1/2}
 \quad\text{for }t\ge1.
\]
Its transition operator has $L^p(\Hh^d,\dd V)\to L^\infty$ norm at most
$C_{p,\lambda}t^{-1/(2p)}$ for $t\ge1$ and $1\le p<\infty$.
\end{lemma}

\begin{proof}
Let $X_t^\lambda=(Z_t,Y_t)$ have generator
$\mathcal Q_\lambda$ from \eqref{eq:Doob-full}, and denote
expectation from $x$ by $\mathbb E_x^\lambda$.
In horospherical coordinates its stochastic equations are
\[\label{eq:Doob-SDE}
 \left\{\begin{aligned}
 \dd Z_t&=\sqrt{2D}\dd B_t+D_\lambda(Z_t)\dd t,\\
 \dd Y_t&=\sqrt{2D}e^{-Z_t}\dd\widetilde B_t,
 \end{aligned}\right.
\]
where $B$ is a one-dimensional standard Brownian motion and
$\widetilde B$ is an independent $(d-1)$-dimensional standard Brownian motion. Since $D_\lambda$ is bounded,
$Z$ is nonexplosive, and its continuous path is bounded on every
finite time interval. Thus
$\int_0^t e^{-2Z_v}\dd v<\infty$ almost surely, which also
makes $Y$ well defined without explosion.
For a smooth relatively compact domain $\Omega\subset\Hh^d$,
let $K_a^\Omega(t,x,y)$ be the zero Dirichlet kernel of
$D\Delta+a(b)$ on $\Omega$, relative to $\dd V(y)$, and let
$q_\lambda^\Omega(t,x,y)$ be the transition density of $X^\lambda$
killed upon leaving $\Omega$. The generator calculation and the
stopped Feynman--Kac formula give
$q_\lambda^\Omega(t,x,y)=e^{-k(\lambda)t}
H_\lambda(y)K_a^\Omega(t,x,y)/H_\lambda(x)$
for $t>0$ and $x,y\in\Omega$. As $\Omega$ exhausts $\Hh^d$,
$K_a^\Omega\uparrow K_a$, so the transition density of
$X_t^\lambda$ is $q_\lambda$ in~\eqref{eq:Doob-kernel}.
Nonexplosion then gives conservativity \cite{Hsu,Pinsky}.

The one-Lipschitz property of $b$ and periodic eigenfunction bounds
give $H_\lambda(y)/H_\lambda(x)\le Ce^{|\lambda|\dist(x,y)}$.
Combining this with \eqref{eq:heat-estimate} at time one and
$K_a(1,x,y)\le e^{a_{\max}}p_1(\dist(x,y))$ yields
\[\label{eq:Doob-short}
 q_\lambda(1,x,y)\le C e^{-c_0\dist(x,y)^2}
 \le C e^{-c_0|b(x)-b(y)|^2}
 \quad {\text{for }x,y\in\Hh^d}
\]
for some $c_0>0$. For $t\ge2$, use the semigroup identity
$q_\lambda(t,x,y)=\int_{\Hh^d}q_\lambda(t-1,x,\zeta)q_\lambda(1,\zeta,y)\dd V(\zeta)$.
Applying \eqref{eq:Doob-short} and using the transition density of $Z_{t-1}$ gives
\[\label{eq:projection-kernel}
 \begin{aligned}
 q_\lambda(t,x,y)
 &\le C\mathbb E_x^\lambda e^{-c_0|b(X_{t-1}^\lambda)-b(y)|^2}\\
 &=C\int_\R q_\lambda^{(1)}(t-1,b(x),z)
                         e^{-c_0|z-b(y)|^2}\dd z
 \le Ct^{-1/2} \quad {\text{for }x,y\in\Hh^d,\quad t\ge2}.
 \end{aligned}
\]
Here the final bound follows from Lemma~\ref{lem:periodic-kernel}. The heat-kernel comparison used in~\eqref{eq:Doob-short}
bounds $q_\lambda(t,x,y)$ uniformly for $1\le t\le2$, so \eqref{eq:projection-kernel}
proves \eqref{eq:full-kernel-bound} for all $t\ge1$.
Conservativity and H\"older's inequality convert the kernel bound
into the $L^p\to L^\infty$ estimate with respect to hyperbolic
volume, by the calculation in Lemma~\ref{lem:periodic-kernel}.
\end{proof}

\begin{proof}[Proof of Theorem~\ref{thm:stability}]
Let $u$ solve \eqref{eq:KPP} with reaction $f_b$ and
initial datum $u_0$, and let $W$ be the front of~\eqref{eq:line}
in the statement. Put $v(x,t)=u(x,t)-W(b(x),t)$.
Subtracting the equations gives $v_t=D\Delta v+c_v(x,t)v$,
where
$c_v(x,t)=\int_0^1F_u(b(x),W(b(x),t)+\theta v(x,t))\dd\theta$.
Since $0\le W(b(x),t)+\theta v(x,t)\le1$ for $0\le\theta\le1$,
\eqref{eq:derivative} gives $c_v(x,t)\le a(b(x))$ and
\eqref{eq:F1} bounds $c_v$, including at $v=0$.
For $\varepsilon>0$, define
$q_\varepsilon(x,t)=(v(x,t)^2+\varepsilon^2)^{1/2}$. The chain rule
at positive time gives
\[\label{eq:absolute-value-regularization}
 (\partial_t-D\Delta)q_\varepsilon
 =\frac{c_v(x,t)v^2}{q_\varepsilon}
    -\frac{D\varepsilon^2|\nabla v|_g^2}{q_\varepsilon^3}
 \le a(b(x))q_\varepsilon,
\]
since $a(b(x))>0$ for $x\in\Hh^d$. Testing \eqref{eq:absolute-value-regularization}
on compact cylinders and passing to $\varepsilon\downarrow0$
yields
\[\label{eq:difference-linear}
 \partial_t|v|\le D\Delta|v|+a(b)|v|
\]
in distributions. Linear comparison therefore gives
$|v(\cdot,t)|\le e^{t(D\Delta+a(b))}|v(\cdot,0)|$.

For real $\lambda$, set
\[\label{eq:weighted-difference}
 V(x,t)=e^{-\lambda ct}\frac{|v(x,t)|}{H_\lambda(x)}
 \quad {\text{for }x\in\Hh^d,\quad t>0}.
\]
By the kernel transformation \eqref{eq:Doob-kernel},
\[
 V(x,t)\le e^{(k(\lambda)-c\lambda)t}
       \int_{\Hh^d}q_\lambda(t,x,y)
                   \frac{|v(y,0)|}{H_\lambda(y)}\dd V(y).
\]
Write $h_0=|v(\cdot,0)|/H_\lambda$.
Apply \eqref{eq:Doob-kernel} first to the bounded data
$H_\lambda\min\{h_0,n\}\le |v(\cdot,0)|$. Positivity and
monotone convergence as $n\to\infty$ give the integral bound
for $V$ displayed above.
At $\lambda=\lam$, \eqref{eq:weighted-data} is exactly
$\|h_0\|_{L^p(\Hh^d,\dd V)}=B_p<\infty$.
Applying Lemma~\ref{lem:full-kernel} and
$k(\lam)=\cp\lam$ proves \eqref{eq:weighted-stability}.
For \eqref{eq:line}, use
$H_\lambda(z)=e^{-\lambda z}\phi_\lambda(z)$ and the kernel
$q_\lambda^{(1)}$ from Lemma~\ref{lem:periodic-kernel}; integration
against $\dd z$ and \eqref{eq:line-weighted-data} give the weighted
estimate in~\eqref{eq:line-stability}. On $\Hh^d$,
the preceding calculation uses $\dd V$.

By comparison and \eqref{eq:front-bracket-explicit},
$W^{\eta_2}(b(x),t)\le u(x,t)\le W^{\eta_1}(b(x),t)$
for $x\in\Hh^d$, $t>0$ and $c\ne0$;
for $c=0$, the bounds are $P(b(x)+j_2L)\le u(x,t)
\le P(b(x)+j_1L)$ with the integers in
Theorem~\ref{thm:stability}.
Given $\varepsilon>0$, the uniform end limits of
$W^{\eta_1}$, $W^{\eta_2}$ and $W$ give $R<\infty$ such that
$|v|\le2\varepsilon$ outside $|b(x)-ct|\le R$ when $c\ne0$.
For $c=0$, the end limits of $P(\cdot+j_1L)$,
$P(\cdot+j_2L)$ and $P$ give this bound outside $|b(x)|\le R$.
For $|b(x)-ct|\le R$,
$e^{\lam(b(x)-ct)}/\phi_{\lam}(b(x))
\ge e^{-\lam R}/\max_{[0,L]}\phi_{\lam}$, so
\eqref{eq:weighted-stability} gives
\[\label{eq:fixed-strip-unweighting}
 \sup_{|b(x)-ct|\le R}|v(x,t)|
 \le e^{\lam R}\max_{[0,L]}\phi_{\lam}\,
 C_p t^{-1/(2p)}e^{-\lam(c-\cp)t}B_p\quad {\text{for }R>0,\quad t\ge1}.
\]
Since $p<\infty$ and $c\ge\cp$,
$t^{-1/(2p)}e^{-\lam(c-\cp)t}\to0$ as $t\to\infty$.
Taking $t\to\infty$ in~\eqref{eq:fixed-strip-unweighting} and then
$\varepsilon\downarrow0$ proves \eqref{eq:uniform-stability}.
For the solution $w$ in part~(iii), comparison with
$W^{\eta_1}$ and $W^{\eta_2}$ when $c\ne0$, or with
$P(\cdot+j_1L)$ and $P(\cdot+j_2L)$ when $c=0$,
controls $|w-W|$ on $|z-ct|>R$ by the uniform end limits of
these profiles. On $|z-ct|\le R$, the weight
$e^{\lam(z-ct)}/\phi_{\lam}(z)$ is at least
$e^{-\lam R}/\max_{[0,L]}\phi_{\lam}$.
The weighted estimate in~\eqref{eq:line-stability} therefore gives
uniform convergence and completes the proof of part~(iii).
\end{proof}

\begin{proposition}\label{prop:linear-position}
Under \eqref{eq:F1}, fix $B\in\R$ and let $w$ solve
\eqref{eq:line} with measurable initial datum $w_0$ satisfying
$0\le w_0(z)\le\one_{\{z\le B\}}$ almost everywhere.
Then for each $\theta>0$,
\[\label{eq:linear-position}
 \{z:w(z,t)\ge\theta\}
 \subset\left\{z:z\le\cp t-\frac1{2\lam}\log t+C_\theta\right\}
 \quad\text{for }t\ge1.
\]
\end{proposition}

\begin{proof}
By the initial bound and periodic positivity of $\phi_{\lam}$,
$e^{\lam z}w(z,0)/\phi_{\lam}(z)$ belongs to $L^1(\R)$,
with norm at most
$e^{\lam B}/(\lam\min_{[0,L]}\phi_{\lam})$.
Linear comparison and the identity
$e^{-k(\lambda)t}H_\lambda^{-1}
 e^{t(D\partial_{zz}+D\beta\partial_z+a)}(H_\lambda h)
=e^{tQ_\lambda}h$, for $h\in C_c^\infty(\R)$, with
$H_\lambda(z)=e^{-\lambda z}\phi_\lambda(z)$,
followed by Lemma~\ref{lem:periodic-kernel}, therefore give
$w(z,t)\le Ct^{-1/2}e^{-\lam(z-\cp t)}\phi_{\lam}(z)$.
Taking logarithms at a point where $w(z,t)\ge\theta$, and
using the periodic upper bound for $\phi_{\lam}$, proves
\eqref{eq:linear-position}.
\end{proof}

\subsection{Logarithmic level positions}

Proposition~\ref{prop:linear-position} gives a logarithmic upper
bound with coefficient $1/(2\lam)$ by linear comparison.
We now prove the sharper coefficient $3/(2\lam)$ for the periodic
coefficient $a(z)$ when $\mu>0$. The scalar estimates below use
the eigenvalue $k$ and minimiser $\lam$ defined in~\eqref{eq:cell}--\eqref{eq:parameters}.

\begin{lemma}\label{lem:periodic-logistic-step}
Assume \eqref{eq:F1} and $\mu>0$, and let $S\in\R$.
Let $W$ solve
\[\label{eq:periodic-logistic-step-equation}
 \left\{\begin{aligned}
 W_t&=DW_{zz}+D\beta W_z+a(z)W(1-W)
       &&\text{for }z\in\R,\ t>0,\\
 W(z,0)&=\one_{(-\infty,S]}(z)&&\text{for }z\in\R.
 \end{aligned}\right.
\]
Use $m_*$ from \eqref{eq:periodic-position}.
For every $\eta\in(0,1)$ there are $C_\eta,T_\eta<\infty$ such that
\[\label{eq:periodic-logistic-step-lower}
 W(z,t)\ge\eta
 \quad\text{for }t\ge T_\eta,\ z\le m_*(t)-C_\eta.
\]
\end{lemma}

\begin{proof}
Consider binary branching diffusion on $\R$ whose particles
move with generator $D\partial_{yy}-D\beta\partial_y$ and branch at
rate $a(-y)$. Denote the law and expectation for one initial particle
at $y$ by $\mathbb P_y$ and $\mathbb E_y$, respectively, and write
$M_t$ for the largest particle position at time $t$. The coefficient $a(-y)$ is positive, $C^1$, and
$L$-periodic. The exponential eigenvalue problem for this process is
\[\label{eq:reflected-branching-cell}
 D\psi''+D(2\lambda-\beta)\psi'
 +(D\lambda^2-D\beta\lambda+a(-y))\psi=k(\lambda)\psi,
 \qquad\psi(y)=\phi_\lambda(-y)
 \quad\text{for }y\in\R,\ \lambda\in\R.
\]
Under the change of time $\tau=Dt$, particle motions have
generator $\partial_{yy}-\beta\partial_y$ and branching occurs at
rate $a(-y)/D$. At rescaled time $\tau$, the rightmost particle is at $M_{\tau/D}$,
and the exponential eigenvalue is $k(\lambda)/D$ by
\eqref{eq:reflected-branching-cell}. Keep $\tau=Dt$ and set $\xi=y/L$. In $(\xi,\tau)$ the generator is
$L^{-2}\partial_{\xi\xi}-(\beta/L)\partial_\xi$, and the branching rate
$a(-L\xi)/D$ has period one. These smooth periodic coefficients and binary
offspring satisfy the assumptions of~\cite[Theorem~5.1]{LTZ22}.
The speed and critical exponent in these variables are $\cp/(DL)$ and
$L\lam$, respectively. The speed is positive by $\mu>0$ and
Lemma~\ref{lem:cell}. The logarithmic centring term in $\xi$ is therefore
$3\log\tau/(2L\lam)$. Multiplication by $L$ restores
$3\log\tau/(2\lam)$ in $y$, and a displacement $h$ in $y$ contributes
the factor $e^{-(L\lam)(h/L)}=e^{-\lam h}$.

For this diffusion, write $Z_\infty$ and $\nu$ for the positive
random variable and the positive continuous periodic function in~\cite[Section~4, equation~(4.1)]{LTZ22}, with the extension furnished
by \cite[Theorem~5.1]{LTZ22} and the spatial period restored to $L$.
Set $\nu_-=\min_{[0,L]}\nu$ and $\nu_+=\max_{[0,L]}\nu$,
so $0<\nu_-\le\nu_+<\infty$.
For each fixed $h\in\R$, the formula in~\cite[Theorem~5.1]{LTZ22} bounds the lower and upper limits of
$\mathbb P_0(M_{\tau/D}-(\cp/D)\tau+3\log\tau/(2\lam)\le h)$,
as $\tau\to\infty$, between
$\mathbb E_0[e^{-Z_\infty\nu_+e^{-\lam h}}]$ and
$\mathbb E_0[e^{-Z_\infty\nu_-e^{-\lam h}}]$.
Bounded convergence makes the first expectation tend to one as
$h\to\infty$ and the second tend to zero as $h\to-\infty$.
These two bounds give tightness under $\mathbb P_0$ as
$\tau\to\infty$. Returning to $t$ and using
$\log(Dt)=\log t+\log D$ gives tightness of
$M_t-\cp t+3\log t/(2\lam)$ as $t\to\infty$.

We record why tightness is uniform in the starting point modulo
$L$. For $h\in\R$, let
$v_h(y,t)=\mathbb P_y(M_t\ge h)$. The branching property gives
\[
 (v_h)_t=D(v_h)_{yy}-D\beta(v_h)_y
       +a(-y)v_h(1-v_h)
 \quad\text{for }t>0,\ y\in\R.
\]
Both $v_h$ and $1-v_h$ are positive solutions of linear parabolic
equations with uniformly bounded zero-order coefficients.
By the parabolic Harnack inequality \cite{Friedman}, there is a constant $C$,
independent of $h$, such that, for $t\ge1$ and $0\le y\le L$,
\[
 \begin{aligned}
 v_h(y,t)&\le C v_h(0,t+1),\\
 1-v_h(y,t)&\le C\bigl(1-v_h(0,t+1)\bigr).
 \end{aligned}
\]
Since $m_*(t+1)-m_*(t)$ is bounded for $t\ge1$, tightness under
$\mathbb P_0$ implies tightness of $M_t-y-m_*(t)$ uniformly for
$y\in[0,L]$. Under $\mathbb P_{y+jL}$ the centred maximum $M_t-(y+jL)-m_*(t)$
has the same law as $M_t-y-m_*(t)$ under $\mathbb P_y$ for each $j\in\Z$.
Thus the tightness is uniform over $y\in\R$. In particular,
\[\label{eq:branching-uniform-lower-tail}
 \lim_{K\to\infty}\limsup_{t\to\infty}
 \sup_{y\in\R}
 \mathbb P_y\bigl(M_t-y\le m_*(t)-K\bigr)=0.
\]

The branching representation of~\eqref{eq:periodic-logistic-step-equation}, followed by the reflection
$y=-z$, gives
\[
 W(z,t)=\mathbb P_{-z}(M_t\ge-S)
 \quad\text{for }z\in\R,\ t>0.
\]
If $z\le m_*(t)+S-K$, then $z-S\le m_*(t)-K$. Applying
\eqref{eq:branching-uniform-lower-tail} with the starting point
$y=-z$ therefore proves \eqref{eq:periodic-logistic-step-lower},
after choosing $K$ and then $T_\eta$.
\end{proof}

\begin{lemma}\label{lem:periodic-dirichlet-log}
Assume \eqref{eq:F1} and $\mu>0$. Set $c=\cp$, $\lambda=\lam$
and $\kappa=3/(2\lambda)$. There are $T>\kappa/c$ and a positive
solution $J$ of
\[\label{eq:periodic-dirichlet-log-equation}
 \begin{gathered}
 \widehat m(t)=ct-\kappa\log\frac{t+T}{T}
       \quad\text{for }t\ge0,\\
 \left\{\begin{aligned}
 J_t&=DJ_{zz}+D\beta J_z+a(z)J
       &&\text{for }t>0,\ z>\widehat m(t),\\
 J(\widehat m(t),t)&=0&&\text{for }t>0,
 \end{aligned}\right.
 \end{gathered}
\]
with a nonzero nonnegative initial function in $C_c^\infty((0,\infty))$,
and constants $c_D,C_D,\varepsilon,t_0>0$ such that
\[\label{eq:periodic-dirichlet-log-bounds}
 c_D\,y e^{-\lambda y}\le J(\widehat m(t)+y,t)
 \le C_D\,y e^{-\lambda y}
 \quad\text{for }t\ge t_0,\ 0<y\le\varepsilon\sqrt t.
\]
The initial function may be prescribed before choosing the constants.
\end{lemma}

\begin{proof}
By Lemma~\ref{lem:cell}, $k'(\lambda)=c>0$ and $k'(\beta/2)=0$.
Strict convexity of $k$ therefore gives $q_0:=\lambda-\beta/2>0$.
Define
\[\label{eq:auxiliary-periodic-potential}
 A(z)=a(z)-\frac{D\beta^2}{4}-\frac{\beta c}{2},\qquad
 \gamma(\eta)=k(\eta+\beta/2)-\frac{\beta c}{2}
 \quad\text{for }z,\eta\in\R.
\]
The principal periodic eigenvalue of
$D\partial_{zz}-2D\eta\partial_z+D\eta^2+A(z)$ is $\gamma(\eta)$.
In view of $k(\lambda)=c\lambda$ and $k'(\lambda)=c$,
\[\label{eq:auxiliary-periodic-tangency}
 \gamma(q_0)=cq_0,\qquad \gamma'(q_0)=c,\qquad
 \gamma(\eta)>c\eta\quad\text{for }\eta\ne q_0.
\]
In particular, $\gamma(0)>0$, and the minimum of
$\gamma(\eta)/\eta$ over $\eta>0$ is $c$, attained at $q_0$.
The corresponding positive periodic eigenfunction is $\phi_\lambda$.
The auxiliary reaction $A(z)v-v^2$ is concave in $v$, has the
linearisation $A(z)$, and is negative for all sufficiently large $v$.
Together with $\gamma(0)>0$, these properties verify the hypotheses
used for the linear estimates in~\cite{HNRR16}.

Let
\[\label{eq:dirichlet-gauge}
 Z(z,t)=e^{\beta(z-ct)/2}J(z,t)
 \quad\text{for }t>0,\ z>\widehat m(t).
\]
By \eqref{eq:periodic-dirichlet-log-equation} and
\eqref{eq:auxiliary-periodic-potential}, $Z_t=DZ_{zz}+A(z)Z$ on
$z>\widehat m(t)$, with zero boundary value at $\widehat m(t)$.
To apply \cite[Proposition~2.3 and equations~(2.13)--(2.18)]{HNRR16}
in its unit-diffusion normalisation, set $\widehat t=Dt$ and
$\widetilde Z(z,\widehat t)=Z(z,\widehat t/D)$. Then
$\widetilde Z_{\widehat t}=\widetilde Z_{zz}+(A(z)/D)\widetilde Z$,
and its boundary is
$\widehat m(\widehat t/D)=(c/D)\widehat t-\kappa\log((\widehat t+DT)/(DT))$.
The auxiliary reaction is $(A(z)v-v^2)/D$, and the periodic spectral
value at zero is $\gamma(0)/D>0$ by
\eqref{eq:auxiliary-periodic-tangency}. In $(z,\widehat t)$ the critical speed is $c/D$ and the critical
exponential decay rate is $q_0$. In the unit-period variables $z/L$ and $\widehat t/L^2$, the critical speed
is $cL/D$, the decay rate is $Lq_0$, and the logarithmic displacement
$\kappa\log((t+T)/T)$ in~\eqref{eq:periodic-dirichlet-log-equation}
becomes $(\kappa/L)\log((t+T)/T)$.
To straighten the boundary, set $c\tau=\widehat m(t)$ and write
$t=\mathfrak t(\tau)$. The two rescaled times are
$\widehat t=D\mathfrak t(\tau)$ and $\widehat\tau=D\tau$.
Their difference is $(D\kappa/c)\log((t+T)/T)=O(\log\tau)$ as
$\tau\to\infty$, since $\tau/t\to1$.
The choice $T>\kappa/c$ makes this change of time increasing. With
\[\label{eq:dirichlet-time-factor}
 Z(z,\mathfrak t(\tau))=e^{-q_0(z-c\tau)}\phi_\lambda(z)
 \alpha(\tau)p(z,\tau),\qquad
 \alpha(\tau)=e^{q_0c(\mathfrak t(\tau)-\tau)}
 \quad\text{for }\tau>0,\ z>c\tau,
\]
direct differentiation yields
\[\label{eq:dirichlet-conjugated-equation}
 \left\{\begin{aligned}
 &\frac{1}{\mathfrak t'(\tau)}p_\tau
 =Dp_{zz}+D\left(-2q_0+2\frac{\phi_\lambda'}{\phi_\lambda}\right)p_z
 \quad\text{for }\tau>0,\ z>c\tau,\\
 &p(c\tau,\tau)=0\quad\text{for }\tau>0.
 \end{aligned}\right.
\]
Moreover, $1/\mathfrak t'(\tau)=1-\kappa/[c(\mathfrak t(\tau)+T)]$.
The difference from one is $O(\tau^{-1})$, and its derivative is
$O(\tau^{-2})$ as $\tau\to\infty$. After setting $\widehat\tau=D\tau$,
the function $\widetilde p(z,\widehat\tau)=p(z,\widehat\tau/D)$ satisfies
$\widetilde p_{\widehat\tau}/\mathfrak t'(\widehat\tau/D)
=\widetilde p_{zz}+(-2q_0+2\phi_\lambda'/\phi_\lambda)\widetilde p_z$
on $z>(c/D)\widehat\tau$, with zero boundary value. In the variables $z/L$ and $\widehat\tau/L^2$, the spatial period is one.
After solving for the time derivative, the diffusion coefficient is
$\mathfrak t'(\tau)=1+O(\tau^{-1})$ as $\tau\to\infty$,
where $\tau=\widehat\tau/D$.
Let $\tilde a_1,\tilde a_2,\tilde\varepsilon,\tilde\tau_0>0$ be the constants
in~\cite[Proposition~2.3]{HNRR16} for
\eqref{eq:dirichlet-conjugated-equation} in the coordinates
$(z/L,\widehat\tau/L^2)$ with zero Dirichlet data on the line
$z/L=(cL/D)(\widehat\tau/L^2)$.
At $z=c\tau+y$, the distance from the boundary in the rescaled spatial
coordinate is $y/L$, and the rescaled time is $D\tau/L^2$.
Thus $(y/L)(D\tau/L^2)^{-3/2}=L^2D^{-3/2}y\tau^{-3/2}$.
Returning to $(z,\tau)$ gives
$a_1 y\tau^{-3/2}\le p(c\tau+y,\tau)\le a_2 y\tau^{-3/2}$
for $\tau\ge\tau_0$ and $0<y\le\varepsilon\sqrt\tau$,
where $a_i=\tilde a_iL^2D^{-3/2}$ for $i=1,2$,
$\varepsilon=\tilde\varepsilon\sqrt D$ and
$\tau_0=\tilde\tau_0L^2/D$.
Furthermore,
$\alpha(\tau)=((\mathfrak t(\tau)+T)/T)^{q_0\kappa}$ and
$\mathfrak t(\tau)/\tau\to1$ as $\tau\to\infty$. Combining these facts with
\eqref{eq:dirichlet-gauge}, at $z=\widehat m(t)+y$, gives
\[
 J(\widehat m(t)+y,t)\asymp
 y e^{-(q_0+\beta/2)y}
 (t+T)^{(q_0+\beta/2)\kappa}t^{-3/2}
 \quad\text{for }t\ge t_0,\ 0<y\le\varepsilon\sqrt t.
\]
Since $q_0+\beta/2=\lambda$ and $\lambda\kappa=3/2$, this proves
\eqref{eq:periodic-dirichlet-log-bounds}.
\end{proof}

\begin{lemma}\label{lem:periodic-busemann-upper}
Assume \eqref{eq:F1} and $\mu>0$, and let $u$ solve \eqref{eq:KPP}
with $f=f_b$. Suppose that $0\le u_0\le1$ and that $u_0(x)=0$ for
almost every $x$ satisfying $b(x)>R$, for some $R\in\R$.
For every $\theta\in(0,1)$ there are $C_\theta,T_\theta<\infty$ such that
\[\label{eq:periodic-busemann-upper}
 u(x,t)<\theta\quad\text{for }t\ge T_\theta,\ x\in\Hh^d
 \text{ with }b(x)\ge m_*(t)+C_\theta.
\]
\end{lemma}

\begin{proof}
Use $J,\widehat m,c,\lambda,\kappa$ from
Lemma~\ref{lem:periodic-dirichlet-log}. Increase $t_0$ so that
$\varepsilon\sqrt{t_0}>1$, and choose $M$ with
$M c_De^{-\lambda}\ge2$. Choose an integer $j$ so that
$\widehat m(t_0)+jL+1>R$. On the moving domain
$z>\widehat m(t+t_0)+jL+1$, the function
$M J(z-jL,t+t_0)$ solves the linear equation with coefficient $a(z)$,
by periodicity. Its value is at least one on
the lateral boundary, by \eqref{eq:periodic-dirichlet-log-bounds}.
At $t=0$, the assumed support bound gives $u_0=0$ on this domain.
By Lemma~\ref{lem:coordinates},
$v(x,t)=MJ(b(x)-jL,t+t_0)$ solves
$v_t=D\Delta v+a(b(x))v$ on this moving domain.
On each finite time interval, $v-u$ is bounded and has nonnegative
initial and lateral values. Its residual under
$\partial_t-D\Delta-a(b(x))$ is
$a(b(x))u-F(b(x),u)\ge0$.
The stopped calculation in Lemma~\ref{lem:moving-domain-FK}, with
$b_0(t)=\widehat m(t+t_0)+jL+1$, therefore gives
\[\label{eq:periodic-upper-dirichlet-comparison}
 u(x,t)\le M J(b(x)-jL,t+t_0)
 \quad\text{for }t\ge0,
 \quad b(x)>\widehat m(t+t_0)+jL+1.
\]
By linear comparison with the positive eigenfunction,
\[\label{eq:periodic-upper-free-exponential}
 u(x,t)\le C e^{-\lambda(b(x)-ct)}\phi_\lambda(b(x))
 \quad\text{for }t\ge0,\ x\in\Hh^d.
\]
Let $y=b(x)-jL-\widehat m(t+t_0)$. For
$1<y\le\varepsilon\sqrt{t+t_0}$, equations
\eqref{eq:periodic-dirichlet-log-bounds} and
\eqref{eq:periodic-upper-dirichlet-comparison} give
$u(x,t)\le MC_D\,y e^{-\lambda y}$. For
$y>\varepsilon\sqrt{t+t_0}$, periodic boundedness of
$\phi_\lambda$ and \eqref{eq:periodic-upper-free-exponential} give
\[
 u(x,t)\le C(1+t)^{3/2}
 e^{-\lambda\varepsilon\sqrt{t+t_0}}
 \quad\text{for }t\ge1.
\]
Given $0<\theta<1$, choose $Y>\max\{1,1/\lambda\}$ such that
$MC_D Y e^{-\lambda Y}<\theta$, and take $T_\theta\ge1$ so that
$C(1+t)^{3/2}e^{-\lambda\varepsilon\sqrt{t+t_0}}<\theta$ for $t\ge T_\theta$.
Since $y e^{-\lambda y}$ decreases for $y\ge Y$, both estimates give
$u(x,t)<\theta$ when $y\ge Y$ and $t\ge T_\theta$.
The difference $\widehat m(t+t_0)+jL-m_*(t)$ is bounded there.
Thus \eqref{eq:periodic-busemann-upper} holds with
$C_\theta=Y+\max\{0,\sup_{t\ge T_\theta}
[\widehat m(t+t_0)+jL-m_*(t)]\}$.
\end{proof}

\begin{lemma}\label{lem:periodic-distance-transfer}
Assume \eqref{eq:F1}, \eqref{eq:axis-regularity}, and $\mu>0$.
Let $U$ solve either
\[
 U_t=DU_{rr}+D\beta\coth r\,U_r+F(r,U)
 \quad\text{for }r>0,\quad t>0,
\]
or
\[
 U_t=DU_{rr}+DB_d(r)U_r+F(r,U)
 \quad\text{for }r>0,\quad t>0,
\]
with $U_r(0,t)=0$ for $t>0$, interpreted by composition with
$r(x)$ or $s(x)$, respectively. Assume that there is $R_0>0$
such that $0\le U(r,0)\le1$, $U(r,0)=0$ for $r>R_0$, and
$\int_0^{R_0}U(r,0)\,\dd r>0$. For every $0<\theta<1$, there are $C_\theta,T_\theta>0$ such that
\[
 [0,m_*(t)-C_\theta]
 \subset\{r\ge0:U(r,t)\ge\theta\}
 \subset[0,m_*(t)+C_\theta]
 \quad\text{for }t\ge T_\theta.
\]
\end{lemma}

\begin{proof}
Write $B(r)=\beta\coth r$ or $B(r)=B_d(r)$, according to the
equation for $U$. The geometric formulas give a constant $C_B>0$
such that
\[\label{eq:periodic-distance-drift}
 |B(r)-\beta|+|B'(r)|\le C_Be^{-2r}
 \quad\text{for }r\ge1.
\]
Define
\[\label{eq:periodic-distance-gauge}
 h_B(r)=\exp\left(\frac12\int_r^\infty[B(y)-\beta]\,\dd y\right),
 \qquad
 E_B(r)=\frac{\beta^2-B(r)^2}{4}-\frac{B'(r)}2
 \quad\text{for }r\ge1.
\]
By \eqref{eq:periodic-distance-drift}, there is $C_E>0$ such that
\[\label{eq:periodic-distance-gauge-bounds}
 |h_B(r)-1|+|E_B(r)|\le C_Ee^{-2r}
 \quad\text{for }r\ge1.
\]
Since $2h_B'/h_B+B=\beta$, direct differentiation gives
\[\label{eq:periodic-distance-gauge-identity}
 h_B^{-1}(\partial_{rr}+B\partial_r)(h_Bv)
 =v_{rr}+\beta v_r+E_B(r)v
 \quad\text{for }r>1,\quad v\in C^2((1,\infty)).
\]

Fix $0<c_1<\cp$. By Theorem~\ref{thm:periodic-speed},
\[\label{eq:periodic-distance-interior}
 \sup_{0\le r\le c_1t}|U(r,t)-1|\to0
 \quad\text{as }t\to\infty.
\]
In the geodesic case, a general nonzero profile has a positive
minimum on every bounded interval at time one. Fix $R_1>0$ and set $\varepsilon=\min_{[0,R_1]}U(\cdot,1)>0$.
Comparison with the datum $\varepsilon\one_{[0,R_1]}$ at time one and Theorem~\ref{thm:periodic-speed}(iii), applied at a speed
$c_2\in(c_1,\cp)$, give \eqref{eq:periodic-distance-interior}, since
$c_1t\le c_2(t-1)$ for all sufficiently large $t$.
Choose $T>0$ such that $c_1T\ge1$ and
$U(c_1t,t)\ge1/2$ for $t\ge T$.
Let $h_{B,-},h_{B,+}>0$ be lower and upper bounds for $h_B$ on
$[c_1T,\infty)$, and choose $q_0>0$ with $q_0h_{B,+}\le1$.
Choose $0<\varepsilon\le1/2$ such that
$C_2\varepsilon\le a_{\min}$, and let $V$ solve
\[
 \left\{\begin{aligned}
 V_t&=DV_{rr}+D\beta V_r+a(r)V(1-V)
       &&\text{for }r\in\R,\ t>0,\\
 V(r,0)&=\one_{(-\infty,c_1T]}(r)&&\text{for }r\in\R.
 \end{aligned}\right.
\]
Define
\[\label{eq:periodic-distance-lower-amplitude}
 q(t)=q_0\exp\left(-DC_E\int_T^t e^{-2c_1\tau}\,\dd\tau\right)
 \quad\text{for }t\ge T.
\]
The integral in~\eqref{eq:periodic-distance-lower-amplitude} is finite,
so $q_*:=\inf_{t\ge T}q(t)>0$.

We claim that
\[\label{eq:periodic-distance-lower-comparison}
 \varepsilon q(t)h_B(r)V(r,t-T)\le U(r,t)
 \quad\text{for }r\ge c_1t,\quad t\ge T.
\]
Put $w(r,t)=\varepsilon V(r,t-T)$ for $r\in\R$ and $t\ge T$.
By \eqref{eq:periodic-distance-gauge-identity}, the residual of
$q(t)h_B(r)w$ in the equation for $U$, evaluated at $r>c_1t$ and $t>T$, is
\[
 \begin{aligned}
 &q(t)h_B(r)\left(\frac{q'}q-DE_B(r)\right)w\\
 &\quad+q(t)h_B(r)\left(a(r)w-\frac{a(r)}\varepsilon w^2\right)
 -F(r,q(t)h_B(r)w).
 \end{aligned}
\]
By \eqref{eq:periodic-distance-gauge-bounds} and
\eqref{eq:periodic-distance-lower-amplitude}, the first term is
nonpositive. Since $0<q(t)h_B(r)\le1$, the Taylor lower bound in~\eqref{eq:uniform-taylor} gives an upper bound for the other terms
of
\[
 -q(t)h_B(r)\left(\frac{a(r)}\varepsilon-C_2q(t)h_B(r)\right)w^2\le0
 \quad\text{for }r>c_1t,\quad t>T.
\]
Thus $q(t)h_B(r)w$ is a subsolution. It is zero initially on
$r>c_1T$, and it is at most $1/2$ on $r=c_1t$.
The choice of $T$ therefore orders it below $U$ on the initial
and lateral boundaries. Comparison proves
\eqref{eq:periodic-distance-lower-comparison}.

We next prove the exponential estimate
\[\label{eq:periodic-distance-sharp-exponential}
 U(r,t)\le C\exp[-\lam(r-\cp t)]
 \quad\text{for }r\ge0,\quad t\ge0.
\]
Take $R>1$ larger than the support of $U(\cdot,0)$, and put
\[
 Q_0(t)=\exp\left(DC_E\int_0^t e^{-2(R+c_1\tau)}\,\dd\tau\right)
 \quad\text{for }t\ge0.
\]
By \eqref{eq:periodic-distance-gauge-identity} and the cell
identity $k(\lam)=\cp\lam$, the function
\[
 Z(r,t)=A Q_0(t)h_B(r)e^{-\lam(r-\cp t)}\phi_{\lam}(r)
 \quad\text{for }r>R+c_1t,\quad t\ge0
\]
satisfies
$Z_t-DZ_{rr}-DB(r)Z_r-a(r)Z=(Q_0'/Q_0-DE_B(r))Z\ge0$.
Since $\cp>c_1$ and $h_B\phi_{\lam}$ has a positive lower
bound on $r\ge R$, choose $A$ so that
$Z(R+c_1t,t)\ge1$ for all $t\ge0$.
Initially $U=0$ on $r>R$. Linear comparison and $F(r,U)\le a(r)U$
therefore give $U\le Z$ on this moving exterior domain.
The bounds for $Q_0,h_B,\phi_{\lam}$ give
\eqref{eq:periodic-distance-sharp-exponential} there.
On $0\le r\le R+c_1t$, it follows from $U\le1$ after
increasing $C$ to at least $e^{\lam R}$.

For the logarithmic upper bound, let $J(z,t)$ denote the linear
Dirichlet solution supplied by Lemma~\ref{lem:periodic-dirichlet-log},
posed on $z>\widehat m(t)$ with $J(\widehat m(t),t)=0$.
Thus, with $\kappa_*=3/(2\lam)$ and some $T_*>0$,
\[
 \widehat m(t)=\cp t-\kappa_*\log((t+T_*)/T_*)
 \quad\text{for }t\ge0,
\]
and there are $\delta,C_J,\epsilon_J>0$ such that
\[\label{eq:periodic-distance-dirichlet-input}
 \begin{gathered}
 J(\widehat m(t)+1,t)\ge\delta,\\
 J(\widehat m(t)+y,t)\le C_J y e^{-\lam y}
 \quad\text{for }1\le y\le\epsilon_J\sqrt t
 \end{gathered}
 \quad\text{for all sufficiently large }t.
\]
Choose $t_0>0$ so that both bounds in~\eqref{eq:periodic-distance-dirichlet-input}
and the inequality $\widehat m'(t)\ge\cp/2$ hold for every $t\ge t_0$.
Choose an integer $j$ so that
$\xi(0):=\widehat m(t_0)+jL+1$ exceeds both one and the support radius of
$U(\cdot,0)$, and define
\[
 \xi(t)=\widehat m(t+t_0)+jL+1,\qquad
 Q(t)=\exp\left(DC_E\int_0^t e^{-2\xi(\tau)}\,\dd\tau\right)
 \quad\text{for }t\ge0.
\]
As $\xi(t)\ge\xi(0)+\cp t/2$, the function $Q$ is bounded.
The periodicity identity $a(r-jL)=a(r)$, together with
\eqref{eq:periodic-distance-gauge-identity}, shows that
\[
 Z_J(r,t)=A_JQ(t)h_B(r)J(r-jL,t+t_0)
 \quad\text{for }r>\xi(t),\quad t\ge0
\]
has residual $(Q'/Q-DE_B(r))Z_J\ge0$ in the linearised equation.
By \eqref{eq:periodic-distance-dirichlet-input}, choose $A_J$
so that $Z_J(\xi(t),t)\ge1$ for all $t\ge0$.
Since $U(r,0)=0$ on $r>\xi(0)$, comparison yields
\[\label{eq:periodic-distance-dirichlet-comparison}
 U(r,t)\le Z_J(r,t)
 \quad\text{for }r\ge\xi(t),\quad t\ge0.
\]

Write $\tau=t+t_0$ and $y=r-jL-\widehat m(\tau)$.
On $1\le y\le\epsilon_J\sqrt\tau$,
\eqref{eq:periodic-distance-dirichlet-input} and
\eqref{eq:periodic-distance-dirichlet-comparison} give
$U(r,t)\le C y e^{-\lam y}$.
On $y\ge\epsilon_J\sqrt\tau$,
\eqref{eq:periodic-distance-sharp-exponential} gives
\[
 U(r,t)\le C(\tau+T_*)^{3/2}e^{-\lam y}
 \le C(\tau+T_*)^{3/2}e^{-\lam\epsilon_J\sqrt\tau}
 \to0
 \quad\text{as }t\to\infty.
\]
The bounded difference $\widehat m(t+t_0)+jL-m_*(t)$ therefore proves the
outer inclusion at every level $0<\theta<1$.

Apply \eqref{eq:periodic-logistic-step-lower} with $\eta=1/2$ at time $t-T$.
Since $m_*(t)-m_*(t-T)\to\cp T$ as $t\to\infty$, there are
$C_0,T_0>0$ such that $V(r,t-T)\ge1/2$ for
$r\le m_*(t)-C_0$ and $t\ge T_0$.
Combining this estimate with
\eqref{eq:periodic-distance-lower-comparison}, the positive bounds
for $q,h_B$, and \eqref{eq:periodic-distance-interior} gives
\[\label{eq:periodic-distance-small-level}
 U(r,t)\ge\alpha
 \quad\text{for }0\le r\le m_*(t)-C_0,\quad t\ge T_0,
\]
for some $\alpha>0$.
Apply Lemma~\ref{lem:plateau-amplification} with initial lower bound
$\alpha$ and error $1-\theta$. Let $R_\theta,S_\theta>0$ denote
the radius and time furnished by Lemma~\ref{lem:plateau-amplification}. Both $r(x)$ and
$s(x)$ are one-Lipschitz,
so the ball of radius $R_\theta$ centred at a point of distance
at most $m_*(t)-C_0-R_\theta$ lies in the region covered by
\eqref{eq:periodic-distance-small-level}.
The amplification lemma yields
\[
 U(r,t+S_\theta)\ge\theta
 \quad\text{for }0\le r\le m_*(t)-C_0-R_\theta,\quad t\ge T_0.
\]
Since $m_*(t+S_\theta)-m_*(t)$ remains bounded, this proves the
inner inclusion.

All comparisons on moving half-lines can be justified on each
finite time interval by putting $y=r-\zeta(t)$, where $\zeta$
is the relevant boundary, and using a barrier
$\eta e^{Ct}(1+y^2)$. The coefficients are bounded on the resulting
half-strip. The maximum principle on bounded intervals, followed
by $\eta\downarrow0$, gives the asserted comparison on the whole
half-line.
\end{proof}

\begin{lemma}\label{lem:horo-subcritical-lower}
Assume \eqref{eq:F1} and $\mu>0$, and let $u$ solve
\eqref{eq:KPP} with $f=f_b$ and a nonzero initial datum
$0\le u_0\le1$. Fix $c\in(0,\cp)$. In the coordinates
\eqref{eq:horo-metric}, let
$\mathcal B=\{y\in\R^{d-1}:|y|<1\}$, and let $\psi$ be the positive
Dirichlet eigenfunction of $-\Delta_y$ on $\mathcal B$, normalised by
$\max_{\overline{\mathcal B}}\psi=1$. There are constants $S,A_0>0$
and $T_0>0$ such that
\[
 u(S+ct,y,T_0+t)\ge A_0\psi(y)
 \quad\text{for }t\ge0\text{ and }y\in \mathcal B.
\]
\end{lemma}

\begin{proof}
Write $\nu>0$ for the eigenvalue of $\psi$. Since
$k(\lambda)-c\lambda>0$ for $\lambda>0$, while
$k(\lambda)-c\lambda\ge\mu>0$ for $\lambda\le0$,
Lemma~\ref{lem:periodic-moving-bumps} applies with $I=\{c\}$.
Let $\varphi(z,t;S)$, $A$, $M$, $\kappa$, and $E>0$ be the
corresponding function and constants. For $t\ge0$ and $z\in\R$,
\[
 \begin{gathered}
 \supp\varphi(\cdot,t;S)\subset[S+ct-A,S+ct+A],\qquad
 0\le\varphi\le M,\\
 \varphi(S+ct,t;S)\ge\kappa,
 \qquad
 (\partial_t-D\partial_{zz}-D\beta\partial_z-a)\varphi
 =-E\varphi\quad\text{where }\varphi(z,t;S)>0.
 \end{gathered}
\]
Choose $p>1$ with $(p-1)a_{\max}\le E/4$, and then choose
$S>A$ so that $D\nu e^{-2(S-A)}\le E/4$. For $\delta>0$,
put $v(z,y,t)=\delta\varphi(z,t;S)^p\psi(y)$ on
$\R\times \mathcal B\times[0,\infty)$. At points where $t>0$ and
$v(z,y,t)>0$, \eqref{eq:horo-metric} and \eqref{eq:uniform-taylor} give
\[
 \begin{aligned}
 v_t-D\Delta v-F(z,v)
 &\le \bigl((p-1)a(z)-pE+D\nu e^{-2z}+C_2v\bigr)v\\
 &\hspace{1cm}
 -D\delta p(p-1)\psi\varphi^{p-2}|\varphi_z|^2
 \le-\frac E4v, 
 \end{aligned}
\]
provided
$\delta M^p\le\min\{1,E/(4(C_2+1))\}$.
Here $z\ge S+ct-A\ge S-A$ on the support of $v$.
The power $p>1$ makes the first derivative in $z$ continuous
at the zeros of $\varphi$, and its second derivative is locally
integrable. Thus the inequality holds weakly throughout
$\R\times \mathcal B$, with $v=0$ on $\R\times\partial \mathcal B$.

Fix $T_0>0$. Strict positivity and continuity of
$u(\cdot,T_0)$ on the compact set
$[S-A,S+A]\times\overline{\mathcal B}$ allow a further decrease of
$\delta$ such that $v(z,y,0)\le u(z,y,T_0)$ there.
Outside this set $v(z,y,0)=0$. Parabolic comparison on
$\R\times \mathcal B$ therefore gives $u(z,y,T_0+t)\ge v(z,y,t)$.
One may first compare on $(-N,N)\times \mathcal B$ for each fixed
finite time interval, taking $N$ large enough to contain the
support of $v$ during that interval. Evaluating at $z=S+ct$
proves the assertion with $A_0=\delta\kappa^p$.
\end{proof}

\begin{lemma}\label{lem:horo-periodic-step-comparison}
Assume the hypotheses of Lemma~\ref{lem:horo-subcritical-lower},
and retain $c,S,T_0,A_0,\mathcal B,\psi,\nu$ from its proof.
Choose $0<\varepsilon\le\min\{1,A_0\}$ so that
$C_2\varepsilon\le a_{\min}$. Let $V$ solve
\[
 \left\{\begin{aligned}
 V_t&=DV_{zz}+D\beta V_z+a(z)V(1-V)
       &&\text{for }z\in\R,\ t>0,\\
 V(z,0)&=\one_{(-\infty,S]}(z)&&\text{for }z\in\R.
 \end{aligned}\right.
\]
Define
\[
 q(t)=\exp\left(-D\nu\int_0^t e^{-2(S+c\tau)}\,\dd\tau\right)
 \quad\text{for }t\ge0.
\]
Then
\[
 u(z,y,T_0+t)\ge\varepsilon q(t)\psi(y)V(z,t)
 \quad\text{for }t\ge0,\ z>S+ct,\text{ and }y\in \mathcal B,
\]
and $q(t)\ge q_\infty:=\exp(-D\nu e^{-2S}/(2c))>0$
for $t\ge0$.
\end{lemma}

\begin{proof}
Put $w(z,t)=\varepsilon V(z,t)$ and
$v(z,y,t)=q(t)\psi(y)w(z,t)$. Comparison gives
$0\le w\le\varepsilon$, and the equation for $V$ yields
\[
 w_t-Dw_{zz}-D\beta w_z=a(z)w-\frac{a(z)}{\varepsilon}w^2
 \quad\text{for }z\in\R\text{ and }t>0.
\]
Since $q'/q=-D\nu e^{-2(S+ct)}$, it follows from
\eqref{eq:horo-metric} and \eqref{eq:uniform-taylor} that, for $t>0$, $z>S+ct$ and $y\in\mathcal B$,
\[
 \begin{aligned}
 v_t-D\Delta v-F(z,v)
 &\le D\nu\bigl(e^{-2z}-e^{-2(S+ct)}\bigr)v\\
 &\hspace{1cm}
 -q\psi\left(\frac{a(z)}{\varepsilon}-C_2q\psi\right)w^2
 \le0.
 \end{aligned}
\]
On $z=S+ct$,
Lemma~\ref{lem:horo-subcritical-lower} gives
$v(S+ct,y,t)\le\varepsilon\psi(y)\le A_0\psi(y)
\le u(S+ct,y,T_0+t)$ for $t>0$ and $y\in \mathcal B$.
Moreover, $v=0$ on $y\in\partial \mathcal B$, and $v(z,y,0)=0$
when $z>S$. Comparison on the moving half-cylinder therefore
gives $u(z,y,T_0+t)\ge v(z,y,t)$ for $t>0$, $z>S+ct$, and $y\in\mathcal B$. To justify comparison at infinity,
put $\zeta=z-S-ct$ and let $K_F$ be a Lipschitz constant for
$F(z,\cdot)$ on $[0,1]$, uniform in $z$. On each finite
time interval, apply the maximum principle to
$v-u-\rho e^{Kt}(1+\zeta^2)$ on $0<\zeta<N$, where
$\rho>0$ and $K>K_F+2D+|D\beta+c|$.
Indeed,
\[
 \begin{aligned}
 &(\partial_t-D\partial_{zz}-D\beta\partial_z-K_F)
       \bigl[e^{Kt}(1+\zeta^2)\bigr]\\
 &\qquad=e^{Kt}\bigl[(K-K_F)(1+\zeta^2)
                    -2(D\beta+c)\zeta-2D\bigr]\ge0
 \quad\text{for }\zeta\ge0.
 \end{aligned}
\]
The added function dominates the bounded difference at $\zeta=N$
for large $N$. Letting $N\to\infty$ and then $\rho\downarrow0$
gives the comparison on the full half-cylinder.
Finally, integrating the definition of $q$ gives
$q(t)\ge\exp(-D\nu e^{-2S}/(2c))$.
\end{proof}

\begin{proof}[Proof of Theorem~\ref{thm:horo-compact-periodic}]
Fix $\theta\in(0,1)$. Linear comparison gives
$u(\cdot,t)\le e^{a_{\max}t}e^{tD\Delta}u_0$.
For each $t>0$, the hyperbolic heat-kernel bound and compact
support of $u_0$ imply $u(x,t)\to0$ as
$\dist(o,x)\to\infty$. Since $u(\cdot,t)$ is continuous,
its $\theta$-superlevel set is compact. By
Theorem~\ref{thm:threshold}, $u(o,t)\to1$ as $t\to\infty$,
so the set is nonempty for all sufficiently large $t$.

Choose $B\in\R$ so that $\supp u_0\subset\{b\le B\}$.
By Lemma~\ref{lem:periodic-busemann-upper},
\[
 u(x,t)<\theta
 \quad\text{for }b(x)\ge m_*(t)+C_\theta^+
 \text{ and }t\ge T_\theta^+,
\]
after choosing $C_\theta^+,T_\theta^+$. This proves the upper
bound for $H_\theta^b(t)$.

For the lower bound, fix $c\in(0,\cp)$ and use
Lemma~\ref{lem:horo-periodic-step-comparison}. Since
$\psi>0$ on $\mathcal B$, the constant
$\psi_0=\min_{|y|\le1/2}\psi(y)$ is positive. By
Lemma~\ref{lem:periodic-logistic-step},
there are $C_0,T_1>0$ such that
\[
 V(z,t)\ge\frac12
 \quad\text{for }z\le m_*(t)-C_0\text{ and }t\ge T_1.
\]
Consequently, with
$\eta=\varepsilon q_\infty\psi_0/2>0$,
\[\label{eq:horo-periodic-strip-lower}
 u(z,y,T_0+t)\ge\eta
 \quad\text{for }t\ge T_1,
 \ S+ct<z\le m_*(t)-C_0,\text{ and }|y|\le\frac12.
\]

Apply Lemma~\ref{lem:plateau-amplification} with initial lower
bound $\eta$ and error $1-\theta$, and denote the resulting
radius and time by $R$ and $T$. Put
$z_t=m_*(t)-C_0-R$ and $x_t=(z_t,0)$.
The upper-half-space distance formula gives
\[
 B(x_t,R)\subset
 \{(z,y):|z-z_t|\le R,\ |y|\le e^{-z_t}\sinh R\}
 \quad\text{for }t\ge0.
\]
As $\cp>c>0$, we have $z_t-R>S+ct$ and
$e^{-z_t}\sinh R\le1/2$ for all sufficiently large $t$.
The lower bound \eqref{eq:horo-periodic-strip-lower} therefore
holds throughout $B(x_t,R)$, and
Lemma~\ref{lem:plateau-amplification} yields
\[\label{eq:busemann-amplified-point}
 u(x_t,T_0+t+T)\ge\theta,
 \qquad
 H_\theta^b(T_0+t+T)\ge m_*(t)-C_0-R
 \quad\text{for all sufficiently large }t.
\]
For the fixed number $D_0=T_0+T$, the definition of $m_*$
gives
\[
 |m_*(t+D_0)-m_*(t)|
 \le\cp D_0+\frac{3}{2\lam}\log(1+D_0)
 \quad\text{for }t\ge0.
\]
Together with the upper bound, this gives the constants $\tilde{C}_3,\tilde{T}_3$ for $H_\theta^b$. We use the pointwise lower bound
at $x_t$ to prove~\eqref{eq:busemann-window}.

Fix $A,R$ and choose a nonzero smooth compactly supported initial function
$0\le f_0\le1/2$. Let $v$ solve~\eqref{eq:KPP} with $f=f_b$ and datum $f_0$.
Applying~\eqref{eq:busemann-amplified-point} to $v$ with $d_0=C_0+R_0$ and
$H=T_0+T$ (here $R_0$ denotes the amplification radius for $v$) gives $s_0>0$ such that
\[\label{eq:busemann-axis-point}
 v(X(m_*(s)-d_0,0),s+H)\ge\theta\qquad(s\ge s_0).
\]
For $a\in\R^{d-1}$ the isometry $P_aX(z,y)=X(z,y+a)$ preserves $b$ and
the reaction. The union of $P_a(\supp f_0)$ over $|a|\le R$ is compact,
as the image of $\supp f_0\times\{|a|\le R\}$ under a continuous map.
By Theorem~\ref{thm:threshold}, $u(\cdot,\tau)\ge1/2$ on that compact
set for every $\tau\ge\tau_0$, for some $\tau_0$. Thus
$u(\cdot,\tau)\ge f_0\circ P_a^{-1}$ for every such $a,\tau$.
By comparison and~\eqref{eq:busemann-axis-point},
\[
 u(X(m_*(s)-d_0,a),\tau+s+H)\ge\theta
 \qquad(s\ge s_0,\ \tau\ge\tau_0,\ |a|\le R).
\]
Increase $s_0$ until $m_*$ is increasing on $[s_0,\infty)$ and
$Z=m_*(s_0)-d_0\ge A$; this is possible because $\cp>0$.
At time $t$ take $s\in[s_0,t-H-\tau_0]$ and $\tau=t-H-s$.
The attained heights fill $[Z,m_*(t-H-\tau_0)-d_0]$.
Since $m_*(t)-m_*(t-H-\tau_0)\le\cp(H+\tau_0)$,
this interval contains $[Z,m_*(t)-d_0-\cp(H+\tau_0)]$ for large $t$.
On the compact set $A\le z\le Z$, $|y|\le R$, local convergence to one
gives $u\ge\theta$ for all sufficiently large times. Together with
\eqref{eq:periodic-busemann-upper} this proves~\eqref{eq:busemann-window}.

For $x\ne o$, write $x=\exp_o(r\omega)$ with $r=\dist(o,x)$ and $g_o(\omega,\omega)=1$. Put $\alpha=\arccos g_o(\omega,e)$.
In the hyperboloid model the time and axial components of $x$ are
$\cosh r$ and $\sinh r\cos\alpha$. Passing to the upper half-space gives
\[
 q=(\cosh r-\sinh r\cos\alpha)^{-1},\qquad
 |y|=q\sinh r\sin\alpha .
\]
Thus $b(x)=r+\log(\sin^2(\alpha/2)+e^{-2r}\cos^2(\alpha/2))$.
For $x\in\mathcal C_\delta$ this yields
$r+A_\delta\le b(x)\le r$ and $|y|\le R_\delta$, where
$A_\delta=\log\sin^2(\delta/2)<0$ and $R_\delta=\cot(\delta/2)>0$:
for $r>0$ use $|y|\le\sin\alpha/(1-\cos\alpha)$, and at $o$ use $b(o)=0$ and $y=0$.
Apply~\eqref{eq:busemann-window} with $A=A_\delta$, $R=R_\delta$ and
denote its position constant by $C_W$. On the sector,
$r\le m_*(t)-C_W$ implies $b(x)\le m_*(t)-C_W$, and
$b(x)\le m_*(t)+C_W$ implies $r\le m_*(t)+C_W-A_\delta$.
Taking $\tilde{C}_2=C_W-A_\delta$ gives~\eqref{eq:busemann-sector}.

For part~(ii), compare $u$ with the solutions of~\eqref{eq:line} from
$\eta\one_{(-\infty,A_-]}$ and $\one_{(-\infty,A_+]}$, lifted by $z=b(x)$.
Denote these scalar solutions by $u_-$ and $u_+$.
By~\eqref{eq:periodic-busemann-upper} the upper comparison gives the outer inclusion.
For the lower comparison let $W$ solve~\eqref{eq:periodic-logistic-step-equation}
with $S=A_-$. Choose $0<\varepsilon\le\min\{\eta,1/2\}$
with $C_2\varepsilon\le a_{\min}$.
By~\eqref{eq:uniform-taylor}, $\varepsilon W$ is a subsolution since
$a(z)\varepsilon W(1-W)\le F(z,\varepsilon W)$.
Hence~\eqref{eq:periodic-logistic-step-lower} gives
$u_-(z,t)\ge\varepsilon/2$ for $z\le m_*(t)-C_0$ and all large $t$.
Let $r_0,\tau_a$ be the interval radius and time in
Lemma~\ref{lem:plateau-amplification} for initial level $\varepsilon/2$
and target $\theta$. The lemma gives $u_-(z,t+\tau_a)\ge\theta$ for
$z\le m_*(t)-C_0-r_0$. Since $m_*(t+\tau_a)-m_*(t)\le\cp\tau_a$,
this proves the inner inclusion in~\eqref{eq:busemann-global}.
Both scalar comparison data are fixed by $\eta,A_-,A_+$.
\end{proof}

\begin{proof}[Proof of Theorem~\ref{thm:periodic-bramson-distance}]
For assertion \textup{(i)}, choose $R_+>0$ so that
$u_0\le\one_{B(o,R_+)}$ almost everywhere, and let $U_+$ solve
the radial equation with $U_+(r,0)=\one_{[0,R_+]}(r)$.
By comparison,
\[
 u(x,t)\le U_+(r(x),t)
 \quad\text{for }x\in\Hh^d,\quad t>0.
\]
Since $u_0$ is nonzero, strict positivity at time one and
continuity on the compact ball $B(o,1)$ give $\eta_->0$ such
that $u(x,1)\ge\eta_-$ for $x\in B(o,1)$.
Let $U_-$ solve the same radial equation from
$U_-(r,0)=\eta_-\one_{[0,1]}(r)$.
The comparison principle then gives
\[
 U_-(r(x),t-1)\le u(x,t)\le U_+(r(x),t)
 \quad\text{for }x\in\Hh^d,\quad t\ge1.
\]
Apply Lemma~\ref{lem:periodic-distance-transfer} to $U_-$ and
$U_+$. Since $m_*(t)-m_*(t-1)$ is bounded for $t\ge2$, the
resulting inner and outer bounds give
\eqref{eq:periodic-bramson-balls}, after increasing
$\tilde{C}_{\theta,1},\tilde{T}_{\theta,1}$.

For assertion \textup{(ii)}, let $U_-$ and $U_+$ solve the
distance-from-$\G$ equation with respective initial profiles
$\eta\one_{[0,R_1]}$ and $\one_{[0,R_2]}$.
The condition $R_1>0$ makes the lower datum nonzero on a set of
positive hyperbolic volume. By the geometric reduction and the
assumed initial inequalities,
\[
 U_-(s(x),t)\le u(x,t)\le U_+(s(x),t)
 \quad\text{for }x\in\Hh^d,\quad t>0.
\]
Applying Lemma~\ref{lem:periodic-distance-transfer} to both
profiles proves \eqref{eq:periodic-bramson-tubes}.
\end{proof}

\paragraph{A common derivative at zero.}
We now consider reactions $f$ satisfying
\eqref{eq:common-assumptions}, the spatial H\"older condition
\eqref{eq:common-spatial-regularity}, and
$\sup_{x\in\Hh^d}\|f(x,\cdot)\|_{C^2([0,1])}<\infty$.
The common derivative is $f_u(x,0)=a_0>0$.
We compare \eqref{eq:KPP} with reactions $g$ depending only on
the population variable and satisfying $g'(0)=a_0$.

\begin{lemma}\label{lem:envelopes}
Assume \eqref{eq:common-assumptions},
\eqref{eq:common-spatial-regularity}, and
$\sup_{x\in\Hh^d}\|f(x,\cdot)\|_{C^2([0,1])}<\infty$.
Then there are $C^2$ functions $g_-,g_+:[0,1]\to\R$ such that,
for either choice of sign,
\[\label{eq:envelope-properties}
 \begin{gathered}
 g_\pm(0)=g_\pm(1)=0,\quad g_\pm'(0)=a_0,\quad g_\pm'(1)<0,\\
 0<g_\pm(v)\le a_0v\ (0<v<1),\quad g_\pm'(v)\le a_0\ (0\le v\le1),\\
 g_-(v)\le f(x,v)\le g_+(v)\quad\text{for }x\in\Hh^d,\ 0\le v\le1.
 \end{gathered}
\]
\end{lemma}

\begin{proof}
Uniformly in $x$, Taylor's formula gives $f(x,v)=a_0v+O(v^2)$ near zero and
$f(x,v)=-f_u(x,1)(1-v)+O((1-v)^2)$ near one.
Since $a_0>0$, $\sup_x f_u(x,1)<0$, and $f$ is uniformly positive on each closed
subinterval of $(0,1)$ by~\eqref{eq:common-assumptions}, there are $c_0,C>0$ and $\delta\in(0,1)$ such that
\[\label{eq:reaction-ratio-bounds}
 \begin{aligned}
 \frac{f(x,v)}{v(1-v)}&\ge c_0
       &&\text{for }x\in\Hh^d,\ 0<v<1,\\
 \frac{f(x,v)}{v(1-v)}&\ge a_0-Cv
       &&\text{for }x\in\Hh^d,\ 0<v\le\delta.
 \end{aligned}
\]
Put $g_-(v)=a_0v(1-v)e^{-Mv}$ and decrease $\delta$ so that
$Cv/a_0\le1/2$ on $[0,\delta]$.
Since $\log(1-Cv/a_0)\ge-2Cv/a_0$, taking $M\ge2C/a_0$ gives
$a_0e^{-Mv}\le a_0-Cv$ there. Increasing $M$ so that
$a_0e^{-M\delta}\le c_0$ treats $v\ge\delta$.
Equation \eqref{eq:reaction-ratio-bounds} gives
$g_-(v)\le f(x,v)$ for $x\in\Hh^d$ and $0\le v\le1$.
Differentiation gives
\[
 g_-'(v)=a_0e^{-Mv}[1-2v-Mv(1-v)]\le a_0
 \quad\text{for }0\le v\le1.
\]
The formula for $g_-$ also gives $g_-'(0)=a_0$, $g_-'(1)<0$
and $g_-(v)>0$ for $0<v<1$.

To construct $g_+$, uniform Lipschitz continuity and
$f(x,1)=0$ give $f(x,v)\le K_0(1-v)$.
Choose $K>K_0$ and define the concave piecewise affine function
$g(v)=\min\{a_0v,K(1-v)\}$ on $\R$.
The bound $f(x,v)\le a_0v$ in~\eqref{eq:common-assumptions}
and the choice $K>K_0$ give $g(v)\ge f(x,v)$ for
$x\in\Hh^d$ and $v\in[0,1]$.
The corner of $g$ is $v_*=K/(a_0+K)\in(0,1)$.
Since $g(v_*)-K_0(1-v_*)>0$, choose a closed interval about
$v_*$ contained in $(0,1)$ on which this difference is at
least $d_*>0$.
Choose a nonnegative, smooth, even function $\chi_h$ supported
in $[-h,h]$, with $\int_\R\chi_h{(z)\dd z}=1$, and form $g*\chi_h$.
Take $h$ so that $[v_*-h,v_*+h]$ lies in the interval where
$g(v)-K_0(1-v)\ge d_*$ and $h\max\{a_0,K\}<d_*/2$.
The convolution is concave, retains the affine pieces away from
the corner, and differs from $g$ by at most $h\max\{a_0,K\}$;
hence it remains above $f$ on $[0,1]$.
Set $g_+=g*\chi_h$. Then $g_+=a_0v$ near zero and
$g_+=K(1-v)$ near one. Concavity gives $g_+'\le a_0$ and
$g_+\le a_0v$, while a sufficiently small smoothing interval
preserves positivity throughout $(0,1)$. These properties prove
\eqref{eq:envelope-properties}.
\end{proof}

\begin{lemma}\label{lem:homogeneous-input}
Let $g\in C^2([0,1])$ satisfy the first two lines of~\eqref{eq:envelope-properties}. Let $W$ be the solution of
\[\label{eq:step-comparison-equation}
 \left\{\begin{aligned}
 W_t&=DW_{zz}+D\beta W_z+g(W)&&\text{for }z\in\R,\ t>0,\\
 W(z,0)&=\one_{(-\infty,0]}(z)&&\text{for }z\in\R.
 \end{aligned}\right.
\]
For every $t>0$, $W(\cdot,t)$ decreases strictly from one to zero.
For $0<\theta<1$, its unique level position $X_\theta(t)$, defined by
$W(X_\theta(t),t)=\theta$, satisfies
\[\label{eq:step-comparison-position}
 X_\theta(t)=m_{a_0}(t)+O_\theta(1)\quad\text{as }t\to\infty,
\]
where $m_{a_0}$ is defined in~\eqref{eq:common-position}. The estimate holds for every $a_0>0$, including
$c_{a_0}>0$, $c_{a_0}=0$ and $c_{a_0}<0$.
\end{lemma}

\begin{proof}
By translation comparison in~\eqref{eq:step-comparison-equation},
$W_z\le0$. The equation
$(W_z)_t=D(W_z)_{zz}+D\beta(W_z)_z+g'(W)W_z$ and the strong
maximum principle give $W_z<0$ at positive times. Heat comparison
for $W$ and $1-W$ gives the end limits zero and one, so every
level in $(0,1)$ is attained uniquely.

Set $V(y,t)=W(y-D\beta t,t)$.
By \eqref{eq:step-comparison-equation},
$V_t=DV_{yy}+g(V)$ on $\R\times(0,\infty)$ and
$V(y,0)=\one_{(-\infty,0]}(y)$.
With $\tau=Dt$ and $\widetilde V(y,\tau)=V(y,\tau/D)$, the
equation becomes $\widetilde V_\tau=\widetilde V_{yy}+\widetilde g(\widetilde V)$,
where $\widetilde g=g/D$ and $\widetilde g'(0)=a_0/D$.
Comparison with a solution of $v_t=Dv_{yy}+g(v)$ whose initial
datum $v_0\not\equiv0$ satisfies
$0\le v_0\le\one_{(-\infty,0]}$ and has compact support,
together with homogeneous spreading \cite{AW75} in this normalised
time variable, gives
$V(0,t)\to1$ as $t\to\infty$.
For $0<\theta<1$, let $\widehat X_\theta(t)$ be the unique point
satisfying $V(\widehat X_\theta(t),t)=\theta$. Since $V(\cdot,t)$ is
strictly decreasing and $V(0,t)\to1$, there is $T_\theta$ such that
$\widehat X_\theta(t)>0$ for every $t\ge T_\theta$.
Applying \cite[Theorem~1.1]{HNRR13} to $\widetilde V$ with the
$C^2$ reaction $\widetilde g$, satisfying
$\widetilde g(0)=\widetilde g(1)=0$,
$\widetilde g'(1)<0<\widetilde g'(0)=a_0/D$ and
$0<\widetilde g(v)\le(a_0/D)v$, gives the level position
$2\sqrt{a_0/D}\,\tau-3\log\tau/(2\sqrt{a_0/D})+O_\theta(1)$;
see also \cite{Bramson83}. The initial datum of $\widetilde V$
is bounded, nonzero and zero for $y>0$, and
$\widetilde g'\le a_0/D$.
Returning to $t$, this position is
$2\sqrt{Da_0}\,t-3\log t/(2\lambda_{a_0})+O_\theta(1)$,
since $\lambda_{a_0}=\sqrt{a_0/D}$ and $D$ is fixed.
The identity $V(y,t)=W(y-D\beta t,t)$ gives $X_\theta(t)=\widehat X_\theta(t)-D\beta t$.
Since $0\le\log(1+t)-\log t\le1/t$ for $t\ge1$,
these identities prove \eqref{eq:step-comparison-position}.

For \eqref{eq:line}, this change of variables replaces the
reaction by $F(y-D\beta t,V)$.
For a distance problem in~\eqref{eq:three-reductions}, write
$U_t=DU_{rr}+DB(r)U_r+F(r,U)$ on $r>0$ and set
$\widetilde U(y,t)=U(y-D\beta t,t)$.
Direct differentiation gives
\[\label{eq:translated-distance-domain}
 \widetilde U_t=D\widetilde U_{yy}
   +D[B(y-D\beta t)-\beta]\widetilde U_y
   +F(y-D\beta t,\widetilde U),\quad\text{for }t>0,\ y>D\beta t.
\]
The condition at $r=0$ becomes
$\widetilde U_y(D\beta t,t)=0$.
The distance comparisons below are formulated in the original
variable $r>0$.
\end{proof}

Throughout the comparisons below, $g\in C^2([0,1])$ satisfies
the first two lines of~\eqref{eq:envelope-properties}.
Let $\Phi_g$ denote the decreasing travelling-wave profile of speed
$2\sqrt{Da_0}$ for $V_t=DV_{yy}+g(V)$, normalised by $\Phi_g(0)=1/2$, so that
\[\label{eq:homogeneous-wave}
 \left\{\begin{aligned}
 &D\Phi_g''+2\sqrt{Da_0}\,\Phi_g'+g(\Phi_g)=0
       &&\text{for }z\in\R,\\
 &\Phi_g(-\infty)=1,\quad\Phi_g(+\infty)=0.
 \end{aligned}\right.
\]
Substituting $e^{-\lambda z}$ into the linearisation of~\eqref{eq:homogeneous-wave} at zero gives the polynomial
$D\lambda^2-2\sqrt{Da_0}\lambda+a_0=D(\lambda-\lambda_{a_0})^2$.
Thus its double root is the exponent $\lambda_{a_0}$ used in~\eqref{eq:common-position}.
For an even, compactly supported measurable function
$U_{h,0}:\R\to[0,1]$, let $U_h$ solve
\[\label{eq:hyperplane-equation}
 (U_h)_t=D(U_h)_{\rho\rho}+D\beta\tanh\rho\,(U_h)_\rho+g(U_h),
 \quad\text{for }\rho\in\R,\ t>0,
\]
with $U_h(\rho,0)=U_{h,0}(\rho)$. Evenness gives
$(U_h)_\rho(0,t)=0$ for $t>0$. For a prescribed measurable
$U_{e,0}:[0,\infty)\to[0,1]$, let $U_e(r,t)$ solve the radial
problem in~\eqref{eq:three-reductions} with $F(r,v)=g(v)$,
$(U_e)_r(0,t)=0$ and $U_e(r,0)=U_{e,0}(r)$. In
Lemma~\ref{lem:tube-sandwich}, $U_e$ uses $g_-$ and $U_h$ uses
$g_+$ from Lemma~\ref{lem:envelopes}. The composition $U_e(s(x),t)$
evaluates this radial profile at the distance $s(x)$ from $\G$.
Lemma~\ref{lem:tube-sandwich} uses the drift ordering
$B_d(s)\le\beta\coth s$ to make this composition a subsolution.

For $B(r)=\beta\coth r$, $B(r)=\beta\tanh r$, or $B(r)=B_d(r)$,
$|B(r)-\beta|+|B'(r)|\le Ce^{-2r}$ on $r\ge1$. If $c>0$,
$\kappa\in\R$, $0<\gamma<1$, and $m(t)=ct-\kappa\log(1+t)$, then
$m(t)-t^\gamma\ge ct/2\ge1$ for all sufficiently large $t$.
Consequently
\[\label{eq:geometric-error-front-strip}
 \sup_{|\xi|\le t^\gamma}
 \bigl(|B(m(t)+\xi)-\beta|+|B'(m(t)+\xi)|\bigr)
 \le Ce^{-ct}\quad\text{for }t\ge T.
\]
In particular, the left-hand side of~\eqref{eq:geometric-error-front-strip} is integrable over $t\in[T,\infty)$; see also
\cite[formula~(5.3) and Lemma~6.1]{GGQ}. For the Busemann coordinate,
the drift equals $D\beta$ at every speed.

\begin{lemma}\label{lem:tube-sandwich}
Assume that $U_{e,0}$ and $U_{h,0}|_{[0,\infty)}$ are nonincreasing,
with $U_{h,0}$ even. Let $U_e$ solve the radial problem in~\eqref{eq:three-reductions} with $F(r,v)=g_-(v)$ and
$U_e(r,0)=U_{e,0}(r)$, and let $U_h$ solve \eqref{eq:hyperplane-equation}
with $g=g_+$ and $U_h(\rho,0)=U_{h,0}(\rho)$. Then
$U_e(s(x),t)$ is a subsolution and $U_h(s(x),t)$ is a supersolution of~\eqref{eq:KPP} with any reaction satisfying
$g_-\le f(x,\cdot)\le g_+$.
\end{lemma}

\begin{proof}
First take smooth nonincreasing initial profiles with bounded
derivatives and smooth reactions. The general bounded monotone
profiles and $C^2$ reactions follow by approximation preserving
the reaction endpoint values and a common Lipschitz bound.
For the radial solution $U_e$ in~\eqref{eq:three-reductions}
with reaction $g_-$, put $Z(r,t)=(U_e)_r(r,t)/\sinh r$ for
$r>0$. Radial regularity gives its even extension with
$Z(0,t)=(U_e)_{rr}(0,t)$, and differentiation yields
\[\label{eq:radial-derivative-lift}
 Z_t=DZ_{rr}+D(d+1)\coth r\,Z_r+
                  \bigl(Dd+g_-'(U_e(r,t))\bigr)Z
 \quad {\text{for }r>0,\quad t>0}.
\]
By \eqref{eq:polar-metric}, this is the radial form of a linear
equation on $\Hh^{d+2}$ with bounded zeroth-order coefficient.
On finite time intervals, radial regularity bounds $Z(r,t)$ for
$r$ near $0$, and boundedness of $(U_e)_r$ bounds $Z$ for $r\ge1$. Since its initial value is nonpositive, the comparison
argument of Lemma~\ref{lem:wellposed} in dimension $d+2$ gives
$Z\le0$, and hence $(U_e)_r\le0$.
For the even solution $U_h$ of~\eqref{eq:hyperplane-equation}
with $g=g_+$, put $Y=(U_h)_\rho$. Differentiation gives
\[\label{eq:hyperplane-derivative-equation}
 \left\{\begin{aligned}
 Y_t&=DY_{\rho\rho}+D\beta\tanh\rho\,Y_\rho\\
 &\quad+\bigl(g_+'(U_h(\rho,t))+D\beta\sech^2\rho\bigr)Y
       &&\text{for }\rho>0,\ t>0,\\
 Y(0,t)&=0&&\text{for }t>0.
 \end{aligned}\right.
\]
Half-line comparison in~\eqref{eq:hyperplane-derivative-equation}
gives $Y\le0$. Positive-time local convergence preserves both
monotonicity assertions under approximation.

By Lemma~\ref{lem:drifts},
$\beta\tanh s\le B_d(s)\le\beta\coth s$ for $s>0$.
In the next display, $(U_e)_s$ denotes $(U_e)_r(s,t)$.
The radial equation, $(U_e)_r\le0$ and $g_-\le f$ give the
following expression for
$(\partial_t-D\Delta)U_e(s(x),t)-f(x,U_e(s(x),t))$:
\[
 D[\beta\coth s-B_d(s)](U_e)_s+g_-(U_e)-f(x,U_e)\le0.
\]
By the equation for $U_h$ and the inequalities $g_+\ge f$ and
$\beta\tanh s\le B_d(s)$, we likewise obtain
\[
 (\partial_t-D\Delta)U_h(s(x),t)-f(x,U_h(s(x),t))\ge0
 \quad\text{for }s(x)>0,\ t>0.
\]
At $s=0$, both first derivatives vanish.
By \eqref{eq:three-reductions}, for a $C^2$ function $x\mapsto U(s(x),t)$,
the Laplacian at $s=0$ equals $(d-1)U_{ss}(0,t)$. The diffusion
terms there are therefore $D(d-1)U_{ss}(0,t)$, $DdU_{ss}(0,t)$
and $DU_{ss}(0,t)$ in the tube, radial and hyperplane equations, respectively. Monotonicity makes these second
derivatives nonpositive, so the same residual inequalities hold
on $\G$. Thus $U_e(s(x),t)$ and $U_h(s(x),t)$ are, respectively, a weak
subsolution and a weak supersolution on $\Hh^d\times(0,\infty)$.
\end{proof}

\begin{lemma}\label{lem:step-slope}
Let $W$ solve \eqref{eq:step-comparison-equation}, and let
$X(t)=X_{1/2}(t)$. For each $0<\delta<1$, there is $\kappa_\delta>0$
such that
\[\label{eq:step-slope}
 -W_z(z,t)\ge\kappa_\delta W(z,t)
 \quad\hbox{if }W(z,t)\le1-\delta,\qquad t>0.
\]
There are also $C,\omega>0$, depending on $D$ and $g$, such that
\[\label{eq:step-left-tail}
 1-W(z,t)\le C e^{\omega(z-X(t))}
 \quad\text{for }z\in\R,\ t>0.
\]
\end{lemma}

\begin{proof}
The substitution $V(y,t)=W(y-D\beta t,t)$ gives
$W(z,t)=V(z+D\beta t,t)$ and $X(t)=\widehat X(t)-D\beta t$,
where $V(\widehat X(t),t)=1/2$. In particular, $W(X(t)+z,t)=V(\widehat X(t)+z,t)$ and
$W_z(X(t)+z,t)=V_y(\widehat X(t)+z,t)$. It therefore suffices
to prove both estimates for $\beta=0$, with constants depending
only on $D$ and $g$.
Choose $c_0>\max\{2\sqrt{Da_0},D\beta\}$ and a decreasing travelling-wave
profile $\Psi:\R\to(0,1)$ of speed $c_0$ for $V_t=DV_{yy}+g(V)$,
normalised by $\Psi(0)=1/2$. Thus
$D\Psi''+c_0\Psi'+g(\Psi)=0$, with $\Psi(-\infty)=1$ and $\Psi(+\infty)=0$.
We first establish its endpoint asymptotics
\[\label{eq:comparison-wave-tails}
 \begin{gathered}
 \Psi(z)\sim C_+e^{-\lambda z}\quad\text{as }z\to\infty,\qquad
 1-\Psi(z)\sim C_-e^{\omega z}\quad\text{as }z\to-\infty,\qquad C_\pm>0,\\
 \lambda=\frac{c_0-\sqrt{c_0^2-4Da_0}}{2D}>0,\qquad
 \omega=\frac{\sqrt{c_0^2-4Dg'(1)}-c_0}{2D}>0.
 \end{gathered}
\]
By Lemmas~\ref{lem:front-barriers}--\ref{lem:supercritical-front},
applied with $F=g$ and coordinate speed $c_0-D\beta>0$,
the construction is invariant under
$(z,t)\mapsto(z+h,t+h/(c_0-D\beta))$ for every real $h$.
The entire solution obtained in~\eqref{eq:entire-front-limit}
therefore has the form $\Psi(z-(c_0-D\beta)t)$. As $z\to\infty$,
$\Psi(z)$ lies
between $C_+e^{-\lambda z}-C'e^{-\lambda' z}$ and
$C_+e^{-\lambda z}$ for some $\lambda'>\lambda$.
Putting $h=\Psi'+\lambda\Psi$ in the wave equation gives
$Dh'+(c_0-D\lambda)h=a_0\Psi-g(\Psi)=O(\Psi^2)$.
Variation of constants and $c_0-D\lambda>D\lambda$ yield
$h=o(e^{-\lambda z})$, so
$-\Psi'/\Psi\to\lambda$ as $z\to\infty$.

For the left tail, put $q=1-\Psi>0$, $P=q'/q>0$ and
$\kappa_1=-g'(1)>0$. The wave equation gives
\[\label{eq:wave-left-riccati}
 \begin{gathered}
 DP'=A(q)-c_0P-DP^2\quad\text{for }z\in\R,\\
 A(v)=\frac{g(1-v)}v\quad\text{for }0<v<1,
 \qquad A(v)=\kappa_1+O(v)\quad\text{as }v\downarrow0.
 \end{gathered}
\]
The number $\omega$ in~\eqref{eq:comparison-wave-tails} is the positive
root of $Dr^2+c_0r=\kappa_1$. For $0<\varepsilon<\kappa_1$ sufficiently small, set
$\omega_\pm(\varepsilon)=(\sqrt{c_0^2+4D(\kappa_1\pm\varepsilon)}-c_0)/(2D)$.
On a sufficiently negative half-line,
$\kappa_1-\varepsilon\le A(q)\le\kappa_1+\varepsilon$.
By \eqref{eq:wave-left-riccati}, a value $P(z)>\omega_+(\varepsilon)$
forces blow-up at a finite point to the left of $z$, while
$0<P(z)<\omega_-(\varepsilon)$ forces $P$ to reach zero at a finite point to the left of $z$.
Since $q,q'>0$ on $\R$, $P=q'/q$ is finite and positive at every finite point.
Hence $\omega_-(\varepsilon)\le P(z)\le\omega_+(\varepsilon)$ on that half-line.
Letting $\varepsilon\downarrow0$ therefore gives
$P(z)\to\omega$ as $z\to-\infty$ and, in particular,
$q(z)\le Ce^{\omega z/2}$ on a left half-line.
Set $R=P-\omega$ and $b_0=c_0/D+2\omega$.
By \eqref{eq:wave-left-riccati} and boundedness of $R$ at
$-\infty$,
\[\label{eq:wave-left-integral}
 R(z)=\int_{-\infty}^z e^{-b_0(z-y)}
                     \bigl[D^{-1}(A(q(y))-\kappa_1)-R(y)^2\bigr]\dd y.
\]
Choose $z_0$ with $|R(y)|\le b_0/2$ for $y\le z_0$.
For $z\le z_0$, set $S(z)=\sup_{y\le z}|R(y)|$. Since $q$ is increasing
and $|A(q)-\kappa_1|\le Cq$ near $q=0$, \eqref{eq:wave-left-integral}
gives $S(z)\le Cq(z)/(Db_0)+S(z)/2$, after decreasing $z_0$ if needed.
Hence $|R(z)|\le 2Cq(z)/(Db_0)$. Hence $R$ is integrable at $-\infty$,
and $(\log q-\omega z)'=R$ shows that
$e^{-\omega z}q(z)$ has a finite positive limit.
This proves the left asymptotic in~\eqref{eq:comparison-wave-tails},
with the derivative asymptotic following from $q'=Pq$;
see also the classical homogeneous wave theory in~\cite{AW75}.

For $n>1$, let $W_n$ solve
$(W_n)_t=D(W_n)_{zz}+D\beta(W_n)_z+g(W_n)$ on $\R\times(0,\infty)$
with $W_n(z,0)=\Psi(nz)$, and compare it with
$\Psi(z-(c_0-D\beta)t-\xi)$ for fixed $\xi\in\R$.
At time zero the difference has exactly one sign change, from
positive to negative. It solves a scalar linear parabolic equation
with diffusion coefficient $D$, drift $D\beta$ and reaction
coefficient bounded by $\max_{[0,1]}|g'|$.
The zero-number theorem \cite{Angenent} therefore preserves at
most one sign change on a bounded interval with fixed nonzero
endpoint signs.

To obtain such an interval, fix $T<\infty$.
By \eqref{eq:comparison-wave-tails} and linear heat comparison,
\[
 W_n(z,t)\le C_{n,T}e^{-n\lambda z},\qquad
 1-W_n(z,t)\le C_{n,T}e^{n\omega z}
 \quad\text{for }0\le t\le T.
\]
For the bound on $1-W_n$, the equation for $W_n$ gives
$(1-W_n)_t=D(1-W_n)_{zz}+D\beta(1-W_n)_z-g(W_n)$;
the function $v\mapsto-g(1-v)$ vanishes at $0$ and has bounded
difference quotient there. By \eqref{eq:comparison-wave-tails}, $W_n$ decays faster at
$+\infty$ than the translated $\Psi$, and $1-W_n$ decays faster
at $-\infty$ than the translated $1-\Psi$.
Thus $W_n(z,t)-\Psi(z-(c_0-D\beta)t-\xi)$ is positive for
sufficiently negative $z$ and negative for sufficiently positive $z$,
uniformly for $0\le t\le T$.
The interval zero-number theorem with these endpoint signs then gives
the whole-line sign-change bound.

As $n\to\infty$, $\Psi(nz)\to\one_{(-\infty,0]}(z)$
for almost every $z\in\R$. Equation \eqref{eq:line-data-continuity} and local
parabolic estimates consequently give convergence of $W_n$ to
$W$, including first spatial derivatives, at positive times.
For each $t>0$ and $\xi\in\R$, the function
$z\mapsto W(z,t)-\Psi(z-\xi)$ has at most one sign change, from
positive to negative: two strict changes would
persist at three fixed spatial points for large $n$, contradicting
the one-sign-change bound for $W_n(z,t)-\Psi(z-(c_0-D\beta)t-\eta)$ with $\eta=\xi-(c_0-D\beta)t$. At every fixed $t>0$, Gaussian bounds for $W$ and $1-W$,
compared with the exponential bounds in~\eqref{eq:comparison-wave-tails},
give $W-\Psi>0$ far to the left and $W-\Psi<0$ far to the right.
At a point where $W(z,t)=\Psi(z-\xi)$, a positive derivative
of the difference would create a negative-to-positive crossing,
contrary to these endpoint signs and the fact that
$W-\Psi$ has at most one sign change.
Thus $W_z(z,t)\le\Psi'(z-\xi)$, including tangential contact.
Since $\Psi^{-1}:(0,1)\to\R$ exists and $\xi$ is arbitrary,
\[\label{eq:step-profile-slope}
 -W_z(z,t)\ge-\Psi'\bigl(\Psi^{-1}(W(z,t))\bigr)\quad {\text{for }z\in\R,\quad t>0}.
\]
By \eqref{eq:comparison-wave-tails}, $-\Psi'/\Psi$ tends to
$\lambda>0$ as $z\to\infty$; its positivity on compact intervals
therefore gives a positive infimum where $0<\Psi\le1-\delta$.
Together with \eqref{eq:step-profile-slope}, this proves
\eqref{eq:step-slope}.

For the left-tail estimate, put
$\Xi(z,t)=\Psi^{-1}(W(z,t))$. Dividing
\eqref{eq:step-profile-slope} by $\Psi'(\Xi)<0$ gives
\[\label{eq:step-inverse-coordinate}
 \Xi_z(z,t)=\frac{W_z(z,t)}{\Psi'(\Xi(z,t))}\ge1,
 \qquad \Xi(X(t),t)=0\quad {\text{for }z\in\R,\quad t>0}.
\]
Integrating \eqref{eq:step-inverse-coordinate} from $z$ to
$X(t)$ shows that $\Xi(z,t)\le z-X(t)$ for $z\le X(t)$,
and hence $W(z,t)\ge\Psi(z-X(t))$.
The left tail in~\eqref{eq:comparison-wave-tails} now gives
\eqref{eq:step-left-tail} on this half-line. Increasing $C$
to at least one covers $z\ge X(t)$ and completes the proof.
\end{proof}

\begin{lemma}\label{lem:normal-fixed-time}
Let $U_e$ solve the radial problem in~\eqref{eq:three-reductions}
with $F(r,v)=g(v)$, and let $U_h$ solve
\eqref{eq:hyperplane-equation} with the same reaction $g$. Suppose their
respective initial profiles take values in $[0,1]$, are
compactly supported, and satisfy
$\int_0^\infty U_e(r,0)\sinh^\beta r\dd r>0$ and
$\int_\R U_h(\rho,0)\dd\rho>0$, with $U_h(\cdot,0)$ even. For every fixed
$T>0$ and for $U=U_e$ or $U=U_h|_{[0,\infty)}$,
\[\label{eq:normal-fixed-time}
 \lim_{r\to\infty}\frac{\log U(r,T)}{r^2}=-\frac1{4DT}
 \quad {\text{for }T>0}.
\]
The step solution $W$ in~\eqref{eq:step-comparison-equation} satisfies
the same limit with $U(r,T)$ replaced by $W(r,T)$.
\end{lemma}

\begin{proof}
For $U=U_e$, let $U^{\mathrm{heat}}$ solve the radial problem
in~\eqref{eq:three-reductions} with $F=0$ and initial profile
$U_e(r,0)$. For $U=U_h$, define $U^{\mathrm{heat}}$ by
\eqref{eq:hyperplane-equation} with $g=0$ and initial profile
$U_h(\rho,0)$.
The bound $0\le g(v)\le a_0v$ gives
$U^{\mathrm{heat}}(\cdot,T)\le U(\cdot,T)
\le e^{a_0T}U^{\mathrm{heat}}(\cdot,T)$ by comparison.
For the radial datum, choose $R$ so that
$\supp[U_e(r(\cdot),0)]\subset B(o,R)$.
If $r(x)=r>R$ and $y\in\supp[U_e(r(\cdot),0)]$, then
$r-R\le\dist(x,y)\le r+R$. Integrating \eqref{eq:heat-estimate}
against the nonnegative function $y\mapsto U_e(r(y),0)$, whose integral
is finite and positive, therefore yields
$\log U_e(r,T)=-r^2/(4DT)+O_T(r+\log(2+r))$.
Division by $r^2$ proves \eqref{eq:normal-fixed-time} for $U_e$.

For the hyperplane heat comparison, put $m(r)=\cosh^\beta r$
and $B(r)=\beta\tanh r$ on $\R$.
Lemma~\ref{lem:liouville} gives
\[
 m^{1/2}(\partial_{rr}+B\partial_r)m^{-1/2}
 =\partial_{rr}-V_B(r),\qquad
 V_B(r)=\frac12B'(r)+\frac14B(r)^2.
\]
The potential $V_B$ is bounded; consequently the Feynman--Kac
formula bounds the Lebesgue kernel of
$e^{DT(\partial_{rr}-V_B)}$ between
$e^{-DT\norm{V_B}_\infty}$ and $e^{DT\norm{V_B}_\infty}$
times the Gaussian kernel.
Multiplying by $m(r)^{-1/2}$ and integrating against the
nonzero compactly supported function $m^{1/2}U_h(\cdot,0)$
preserves the quadratic logarithmic term $-r^2/(4DT)$, since
$\log m(r)=\beta r+O(1)$ as $r\to\infty$.
The bounds $U^{\mathrm{heat}}\le U_h\le e^{a_0T}U^{\mathrm{heat}}$
at time $T$ therefore prove
\eqref{eq:normal-fixed-time} for $U_h$.

Finally, let $w$ solve $w_t=Dw_{zz}+D\beta w_z$ on
$\R\times(0,\infty)$ with $w(z,0)=\one_{(-\infty,0]}(z)$.
At $(r,T)$, $w(r,T)$ is
$(4\pi DT)^{-1/2}\int_{-\infty}^0
 e^{-(r+D\beta T-y)^2/(4DT)}\dd y$, whose logarithm is
$-(r+D\beta T)^2/(4DT)+O_T(\log r)$.
The solution $W(\cdot,T)$ lies between this Gaussian step
integral and $e^{a_0T}$ times that integral, proving
$\log W(r,T)/r^2\to-1/(4DT)$ as $r\to\infty$.
\end{proof}

\begin{lemma}\label{lem:normal-interior}
Assume $a_0>D\rho_d$. Let $U$ be either of the two solutions in
Lemma~\ref{lem:normal-fixed-time}, and assume additionally that its
initial profile is nonincreasing on $[0,\infty)$. For every
$0<c<c_{a_0}$,
\[\label{eq:normal-interior-convergence}
 \lim_{t\to\infty}\sup_{0\le r\le ct}|U(r,t)-1|=0
 \quad {\text{for }0<c<c_{a_0}}.
\]
For each fixed $0<c_1<c_{a_0}$, there are $C,\omega_1,T_1>0$ such that
\[\label{eq:normal-boundary-deficit}
 0\le1-U(c_1t,t)\le C e^{-\omega_1t}\quad\text{for }t\ge T_1.
\]
\end{lemma}

\begin{proof}
For $U_e$, Theorem~\ref{thm:periodic-speed} applied to the
reaction $F(z,v)=g(v)$ gives
\eqref{eq:normal-interior-convergence},
since $k(\lambda)=D\lambda^2-D\beta\lambda+a_0$ and
$\cp=c_{a_0}$.
For $U_h$, let $\widehat U_e$ solve the radial problem in~\eqref{eq:three-reductions} with $F(r,v)=g(v)$ and
$\widehat U_e(r,0)=U_h(r,0)$ for $r\ge0$.
By the monotonicity argument of Lemma~\ref{lem:tube-sandwich},
$(\widehat U_e)_r\le0$ for $t>0$; the inequality
$\beta\tanh r\le\beta\coth r$ makes it a subsolution of~\eqref{eq:hyperplane-equation} on the positive half-line.
At $r=0$, the subsolution inequality follows from
$(\widehat U_e)_r(0,t)=0$ and $(\widehat U_e)_{rr}(0,t)\le0$,
as in Lemma~\ref{lem:tube-sandwich}.
Thus comparison gives $\widehat U_e\le U_h$, proving
\eqref{eq:normal-interior-convergence} in both cases.

Choose $0<c_0<c_1<c_2<c_{a_0}$.
Since $g'(1)<0$, there are $\kappa,\delta>0$ with
$g(1-v)\ge\kappa v$ for $0\le v\le\delta$.
Set $z(r,t)=1-U(r,t)$, with $B(r)=\beta\coth r$ for $U=U_e$
and $B(r)=\beta\tanh r$ for $U=U_h|_{[0,\infty)}$. By
\eqref{eq:normal-interior-convergence}, $0\le z(r,t)\le\delta$ on
$c_0t\le r\le c_2t$ for all sufficiently large $t$, and hence
\[
 z_t-Dz_{rr}-DB(r)z_r+\kappa z\le0,
 \quad\text{for }c_0t<r<c_2t\text{ and all sufficiently large }t.
\]
Since $c_0>0$, choose $B_*>0$ such that
$|B(r)|\le B_*$ on $c_0t\le r\le c_2t$ for all sufficiently large $t$.
Choose $0<\alpha<\kappa$ and $\eta>0$ sufficiently small
that $\kappa-\eta(c_2+DB_*)-D\eta^2>0$.
Choose $T_0$ so that these bounds hold for all later times, and set
\[\label{eq:interior-deficit-barrier}
 Z(r,t)=e^{-\alpha(t-T_0)}
       +e^{-\eta(r-c_0t)}+e^{-\eta(c_2t-r)}
\]
Applying $\partial_t-D\partial_{rr}-DB(r)\partial_r+\kappa$
to each term in~\eqref{eq:interior-deficit-barrier} and dividing
by that term gives $\kappa-\alpha$,
$\kappa+\eta c_0-D\eta^2+D\eta B$ and
$\kappa-\eta c_2-D\eta^2-D\eta B$, respectively, all positive
by the choices above. The initial and lateral boundary inequalities are
$Z(r,T_0)\ge1\ge z(r,T_0)$ for $c_0T_0\le r\le c_2T_0$,
and $Z(c_jt,t)\ge1\ge z(c_jt,t)$ for $j\in\{0,2\}$
and $t\ge T_0$.
Comparison gives $z\le Z$, and evaluating
\eqref{eq:interior-deficit-barrier} at $r=c_1t$ proves
\eqref{eq:normal-boundary-deficit} with any positive rate at
most $\min\{\alpha,\eta(c_1-c_0),\eta(c_2-c_1)\}$.
\end{proof}

\begin{lemma}\label{lem:normal-shift-comparison}
Assume $a_0>D\rho_d$, and let $U$ satisfy the assumptions of
Lemma~\ref{lem:normal-interior}. Fix $0<c_1<c_{a_0}$. There are $T>0$,
$0<\tau_-<T<\tau_+$, $0<\varepsilon_0<1$, $\nu>0$, and bounded
nondecreasing $C^1$ functions $\zeta_\pm:[T,\infty)\to[0,\infty)$
such that, with $\varepsilon(t)=\varepsilon_0e^{-\nu(t-T)}$,
\[\label{eq:normal-shift-comparison}
 \begin{split}
 &(1-\varepsilon(t))
       W(r+\zeta_-(t),t-T+\tau_-)\le U(r,t)\\
 &\hspace{12mm}\le
 \min\{1,(1+\varepsilon(t))
       W(r-\zeta_+(t),t-T+\tau_+)\},
 \quad\text{for }r\ge c_1t,\ t\ge T.
 \end{split}
\]
Here $W$ is the solution with initial datum
$\one_{(-\infty,0]}$ of~\eqref{eq:step-comparison-equation} with the same reaction $g$ as $U$.
\end{lemma}

\begin{proof}
Take $B(r)=\beta\coth r$ for the radial solution of~\eqref{eq:three-reductions}, and $B(r)=\beta\tanh r$ for the
restriction of~\eqref{eq:hyperplane-equation} to $r\ge0$.
Thus $U_t=DU_{rr}+DB(r)U_r+g(U)$, and for $r\ge1$,
\[\label{eq:normal-drift-error}
 |B(r)-\beta|\le C_Be^{-2r}
 \quad {\text{for }r\ge1}.
\]
One may take $C_B=2\beta/(1-e^{-2})$ in both cases.
Choose $\kappa>0$ and $0<\delta<1/4$ so that
$g'(v)\le-\kappa$ on $[1-2\delta,1]$, and fix
$0<\varepsilon_0<\delta/2$.
Let $\omega_1$ and $\omega$ be the constants in~\eqref{eq:normal-boundary-deficit} and \eqref{eq:step-left-tail},
respectively, and choose
\[
 0<\nu<\min\{\kappa,\omega_1,\omega(c_{a_0}-c_1)\}.
\]
With $K_g=\sup_{[0,1]}|g'|$, take $K$ sufficiently large that
\[\label{eq:shift-gain-choice}
 (1-\varepsilon_0)K\kappa_\delta>\nu+2K_g.
\]
For initial shifts to be chosen below, define
\[\label{eq:shift-equations}
 \zeta_\pm'(t)=K\varepsilon(t)+DC_Be^{-2c_1t}
 \quad {\text{for }t\ge T}.
\]
Integrating \eqref{eq:shift-equations} gives
\[\label{eq:integrated-shifts}
 \begin{split}
 \zeta_\pm(t)-\zeta_\pm(T)
 ={}&\frac{K\varepsilon_0}{\nu}
                   \bigl(1-e^{-\nu(t-T)}\bigr)\\
 &+\frac{DC_B}{2c_1}
                   \bigl(e^{-2c_1T}-e^{-2c_1t}\bigr),\qquad t\ge T.
 \end{split}
\]
Hence, for each sign,
\[\label{eq:individual-shift-bound}
 0\le\zeta_\pm(t)-\zeta_\pm(T)
 \le \frac{K\varepsilon_0}{\nu}
             +\frac{DC_Be^{-2c_1T}}{2c_1}\quad\text{for }t\ge T.
\]
By \eqref{eq:shift-equations} and \eqref{eq:individual-shift-bound},
$\zeta_+(t)$ and $\zeta_-(t)$ converge to finite limits as
$t\to\infty$.

Take $T$ with $c_1T\ge1$ and first verify the differential
inequalities on $r>c_1t$.
For the upper bound in~\eqref{eq:normal-shift-comparison},
where $(1+\varepsilon(t))
W(r-\zeta_+(t),t-T+\tau_+)<1$, set
$w=W(r-\zeta_+(t),t-T+\tau_+)$ and $q=1+\varepsilon(t)$.
Using \eqref{eq:step-comparison-equation}, the expression $(\partial_t-D\partial_{rr}-DB(r)\partial_r)(qw)
-g(qw)$ is
\[\label{eq:upper-shift-residual}
 -\nu\varepsilon w
 -q\bigl[\zeta_+'+D(B(r)-\beta)\bigr]W_z
 +qg(w)-g(qw).
\]
For the lower bound in~\eqref{eq:normal-shift-comparison}, set
$w=W(r+\zeta_-(t),t-T+\tau_-)$ and $q=1-\varepsilon(t)$;
the same computation gives
\[\label{eq:lower-shift-residual}
 \nu\varepsilon w
 +q\bigl[\zeta_-'-D(B(r)-\beta)\bigr]W_z
 +qg(w)-g(qw).
\]
In~\eqref{eq:upper-shift-residual}, $W_z$ is evaluated at
$(r-\zeta_+(t),t-T+\tau_+)$; in~\eqref{eq:lower-shift-residual},
it is evaluated at $(r+\zeta_-(t),t-T+\tau_-)$.
By \eqref{eq:normal-drift-error} and \eqref{eq:shift-equations},
both $\zeta_+'+D(B(r)-\beta)$ and $\zeta_-'-D(B(r)-\beta)$
are at least $K\varepsilon(t)$ on $r\ge c_1t$.
When $w\le1-\delta$, the Lipschitz bound gives
$|qg(w)-g(qw)|\le2K_g\varepsilon w$, and
Lemma~\ref{lem:step-slope} gives $-W_z\ge\kappa_\delta w$.
Thus \eqref{eq:shift-gain-choice} makes the expression in~\eqref{eq:upper-shift-residual} nonnegative and that in~\eqref{eq:lower-shift-residual} nonpositive.
When $w\ge1-\delta$, the interval between $w$ and $qw$
lies in $[1-2\delta,1]$ for $q=1-\varepsilon(t)$ and, on
$qw<1$, also for $q=1+\varepsilon(t)$;
from $g'\le-\kappa$ and $g(w)\ge0$ we obtain
\[
 qg(w)-g(qw)\ge\kappa\varepsilon w\quad\text{for }q=1+\varepsilon,
 \qquad
 qg(w)-g(qw)\le-\kappa\varepsilon w\quad\text{for }q=1-\varepsilon.
\]
Together with $\nu<\kappa$ and $W_z<0$, these inequalities
give a nonnegative expression in~\eqref{eq:upper-shift-residual}
and a nonpositive expression in~\eqref{eq:lower-shift-residual}. Clipping the upper function at
one preserves the supersolution inequality, since its spatial
derivative jumps downward from zero to $qW_z<0$ and hence
contributes a nonnegative measure to the negative Laplacian.

With $\delta,\varepsilon_0,\nu,K$ fixed independently of the initial
shifts, enlarge $T$ so that \eqref{eq:normal-boundary-deficit} gives
\[
 1-U(c_1t,t)\le\varepsilon(t)\quad\text{for }t\ge T.
\]
This is possible because $\omega_1>\nu$. Since $W\le1$,
$(1-\varepsilon(t))W\le1-\varepsilon(t)\le U(c_1t,t)$,
which orders the lower comparison and $U$ at $r=c_1t$.
Fix $\tau_-=T/2$ and $\tau_+=2T$.
For the upper comparison, putting $\tau=t-T+\tau_+$ in~\eqref{eq:step-comparison-position} and \eqref{eq:step-left-tail}
gives
\[
 1-W(c_1t-\zeta_+(t),\tau)
 \le C\exp\bigl(\omega[c_1t-\zeta_+(T)-c_{a_0}\tau]
                  +\tfrac{3\omega}{2\lambda_{a_0}}\log(1+\tau)\bigr).
\]
Here $C$ absorbs the bounded difference
$X_{1/2}(\tau)-m_{a_0}(\tau)$ from
\eqref{eq:step-comparison-position} and
any initial bounded time interval. After division by
$\varepsilon(t)$, the factor depending on $t$ is bounded on $[T,\infty)$ because
$\nu<\omega(c_{a_0}-c_1)$.
Increasing $\zeta_+(T)$ therefore makes the ratio at most
$1/(1+\varepsilon_0)$, so $(1+\varepsilon(t))W\ge1$ at
$r=c_1t$. The clipped upper comparison consequently equals one
there and dominates $U$.

At the fixed initial time $T$,
Lemma~\ref{lem:normal-fixed-time} gives
\[\label{eq:initial-gaussian-ratio-rates}
 \begin{aligned}
 \lim_{r\to\infty}\frac1{r^2}
       \log\frac{W(r,\tau_-)}{U(r,T)}
     &=\frac1{4DT}-\frac1{4D\tau_-}<0,\\
 \lim_{r\to\infty}\frac1{r^2}
       \log\frac{U(r,T)}{W(r,\tau_+)}
     &=\frac1{4D\tau_+}-\frac1{4DT}<0.
 \end{aligned}
\]
By \eqref{eq:initial-gaussian-ratio-rates} and $U\le1$, choose $R_0>c_1T$ such that
$(1-\varepsilon_0)W(r,\tau_-)\le U(r,T)\le\min\{1,(1+\varepsilon_0)W(r,\tau_+)\}$
for $r\ge R_0$. Monotonicity of $W$ shows that nonnegative shifts
$\zeta_-(T)$ and $\zeta_+(T)$ preserve these inequalities on $[R_0,\infty)$.
Set $m_T=\min_{[c_1T,R_0]}U(\cdot,T)>0$. Choose nonnegative shifts
so that $W(c_1T+\zeta_-(T),\tau_-)\le m_T$ and
$W(R_0-\zeta_+(T),\tau_+)\ge(1+\varepsilon_0)^{-1}$.
Both choices are possible because $W(-\infty,\tau_+)=1$ and
$W(+\infty,\tau_-)=0$. Monotonicity of $W$ gives the two initial
inequalities on $[c_1T,R_0]$. Increasing $\zeta_+(T)$ also
preserves the upper inequality on $r=c_1t$.

To conclude comparison on the unbounded moving domain, fix
$T_2>T$ and use $y=r-c_1t\ge0$.
The drift $DB(y+c_1t)+c_1$ is bounded on this half-strip.
Let $\underline U(r,t)$ and $\overline U(r,t)$ be the lower
and upper functions in~\eqref{eq:normal-shift-comparison}, and set
$E_-(y,t)=(\underline U(y+c_1t,t)-U(y+c_1t,t))_+$ and
$E_+(y,t)=(U(y+c_1t,t)-\overline U(y+c_1t,t))_+$.
The differential inequalities in~\eqref{eq:upper-shift-residual}
and \eqref{eq:lower-shift-residual} imply, for $E=E_-$ or $E=E_+$,
$E_t-DE_{yy}-[DB(y+c_1t)+c_1]E_y\le K_gE$ weakly on
$y>0$, $T<t\le T_2$, with zero initial and lateral trace. Subtract
$\eta e^{C(t-T)}(1+y^2)$, where $\eta>0$ and
$C>K_g+2D+1+\sup_{r\ge c_1T}|DB(r)+c_1|$.
Applying $\partial_t-D\partial_{yy}-(DB+c_1)\partial_y-K_g$
to $\eta e^{C(t-T)}(1+y^2)$ and dividing by
$\eta e^{C(t-T)}$ gives
$(C-K_g)(1+y^2)-2D-2(DB+c_1)y>0$, using
$2|y|\le1+y^2$. Because $E$ is bounded, any positive maximum of
$E(y,t)-\eta e^{C(t-T)}(1+y^2)$ lies in a bounded region,
where the maximum principle excludes such a maximum. Letting $\eta\downarrow0$ proves the
order on $[T,T_2]$, and arbitrariness of $T_2$ gives
\eqref{eq:normal-shift-comparison} for every $t\ge T$.
\end{proof}

\begin{lemma}\label{lem:normal-logarithmic-position}
Under the assumptions of Lemma~\ref{lem:normal-interior}, every
$0<\theta<1$ has constants $C_\theta,T_\theta<\infty$ such that
\[\label{eq:normal-logarithmic-position}
 [0,m_{a_0}(t)-C_\theta]
 \subset\{r\ge0:U(r,t)\ge\theta\}
 \subset[0,m_{a_0}(t)+C_\theta]\quad\text{for }t\ge T_\theta.
\]
\end{lemma}

\begin{proof}
Fix $0<c_1<c_{a_0}$ and apply Lemma~\ref{lem:normal-shift-comparison}
on $r\ge c_1t$.
If $U(r,t)\ge\theta$, the upper inequality in~\eqref{eq:normal-shift-comparison} gives
$W(r-\zeta_+(t),t-T+\tau_+)\ge\theta/(1+\varepsilon_0)$.
By Lemma~\ref{lem:homogeneous-input} and boundedness of
$\zeta_+$, this implies $r\le m_{a_0}(t)+C_\theta$;
the replacement of $t$ by $t-T+\tau_+$ changes $m_{a_0}(t)$ by a bounded amount.
Indeed, for each fixed $h\in\R$ and $t\ge\max\{0,-h\}$,
\[\label{eq:fixed-age-position-difference}
 m_{a_0}(t+h)-m_{a_0}(t)
 =c_{a_0} h-\frac{3}{2\lambda_{a_0}}
                  \log\frac{1+t+h}{1+t}\to c_{a_0} h
                  \quad\text{as }t\to\infty.
\]
For the inner inclusion, choose $\theta'\in(\theta,1)$.
At sufficiently large times,
$1-\varepsilon(t)\ge\theta/\theta'$.
The lower comparison and \eqref{eq:step-comparison-position}
at level $\theta'$, together with
\eqref{eq:fixed-age-position-difference}, yield
$U(r,t)\ge\theta$ for
$c_1t\le r\le m_{a_0}(t)-C_\theta$ after increasing the constant.
Since $c_1<c_{a_0}$, this interval is eventually nonempty.
Equation \eqref{eq:normal-interior-convergence} covers
$0\le r\le c_1t$, and this interval also lies below
$m_{a_0}(t)+C_\theta$ for large $t$.
Combining \eqref{eq:normal-interior-convergence} with
\eqref{eq:normal-shift-comparison} proves both inclusions in~\eqref{eq:normal-logarithmic-position}.
\end{proof}

\begin{proof}[Proof of Theorem~\ref{thm:common}]
The inequality $f(x,v)\le a_0v$ from
\eqref{eq:common-assumptions} and linear comparison give
$u(\cdot,t)\le e^{a_0t}e^{tD\Delta}u_0$.
Using \eqref{eq:kernel-L2} and \eqref{eq:heat-diagonal},
\[
 \norm{u(\cdot,t)}_{L^{\infty}(\Hh^d)}
 \le e^{a_0t}p_{2t}(0)^{1/2}\norm{u_0}_{L^{2}(\Hh^d)}
 \le Ct^{-3/4}e^{(a_0-D\rho_d)t}\norm{u_0}_{L^{2}(\Hh^d)}\quad\text{for }t\ge1.
\]
This proves \eqref{eq:common-decay}, including the critical
case $a_0=D\rho_d$.

Now suppose $a_0>D\rho_d$ and choose $g_-$ and $g_+$ as in
Lemma~\ref{lem:envelopes}.
For compactly supported $u_0$, choose $R_2$ with
$u_0\le\one_{B(o,R_2)}$, and let $U_e^+$ solve the radial
problem in~\eqref{eq:three-reductions} with
$F(r,v)=g_+(v)$ and $U_e^+(r,0)=\one_{[0,R_2]}(r)$.
Comparison gives $u(x,t;u_0,f)\le U_e^+(r(x),t)$.
Since $u_0\not\equiv0$ and $f(x,v)\ge0$ for $0\le v\le1$
by \eqref{eq:common-assumptions}, the formula
\eqref{eq:duhamel} gives $u(x,1)\ge(P_1u_0)(x)>0$.
Continuity on a compact ball then gives
$u(x,1)\ge\eta\one_{B(o,R_1)}(x)$ for some $R_1,\eta>0$.
Let $U_e^-$ solve the radial problem in~\eqref{eq:three-reductions}
with $F(r,v)=g_-(v)$ and
$U_e^-(r,0)=\eta\one_{[0,R_1]}(r)$.
Comparison gives $U_e^-(r(x),t-1)\le u(x,t;u_0,f)$ for
$t\ge1$.
Both comparison data are nonincreasing and compactly supported,
so Lemma~\ref{lem:normal-logarithmic-position} gives the inclusions in~\eqref{eq:normal-logarithmic-position}
for $U_e^-$ and $U_e^+$, with the same $m_{a_0}(t)$,
since $g_-'(0)=g_+'(0)=a_0$.
The difference $m_{a_0}(t)-m_{a_0}(t-1)$ is bounded for large $t$ by
\eqref{eq:fixed-age-position-difference}; hence these two
comparisons prove \eqref{eq:common-balls}.

For tube data, take $U_e$ from \eqref{eq:three-reductions}
with $F(r,v)=g_-(v)$ and
$U_e(r,0)=\eta\one_{[0,R_1]}(r)$, and take $U_h$ from
\eqref{eq:hyperplane-equation} with $g=g_+$ and
$U_h(\rho,0)=\one_{[-R_2,R_2]}(\rho)$.
By Lemma~\ref{lem:tube-sandwich},
$U_e(s(x),t)\le u(x,t)\le U_h(s(x),t)$ for $t\ge0$. Applying
Lemma~\ref{lem:normal-logarithmic-position} to each gives the
level-set inclusions in~\eqref{eq:normal-logarithmic-position}
with the same function $m_{a_0}(t)$. Taking the larger of the two
position constants and time thresholds gives
$\hat{C}_{\theta,2}$ and $\hat{T}_{\theta,2}$ in
\eqref{eq:common-tubes}.
The assumption $R_1>0$ ensures that the lower initial profile is
nonzero on a set of positive hyperbolic volume.

For horoball data, let $W_-$ and $W_+$ solve
\eqref{eq:step-comparison-equation} with $g_-$ and $g_+$,
respectively. The functions $W_-(b(x)-A,t)$ and
$W_+(b(x)-B,t)$ have initial values
$\one_{\mathcal H_A}$ and $\one_{\mathcal H_B}$, and
comparison encloses $u$ between them.
By Lemma~\ref{lem:homogeneous-input}, the level positions of
$W_-(b-A,t)$ and $W_+(b-B,t)$ are $m_{a_0}(t)+O(1)$ as $t\to\infty$,
for every $a_0>0$, whether $c_{a_0}$ is positive, zero or negative.
The inequalities $W_-(b(x)-A,t)\le u(x,t)
\le W_+(b(x)-B,t)$ therefore prove \eqref{eq:common-horoballs}.
\end{proof}

The identity $g_-'(0)=g_+'(0)=a_0$ fixes the linear and
logarithmic terms in $m_{a_0}(t)$ for each comparison pair.
The constants in Theorem~\ref{thm:common} absorb the separate
bounded position errors of $U_e^-$ and $U_e^+$ for ball data,
$U_e$ and $U_h$ for tube data, and
$W_-(\cdot-A,t)$ and $W_+(\cdot-B,t)$ for horoball data,
as supplied by Lemmas~\ref{lem:normal-logarithmic-position} and
\ref{lem:homogeneous-input}.

\section{Solutions near a geodesic}\label{sec:localized}

Throughout this section $d=2$, $\beta=1$, $\sigma=\sigma_2$,
and $f(x,v)=F(s(x),v)$. Constants may depend on the fixed diffusion
coefficient $D>0$. We use signed normal distance $\rho\in\R$ when integrating
on the full plane and write $s=|\rho|$ for an even profile. Thus
\[\label{eq:plane-Fermi}
 \begin{gathered}
 g=\dd\rho^2+\cosh^2\rho\dd\ell^2,\qquad
 \dd V=\cosh\rho\dd\rho\dd\ell,\\
 \Delta=\partial_{\rho\rho}+\tanh\rho\partial_\rho
                      +\sech^2\rho\partial_{\ell\ell}
 \quad {\text{for }(\rho,\ell)\in\R^2}.
 \end{gathered}
\]
In this section $u(\rho,\ell,t)$ denotes the solution of~\eqref{eq:KPP} expressed in these coordinates, with reaction
$F(|\rho|,u)$. For a function of $s$, a prime denotes differentiation
with respect to $s$. Derivatives of $\Lambda(\lambda)$ are written
$\Lambda^\prime$ and $\Lambda^{\prime\prime}$; derivatives of
parameter-dependent profiles in $\lambda$ are written $\partial_\lambda$. For an even $C^2$ function $\rho\mapsto P(\rho)$, the
restriction to $s\ge0$ satisfies $P'(0)=0$.

\subsection{Stationary solutions}

\begin{lemma}\label{lem:plane-operator}
The map $v\mapsto\sqrt{\cosh s}\,v$ is unitary from
$L^2((0,\infty),\cosh s\dd s)$ to $L^2(0,\infty)$. It conjugates the self-adjoint operator associated with the form
in~\eqref{eq:sigma}, whose differential expression is
$D(v''+\tanh s\,v')+a(s)v$ and whose endpoint condition is
$v'(0)=0$, to
\[\label{eq:plane-Schrodinger}
 \J=D\partial_{ss}+a(s)-\frac D4-\frac D4\sech^2s,
 \qquad \operatorname{Dom}(\J)=\{w\in H^2(0,\infty):w'(0)=0\}.
\]
The upper endpoint of the essential spectrum of $\J$ is $\mu$. If $\sigma>\mu$, then
$\sigma$ is an isolated simple eigenvalue with a strictly positive
eigenfunction satisfying the Neumann condition at zero.
\end{lemma}

\begin{proof}
The identity $w=\sqrt{\cosh s}\,v$ gives
$\int_0^\infty|w|^2\dd s=\int_0^\infty|v|^2\cosh s\dd s$, and Lemma~\ref{lem:liouville}, with $B=\tanh s$,
gives the expression in~\eqref{eq:plane-Schrodinger}.
Since $(\cosh s)'=0$ at zero and the potential is bounded, $\J$ has the operator domain in~\eqref{eq:plane-Schrodinger} and form domain $H^1(0,\infty)$. By Lemma~\ref{lem:far-tests}, the upper endpoint of its essential
spectrum is $\mu$; thus $\sigma>\mu$ is an isolated eigenvalue of
finite multiplicity \cite{ReedSimon}.

Replacing a real eigenfunction of $\J$ at $\sigma$ by its absolute value
preserves the quadratic form. Elliptic regularity and the strong
maximum principle therefore give a positive eigenfunction on
$(0,\infty)$, which is also positive at zero: otherwise its
zero value and Neumann derivative would force the identically
zero solution by uniqueness for the ODE in~\eqref{eq:plane-Schrodinger}.
For a fixed eigenvalue, uniqueness for the initial-value problem
determines a Neumann eigenfunction from its value at zero.
Any two eigenfunctions at that eigenvalue are therefore proportional,
proving simplicity.
\end{proof}

\begin{lemma}\label{lem:cap}
Suppose $\mu\le0$. For all sufficiently large $R$, the function
\[\label{eq:cap-definition}
 U_R(s)=
 \begin{cases}
 1,&0\le s\le R,\\
 h(s)/h(R),&s>R,
 \end{cases}
 \qquad h(s)=(\cosh s)^{-1/2}\phiZ(s),
\]
gives, by composition $x\mapsto U_R(s(x))$, a positive weak
stationary supersolution of~\eqref{eq:KPP} with $f=f_{\G}$ on
$\Hh^2$, and $U_R\le1$. For every compactly supported $u_0$ with $0\le u_0\le1$, there is
such an $R$ for which $u_0(x)\le U_R(s(x))$. More generally, if
$u_0(x)\le\min\{1,C_0h(s(x))\}$, one can choose $R$ so that $u_0(x)\le U_R(s(x))$.
\end{lemma}

\begin{proof}
By Lemma~\ref{lem:plane-operator} and the cell equation
$\K_{\lambda_0}\phiZ=\mu\phiZ$,
\[\label{eq:h-super}
 \bigl[D(\partial_{ss}+\tanh s\partial_s)+a(s)\bigr]h
       =\left(\mu-\frac D4\sech^2s\right)h\le0
 \quad {\text{for }s>0}.
\]
The definitions of $h$ and $H_{\mathrm{sym}}$ also give
\[\label{eq:h-log-derivative}
 \frac{h'}h
 =\frac{H_{\mathrm{sym}}'}{H_{\mathrm{sym}}}+\frac12(1-\tanh s),
 \qquad H_{\mathrm{sym}}=e^{-s/2}\phiZ(s)\quad {\text{for }s\ge0}.
\]
Lemma~\ref{lem:decreasing-ground} and
\eqref{eq:h-log-derivative} allow us to fix $R>0$ such that
$h'(s)<0$ for every $s\ge R$.
On $s>R$, \eqref{eq:h-super} and the KPP bound in~\eqref{eq:F1} make $h(s)/h(R)$ a stationary supersolution;
on $s<R$, the constant one solves the stationary equation.
At $s=R$, the outward derivative jumps from zero to
$h'(R)/h(R)<0$.
More precisely, if $\Sigma_\pm=\{\rho=\pm R\}$ and
$\delta_{\Sigma_\pm}$ denotes induced length measure, the
singular part of $D\Delta[U_R(s(x))]$ is
$D[h'(R)/h(R)](\delta_{\Sigma_+}+\delta_{\Sigma_-})$.
The nonpositive singular part of $D\Delta[U_R(s(x))]$ gives
the weak supersolution inequality on $\Sigma_+$ and $\Sigma_-$; $U_R(s(x))=1$ on $\{s(x)<R\}$ gives regularity across $\G$.

Periodic positivity of $\phiZ$ gives $h>0$ and $h(s)\to0$.
For compactly supported $u_0$, take $R$ with
$\supp u_0\subset\{x\in\Hh^2:s(x)<R\}$.
If $u_0\le\min\{1,C_0h\}$, choose $R$ with $1/h(R)\ge C_0$.
In either case, the two branches of~\eqref{eq:cap-definition}
give $u_0(x)\le U_R(s(x))$.
\end{proof}

\begin{lemma}\label{lem:planar-minimal}
If $\mu\le0<\sigma$, the smallest positive stationary solution has the
form $p_*(x)=p(s(x))$ and satisfies \eqref{eq:localized-bounds}.
The confinement estimate \eqref{eq:confinement} holds.
\end{lemma}

\begin{proof}
By Theorem~\ref{thm:dimension}, the smallest positive stationary solution
has the form $p_*(x)=p(s(x))$.
For each $U_R$ in Lemma~\ref{lem:cap}, take
$\varepsilon\widetilde\varphi$ from Lemma~\ref{lem:minimal-state}
small enough that
$\varepsilon\widetilde\varphi(x)\le U_R(s(x))$ on $\Hh^2$. Let $u_\varepsilon$ solve \eqref{eq:KPP} with $f=f_{\G}$
and $u_\varepsilon(x,0)=\varepsilon\widetilde\varphi(x)$. Comparison gives
$u_\varepsilon(x,t)\le U_R(s(x))$; letting $t\to\infty$ gives $p_*(x)\le U_R(s(x))$.
Since $\min_{[0,R]}h>0$, we have
$p(s)\le U_R(s)\le Ch(s)$ for $s\ge0$.
As $h(s)\to0$ when $s\to\infty$, it follows that $p(s)\to0$
as $s\to\infty$.
The strong maximum principle gives $0<p<1$, and even regularity
at the geodesic yields
\[\label{eq:planar-stationary}
 \left\{\begin{aligned}
 &D(p''+\tanh s\,p')+F(s,p)=0&&\text{for }s>0,\\
 &p'(0)=0.
 \end{aligned}\right.
\]
Multiplying \eqref{eq:planar-stationary} by $\cosh s$ and
using $p'(0)=0$ gives
\[\label{eq:planar-stationary-flux}
 D\cosh s\,p'(s)=-\int_0^s\cosh v\,F(v,p(v))\dd v<0
 \quad\text{for }s>0.
\]
Finally, for compactly supported $u_0$, Lemma~\ref{lem:cap}
gives $u_0(x)\le U_R(s(x))$. Comparison yields
$u(x,t;u_0,f_{\G})\le U_R(s(x))\le Ch(s(x))$ for every
$t\ge0$. Combining this with $u\le1$ proves
\eqref{eq:confinement}.
\end{proof}

\paragraph{Uniqueness and convergence.}

For uniqueness among stationary solutions depending on both $\rho$
and $\ell$, we first construct cutoffs $\chi_R$ for which
$\int_{\Hh^2}q^2|\nabla\chi_R|_g^2\dd V\to0$ as $R\to\infty$
whenever $0<q(x)\le Ch(s(x))$. 

\begin{lemma}\label{lem:weighted-cutoffs}
Let $q:\Hh^2\to(0,\infty)$ be measurable and satisfy
$q(x)\le Ch(s(x))$. There are smooth
compactly supported functions $0\le\chi_R\le1$, tending to one on each
compact set as $R\to\infty$, such that
\[\label{eq:cutoff-energy}
 \int_{\Hh^2}q^2|\nabla\chi_R|_g^2\dd V\le\frac{C'}R
 \quad\text{for }R\ge1.
\]
\end{lemma}

\begin{proof}
Choose a smooth even function $\chi$ with $0\le\chi\le1$,
equal to one on $[-1,1]$ and zero outside $[-2,2]$, and set
\[\label{eq:anisotropic-cutoff}
 \chi_R(\rho,\ell)=\chi(\rho/R^2)\chi(\ell/R)
 \quad {\text{for }(\rho,\ell)\in\R^2,\quad R\ge1}.
\]
By the assumed bound for $q$ and \eqref{eq:h-weight},
$q^2\cosh\rho\le C^2\phiZ(|\rho|)^2\le C_1$.
The term containing $\partial_\rho\chi_R$ in
$\int_{\Hh^2}q^2|\nabla\chi_R|_g^2\dd V$ is therefore at most
$CR^{-4}R^2R=C/R$.
For the longitudinal derivative, \eqref{eq:plane-Fermi} gives
the additional factor $\sech^2\rho$, and its contribution is
bounded by
\[
 C R^{-2}\int_{-2R}^{2R}\dd\ell
             \int_\R\sech^2\rho\dd\rho\le C/R.
\]
Combining the bounds for the two cutoff derivatives proves
\eqref{eq:cutoff-energy}.
In the global Fermi coordinates,
$\supp\chi_R\subset\{|\rho|\le2R^2,\ |\ell|\le2R\}$,
so each support is compact, and $\chi_R\to1$ locally uniformly.
\end{proof}

\begin{lemma}\label{lem:localized-uniqueness}
Assume \eqref{eq:quotient} and $\mu\le0<\sigma$. If a stationary
solution $q$ of~\eqref{eq:KPP} with $f=f_{\G}$ satisfies
$0<q\le1$ and $q(x)\le Ch(s(x))$, then $q=p_*$.
\end{lemma}

\begin{proof}
By minimality, $0<p_*\le q$. Put $\theta=p_*/q$, so
$0<\theta\le1$. By \eqref{eq:quotient}, the two stationary equations
and $q^2\nabla\theta=q\nabla p_*-p_*\nabla q$ give
\[\label{eq:ratio-equation}
 -D\operatorname{div}(q^2\nabla\theta)
 =p_*q\left(\frac{f(x,p_*)}{p_*}-\frac{f(x,q)}q\right)=:G\ge0
 \quad {\text{for }x\in\Hh^2}.
\]
For fixed $\varepsilon>0$, test \eqref{eq:ratio-equation}
with $\chi_R^2/(\theta+\varepsilon)$, using the cutoffs of
Lemma~\ref{lem:weighted-cutoffs}.
The test is admissible by local Sobolev approximation, since its
denominator is positive and $p_*$ and $q$ are $C^2$ on the compact support of $\chi_R$. Integration against $\dd V$ gives
\[
 0\le{\int_{\Hh^2}}\frac{G\chi_R^2}{\theta+\varepsilon}\dd V
 =2D{\int_{\Hh^2}}\frac{q^2\chi_R}{\theta+\varepsilon}
                    \nabla\theta\cdot\nabla\chi_R\dd V
   -D{\int_{\Hh^2}}\frac{q^2\chi_R^2}{(\theta+\varepsilon)^2}
                           |\nabla\theta|_g^2\dd V.
\]
Dividing by $D$ and applying $2\xi\eta\le \xi^2/2+2\eta^2$
for $\xi,\eta\in\R$, we obtain
\[\label{eq:ratio-energy}
 {\int_{\Hh^2}}\frac{q^2\chi_R^2}{(\theta+\varepsilon)^2}
                           |\nabla\theta|_g^2\dd V
 \le4{\int_{\Hh^2}} q^2|\nabla\chi_R|_g^2\dd V.
\]
By Lemma~\ref{lem:weighted-cutoffs},
$4\int_{\Hh^2}q^2|\nabla\chi_R|_g^2\dd V\to0$ as $R\to\infty$.
On each fixed compact set, $\chi_R=1$ for large $R$ and
$q^2/(\theta+\varepsilon)^2$ is bounded below by a positive constant.
Hence \eqref{eq:ratio-energy}, with $\varepsilon$ fixed, gives
$\nabla\theta=0$; connectedness makes $\theta$ constant.
If $\theta<1$, \eqref{eq:quotient} would
make the right-hand side of~\eqref{eq:ratio-equation} strictly
positive, whereas its left-hand side is zero. Hence
$\theta=1$ and $q=p_*$.
\end{proof}

\begin{lemma}\label{lem:localized-convergence}
Assume \eqref{eq:quotient} and $\mu\le0<\sigma$. Every
nonzero datum with $0\le u_0\le\min\{1,C_0h(s(x))\}$ satisfies
\eqref{eq:strip-convergence}.
\end{lemma}

\begin{proof}
By Lemma~\ref{lem:cap}, choose $U_R$ with
$u_0(x)\le U_R(s(x))$, and let $v$ solve \eqref{eq:KPP}
with $f=f_{\G}$ and $v(x,0)=U_R(s(x))$.
Since $x\mapsto U_R(s(x))$ is a stationary supersolution,
comparison gives $p_*(x)\le v(x,t)\le U_R(s(x))$ and
$v(x,t+\tau)\le v(x,t)$ for $t,\tau\ge0$.
Lemma~\ref{lem:compactness} therefore gives a locally uniform
stationary limit $q$, with $q>0$ and $q(x)\le Ch(s(x))$.
By Lemma~\ref{lem:localized-uniqueness}, $q=p_*$.
Combining $u\le v$ with the locally uniform lower bound of
Lemma~\ref{lem:minimal-state} proves $u\to p_*$ locally
uniformly as $t\to\infty$.

Fix $M\ge0$ and $\varepsilon>0$. The bounds
$u(x,t)\le Ch(s(x))$ and $p_*(x)\le Ch(s(x))$, with
$h(s)\to0$ as $s\to\infty$, give
$|u(x,t)-p_*(x)|<\varepsilon$ on $\{s(x)\ge S\}$
for sufficiently large $S$, uniformly in $t\ge0$. On the compact rectangle
$\{(\rho,\ell):|\rho|\le S,\ |\ell|\le M\}$, local uniform
convergence gives $|u(\rho,\ell,t)-p(|\rho|)|<\varepsilon$
for all sufficiently large $t$. Letting $\varepsilon\downarrow0$ proves
\eqref{eq:strip-convergence}.
\end{proof}

\paragraph{The critical stationary tail.}

\begin{lemma}\label{lem:critical-tail}
Assume $\mu=0<\sigma$, \eqref{eq:F1}, and
\eqref{eq:axis-regularity} for the geodesic-distance reaction. The profile $p$ of
Lemma~\ref{lem:planar-minimal} satisfies \eqref{eq:critical-tail} with
$A>0$.
\end{lemma}

\begin{proof}
Define $y(s)=p(s)/h(s)$ for $s\ge0$, with $h$ from
\eqref{eq:cap-definition}. By Lemma~\ref{lem:planar-minimal}, $y$ is
positive and bounded. Substituting in~\eqref{eq:planar-stationary}, using \eqref{eq:h-super} at
$\mu=0$ and $R_F(s,p(s))=a(s)p(s)-F(s,p(s))$ from
\eqref{eq:uniform-taylor}, gives
\[\label{eq:critical-y}
 (\phiZ^2 y')'
 =\frac14\sech^2s\,\phiZ^2y
       +\frac1D\phiZ\sqrt{\cosh s}\,R_F(s,h y)=:G(s)\ge0
 \quad {\text{for }s>0}.
\]
The bounds for $\phiZ$, its reciprocal and $y$, together with
$0\le R_F(s,p)\le C_2p^2$, imply
\[\label{eq:critical-G}
 0\le G(s)\le C\left(e^{-2s}y(s)+e^{-s/2}y(s)^2\right)
 \quad {\text{for }s\ge0}.
\]
Thus $\int_0^\infty(1+s)G(s)\dd s<\infty$, and
\eqref{eq:critical-y} gives a finite limit for $\phiZ(s)^2y'(s)$ as $s\to\infty$.
Since $0<y\le C$ and $\phiZ$ is positive and periodic,
this limit is zero: a positive limit would make $y$ unbounded,
whereas a negative limit would force $y(s)<0$ for large $s$.
Integrating \eqref{eq:critical-y} from infinity therefore yields
\[\label{eq:y-derivative-tail}
 y'(s)=-\phiZ(s)^{-2}\int_s^\infty G(v)\dd v\le0
 \quad {\text{for }s\ge0}.
\]
Thus $y(s)\to A_0\ge0$, and
\eqref{eq:critical-G}--\eqref{eq:y-derivative-tail} give
$y(s)=A_0+O(e^{-s/2})$ as $s\to\infty$.
Integrating once more, using
$\int_0^\infty sG(s)\dd s<\infty$, gives
\[\label{eq:critical-amplitude-integral}
 A_0=\frac{p(0)}{\phiZ(0)}
 -\int_0^\infty\phiZ(r)^{-2}
                   \int_r^\infty G(v)\dd v\dd r.
\]
To prove $A_0>0$, suppose $A_0=0$. A further integration of~\eqref{eq:y-derivative-tail} then gives
\[\label{eq:Volterra-tail}
 y(s)=\int_s^\infty\phiZ(r)^{-2}
                         \int_r^\infty G(v)\dd v\dd r\quad {\text{for }s\ge0}.
\]
Let $Y=\sup_{v\ge0}y(v)<\infty$.
By monotonicity of $y$ and \eqref{eq:critical-G}, the
right-hand side of~\eqref{eq:Volterra-tail} is bounded by
$C(e^{-2s}+Ye^{-s/2})y(s)$; here
$\int_s^\infty\int_r^\infty e^{-\alpha_*v}\dd v\dd r=e^{-\alpha_*s}/\alpha_*^2$
for $\alpha_*>0$.
Since $y(s)>0$, this contradicts \eqref{eq:Volterra-tail} for large
$s$. Hence $A_0>0$.
Finally, multiplying the expansion for $y$ by
$h=(\cosh s)^{-1/2}\phiZ(s)$ and using
$(\cosh s)^{-1/2}=\sqrt2e^{-s/2}(1+O(e^{-2s}))$
gives $p(s)=\sqrt2A_0e^{-s/2}\phiZ(s)+O(e^{-s})$.
Taking $A=\sqrt2A_0$ proves \eqref{eq:critical-tail}.
\end{proof}

\begin{proof}[Proof of Theorem~\ref{thm:localized}]
Lemma~\ref{lem:planar-minimal} gives the bounds in~\eqref{eq:localized-bounds} and the confinement estimate
\eqref{eq:confinement}. Under \eqref{eq:quotient},
Lemmas~\ref{lem:localized-uniqueness} and
\ref{lem:localized-convergence} prove uniqueness in the stated
class $0<q(x)\le1$, $q(x)\le Ch(s(x))$, and convergence
uniformly on $\{x\in\Hh^2:|\ell(x)|\le M\}$ for each fixed $M\ge0$.
At $\mu=0$, Lemma~\ref{lem:critical-tail} gives the asserted
expansion \eqref{eq:critical-tail} with $A>0$, completing the proof.
\end{proof}

\begin{lemma}\label{lem:shooting-family}
Assume \eqref{eq:quotient} and $\mu\le0<\sigma$. For each
$p(0)<\alpha<1$, the initial-value problem
\[\label{eq:shooting-IVP}
 \left\{\begin{aligned}
 &D(P_\alpha''+\tanh s\,P_\alpha')+F(s,P_\alpha)=0
 \quad {\text{for }s>0},\\
 &P_\alpha(0)=\alpha,\qquad P_\alpha'(0)=0,
 \qquad {p(0)<\alpha<1}
 \end{aligned}\right.
\]
has a global solution on $[0,\infty)$ satisfying
$P_\alpha(s)>p(s)$, $P_\alpha'(s)<0$ for $s>0$, and
$P_\alpha(s)\to0$ as $s\to\infty$.
\end{lemma}

\begin{proof}
The initial-value problem \eqref{eq:shooting-IVP} is regular
at zero. While $0<P_\alpha<1$, its weighted derivative obeys
$D(\cosh s\,P_\alpha')'=-\cosh s\,F(s,P_\alpha)<0$
by \eqref{eq:F1}. Since $P_\alpha'(0)=0$, this gives
$P_\alpha'<0$ for $s>0$ on the interval of existence where $P_\alpha>0$.
On the interval where $P_\alpha>0$, define
$q(s)=P_\alpha(s)/p(s)$. Subtracting the equations
\eqref{eq:shooting-IVP} and \eqref{eq:planar-stationary}
gives
\[\label{eq:shooting-ratio}
 D(\cosh s\,p^2q')'
 =\cosh s\,pP_\alpha
      \left(\frac{F(s,p)}p-\frac{F(s,P_\alpha)}{P_\alpha}\right)
 \quad {\text{for }s>0\text{ with }P_\alpha(s)>0}.
\]
At $s=0$, $q(0)=\alpha/p(0)>1$ and $q'(0)=0$.
By \eqref{eq:quotient}, the right-hand side of~\eqref{eq:shooting-ratio} is strictly positive while $q>1$,
so $q'>0$ for $s>0$ up to any possible first point at which $q=1$.
Since $q(0)>1$ and $q'>0$, no such point exists; hence $P_\alpha>p>0$
throughout its existence interval.
Since $P_\alpha$ is decreasing, it remains below one. With
$K_F=\max_{[0,L]\times[0,1]}F$, integration of
$D(\cosh s\,P_\alpha')'=-\cosh s\,F(s,P_\alpha)$ gives
$|P_\alpha'(s)|\le (K_F/D)\tanh s\le K_F/D$.
At any finite endpoint $s_*$ of the interval of existence,
$P_\alpha(s)\ge p(s)$ and $p(s_*)>0$ bound $P_\alpha$ away from zero. ODE continuation therefore gives a global
solution on $[0,\infty)$.

Since $P_\alpha$ is decreasing and positive,
$l=\lim_{s\to\infty}P_\alpha(s)$ exists and lies in $[0,\alpha)$.
If $l>0$, periodicity and \eqref{eq:F1} give
$F(s,P_\alpha(s))\ge\eta>0$ for large $s$.
The integrated flux equation then yields
$P_\alpha'(s)\le-\eta/(2D)$ for sufficiently large $s$,
contrary to positivity of $P_\alpha$ on the half-line. Hence $l=0$.
\end{proof}

\begin{proof}[Proof of Theorem~\ref{thm:other-states}]
Take $p_\alpha=P_\alpha$ from Lemma~\ref{lem:shooting-family}
and extend it evenly across the geodesic.
The condition $P_\alpha'(0)=0$ and
\eqref{eq:plane-Fermi} make $x\mapsto P_\alpha(s(x))$ a
classical stationary solution of~\eqref{eq:KPP} with $f=f_{\G}$,
with $P_\alpha(s)>p(s)$, $P_\alpha'(s)<0$ for $s>0$, and $\lim_{s\to\infty}P_\alpha(s)=0$.
If $P_\alpha(s)\le Ch(s)$ on $[0,\infty)$ for some $C>0$,
Lemma~\ref{lem:localized-uniqueness} would force it to equal
$p$, contradicting $P_\alpha(0)=\alpha>p(0)$.
Thus $\sup_{s\ge0}P_\alpha(s)/h(s)=\infty$.
Since $P_\alpha(0)=\alpha$, distinct values of $\alpha$
give distinct stationary solutions.
\end{proof}

\subsection{A periodic example}

Define the smooth nondecreasing function $\varsigma:\R\to[0,1]$ by
\[\label{eq:smooth-step}
 E(v)=\begin{cases}e^{-1/v},&v>0,\\0,&v\le0,\end{cases}
 \qquad \varsigma(v)=\frac{E(v)}{E(v)+E(1-v)}
 \quad {\text{for }v\in\R}.
\]
Thus $\varsigma=0$ on $(-\infty,0]$ and $\varsigma=1$ on $[1,\infty)$, with all
one-sided derivatives matching at the endpoints. On one period define
\[\label{eq:example-coefficient}
 a(s)=D\left(\frac1{100}+\frac{99}{100}\chi(s)\right),\qquad
 \chi(s)=
 \begin{cases}
 1,&0\le s\le1,\\
 1-\varsigma(25(s-1)),&1\le s\le26/25,\\
 0,&26/25\le s\le199/25,\\
 \varsigma(50(s-199/25)),&199/25\le s\le8.
 \end{cases}
\]
In~\eqref{eq:example-coefficient}, $\chi(s)=1$ for
$399/50\le s\le8$. The periodic extension of $a$ is therefore
smooth and equals $D$ in a neighbourhood of $0$ modulo $8$. Consequently $F(s,v)=a(s)v(1-v)$ satisfies
\eqref{eq:F1}--\eqref{eq:axis-regularity} and
\eqref{eq:quotient}.

\begin{lemma}\label{lem:example-bulk}
The periodic eigenvalue $\mu$ for \eqref{eq:example-coefficient} is
strictly negative.
\end{lemma}

\begin{proof}
By \eqref{eq:example-coefficient}, $a$ is bounded above by the
$8$-periodic step potential $a^{\rm step}$ equal to $D$ on
$[-1/25,26/25]$ modulo eight and to $D/100$ elsewhere.
The interval where $a^{\rm step}=D$ has length $W=27/25$;
its complement in one period has length $G=173/25$. After the change of coordinate $t=s-1/2$ modulo $8$,
write $a^{\rm step}(t)$ for the translated potential. Then
$\{t\in[-4,4]:a^{\rm step}(t)=D\}=[-W/2,W/2]$.
For the equation $D\varphi''+(a^{\rm step}-D/4)\varphi=0$,
put
\[\label{eq:example-kappa}
 k=\frac{\sqrt3}{2},\qquad \kappa=\frac{\sqrt6}{5}.
\]
Define the positive continuous even periodic function
\[\label{eq:piecewise-test}
 \varphi(t)=
 \begin{cases}
 \cos(kt),&0\le |t|\le W/2,\\
 \displaystyle
 \frac{\cos(kW/2)}{\cosh(\kappa G/2)}
                    \cosh(\kappa(4-|t|)),&W/2\le|t|\le4.
 \end{cases}
\]
By \eqref{eq:example-kappa}, it solves the equation on each
constant-potential interval, and its derivatives vanish at the
identified endpoints $-4$ and $4$.
At $t=-W/2$ and $t=W/2$, the jump in $\varphi'$ is negative provided
\[\label{eq:step-sign-inequality}
 k\tan(kW/2)<\kappa\tanh(\kappa G/2).
\]
To verify \eqref{eq:step-sign-inequality}, use
$k<867/1000$, $kW/2<468/1000$,
$\kappa>489/1000$ and $\kappa G/2>169/100$.
The bounds $\sin x\le x$ and $\cos x\ge1-x^2/2$ give
\[\label{eq:step-rational-bounds}
 k\tan(kW/2)
 <\frac{(867/1000)(468/1000)}{1-(468/1000)^2/2}
 <\frac{456}{1000}
 <\frac{489}{1000}\frac{934}{1000}
 <\kappa\tanh(\kappa G/2).
\]
The final comparison in~\eqref{eq:step-rational-bounds} follows from
$\tanh(169/100)>934/1000$, equivalently
$e^{169/50}>967/33$; the nonnegative exponential series through
degree twelve gives
$\sum_{j=0}^{12}(169/50)^j/j!>967/33$.
These rational bounds prove \eqref{eq:step-sign-inequality}.

It follows that
$(D\partial_{tt}+a^{\rm step}-D/4)\varphi$ is a nonzero
negative sum of point masses at $t=-W/2$ and $t=W/2$.
Let $\psi>0$ be an $8$-periodic principal eigenfunction satisfying
$(D\partial_{tt}+a^{\rm step}-D/4)\psi=\mu^{\rm step}\psi$.
Distributional integration by parts on the circle gives
\[
 \mu^{\rm step}\int_{-4}^4\psi\varphi\dd t
 =\int_{-4}^4\psi
      (D\partial_{tt}+a^{\rm step}-D/4)\varphi\dd t<0.
\]
By positivity of $\psi\varphi$ and the periodic Rayleigh principle,
$\mu\le\mu^{\rm step}<0$.
\end{proof}

\begin{lemma}\label{lem:example-tube}
For $a$ in~\eqref{eq:example-coefficient}, the transverse
spectral endpoint in~\eqref{eq:sigma} satisfies $\sigma>D/200$.
\end{lemma}

\begin{proof}
Use $v(s)=\sech s$ in the transverse form \eqref{eq:sigma}.
The function $v$ is the limit of compactly supported tests in
$H^1((0,\infty),\cosh s\dd s)$ and satisfies
$v'(0)=0$. Direct integration gives
\[\label{eq:sech-integrals}
 \int_0^\infty v^2\cosh s\dd s=\frac\pi2,
 \qquad
 \int_0^\infty(v')^2\cosh s\dd s=\frac\pi4.
\]
By \eqref{eq:example-coefficient}, $a=D$ on $[0,1]$ and
$a\ge D/100$ elsewhere. Substituting \eqref{eq:sech-integrals}
into the Rayleigh quotient therefore yields
\[
 \sigma\ge\frac D{100}-\frac D2
       +\frac{99D}{100}\frac{\int_0^1\sech s\dd s}{\pi/2}.
\]
Since $\int_0^1\sech s\dd s=\arctan(\sinh1)>\pi/4$,
the right-hand side is strictly greater than
$D/100-D/2+99D/200=D/200$, as required.
\end{proof}

\begin{proof}[Proof of Example~\ref{ex:periodic-reaction}]
The coefficient $a$ in~\eqref{eq:example-coefficient} has the
smoothness, periodicity and range stated in the theorem.
By Lemmas~\ref{lem:example-bulk} and \ref{lem:example-tube},
it satisfies $\mu<0<\sigma_2$ and $\sigma_2>D/200$.
Theorems~\ref{thm:dimension} and \ref{thm:localized} consequently
give the smallest positive stationary solution
$p_*(x)=p(s(x))$, with $p(s)\to0$ as $s\to\infty$, and
\eqref{eq:strip-convergence} for every nonzero compactly supported datum.
For $f=f_b$ or $f=f_o$, Theorem~\ref{thm:threshold} gives
uniform extinction for $u_0\in L^2(\Hh^2)$.
\end{proof}

\begin{remark}
Lemma~\ref{lem:example-bulk} obtains $\mu\le\mu^{\rm step}<0$
from $a\le a^{\rm step}$ and the negative derivative jumps of the
piecewise smooth periodic function $\varphi$ in~\eqref{eq:piecewise-test}.
Lemma~\ref{lem:example-tube} gives $\sigma>D/200>0$ by inserting
$v(s)=\sech s$ into the Rayleigh quotient \eqref{eq:sigma}.
\end{remark}

For the study of propagation along the geodesic, assume
$\mu\le0<\sigma$ throughout the rest of this section and define
$\mathcal V(s)=\sech^2s$. By \eqref{eq:plane-Fermi}, $\mathcal V=|\nabla\ell|_g^2$; substitution of $e^{-\lambda\ell}H(s)$ into
$D\Delta+a(s)$, followed by $H=\psi/\sqrt{\cosh s}$, gives
$\J+D\lambda^2\mathcal V(s)$.

\subsection{Propagation along the geodesic}

\begin{lemma}\label{lem:longitudinal-branch}
For every real $\lambda$, the operator
$\J_\parallel(\lambda)=\J+D\lambda^2\mathcal V$ on $L^2(0,\infty)$,
with domain $\operatorname{Dom}(\J)$, has essential spectrum with upper endpoint
$\mu$ and a simple isolated largest eigenvalue $\Lambda(\lambda)\ge\sigma>0$.
Let $\psi_\lambda>0$ be its eigenfunction normalised by
$\int_0^\infty\psi_\lambda^2\dd s=1$. The function $\lambda\mapsto\Lambda(\lambda)$ is real analytic,
as is $\lambda\mapsto\psi_\lambda$ with values in $\operatorname{Dom}(\J)$
equipped with its graph norm, and
\[\label{eq:longitudinal-HF}
 \Lambda'(\lambda)=2D\lambda\int_0^\infty \mathcal V\psi_\lambda^2\dd s,
 \qquad \Lambda''(\lambda)>0
 \quad {\text{for }\lambda\in\R}.
\]
The function $\Lambda$ is even. For $\lambda\ge1$,
\[\label{eq:longitudinal-growth}
 D\lambda^2-C\lambda-C\le\Lambda(\lambda)\le D\lambda^2+a_{\max}
 \quad {\text{for }\lambda\ge1}.
\]
The quotient $\Lambda(\lambda)/\lambda$, $\lambda>0$, has a unique
minimiser $\lambda_\parallel$, and its positive minimum is
$c_\parallel$ in~\eqref{eq:longitudinal-intro}.
\end{lemma}

\begin{proof}
Since $\mathcal V(s)\to0$ as $s\to\infty$, multiplication by
$D\lambda^2\mathcal V$ is bounded and relatively compact in the sense
proved in Lemma~\ref{lem:far-tests}; thus the upper endpoint of the essential spectrum remains $\mu$.
The Rayleigh principle gives $\Lambda(\lambda)\ge\sigma>0\ge\mu$.
The spectral theorem and the positivity and ODE argument of
Lemma~\ref{lem:plane-operator} give an isolated simple eigenvalue
with a positive eigenfunction. Bounded analytic
perturbation on $\operatorname{Dom}(\J)$ in~\eqref{eq:plane-Schrodinger} gives real analytic
dependence under positive normalisation; the potential $D\lambda^2\mathcal V$
makes $\Lambda$ even.

Differentiate the normalised eigenvalue equation and take its
$L^2(0,\infty)$ inner product with $\psi_\lambda$.
This proves the formula for $\Lambda'$ in~\eqref{eq:longitudinal-HF}, while differentiation
of the normalisation gives
$\langle\psi_\lambda,\partial_\lambda\psi_\lambda\rangle=0$.
Put $\overline{\mathcal V}_\lambda=\int_0^\infty \mathcal V(s)\psi_\lambda(s)^2\dd s$ and
$g_\lambda=(\mathcal V-\overline{\mathcal V}_\lambda)\psi_\lambda$.
The inverse operator
\[
 \mathscr R_\lambda=
 (\Lambda(\lambda)-\J_\parallel(\lambda))^{-1}
       \quad\hbox{on }\{\psi_\lambda\}^{\perp}
\]
is bounded and self-adjoint. Indeed,
$\delta_\lambda=\Lambda(\lambda)-
\sup\sigma(\J_\parallel(\lambda)|_{\{\psi_\lambda\}^\perp})>0$,
so the spectral theorem gives
$\langle g,\mathscr R_\lambda g\rangle\ge0$ for
$g\perp\psi_\lambda$ and
$\|\mathscr R_\lambda\|\le\delta_\lambda^{-1}$.
The differentiated equation therefore gives
$\partial_\lambda\psi_\lambda=2D\lambda\mathscr R_\lambda g_\lambda$,
and a second differentiation of~\eqref{eq:longitudinal-HF} yields
\[\label{eq:longitudinal-convexity}
 \Lambda''(\lambda)
 =2D\overline{\mathcal V}_\lambda
       +8D^2\lambda^2\langle g_\lambda,
                             \mathscr R_\lambda g_\lambda\rangle>0\quad {\text{for }\lambda\in\R}.
\]
The inner products are taken in $L^2(0,\infty)$, and
$\overline{\mathcal V}_\lambda>0$ makes the inequality strict.

For the lower large-parameter bound, choose smooth $\eta$
supported in $[0,1)$, with $\eta'(0)=0$ and $\int_0^\infty\eta(s)^2\dd s=1$,
and test with $w_\lambda(s)=\lambda^{1/4}\eta(\sqrt\lambda s)$.
A change of variable gives
$\int_0^\infty|w_\lambda'|^2\dd s
=\lambda\int_0^\infty|\eta'|^2\dd s$.
Since $\mathcal V(s)=1-\tanh^2s\ge1-s^2$,
$D\lambda^2\int_0^\infty \mathcal V|w_\lambda|^2\dd s
\ge D\lambda^2-D\lambda\int_0^\infty s^2\eta(s)^2\dd s$.
The function $a(s)-D/4-D\mathcal V(s)/4$ is bounded below, proving the first
inequality in~\eqref{eq:longitudinal-growth}; the second follows
by omitting the nonpositive term
$-D\int_0^\infty|w'|^2\dd s$ from the Rayleigh quotient
and using $\mathcal V\le1$.

By \eqref{eq:longitudinal-growth} and $\Lambda(0)=\sigma>0$,
the quotient $\Lambda(\lambda)/\lambda$ tends to infinity
at both ends of $(0,\infty)$ and attains an interior minimum.
For $G(\lambda)=\lambda\Lambda'(\lambda)-\Lambda(\lambda)$,
\eqref{eq:longitudinal-convexity} gives
$G'(\lambda)=\lambda\Lambda''(\lambda)>0$ and
$G(0+)=-\sigma<0$. Thus $G$ has exactly one zero, at the minimiser of
$\Lambda(\lambda)/\lambda$ on $(0,\infty)$. Since $\Lambda>0$, the attained
minimum is positive, proving the characterisation of
$c_\parallel$.
\end{proof}

\begin{lemma}\label{lem:longitudinal-upper}
For every $\lambda>0$ and compactly supported $0\le u_0\le1$, the
solution of~\eqref{eq:KPP} with $f=f_\G$ satisfies
\[\label{eq:longitudinal-exponential}
 \begin{gathered}
 u(x,t)\le C_{\lambda,u_0}
       e^{\Lambda(\lambda)t-\lambda|\ell(x)|}h_\lambda(s(x))
 \quad {\text{for }x\in\Hh^2,\quad t>0},\\
 h_\lambda(s)=\frac{\psi_\lambda(s)}{\sqrt{\cosh s}}
 \quad {\text{for }\lambda>0,\quad s\ge0}.
 \end{gathered}
\]
In particular, \eqref{eq:longitudinal-outer} holds, with exponential
decay for every fixed $c>c_\parallel$.
\end{lemma}

\begin{proof}
By Lemma~\ref{lem:longitudinal-branch},
$\psi_\lambda\in H^2(0,\infty)$ solves
$\J_\parallel(\lambda)\psi_\lambda=\Lambda(\lambda)\psi_\lambda$.
Thus $\psi_\lambda$ and $h_\lambda$ are bounded, the compositions $\psi_\lambda(|\rho|)$ and $h_\lambda(|\rho|)$
are $C^2$ at $\rho=0$ by the Neumann condition, and
\eqref{eq:plane-Fermi} gives
\[\label{eq:separated-eigenfunctions}
 (D\Delta+a(s))h_\pm=\Lambda(\lambda)h_\pm,
 \qquad h_\pm(x)=e^{\mp\lambda\ell(x)}h_\lambda(s(x))\quad {\text{for }x\in\Hh^2,\quad\lambda>0}.
\]
If $K$ is a compact set containing the support of $u_0$,
positivity gives $\min_Kh_\pm>0$.
Taking
$C=\max\{(\min_Kh_+)^{-1},(\min_Kh_-)^{-1}\}$ therefore
ensures $u_0\le Ch_\pm$.

For Brownian motion $X_t$ with generator $D\Delta$ and $X_0=x$,
It\^o's formula
and \eqref{eq:separated-eigenfunctions} make
\[
 \exp\left(\int_0^t a(s(X_v))\dd v-\Lambda(\lambda)t\right)h_\pm(X_t)
\]
a positive local martingale.
Stopping on increasing compact sets and applying Fatou's lemma
gives
$e^{t(D\Delta+a(s))}h_\pm\le e^{\Lambda(\lambda)t}h_\pm$.
By KPP linear comparison and $u_0\le Ch_\pm$, the minimum of the
two bounds, using
$\min\{e^{-\lambda\ell},e^{\lambda\ell}\}=e^{-\lambda|\ell|}$,
is \eqref{eq:longitudinal-exponential}.
Setting $\lambda=\lambda_\parallel$ and using
$\Lambda(\lambda_\parallel)=c_\parallel\lambda_\parallel$
then bounds the solution on $|\ell(x)|\ge ct$ by
$C\norm{h_{\lambda_\parallel}}_\infty
 e^{-\lambda_\parallel(c-c_\parallel)t}$.
This proves \eqref{eq:longitudinal-outer} and its exponential rate.
\end{proof}

\paragraph{Finite tubes and lower solutions.}

For $R>0$, let $\J_{R,\parallel}(\lambda)$ be the realisation
on $L^2(0,R;\dd s)$ of $\J+D\lambda^2\sech^2s$ from
\eqref{eq:longitudinal-intro}. Its operator domain is
$\{w\in H^2(0,R):w'(0)=0,\ w(R)=0\}$, so the boundary
conditions are
\[\label{eq:finite-tube-boundary}
 w'(0)=0,\qquad w(R)=0.
\]
Denote its principal eigenvalue by $\Lambda_R(\lambda)$ and
its positive eigenfunction by $\psi_{R,\lambda}$, normalised by
$\int_0^R\psi_{R,\lambda}^2\dd s=1$. Thus
$\J_{R,\parallel}(\lambda)\psi_{R,\lambda}
=\Lambda_R(\lambda)\psi_{R,\lambda}$ with
\eqref{eq:finite-tube-boundary}. Put $\sigma_R=\Lambda_R(0)$. Here $\sigma>0$ by the hypothesis of Theorem~\ref{thm:longitudinal},
and Lemma~\ref{lem:finite-tube-speeds} proves $\sigma_R\uparrow\sigma$.
For the speed estimates fix $R_0\ge1$ with $\sigma_{R_0}>0$ and take
$R\ge R_0$, so $\sigma_R\ge\sigma_{R_0}>0$.

\begin{lemma}\label{lem:finite-tube-speeds}
As $R\to\infty$, $\sigma_R\uparrow\sigma$ and
$\Lambda_R(\lambda)\uparrow\Lambda(\lambda)$ for every real $\lambda$.
The convergence $\Lambda_R\to\Lambda$ is uniform on compact
intervals of $\lambda\in\R$. For $R\ge R_0$, define
\[\label{eq:finite-tube-speed}
 c_R=\min_{\lambda>0}\frac{\Lambda_R(\lambda)}\lambda.
\]
Then $0<c_R\uparrow c_\parallel$ as $R\to\infty$.
For such a fixed $R$, if $0<c<c_R$, then
\[\label{eq:finite-moving-gap}
 \min_{\lambda\in\R}[\Lambda_R(\lambda)-v\lambda]>0
 \quad\text{for }-c\le v\le c,
\]
with a lower bound independent of $v$ in this interval.
\end{lemma}

\begin{proof}
The finite-interval form domain is
$\{w\in H^1(0,R):w(R)=0\}$.
For $w\in H^1(0,R)$ with $w(R)=0$, extend $w$ by $0$ on
$[R,\infty)$. This extension belongs to $H^1(0,\infty)$, and
the resulting subspaces increase with $R$. Since compactly supported form tests are dense on
the half-line, the Rayleigh principle gives
$\Lambda_R(\lambda)\uparrow\Lambda(\lambda)$ and
$\sigma_R\uparrow\sigma$.
Since $0<\mathcal V\le1$,
$|\Lambda_R(\lambda)-\Lambda_R(\nu)|\le D|\lambda^2-\nu^2|$;
hence the convergence is uniform on compact parameter intervals.

We first bound the minimising values of $\lambda$ for $\Lambda_R(\lambda)/\lambda$ in~\eqref{eq:finite-tube-speed} independently of $R$. For $R\ge R_0$
and $\lambda>0$, the Rayleigh principle gives
$\Lambda_R(\lambda)\ge\sigma_R\ge\sigma_{R_0}>0$.
For $\lambda\ge1$ the tests in~\eqref{eq:longitudinal-growth} are supported
in $[0,R_0)$ and give $\Lambda_R(\lambda)\ge D\lambda^2-C\lambda-C$
with $C$ independent of $R$. Hence $\Lambda_R(\lambda)/\lambda$
tends to infinity at both ends of $(0,\infty)$ uniformly in $R\ge R_0$,
whereas $\Lambda_R(1)\le\Lambda(1)$. Choose
$J=[\underline\lambda,\overline\lambda]\subset(0,\infty)$ containing
$\lambda_\parallel$ so that $\Lambda_R(\lambda)/\lambda>\Lambda(1)$
for $\lambda\notin J$, $\lambda>0$, and $R\ge R_0$.
The quotient therefore attains its minimum on $J$ for every $R\ge R_0$.
Write $\lambda_R$ for any minimiser.

We prove \eqref{eq:finite-tube-eigenvalue-rate} for the difference $\Lambda(\lambda)-\Lambda_R(\lambda)$ on an arbitrary compact interval
$I\subset\R$, so that it can in particular be applied with $I=J$.
Put $K=\max_{\lambda\in I}|\lambda|$
and $M_0=\|\phiZ'/\phiZ\|_\infty$, and choose
$\kappa>0$, $S>1$ so that
\[\label{eq:tube-tail-choices}
 D(\kappa^2+2\kappa M_0)\le\frac\sigma4,\qquad
 D(K^2+1/4)\sech^2S\le\frac\sigma4.
\]
The equation
$\J_\parallel(\lambda)\psi_\lambda=\Lambda(\lambda)\psi_\lambda$
and $\int_0^\infty\psi_\lambda^2\dd s=1$ give
$\sup_{\lambda\in I}\|\psi_\lambda\|_{H^2(0,\infty)}<\infty$:
the potential and eigenvalues are bounded, so the equation
bounds $\psi_\lambda''$ in $L^2$, and the quadratic form
bounds $\psi_\lambda'$.
In particular, these eigenfunctions are uniformly bounded, and
each tends to zero at infinity.
For $Y_\lambda=\psi_\lambda/\phiZ$, the eigenvalue and
periodic cell equations give, on $s\ge S$,
\[\label{eq:transverse-tail-quotient}
 \begin{gathered}
 D\left(Y_\lambda''+2\frac{\phiZ'}{\phiZ}Y_\lambda'\right)
       -d_\lambda(s)Y_\lambda=0,\\
 d_\lambda(s)=\Lambda(\lambda)-\mu-D(\lambda^2-1/4)\mathcal V(s)
                 \ge\frac{3\sigma}4
 \quad {\text{for }s\ge S,\quad\lambda\in I}.
 \end{gathered}
\]
By \eqref{eq:tube-tail-choices}, applying
$D\partial_{ss}+2D(\phiZ'/\phiZ)\partial_s-d_\lambda(s)$
to $e^{-\kappa s}$ gives at most
$-(\sigma/2)e^{-\kappa s}$.
Choose $C_I$ so that $C_Ie^{-\kappa S}\ge Y_\lambda(S)$
for all $\lambda\in I$. Since $Y_\lambda(s)$ and
$C_Ie^{-\kappa s}$ tend to $0$ as $s\to\infty$, the negative
zeroth-order coefficient in~\eqref{eq:transverse-tail-quotient}
excludes a positive maximum of $Y_\lambda-C_Ie^{-\kappa s}$.
Hence $\psi_\lambda(s)\le C_Ie^{-\kappa s}$ for $s\ge S$.
The eigenvalue equation gives
$|\psi_\lambda''(s)|\le C_Ie^{-\kappa s}$ on $s\ge S$,
uniformly for $\lambda\in I$. The mean value theorem on unit
intervals and integration then give
$|\psi_\lambda'(s)|\le C_Ie^{-\kappa s}$ for $s\ge S+1$. Enlarging $C_I$ on $[0,S+1]$ gives
\[\label{eq:transverse-uniform-tail}
 \sup_{\lambda\in I}
 \bigl(|\psi_\lambda(s)|+|\psi_\lambda'(s)|\bigr)
       \le C_Ie^{-\kappa s}\quad\text{for }s\ge0.
\]

For $R>S+2$, choose a smooth $0\le\chi_R\le1$, equal to
one on $[0,R-1]$, zero at $R$, with $|\chi_R'|\le C$.
For $0\ne w\in H^1(0,R)$ with $w(R)=0$, define
\[\label{eq:finite-tube-Rayleigh-quotient}
 \mathscr Q_{R,\lambda}(w)
 =\frac{\int_0^R\bigl([a-D/4+D(\lambda^2-1/4)\mathcal V]|w|^2-D|w'|^2\bigr)\dd s}
        {\int_0^R|w|^2\dd s}.
\]
Multiplying the eigenvalue equation by $\chi_R^2\psi_\lambda$
and integrating by parts over $(0,R)$ gives
\[\label{eq:tube-cutoff-Rayleigh}
 \Lambda(\lambda)-\mathscr Q_{R,\lambda}(\chi_R\psi_\lambda)
 =\frac{D\int_0^R|\chi_R'|^2\psi_\lambda^2\dd s}
        {\int_0^R\chi_R^2\psi_\lambda^2\dd s}.
\]
Since $\chi_R=1$ on $[0,R-1]$ and $\|\psi_\lambda\|_2=1$,
\eqref{eq:transverse-uniform-tail} gives
$\int_0^R\chi_R^2\psi_\lambda^2\dd s
\ge1-(C_I^2/(2\kappa))e^{-2\kappa(R-1)}$.
Fix $R_I>S+2$ with $(C_I^2/(2\kappa))e^{-2\kappa(R_I-1)}\le1/2$.
The denominator of~\eqref{eq:tube-cutoff-Rayleigh} is then at least
$1/2$ for every $R\ge R_I$ and $\lambda\in I$.
Since $\chi_R'$ is supported in $[R-1,R]$ and is uniformly bounded,
the numerator is at most $C e^{-2\kappa R}$ uniformly on $I$.
The Rayleigh principle therefore yields
\[\label{eq:finite-tube-eigenvalue-rate}
 0\le\Lambda(\lambda)-\Lambda_R(\lambda)
       \le C_Ie^{-2\kappa R}
       \quad\text{for }\lambda\in I,\ R\ge R_I.
\]
If $\Lambda_R(\lambda)=\Lambda(\lambda)$ for finite $R$,
the extension of $\psi_{R,\lambda}$ by $0$ on $[R,\infty)$
would attain the half-line Rayleigh maximum. By the spectral theorem, it would
then be a nonzero multiple of $\psi_\lambda$ by the simplicity in
Lemma~\ref{lem:longitudinal-branch}. This contradicts its vanishing on
$(R,\infty)$, where $\psi_\lambda>0$.
Thus $\Lambda_R(\lambda)<\Lambda(\lambda)$ for finite $R$
\cite{ReedSimon}.

For $R_0\le R\le R'$, zero extension in the Rayleigh quotient gives
$\Lambda_R(\lambda)\le\Lambda_{R'}(\lambda)\le\Lambda(\lambda)$.
Dividing by $\lambda>0$ and minimising gives $c_R\le c_{R'}\le c_\parallel$.
Since $\lambda_R\in J$ and $\underline\lambda>0$,
$0\le c_\parallel-c_R\le
\underline\lambda^{-1}\sup_{\lambda\in J}|\Lambda(\lambda)-\Lambda_R(\lambda)|$.
For each $\lambda\in J$, the Rayleigh principle gives
$\Lambda_R(\lambda)\uparrow\Lambda(\lambda)$ as shown at the start of the proof.
Since $J$ is compact and
$|\Lambda_R(\lambda)-\Lambda_R(\nu)|\le D|\lambda^2-\nu^2|$ for $\lambda,\nu\in J$,
we have $\sup_{\lambda\in J}|\Lambda(\lambda)-\Lambda_R(\lambda)|\to0$.
Hence $c_R\to c_\parallel$.
In~\eqref{eq:finite-tube-eigenvalue-rate} take $I=J$, with
$K,\kappa,C_I,R_I$ evaluated for this interval, and put $R_0'=\max\{R_0,R_I\}$.
We obtain
\[\label{eq:finite-tube-speed-rate}
 0<c_\parallel-c_R\le Ce^{-2\kappa R},\qquad
 |\lambda_R-\lambda_\parallel|\le Ce^{-\kappa R}
                       \quad\text{for }R\ge R_0'.
\]
Indeed, for $q(\lambda)=\Lambda(\lambda)/\lambda$,
uniqueness of the minimum and
$q''(\lambda_\parallel)=\Lambda''(\lambda_\parallel)/
\lambda_\parallel>0$ imply
$q(\lambda)-c_\parallel\ge\eta|\lambda-\lambda_\parallel|^2$
on $[\underline{\lambda},\overline{\lambda}]$ for some $\eta>0$, while
\[\label{eq:finite-tube-minimizer-bound}
 0\le q(\lambda_R)-c_\parallel
 \le \frac{\Lambda(\lambda_R)-\Lambda_R(\lambda_R)}{\lambda_R}
 \le Ce^{-2\kappa R}.
\]
Combining the quadratic lower bound for $q-c_\parallel$ with
\eqref{eq:finite-tube-minimizer-bound} gives
\eqref{eq:finite-tube-speed-rate}. The strict inequality
$c_R<c_\parallel$ follows by testing at $\lambda_\parallel$.
In~\eqref{eq:finite-tube-speed-rate}, $C$ and $\kappa$ are independent of
$R\ge R_0'$. The tail estimate \eqref{eq:transverse-uniform-tail} uses \eqref{eq:tube-tail-choices}
with $K=\max_{\lambda\in J}|\lambda|$. The estimate for $\lambda_R$ in~\eqref{eq:finite-tube-speed-rate} also uses the quadratic lower bound for
$q-c_\parallel$ on $J$.

For each fixed $R\ge R_0$ and every $\lambda\in\R$, the principal
eigenvalue $\Lambda_R(\lambda)$ on the spatial interval $(0,R)$ with
boundary conditions \eqref{eq:finite-tube-boundary} is simple.
Differentiation as in~\eqref{eq:longitudinal-convexity} gives
$\Lambda_R''(\lambda)\ge2D\sech^2 R>0$. Integration on either side
of zero yields $\Lambda_R'(\R)=\R$.
Every minimiser of $\Lambda_R(\lambda)/\lambda$ on $(0,\infty)$ satisfies
$\lambda\Lambda_R'(\lambda)-\Lambda_R(\lambda)=0$.
The derivative of the left-hand side is $\lambda\Lambda_R''(\lambda)>0$,
so the minimiser $\lambda_R$ is unique.
For $|v|\le c<c_R$ and $\lambda\ne0$, evenness gives
$\Lambda_R(\lambda)-v\lambda
\ge\Lambda_R(|\lambda|)-c|\lambda|>0$;
at zero the expression is $\sigma_R>0$.
Quadratic growth and compactness of the velocity interval then
prove the uniform gap \eqref{eq:finite-moving-gap}.
\end{proof}

\begin{lemma}\label{lem:tangential-bumps}
Fix $R\ge R_0$ and $0<c<c_R$, with $c_R>0$ by Lemma~\ref{lem:finite-tube-speeds}. There are $A,\kappa>0$
and compactly supported weak subsolutions $w_v(x,t)$ of~\eqref{eq:KPP} with $f=f_\G$, one for each $v\in[-c,c]$,
with the following properties. Their supports are contained
in
\[\label{eq:tangential-bump-support}
 \{x\in\Hh^2:|\rho(x)|\le R,\quad |\ell(x)-vt|\le A\},
\]
their values lie in $[0,1]$, and for a fixed $0<R_1<R$,
\[\label{eq:tangential-bump-lower}
 w_v(x,t)\ge\kappa
 \quad\text{for }x\in\Hh^2,\ t\in\R,\ v\in[-c,c]
 \text{ with }|\rho(x)|\le R_1,\ |\ell(x)-vt|\le1.
\]
The constants $A,\kappa,R_1$ are independent of $v\in[-c,c]$
and $t$. For every $0<\theta\le1$, the functions $\theta w_v$
retain these properties with lower bound $\theta\kappa$ in~\eqref{eq:tangential-bump-lower}.
\end{lemma}

\begin{proof}
By Lemma~\ref{lem:finite-tube-speeds},
Lemma~\ref{lem:complex-root} applies to $K=\Lambda_R$ for every
velocity $v\in[-c,c]$, with the positive gap in~\eqref{eq:finite-moving-gap}
uniform in $v$.
For sufficiently small $\eta>0$, continue $\psi_{R,\lambda}$
analytically to $\zeta_v=\alpha(v,\eta)+i\eta$, with
$\alpha$ from Lemma~\ref{lem:complex-root}, using
$\int_0^R\psi_{R,\zeta}(s)^2\dd s=1$ and the branch
positive for real parameters.
The operator domain is the fixed space of $H^2(0,R)$ functions
satisfying \eqref{eq:finite-tube-boundary}.
Let $\lambda_v$ solve $\Lambda_R'(\lambda_v)=v$.
The parameters $\lambda_v$, $v\in[-c,c]$, range over a compact
interval, and the Hopf lemma gives
$\sup_{v\in[-c,c]}\psi_{R,\lambda_v}'(R)<0$.
Together with interior positivity and $C^1$ analytic dependence,
this gives constants $d_R,C_R>0$ such that, for $0\le s\le R$ and $v\in[-c,c]$,
\[\label{eq:finite-eigenfunction-boundary-ratio}
 \begin{aligned}
  \psi_{R,\lambda_v}(s)&\ge d_R(R-s),\\
  |\psi_{R,\zeta_v}(s)-\psi_{R,\lambda_v}(s)|
       &\le C_R\eta(R-s).
 \end{aligned}
\]
Indeed, $C^1$ dependence gives
$\|\psi_{R,\zeta_v}'-\psi_{R,\lambda_v}'\|_\infty\le C_R\eta$
uniformly in $v\in[-c,c]$. Integrating from $s$ to $R$ and using
$\psi_{R,\zeta_v}(R)=\psi_{R,\lambda_v}(R)=0$ gives the second inequality. Thus \eqref{eq:finite-eigenfunction-boundary-ratio}
gives $\psi_{R,\zeta_v}/\psi_{R,\lambda_v}=1+O(\eta)$,
including the continuous extension of this ratio at $R$.
Choose the continuous argument of
$\psi_{R,\zeta_v}/\psi_{R,\lambda_v}$ equal to $0$ at $\eta=0$;
with this choice,
$|\arg(\psi_{R,\zeta_v}(s)/\psi_{R,\lambda_v}(s))|\le C_R\eta$
for $s\in[0,R]$ and $v\in[-c,c]$.
The constants $d_R,C_R$ depend on $R$ and $[-c,c]$;
see~\cite[Chapter~3]{GilbargTrudinger} for the Hopf principle.

On $|\rho|<R$, put
\[\label{eq:complex-tangential-bump}
 Z_v(\rho,\ell,t)=
 e^{-\zeta_v(\ell-vt)}
            \frac{\psi_{R,\zeta_v}(|\rho|)}{\sqrt{\cosh\rho}}\quad {\text{for }|\rho|<R,\quad\ell,t\in\R,\quad v\in[-c,c]}.
\]
The condition $\psi_{R,\zeta_v}'(0)=0$ makes
$Z_v(\rho,\ell,t)$ of class $C^2$ across $\rho=0$,
and $\psi_{R,\zeta_v}(R)=0$ gives $Z_v(\pm R,\ell,t)=0$.
By \eqref{eq:plane-Fermi}, the eigenvalue equation and
\eqref{eq:complex-dispersion},
\[
 (\partial_t-D\Delta-a(|\rho|))Z_v=-E_v Z_v,
 \qquad E_v\ge E_*>0.
\]
Here $E_v=E(v,\eta)$ from \eqref{eq:complex-dispersion}
with $K=\Lambda_R$. Define $\varphi_v=(\operatorname{Re}Z_v)_+$ on
$\{|\rho|<R,\ |\ell-vt|<\pi/\eta\}$ and set $\varphi_v=0$
outside this set.
The $O(\eta)$ bound on the argument gives
$\operatorname{Re}Z_v<0$ near $\ell-vt=\pm\pi/\eta$, as in
Lemma~\ref{lem:periodic-moving-bumps}. At $\rho=\pm R$,
each connected component of $\{\varphi_v>0\}$ has nonpositive
outward normal derivative along this boundary,
so extension by zero adds a nonnegative boundary measure to
$D\Delta\varphi_v$.
Thus $\varphi_v$ is a weak subsolution of
$(\partial_t-D\Delta-a(|\rho|))\varphi_v\le-E_v\varphi_v$.
For sufficiently small fixed $\eta$, $\sup_{v,x,t}\varphi_v\le M$
and $\varphi_v(\rho,\ell,t)\ge\kappa_0>0$ whenever
$|\rho|\le R_1$, $|\ell-vt|\le1$.

By \eqref{eq:uniform-taylor},
$(\partial_t-D\Delta)(\delta\varphi_v)-F(|\rho|,\delta\varphi_v)$
is at most $-E_*\delta\varphi_v+C_2\delta^2\varphi_v^2$.
Choose $\delta>0$ sufficiently small that this is nonpositive
and $\delta M\le1$. Then $w_v=\delta\varphi_v$, $A=\pi/\eta$ and
$\kappa=\delta\kappa_0$ satisfy
\eqref{eq:tangential-bump-support}--\eqref{eq:tangential-bump-lower},
and $w_v$ is a weak subsolution of~\eqref{eq:KPP} with $f=f_\G$.
Decreasing $\delta$ preserves
$(\partial_t-D\Delta)(\delta\varphi_v)-F(|\rho|,\delta\varphi_v)\le0$
and replaces the lower bound by $\delta\kappa_0$ in~\eqref{eq:tangential-bump-lower}.
\end{proof}

\begin{lemma}\label{lem:tangential-positive}
Let $u_0$ be nonzero and compactly supported. For each
$0<c<c_\parallel$, there are $R_1,\eta,T>0$ such that
\[\label{eq:tangential-positive}
 u(x,t;u_0,f_\G)\ge\eta
 \quad\text{for }x\in\Hh^2,\ t\ge T
 \text{ with }|\ell(x)|\le ct,\ s(x)\le R_1.
\]
\end{lemma}

\begin{proof}
By Lemma~\ref{lem:finite-tube-speeds}, choose a finite $R$
and $c_1$ with $c<c_1<c_R$.
Let $A$ be given by Lemma~\ref{lem:tangential-bumps} for
$v\in[-c_1,c_1]$. By \eqref{eq:tangential-bump-support},
all $w_v(\cdot,0)$ are supported in
$\{x\in\Hh^2:|\rho(x)|\le R,\ |\ell(x)|\le A\}$.
Lemma~\ref{lem:wellposed} gives a positive minimum of
$u(x,1;u_0,f_{\G})$ on this compact set.
Decrease the common amplitude so that
$w_v(x,0)\le u(x,1;u_0,f_\G)$ for all $v\in[-c_1,c_1]$.
Comparison then gives $w_v(x,t-1)\le u(x,t;u_0,f_\G)$ for $t\ge1$.
For $t>1$ and $|\ell|\le c_1(t-1)$, choose $v=\ell/(t-1)$ and keep it
fixed on $[1,t]$. Then $v(t-1)=\ell$, so \eqref{eq:tangential-bump-lower} gives
$w_v(\rho,\ell,t-1)\ge\kappa$ for $|\rho|\le R_1$.
Since $ct\le c_1(t-1)$ for large $t$, this proves
\eqref{eq:tangential-positive}.
\end{proof}

\begin{lemma}\label{lem:tangential-lower-profile}
For a nonzero compactly supported datum, \eqref{eq:longitudinal-inner}
holds for each $0<c<c_\parallel$. Under \eqref{eq:quotient},
\eqref{eq:longitudinal-profile} also holds.
\end{lemma}

\begin{proof}
Choose $c<c_1<c_\parallel$ and use
Lemma~\ref{lem:tangential-positive} to obtain a common positive
bound $\eta>0$ on
$\{|\ell|\le c_1t,\ |\rho|\le R_1,\ t\ge T\}$.
Fix a nonzero continuous compactly supported $v_0$ with
$0\le v_0\le\eta$ and support in
$|\rho|<R_1$, $|\ell|<A$, where $A<\infty$ is fixed.
Let $v$ solve \eqref{eq:KPP} with $f=f_{\G}$ and
$v(x,0)=v_0(x)$, and write it in the Fermi coordinates
\eqref{eq:plane-Fermi} as $v(\rho,\ell,t)$.
For any $S<\infty$ and $\varepsilon>0$,
Lemma~\ref{lem:minimal-state} gives a time $T_0$ such that
\[\label{eq:seed-to-profile}
 v(\rho,0,T_0)\ge p(|\rho|)-\varepsilon
                \quad\text{for }|\rho|\le S.
\]
Let $T$ be the starting time in~\eqref{eq:tangential-positive}
for speed $c_1$, and choose $T_\varepsilon$ so that
$t-T_0\ge T$ and $(c_1-c)t>c_1T_0+A$ for
$t\ge T_\varepsilon$.
For $|\ell_0|\le ct$, the translate $v_0(\rho,\ell-\ell_0)$ lies
below $u(\rho,\ell,t-T_0)$ by these choices. Since translations of $\ell$ preserve the reaction and
\eqref{eq:plane-Fermi}, comparison on $[t-T_0,t]$ gives
$u(\rho,\ell,t)\ge v(\rho,\ell-\ell_0,T_0)$,
in particular $u(\rho,\ell_0,t)\ge v(\rho,0,T_0)$.
By \eqref{eq:seed-to-profile}, the lower bound is uniform for
$|\rho|\le S$ and $|\ell_0|\le ct$.
For $|\rho|>S$, the estimate $p\le Ch$ in~\eqref{eq:localized-bounds} gives $p(|\rho|)<\varepsilon$ when
$S$ is large. Since $u\ge0$, we then have
$u(\rho,\ell,t)\ge p(|\rho|)-\varepsilon$ in this region.
Letting $\varepsilon\downarrow0$ proves
\eqref{eq:longitudinal-inner} uniformly on $\{x\in\Hh^2:|\ell(x)|\le ct\}$.

Under \eqref{eq:quotient}, choose $R$ so that
$U_R(s(x))\ge u_0(x)$ by
Lemma~\ref{lem:cap} and let $v_R$ solve \eqref{eq:KPP} with
$f=f_{\G}$ and initial value $U_R(s(x))$.
Since $U_R\le Ch$, Lemma~\ref{lem:localized-convergence} with $M=0$ gives
$\sup_{\rho\in\R}|v_R(\rho,0,t)-p(|\rho|)|\to0$ as $t\to\infty$.
Uniqueness and invariance under translations of $\ell$ give
$v_R(\rho,\ell,t)=v_R(\rho,0,t)$; hence $u\le v_R$ implies
\[
 \limsup_{t\to\infty}\sup_{x\in\Hh^2}
                  [u(x,t)-p_*(x)]\le0.
\]
Combining this upper bound with
\eqref{eq:longitudinal-inner} gives
\eqref{eq:longitudinal-profile}.
\end{proof}

\begin{proof}[Proof of Theorem~\ref{thm:longitudinal}(i)--(iii) and (v)]
By Lemma~\ref{lem:longitudinal-branch}, the function $\Lambda$ has the properties in~\eqref{eq:longitudinal-HF} and the characterisation
$c_\parallel=\min_{\lambda>0}\Lambda(\lambda)/\lambda$.
Lemma~\ref{lem:longitudinal-upper} proves
\eqref{eq:longitudinal-outer}, and
Lemma~\ref{lem:tangential-lower-profile} proves
\eqref{eq:longitudinal-inner} and, under \eqref{eq:quotient},
\eqref{eq:longitudinal-profile}.
Part~(iv) is proved in Section~\ref{sec:axial-delay} using the transverse diffusion constructed below.
\end{proof}

\paragraph{The transverse diffusion.}

We derive a formula for $\Lambda''(\lambda)$ using $\psi_\lambda$
and construct the martingale in
Proposition~\ref{prop:longitudinal-martingale} by the method for
stationary reversible processes \cite{KV86}. For a positive probability
density $\pi$ on $[0,\infty)$, the weighted form domain is
$H^1(\pi\dd s)=\{g\in H^1_{\rm loc}([0,\infty)):
{\int_0^\infty}(|g|^2+|g'|^2)\pi\dd s<\infty\}$; it imposes no zero trace at the
reflecting endpoint.

\begin{lemma}\label{lem:transverse-corrector}
Fix $\lambda\in\R$. Define
\[\label{eq:transverse-generator}
 \pi_\lambda^{\perp}(s)=\psi_\lambda(s)^2,\qquad
 A_\lambda=D\partial_{ss}
                +2D\frac{\psi_\lambda'}{\psi_\lambda}\partial_s,
 \qquad \chi_\lambda=-\partial_\lambda\log\psi_\lambda.
\]
Then ${\int_0^\infty}\pi_\lambda^{\perp}{\dd s}=1$, ${\int_0^\infty}\chi_\lambda\pi_\lambda^{\perp}{\dd s}=0$,
$\chi_\lambda'(0)=0$, and
\[\label{eq:corrector-equation}
 A_\lambda\chi_\lambda=2D\lambda \mathcal V-\Lambda'(\lambda)
 \quad {\text{for }s>0}.
\]
We have $\chi_\lambda\in L^2(\pi_\lambda^{\perp}\dd s)$ and
$\chi_\lambda'\in L^2(\pi_\lambda^{\perp}\dd s)$, with
\[\label{eq:variance-energy}
 \Lambda''(\lambda)=2D\int_0^\infty
                 [\mathcal V+|\chi_\lambda'|^2]\pi_\lambda^{\perp}\dd s>0
 \quad {\text{for }\lambda\in\R}.
\]
There is $\delta_\lambda>0$ such that
\[\label{eq:transverse-Poincare}
 D{\int_0^\infty} |g'|^2\pi_\lambda^{\perp}\dd s
 \ge\delta_\lambda{\int_0^\infty}|g|^2\pi_\lambda^{\perp}\dd s
 \quad\hbox{if }{\int_0^\infty} g\pi_\lambda^{\perp}\dd s=0
\]
for every $g\in H^1(\pi_\lambda^{\perp}\dd s)$.
\end{lemma}

\begin{proof}
The normalisation in Lemma~\ref{lem:longitudinal-branch} and
its derivative give ${\int_0^\infty}\pi_\lambda^{\perp}{\dd s}=1$ and
${\int_0^\infty}\chi_\lambda\pi_\lambda^{\perp}{\dd s}=0$.
Differentiating $\psi_\lambda'(0)=0$ with respect to $\lambda$
and using $\chi_\lambda=-\partial_\lambda\log\psi_\lambda$ gives
$\chi_\lambda'(0)=0$.
For $q=(\partial_\lambda\psi_\lambda)/\psi_\lambda$,
the differentiated eigenvalue equation is
$A_\lambda q=\Lambda'-2D\lambda \mathcal V$.
Thus $\chi_\lambda=-q$ satisfies
\eqref{eq:corrector-equation}, and
$\chi_\lambda\psi_\lambda=-\partial_\lambda\psi_\lambda$
gives $\int_0^\infty|\chi_\lambda|^2\pi_\lambda^{\perp}\dd s
=\|\partial_\lambda\psi_\lambda\|_{L^2(0,\infty)}^2<\infty$.
By analytic perturbation on the fixed domain,
$\partial_\lambda\psi_\lambda\in \operatorname{Dom}(\J_\parallel(\lambda))=\operatorname{Dom}(\J)$.

To identify the form domains, set
$r=\psi_\lambda'/\psi_\lambda$.
The eigenvalue equation gives $r'=b_\lambda-r^2$, $r(0)=0$,
where
$b_\lambda(s)=(\Lambda(\lambda)-a(s))/D+1/4-(\lambda^2-1/4)\mathcal V(s)$
is bounded and real. Choose $Q>\|b_\lambda\|_\infty$.
Comparison with $\sqrt Q$ bounds $r$ from above; if
$r<-\sqrt Q$, the inequality
$(-r)'\ge(-r)^2-Q$ would force $-r(s)\to\infty$ at a finite value of $s$,
contrary to smoothness and positivity of $\psi_\lambda$ on $[0,\infty)$.
Hence $|r|\le\sqrt Q$.
For $v=\psi_\lambda g$, the identities
$v'=\psi_\lambda g'+rv$ and $\psi_\lambda g'=v'-rv$ therefore
show that $g\in H^1(\pi_\lambda^{\perp}\dd s)$ if and only if
$\psi_\lambda g\in H^1(0,\infty)$.
The self-adjoint realisation of $A_\lambda$ is consequently
\[\label{eq:transverse-generator-domain}
 \begin{aligned}
 \operatorname{Dom}(A_\lambda)&=\{g\in L^2(\pi_\lambda^{\perp}\dd s):
           \psi_\lambda g\in \operatorname{Dom}(\J_\parallel(\lambda))\},\\
 A_\lambda g&=\psi_\lambda^{-1}
          (\J_\parallel(\lambda)-\Lambda(\lambda))(\psi_\lambda g).
 \end{aligned}
\]
Since $\psi_\lambda'(0)=0$ and $\psi_\lambda(0)>0$,
the regular endpoint condition is $g'(0)=0$.
Analytic dependence in the fixed operator domain then gives
$\chi_\lambda\in \operatorname{Dom}(A_\lambda)$.

Applying integration by parts on the half-line with multiplier
$\psi_\lambda$, whose Neumann derivative gives zero boundary term at $0$, to
$v=\partial_\lambda\psi_\lambda$, first with compact
cutoffs and then by form convergence, gives
\[
 D{\int_0^\infty}|\chi_\lambda'|^2\pi_\lambda^{\perp}\dd s
 =\langle\partial_\lambda\psi_\lambda,
       (\Lambda-\J_\parallel)\partial_\lambda\psi_\lambda\rangle
 =-2D\lambda{\int_0^\infty} \mathcal V\chi_\lambda\pi_\lambda^{\perp}\dd s<\infty.
\]
Differentiation of~\eqref{eq:longitudinal-HF} now proves
\eqref{eq:variance-energy}.
Finally, the unitary map $g\mapsto\psi_\lambda g$ from
$L^2((0,\infty),\pi_\lambda^{\perp}\dd s)$ to $L^2(0,\infty)$
conjugates $A_\lambda$ to $\J_\parallel(\lambda)-\Lambda(\lambda)$
by \eqref{eq:transverse-generator-domain}.
The gap between $\Lambda(\lambda)$ and the spectrum of
$\J_\parallel(\lambda)$ on $\{\psi_\lambda\}^\perp$,
from Lemma~\ref{lem:longitudinal-branch}, gives the form
inequality on that subspace.
Transforming back yields \eqref{eq:transverse-Poincare}.
\end{proof}

\begin{proposition}\label{prop:longitudinal-martingale}
Fix $\lambda\in\R$ and put
$h_\lambda(s)=\psi_\lambda(s)/\sqrt{\cosh s}$.
On $[0,\infty)\times\R$, the operator
$\varphi\mapsto e^{\lambda\ell}h_\lambda(s)^{-1}
[D\Delta+a(s)-\Lambda(\lambda)]
(e^{-\lambda\ell}h_\lambda(s)\varphi)$, with reflection at $s=0$,
has differential expression
\[\label{eq:longitudinal-Doob}
 A_\lambda+D\mathcal V(s)(\partial_{\ell\ell}-2\lambda\partial_\ell)\quad {\text{for }s>0,\quad\ell\in\R}.
\]
Let $(\mathsf S_t,\ell_t)$ be the reflected diffusion with generator
\eqref{eq:longitudinal-Doob},
with $\mathsf S_0$ of density $\pi_\lambda^{\perp}$ and fixed $\ell_0\in\R$. Then
\[\label{eq:longitudinal-martingale}
 M_t=\ell_t-\ell_0+\Lambda'(\lambda)t
                    +\chi_\lambda(\mathsf S_t)-\chi_\lambda(\mathsf S_0)\quad {\text{for }t\ge0}
\]
is a square-integrable martingale with
\[\label{eq:martingale-variance}
 \mathbb E M_t^2=\Lambda''(\lambda)t
 \quad {\text{for }t\ge0}.
\]
In particular, the mean longitudinal displacement is
$-\Lambda'(\lambda)t$, and
$\operatorname{Var}(\ell_t-\ell_0)/t\to\Lambda''(\lambda)$ as
$t\to\infty$.
\end{proposition}

\begin{proof}
By \eqref{eq:plane-Fermi} and
$h_\lambda=\psi_\lambda/\sqrt{\cosh s}$, the conjugated operator in Proposition~\ref{prop:longitudinal-martingale}
has transverse drift
$D(\tanh s+2h_\lambda'/h_\lambda)=2D\psi_\lambda'/\psi_\lambda$
and longitudinal drift $-2D\lambda \mathcal V$.
This proves \eqref{eq:longitudinal-Doob}.
Put $r_\lambda=\psi_\lambda'/\psi_\lambda$ and choose
\[\label{eq:transverse-Riccati-constant}
 Q=1+\sup_{s\ge0}\left|
 \frac{\Lambda(\lambda)-a(s)}D+\frac14-(\lambda^2-1/4)\mathcal V(s)\right|.
\]
The eigenvalue equation and the Neumann condition give
\[\label{eq:transverse-Riccati-bound}
 \begin{gathered}
 \left\{\begin{aligned}
 r_\lambda'&=\frac{\Lambda(\lambda)-a}D+\tfrac14
                     -(\lambda^2-\tfrac14)\mathcal V-r_\lambda^2,\\
 r_\lambda(0)&=0,
 \end{aligned}\right.\\
 |r_\lambda(s)|\le\sqrt Q,\qquad
 |r_\lambda'(s)|\le2Q\quad\text{for }s\ge0.
 \end{gathered}
\]
The Riccati comparison in Lemma~\ref{lem:transverse-corrector}
gives $|r_\lambda|\le\sqrt Q$, and the equation then gives
$|r_\lambda'|\le2Q$. On a compact parameter interval $I$, replace $Q$ by
$Q_I=\sup_{\lambda\in I}Q(\lambda)<\infty$, where $Q(\lambda)$
is the right-hand side of~\eqref{eq:transverse-Riccati-constant}.
Thus $|r_\lambda(s)|\le\sqrt{Q_I}$ and
$|r_\lambda'(s)|\le2Q_I$ for $s\ge0$ and $\lambda\in I$.
Thus \eqref{eq:transverse-Riccati-bound} gives the bounded Lipschitz drift
$s\mapsto2D\psi_\lambda'(s)/\psi_\lambda(s)$.

Let $K_\lambda^N(t,s,r)$ be the integral kernel of
$e^{t\J_\parallel(\lambda)}$ relative to $\dd r$, with the
Neumann condition at zero.
The semigroup
$g\mapsto e^{-\Lambda(\lambda)t}\psi_\lambda^{-1}
e^{t\J_\parallel(\lambda)}(\psi_\lambda g)$ has transition density
\[\label{eq:transverse-transition-density}
 p_\lambda^\perp(t,s,r)
   =e^{-\Lambda(\lambda)t}K_\lambda^N(t,s,r)
                       \frac{\psi_\lambda(r)}{\psi_\lambda(s)}\quad {\text{for }t>0,\quad s,r\ge0}.
\]
The identity
$e^{t\J_\parallel(\lambda)}\psi_\lambda
=e^{\Lambda(\lambda)t}\psi_\lambda$ gives
$\int_0^\infty p_\lambda^\perp(t,s,r)\dd r=1$, and symmetry
of $K_\lambda^N$ gives
$\pi_\lambda^{\perp}(s)p_\lambda^\perp(t,s,r)
=\pi_\lambda^{\perp}(r)p_\lambda^\perp(t,r,s)$.
Consequently
$\int_0^\infty\pi_\lambda^{\perp}(s)p_\lambda^\perp(t,s,r)\dd s
=\pi_\lambda^{\perp}(r)$, so $\pi_\lambda^{\perp}(s)\dd s$ is invariant.
The coupled coordinates solve
\[\label{eq:longitudinal-SDE}
 \left\{\begin{aligned}
 \dd \mathsf S_t&=\sqrt{2D}\dd B_t+
       2D\frac{\psi_\lambda'(\mathsf S_t)}{\psi_\lambda(\mathsf S_t)}\dd t+\dd L_t^0,\\
 \dd\ell_t&=\sqrt{2D\mathcal V(\mathsf S_t)}\dd\widetilde B_t
                                  -2D\lambda \mathcal V(\mathsf S_t)\dd t,
 \end{aligned}\right.
\]
Take $\mathsf S_0$ with density $\pi_\lambda^{\perp}$ and standard Brownian motions
$B,\widetilde B$, all mutually independent, and fix $\ell_0$. Use the usual augmentation of the filtration
generated by $\mathsf S_0$, $B$ and $\widetilde B$. The regulator $L_t^0$ is continuous, nondecreasing,
starts at zero and increases only when $\mathsf S_t=0$.
For $T>0$, define $\Gamma:C([0,T])\to C([0,T])$ by
$\Gamma(w)(t)=w(t)+\sup_{0\le v\le t}(-w(v))_+$ for $0\le t\le T$.
The triangle inequality shows that $\Gamma$ is $2$-Lipschitz in
the uniform norm. To solve the first equation of~\eqref{eq:longitudinal-SDE}, we apply Picard iteration on continuous
nonnegative paths starting at $\mathsf S_0$ to
\[\label{eq:reflected-Picard-map}
 \mathsf S_t=\Gamma\left(\mathsf S_0+\sqrt{2D}B+
    2D\int_0^{\,\cdot}\frac{\psi_\lambda'(\mathsf S_v)}
                                  {\psi_\lambda(\mathsf S_v)}\dd v\right)(t)
 \quad\text{for }0\le t\le T.
\]
Let $K$ be a Lipschitz constant of the drift
$2D\psi_\lambda'/\psi_\lambda$ in~\eqref{eq:transverse-Riccati-bound}.
Since $\Gamma$ is $2$-Lipschitz, the map in~\eqref{eq:reflected-Picard-map}
has Lipschitz constant at most $2KT\le1/2$ for $0<T\le1/[4(1+K)]$. Repeating the construction gives a unique
adapted continuous solution for all $t\ge0$, and boundedness of the
drift prevents explosion. The running supremum in $\Gamma$ gives
the regulator $L^0$, with $\dd L^0_t$ supported on $\{t:\mathsf S_t=0\}$. This is the
half-line Skorokhod construction; see~\cite{LionsSznitman}.
The stochastic normalisations are those of~\cite{Hsu,Pinsky}.

Applying It\^o's formula with \eqref{eq:corrector-equation}
and $\chi_\lambda'(0)=0$ cancels both the drift and regulator
in~\eqref{eq:longitudinal-martingale}, giving quadratic
variation
\[
 \langle M\rangle_t=2D\int_0^t
           [\mathcal V(\mathsf S_v)+|\chi_\lambda'(\mathsf S_v)|^2]\dd v.
\]
Put $\tau_n=\inf\{v\ge0:\mathsf S_v\ge n\}$.
It\^o's formula identifies $M_{t\wedge\tau_n}$ with the
stochastic integrals, stopped at $t\wedge\tau_n$, in
\[\label{eq:corrected-martingale-integrals}
 M_t=\sqrt{2D}\int_0^t\chi_\lambda'(\mathsf S_v)\dd B_v
             +\int_0^t\sqrt{2D\mathcal V(\mathsf S_v)}\dd\widetilde B_v
 \quad {\text{for }t\ge0}.
\]
By stationarity and Lemma~\ref{lem:transverse-corrector},
$\mathbb E\int_0^t|\chi_\lambda'(\mathsf S_v)|^2\dd v
=t\int_0^\infty|\chi_\lambda'|^2\pi_\lambda^{\perp}\dd s<\infty$,
and $0<\mathcal V\le1$.
The It\^o isometry therefore gives $L^2$ convergence of the
stopped stochastic integrals. By nonexplosion and continuity, the limit is the process in~\eqref{eq:longitudinal-martingale}, which is therefore a
square-integrable martingale.
The two stochastic integrals in~\eqref{eq:corrected-martingale-integrals}
have zero cross variation because $B$ and $\widetilde B$ are independent.
Taking expectations in the quadratic variation and using
\eqref{eq:variance-energy} proves
\eqref{eq:martingale-variance}.

The mean displacement follows from \eqref{eq:longitudinal-SDE}
and \eqref{eq:longitudinal-HF}.
Set $Z_t=\ell_t-\ell_0+\Lambda'(\lambda)t$ and
$C_t=\chi_\lambda(\mathsf S_t)-\chi_\lambda(\mathsf S_0)$.
We have $\mathbb E Z_t=\mathbb E C_t=0$ and $Z_t=M_t-C_t$.
Stationarity bounds
$\|C_t\|_{L^2(\mathbb P)}\le2\|\chi_\lambda\|_{L^2(\pi_\lambda^{\perp}\dd s)}$,
so Cauchy--Schwarz gives
$|\operatorname{Var}(Z_t)-\mathbb EM_t^2|
\le2\|M_t\|_{L^2(\mathbb P)}\|C_t\|_{L^2(\mathbb P)}
+\|C_t\|_{L^2(\mathbb P)}^2$.
By \eqref{eq:martingale-variance},
$|\operatorname{Var}(Z_t)-\mathbb E M_t^2|\le C\sqrt t+C$.
Dividing by $t$ proves
$\operatorname{Var}(\ell_t-\ell_0)/t\to\Lambda''(\lambda)$
as $t\to\infty$.
\end{proof}

Fix $\lambda\in\R$. For $t>0$ and $s,r\ge0$, write
$\mathbb E^{\lambda,\perp}_{s,r;t}$ for expectation under the bridge
of the transverse diffusion in~\eqref{eq:longitudinal-SDE}
from $\mathsf S_0=s$ to $\mathsf S_t=r$, and set
$\mathcal A_t=\int_0^t\mathcal V(\mathsf S_v)\dd v$. A bridge is understood by
disintegration against the endpoint density
$p_\lambda^\perp(t,s,r)\dd r$; see~\cite{ChaumontUribeBravo}. Conditional on the complete transverse path,
\eqref{eq:longitudinal-SDE} gives a Gaussian longitudinal displacement
with mean $-2D\lambda\mathcal A_t$ and variance $2D\mathcal A_t$.
For $t>0$, $s\ge0$ and almost every $(r,y)\in[0,\infty)\times\R$, define
\[\label{eq:joint-clock-density-main}
 \begin{gathered}
 G^\lambda_{s,r;t}(y)
 =\mathbb E^{\lambda,\perp}_{s,r;t}
 \left[\frac1{\sqrt{4\pi D\mathcal A_t}}
   \exp\left(-\frac{(y+2D\lambda\mathcal A_t)^2}{4D\mathcal A_t}\right)\right],\\
 q_\lambda^{\mathrm{joint}}(t;s,r,y)
 =\frac{p_\lambda^\perp(t,s,r)}{\pi_\lambda^{\perp}(r)}G^\lambda_{s,r;t}(y).
 \end{gathered}
\]
Then $q_\lambda^{\mathrm{joint}}$ is a version of the joint transition
density of $(\mathsf S_t,\ell_t-\ell_0)$ conditional on $\mathsf S_0=s$,
relative to $\pi_\lambda^{\perp}(r)\dd r\dd y$.
For each continuous transverse path, $0<\mathcal A_t\le t$, and the
Gaussian factor has integral one in $y$. By conditioning and Tonelli,
the function $q_\lambda^{\mathrm{joint}}$ in~\eqref{eq:joint-clock-density-main} is a nonnegative probability density
relative to $\pi_\lambda^{\perp}(r)\dd r\dd y$, defined for almost every
$(r,y)\in[0,\infty)\times\R$.
Reflection at $0$ corresponds to the coordinate $s=|\rho|$.
For each fixed $t>0$ and $s\ge0$, and almost every $(r,y)$ relative to
$\pi_\lambda^{\perp}(r)\dd r\dd y$ with $0<G^\lambda_{s,r;t}(y)<\infty$,
conditioning also on $\ell_t-\ell_0=y$ gives a probability measure
on transverse paths whose density with respect to the transverse
bridge from $s$ to $r$ over time $t$ is
\[\label{eq:conditioned-transverse-bridge}
 \frac{1}{G^\lambda_{s,r;t}(y)\sqrt{4\pi D\mathcal A_t}}
 \exp\left(-\frac{(y+2D\lambda\mathcal A_t)^2}{4D\mathcal A_t}\right).
\]
The definition of $G^\lambda_{s,r;t}(y)$ in~\eqref{eq:joint-clock-density-main} makes the bridge integral of~\eqref{eq:conditioned-transverse-bridge} equal to one
\cite{ChaumontUribeBravo}.

For exponential deviation estimates under a spectral gap, see~\cite[Theorem~1.1]{Lezaud01}.

\section{Logarithmic delay along a geodesic}\label{sec:axial-delay}
We prove Theorem~\ref{thm:longitudinal}(iv). Throughout this section,
$d=2$, $f=f_{\G}$ and $\mu\le0<\sigma_2$. We use $y=\ell(X)$ for axial arclength,
$\rho$ for signed transverse distance, and $s=|\rho|$. Set $q(s)=D\sech^2s$.
In Fermi coordinates we write $u(\rho,y,t)$ for the solution at the corresponding point of $\Hh^2$.
For solutions even in $\rho$, equation~\eqref{eq:plane-Fermi} gives
\[\label{eq:reflected-pde}
 \left\{
 \begin{aligned}
 u_t&=D(u_{ss}+\tanh s\,u_s)+q(s)u_{yy}+F(s,u)
       &&\text{for }s>0,\ y\in\R,\ t>0,\\
 u_s(0,y,t)&=0&&\text{for }y\in\R,\ t>0.
 \end{aligned}
 \right.
\]

\subsection{Transverse eigenfunction and diffusion estimates}
By Lemma~\ref{lem:longitudinal-branch}, $\J_\parallel(k)$ has a positive
$L^2$-normalised eigenfunction $\psi_k$ at $\Lambda(k)$ for each $k\in\R$.
\begin{lemma}
For each fixed $k$, there are $C,\eta,R>0$ such that
\[\label{eq:eigenfunction-tails}
 \left|\frac{\psi_k'(s)}{\psi_k(s)}\right|\le C
 \text{ for }s\ge0,
 \qquad
 \frac{\psi_k(r)}{\psi_k(s)}\le C e^{-\eta(r-s)}
 \text{ for }r\ge s\ge R.
\]
\end{lemma}
\begin{proof}
Write the eigenfunction equation as
$\psi_k''=Q_k(s)\psi_k$, where $Q_k$ is bounded. Its logarithmic derivative
$r=\psi_k'/\psi_k$ satisfies $r'=Q_k-r^2$ and $r(0)=0$. Fix $r_0>\sqrt{2\|Q_k\|_\infty}$.
Whenever $|r|\ge r_0$, we have $r'\le-r^2/2$. Thus $r$ cannot cross
$r_0$ from below. If $r$ reached $-r_0$, then $-r$ would thereafter
satisfy $(-r)'\ge(-r)^2/2$ and blow up at a finite value of $s$.
This contradicts positivity and regularity of $\psi_k$ on the half-line
and proves the first bound in~\eqref{eq:eigenfunction-tails}.

Set $w=\psi_k/\phiZ$. By~\eqref{eq:bulk-ground},
\[\label{eq:tail-comparison-ode}
 Dw''+2D\frac{\phiZ'}{\phiZ}w'-d_k(s)w=0,
 \qquad
 d_k(s)=\Lambda(k)-\mu-D(k^2-1/4)\sech^2s
 \text{ for }s>0.
\]
Choose $R,d_0>0$ with $d_k\ge d_0$ on $[R,\infty)$, and then choose
$\eta>0$ with $D\eta^2+2D\eta\norm{\phiZ'/\phiZ}_\infty<d_0$.
For each $s\ge R$, the function $r\mapsto e^{-\eta(r-s)}$ is a
supersolution of~\eqref{eq:tail-comparison-ode} on $[s,\infty)$.
Multiply the differential inequality for
$[w(r)-w(s)e^{-\eta(r-s)}]_+$ by that positive part with weight $\phiZ^2$.
Its trace at $s$ is zero, and removal of cutoffs at infinity is justified
by $w\in H^1(0,\infty)$. Integration by parts gives a sum of nonpositive
terms, including $-d_0$ times the weighted squared norm. The positive part
therefore vanishes. Boundedness above and below of $\phiZ$ proves the
second estimate in~\eqref{eq:eigenfunction-tails}.
\end{proof}

Write $c=c_\parallel$, $\lambda=\lambda_\parallel$ and $\psi=\psi_\lambda$
in this section. Define
\[\label{eq:h-hzero}
 h_\lambda(s)=\frac{\psi(s)}{\sqrt{\cosh s}}
 \text{ for }s\ge0.
\]

Define the probability measure $\varpi(\mathrm ds)=\psi(s)^2\mathrm ds$ on
$[0,\infty)$, and write $\varpi(f)=\int f\dd\varpi$. Let $A$ be the reflected
diffusion generator
\[\label{eq:A-generator}
 Af=Df''+2D\frac{\psi'}\psi f'
 \text{ for }s>0,\qquad f'(0)=0,
\]
with semigroup $P_t$. Put $\bar q=\varpi(q)$ and $\mathfrak b(s)=c-2\lambda q(s)$.
By~\eqref{eq:longitudinal-HF}, $c=2\lambda\bar q$ and $\varpi(\mathfrak b)=0$.
For $\nu\in\{-1,1\}$, let $(S_t,Z_t)$ have generator
\[\label{eq:joint-generator}
 \mathcal L^\nu=A+q(s)\partial_{zz}+\nu\mathfrak b(s)\partial_z
 \text{ for }s>0,\ z\in\R,
\]
with reflection at $s=0$. Its law and expectation from $(s,z)$ are denoted
by $\Prb^\nu_{s,z}$ and $\E^\nu_{s,z}$. The transverse law is independent
of $z$ and $\nu$, so we write $\Prb_s$ and $\E_s$ for its probability and
expectation. Subscripts denoting an initial probability measure or density
have the corresponding meaning.

\begin{lemma}\label{lem:diffusion}
The semigroup $P_t$ is reversible with respect to $\varpi$ and has a spectral
gap $g>0$. There is a function $\chi\in C^2([0,\infty))$ satisfying
\[\label{eq:corrector}
 \begin{gathered}
 A\chi=-\mathfrak b\text{ for }s>0,
 \qquad \chi'(0)=0,\quad\varpi(\chi)=0,\\
 |\chi'(s)|\le C,\quad |\chi(s)|\le C(1+s)
 \text{ for }s\ge0.
 \end{gathered}
\]
For $\chi_\nu=\nu\chi$, the process
\[\label{eq:martingale}
 \begin{gathered}
 M_t=Z_t-Z_0+\chi_\nu(S_t)-\chi_\nu(S_0),\\
 \langle M\rangle_t=2\int_0^t(q+D|\chi'|^2)(S_r)\dd r
 \text{ for }t\ge0
 \end{gathered}
\]
is a continuous martingale with bounded quadratic variation rate.

For a fixed compact $I\subset[0,\infty)$, the transverse transition density
$p_t(s,r)$ relative to $\varpi(\mathrm dr)$ satisfies
\[\label{eq:compact-density}
 p_t(s,r)\le C_I
 \text{ for }s\in I,\ r\ge0,\ t\ge1.
\]
The estimates
\[\label{eq:corrector-moments}
 \begin{gathered}
 \Prb\left(\sup_{0\le r\le H}|\chi(S_r)|>v\right)
       \le C(1+H)e^{-\eta v}
       \text{ for }H,v\ge0,\\
 \E\sup_{2^j\le r\le2^{j+1}}|\chi(S_r)|^2\le C(1+j)^2
       \text{ for }j\ge0,\qquad
 \sup_{r\ge0}\E|\chi(S_r)|^2\le C
 \end{gathered}
\]
hold uniformly for $S_0\in I$. They also hold when the initial law is
$f\dd\varpi$, where $f\ge0$, $\varpi(f)=1$, and
$\norm{f}_{L^p(\varpi)}\le M$ for fixed $p>1$ and $M<\infty$.

If $R$ is sufficiently large and
$\zeta_R=\inf\{t\ge0:S_t\le R\}$, then for every sufficiently small
$\varepsilon>0$ there is $\kappa_\varepsilon>0$ such that
\[\label{eq:return-time}
 \E_s\zeta_R\le C(1+s),\qquad
 \E_s e^{\kappa_\varepsilon\zeta_R}\le e^{\varepsilon s}
 \text{ for }s\ge0.
\]
\end{lemma}

\begin{proof}
Multiplication by $\psi$ maps $L^2(\varpi)$ unitarily onto $L^2(\mathrm ds)$
and takes $A$ to $\J_\parallel(\lambda)-\Lambda(\lambda)$. The simple
isolated top eigenvalue from Lemma~\ref{lem:longitudinal-branch} therefore gives
\[\label{eq:gap}
 \norm{P_t f-\varpi(f)}_{L^2(\varpi)}
      \le e^{-gt}\norm{f-\varpi(f)}_{L^2(\varpi)}
 \text{ for }t\ge0,\ f\in L^2(\varpi).
\]
The logarithmic derivative in~\eqref{eq:eigenfunction-tails} is bounded,
so the reflected diffusion is nonexplosive. Its semigroup is Markov;
self-adjointness and $A1=0$ give reversibility and invariance of $\varpi$.

By~\eqref{eq:eigenfunction-tails},
\[\label{eq:tail-ratio}
 T_0(s):=\frac{\int_s^\infty\psi(r)^2\dd r}{\psi(s)^2}\le C
 \text{ for }s\ge0.
\]
Define $\chi$ first by its derivative,
\[\label{eq:corrector-derivative}
 \chi'(s)=\frac1{D\psi(s)^2}
                 \int_s^\infty\mathfrak b(r)\psi(r)^2\dd r
 \text{ for }s\ge0,
\]
and choose its additive constant so that $\varpi(\chi)=0$.
The numerator vanishes at zero because $\varpi(\mathfrak b)=0$.
Thus~\eqref{eq:tail-ratio} proves~\eqref{eq:corrector}, while the exponential
tail of $\psi$ makes its centreing integral finite. In particular
$\varpi(e^{\eta|\chi|})<\infty$ for some $\eta>0$.

With independent standard Brownian motions $B^1,B^2$, the stochastic
equations are
\[\label{eq:sde}
 \left\{
 \begin{aligned}
 \mathrm dS_t&=\sqrt{2D}\,\mathrm dB^1_t
       +2D(\psi'/\psi)(S_t)\mathrm dt+\mathrm d\ell_t^0,\\
 \mathrm dZ_t&=\nu\mathfrak b(S_t)\mathrm dt
                  +\sqrt{2q(S_t)}\,\mathrm dB^2_t
 \end{aligned}
 \right.
 \quad\text{for }t\ge0.
\]
Here $\ell^0$ is the nondecreasing reflection term, increasing only when
$S_t=0$. Applying It\^o's formula~\cite{RevuzYor} to $\chi_\nu(S_t)$, and using
$\chi'(0)=0$, proves~\eqref{eq:martingale}. Its bracket rate is bounded
by~\eqref{eq:corrector}; it is at least $2q(S_t)$.

Let $K_t^N$ be the Neumann heat kernel, relative to Lebesgue measure, of
$\J_\parallel(\lambda)$. Its bounded potential and the reflected Gaussian
kernel give
$K_1^N(s,r)\le C(e^{-(s-r)^2/(4D)}+e^{-(s+r)^2/(4D)})$.
The ground-state identity is
\[\label{eq:conjugated-kernel}
 p_1(s,r)=e^{-\Lambda(\lambda)}
                  \frac{K_1^N(s,r)}{\psi(s)\psi(r)}
 \text{ for }s,r\ge0.
\]
Since $\psi(r)\ge\psi(s)e^{-C|r-s|}$ by
\eqref{eq:eigenfunction-tails}, this proves~\eqref{eq:compact-density} at
time one. The semigroup property and invariance of $\varpi$ prove it for
every $t\ge1$.

The Skorokhod reflection formula and the bounded drift in~\eqref{eq:sde}
give
\[\label{eq:unit-increments}
 \sup_{0\le r\le1}|S_r-S_0|
       \le C+C\sup_{0\le r\le1}|B^1_r|
 \text{ almost surely for }S_0\ge0.
\]
Together with the Lipschitz bound on $\chi$ and
$\varpi(e^{\eta|\chi|})<\infty$, this gives an exponential tail for its
maximum on one time interval, under $\varpi$ and from a compact set.
For later intervals, use~\eqref{eq:compact-density} and the Markov property.
For an initial density $f\in L^p(\varpi)$, path-space H\"older's inequality
and the contraction of $P_t$ on $L^p(\varpi)$ give the same conclusion.
A union bound over unit intervals proves the first estimate in
\eqref{eq:corrector-moments}. On $[2^j,2^{j+1}]$ the union bound over its $2^j$ unit intervals gives
a tail bounded by $\min\{1,C2^j e^{-\eta v}\}$ at height $v$.
Integrating this bound against $2v\dd v$ gives $C(1+j)^2$.
For a fixed time $r$, integration of the tail on the single unit interval
containing $r$ gives $\E|\chi(S_r)|^2\le C$, uniformly in $r$.
These are the two moment bounds in~\eqref{eq:corrector-moments}.

For $s\ge R$, define
$T_R(s)=D^{-1}\int_R^sT_0(r)\dd r$. By~\eqref{eq:tail-ratio},
$0\le T_R(s)\le C(s-R)$ and $T_R'$ is bounded.
Differentiation gives $AT_R=-1$ on $(R,\infty)$.
Stopping It\^o's formula at $t\wedge\zeta_R$, first on compact intervals,
gives $\E_s(t\wedge\zeta_R)\le T_R(s)$ upon removal of the upper endpoint
by Fatou's lemma. Monotone convergence gives $\E_s\zeta_R\le T_R(s)$.
For small $\eta>0$, $V=e^{\eta T_R}$ satisfies
\[\label{eq:return-lyapunov}
 AV=V(-\eta+D\eta^2|T_R'|^2)\le-\eta V/2
 \text{ for }s>R.
\]
The positive local supermartingale obtained from
\eqref{eq:return-lyapunov}, stopped at $\zeta_R$, gives
$\E_s e^{\eta\zeta_R/2}\le e^{\eta T_R(s)}$ by Fatou's lemma.
Choosing $\eta$ with $\eta T_R(s)\le\varepsilon s$ proves
\eqref{eq:return-time}. The convention $\zeta_R=0$ for $s\le R$ completes
the proof.
\end{proof}

\subsection{Transition estimates with an absorbing boundary}
\label{sec:absorbing}

We estimate the diffusion generated by $\mathcal L^\nu$, for
$\nu\in\{-1,1\}$, before it reaches a prescribed boundary in the axial
coordinate. We retain $P_t=e^{tA}$ and $\bar q=\varpi(q)>0$. Write
$\norm{\cdot}_p=\norm{\cdot}_{L^p(\varpi)}$ and
$\ip{f}{g}_\varpi=\int \overline f g\dd\varpi$. We use $H^1(\varpi)$ for the
space of locally absolutely continuous functions $f$ with
$\varpi(|f|^2+|f'|^2)<\infty$. For a probability density $f$
relative to $\varpi$, the notation $\Prb^\nu_{f,z}$ means that $S_0$ has law
$f\varpi$ and $Z_0=z$. Point starts are denoted by $\Prb^\nu_{s,z}$.

Fix $K>0$. On each time interval under consideration, let $B$ be a
$C^1$ function satisfying
\[\label{eq:wall-modulus}
 |B(v)-B(w)|\le K\log(1+|v-w|)
 \text{ for all }v,w\text{ in that interval}.
\]
For $0\le r<t$, define the first hitting time, measured from $r$, and
the corresponding transition operator for $s\ge0$, $z>B(r)$ and $G\ge0$ by
\[\label{eq:killed-evolution}
 \begin{aligned}
 \tau_r^B
   &=\inf\{v\in[0,t-r]:Z_v\le B(r+v)\},\\
 K^{B,\nu}_{r,t}G(s,z)
   &=\E^\nu_{s,z}\bigl[G(S_{t-r},Z_{t-r});\,\tau_r^B>t-r\bigr].
 \end{aligned}
\]
We use $\inf\varnothing=\infty$, set $K^{B,\nu}_{r,r}G=G$, and extend
$K^{B,\nu}_{r,t}G$ by zero to $z\le B(r)$ when $r<t$. When $r=0$,
we write $\tau_B=\tau_0^B$. Restrictions, reversals and time translations
of $B$ satisfy~\eqref{eq:wall-modulus} with the same constant.

The adjoint of $\mathcal L^\nu$ relative to $\varpi(\mathrm ds)\,\mathrm dz$
is $\mathcal L^{-\nu}$. For nonnegative measurable functions $\Phi,G$
on $[0,\infty)\times\R$, let
$\widehat B(v)=B(t-v)$ for $0\le v\le t-r$. Then, for $0\le r<t$,
\[\label{eq:killed-duality}
 \begin{aligned}
 &\int_0^\infty\int_\R
       \Phi(s,z)K^{B,\nu}_{r,t}G(s,z)\,\mathrm dz\,\varpi(\mathrm ds)\\
 &\qquad=
 \int_0^\infty\int_\R
       G(s,z)K^{\widehat B,-\nu}_{0,t-r}\Phi(s,z)
          \,\mathrm dz\,\varpi(\mathrm ds).
 \end{aligned}
\]
To verify this identity, first take bounded, compactly supported $\Phi,G$
and restrict the diffusion to a bounded rectangle. Integration by parts
uses reversibility of $A$, the independence of $q,\mathfrak b$ from $z$, the
Neumann condition at $s=0$, and the zero values at the absorbing boundary.
The product of a forward and a backward solution then has constant
integral; the boundary term from the moving endpoint vanishes. Exhausting
the half-plane gives~\eqref{eq:killed-duality}, and monotone convergence
extends it to nonnegative $\Phi,G$.

\begin{lemma}\label{lem:clock}
There is $\kappa>0$ with the following property. Fix a compact interval
$I\subset[0,\infty)$, $p>1$, and $M\ge1$. There are $\kappa_0,C>0$ such that
\[\label{eq:clock-lower}
 \Prb\left(\int_v^{v+H}q(S_w)\dd w<\kappa H\right)
 \le Ce^{-\kappa_0H}
 \text{ for }v\ge0,\ H\ge1,
\]
uniformly for a point start $S_0=s\in I$, or for an initial density
$f\ge0$ with $\varpi(f)=1$ and $\norm f_{L^{p}(\varpi)}\le M$.
\end{lemma}

\begin{proof}
Let $g>0$ be the spectral gap in~\eqref{eq:gap}. For
$f_*\in H^1(\varpi)$ with $\norm{f_*}_{L^{2}(\varpi)}=1$, write $f_*=m+g_*$, where
$m=\varpi(f_*)$ and $\varpi(g_*)=0$. Since $0<q\le D$, Young's inequality gives, for $\delta>0$ and $f_*\in H^1(\varpi)$ with $\norm{f_*}_{L^{2}(\varpi)}=1$,
\[\label{eq:clock-form}
 \begin{aligned}
 \ip{f_*}{(A-\delta q)f_*}_\varpi
 &\le-g\norm{g_*}_{L^{2}(\varpi)}^2-\delta\bar q|m|^2
                  +2\delta D|m|\norm{g_*}_{L^{2}(\varpi)}\\
 &\le-\frac g2\norm{g_*}_{L^{2}(\varpi)}^2
       -\left(\delta\bar q-\frac{2\delta^2D^2}{g}\right)|m|^2.
 \end{aligned}
\]
For functions in the form domain, the left side denotes the associated
quadratic form. As $|m|^2+\norm{g_*}_{L^{2}(\varpi)}^2=1$, a sufficiently small
$\delta_0>0$ satisfies
\[\label{eq:clock-spectrum}
 \sup\sigma(A-\delta q)\le-\frac34\delta\bar q
 \text{ for }0<\delta\le\delta_0.
\]
For a stationary start, Feynman--Kac and the spectral theorem therefore
give, for $H\ge0$ and $0<\delta\le\delta_0$,
\[\label{eq:clock-stationary}
 \begin{aligned}
 \E_\varpi\exp\left(-\delta\int_0^Hq(S_w)\dd w\right)
 &=\ip{1}{e^{H(A-\delta q)}1}_\varpi
 \le e^{-3\delta\bar q H/4},\\
 \Prb_\varpi\left(\int_0^Hq(S_w)\dd w<\bar qH/2\right)
 &\le e^{-\delta\bar q H/4}.
 \end{aligned}
\]

Put $p'=p/(p-1)$. If $S_0$ has density $f$ as in the statement, then
H\"older's inequality on path space yields
\[\label{eq:clock-density}
 \E_f\exp\left(-\delta\int_0^Hq(S_w)\dd w\right)
 \le M\left[
 \E_\varpi\exp\left(-p'\delta\int_0^Hq(S_w)\dd w\right)
 \right]^{1/p'}
 \text{ for }H\ge0,\ \delta>0.
\]
Choose $p'\delta\le\delta_0$ and apply~\eqref{eq:clock-stationary}.
The density of $S_v$ is $P_vf$, whose $L^p(\varpi)$ norm is at most $M$.
The same estimate consequently holds on every interval $[v,v+H]$.

For $S_0=s\in I$, \eqref{eq:compact-density} bounds the transverse density
by $C_I\varpi$ at every time at least one. If $v<1$ and $H\ge2$, discard
the part of the clock before time one and apply
\eqref{eq:clock-stationary} on $[1,v+H]$. This interval has length at
least $H-1\ge H/2$. Taking $\kappa=\bar q/4$ proves
\eqref{eq:clock-lower} for $H\ge2$. Increasing $C$ covers $1\le H<2$.
\end{proof}

\begin{lemma}\label{lem:survival}
Fix a compact interval $I\subset[0,\infty)$, $p>1$, and $M\ge1$.
There is $C>0$ such that
\[\label{eq:survival}
 \Prb^\nu_{s,B(0)+x}(\tau_B>H)
 \le\frac{C(1+x)}{\sqrt{1+H}}
 \text{ for }s\in I,\ x\ge0,\ H\ge0,\ \nu\in\{-1,1\}.
\]
The same bound holds for $\Prb^\nu_{f,B(0)+x}$ whenever
$f\ge0$, $\varpi(f)=1$, and $\norm f_{L^{p}(\varpi)}\le M$. The constant is uniform
over all boundaries satisfying~\eqref{eq:wall-modulus}.
\end{lemma}

\begin{proof}
Translate $B(0)$ to zero, fix $\nu$, and put $\xi=\nu\chi$.
We write $\Prb,\E$ for either initial law in the statement. All constants
below are uniform over these laws.

The moment estimates in~\eqref{eq:corrector-moments} and the density bound
\eqref{eq:compact-density} give an exponential tail for the maximum of
$|\xi(S_v)|$ on each unit interval. For an initial density $f$, this follows
from H\"older's inequality and $\norm{P_jf}_{L^{p}(\varpi)}\le M$ at the beginning
of the interval. For a compact point start, the first interval is covered
by~\eqref{eq:corrector-moments}, and the later intervals use
\eqref{eq:compact-density}. A union bound and integration of the resulting
tail yield
\[\label{eq:corrector-maxima}
 \begin{gathered}
 \Prb\left(\sup_{0\le v\le H}|\xi(S_v)|>u\right)
       \le C(1+H)e^{-\eta_0u}
       \text{ for }u,H\ge0,\\
 \E\left[\sup_{2^j\le v\le2^{j+1}}|\xi(S_v)|^2\right]
       \le C(1+j)^2
       \text{ for integers }j\ge0,\\
 \sup_{v\ge0}\E|\xi(S_v)|^2\le C.
 \end{gathered}
\]
Here and below, a bound larger than one on a probability may be replaced
by one.

Choose $C_0$ so that the maximum in the first line of
\eqref{eq:corrector-maxima} is at most $C_0\log(2+H)$ except on an event
of probability $C(1+H)^{-2}$. On survival to time $H$, the martingale
in~\eqref{eq:martingale} satisfies
\[\label{eq:martingale-survival-barrier}
 M_v=Z_v-x+\xi(S_v)-\xi(S_0)
 \ge-x-K\log(1+H)-2C_0\log(2+H)
 \text{ for }0\le v\le H
\]
on this event. Its bracket is at least
$2\int_0^Hq(S_v)\dd v$. Lemma~\ref{lem:clock}, the martingale time-change
theorem, and the Brownian reflection principle therefore give
\[\label{eq:survival-rough}
 \Prb(\tau_B>H)
 \le C\frac{1+x+\log(2+H)}{\sqrt{1+H}}
 \text{ for }x\ge0,\ H\ge0.
\]

We next estimate the first moment of the surviving height. Set
$X_H=Z_H-B(H)$. Lemma~\ref{lem:diffusion} bounds the bracket rate of
$M$, so optional stopping at the bounded time $H\wedge\tau_B$ is
applicable. It gives, for $H\ge0$,
\[\label{eq:stopped-survival-identity}
 \begin{aligned}
 \E[X_H;\tau_B>H]
 ={}&x+\E\xi(S_0)-\E[\xi(S_H);\tau_B>H]\\
 &-\E[\xi(S_{\tau_B});\tau_B\le H]
  -\E[B(\tau_B);\tau_B\le H]
  -B(H)\Prb(\tau_B>H).
 \end{aligned}
\]
The first-unit stopped corrector is integrable by
\eqref{eq:corrector-moments}. On each subsequent dyadic interval,
Cauchy--Schwarz, \eqref{eq:corrector-maxima}, and
\eqref{eq:survival-rough} imply, for each integer $j\ge0$ with $2^j<H$,
\[\label{eq:stopped-corrector-bound}
 \begin{aligned}
 &\E\bigl[|\xi(S_{\tau_B})|;
             2^j<\tau_B\le\min(2^{j+1},H)\bigr]\\
 &\qquad\le
 \left(\E\sup_{2^j\le v\le2^{j+1}}|\xi(S_v)|^2\right)^{1/2}
           \Prb(\tau_B>2^j)^{1/2}\\
 &\qquad\le C(1+j)(1+x+j)^{1/2}2^{-j/4}.
 \end{aligned}
\]
The numerical series is bounded by $C\sqrt{1+x}$. The same dyadic
decomposition, using $|B(v)|\le K\log(1+v)$, gives
\[\label{eq:stopped-wall-bound}
 \E[|B(\tau_B)|;\tau_B\le H]
 \le C+C\sum_{j=0}^\infty(1+j)(1+x+j)2^{-j/2}
 \le C(1+x)
 \text{ for }H\ge0,\ x\ge0.
\]
The initial and terminal correctors in
\eqref{eq:stopped-survival-identity} have bounded first moments by
\eqref{eq:corrector-maxima}. Finally, for $H\ge0$ and $x\ge0$,
\[\label{eq:terminal-wall-bound}
 \begin{aligned}
 |B(H)|\Prb(\tau_B>H)
 &\le CK\log(1+H)
       \frac{1+x+\log(2+H)}{\sqrt{1+H}}\\
 &\le C(1+x).
 \end{aligned}
\]
Substitution in~\eqref{eq:stopped-survival-identity} proves
\[\label{eq:survival-first-moment}
 \E[X_H;\tau_B>H]\le C(1+x)
 \text{ for }H\ge0,\ x\ge0.
\]

Let $s(H)=\Prb(\tau_B>H)$ and fix $\delta>0$. By
\eqref{eq:survival-first-moment}, the mass with $X_H>\delta\sqrt H$ is
at most $C(1+x)/(\delta\sqrt H)$. For the remaining mass, retain
survival only to time $H/2$. Conditional on the transverse path and the
axial history up to $H/2$, the last axial increment is Gaussian with
variance $2\int_{H/2}^Hq(S_v)\dd v$. On the event that this clock is at
least $\kappa H/2$, its density is bounded by $C/\sqrt H$. The possible
endpoint interval has length $\delta\sqrt H$. Lemma~\ref{lem:clock}
then gives
\[\label{eq:survival-recursion}
 s(H)\le\frac{C(1+x)}{\delta\sqrt H}
           +C\delta\,s(H/2)+Ce^{-\eta_1H}
 \text{ for }H\ge4,\ x\ge0.
\]
The exceptional clock probability is estimated under the unconditional
transverse marginal. Choose $\delta$ so that $C\delta<1/\sqrt2$.
Iteration of~\eqref{eq:survival-recursion} over successive halves gives
a geometric sum with ratio $C\delta\sqrt2<1$. Its last term is bounded
by $CH^{-1/2}$, and $e^{-\eta_1H}\le C H^{-1/2}$ for $H\ge4$.
This proves~\eqref{eq:survival}; bounded $H$ is covered by increasing $C$.
\end{proof}

\begin{lemma}\label{lem:weighted-transition}
Fix a compact interval $I\subset[0,\infty)$ and $p>1$. Let $V\ge1$
belong to $L^p(\varpi)$, and let $\varrho:[0,\infty)\to[0,\infty)$ be
bounded and measurable, with
\[\label{eq:transition-profile}
 \int_0^\infty(1+w)\varrho(w)\dd w<\infty.
\]
Suppose that the nonnegative terminal test $G$ satisfies
\[\label{eq:transition-test}
 G(s,B(H)+w)\le V(s)\varrho(w)
 \text{ for a.e. }s\ge0,\ w\ge0.
\]
There is $C>0$, depending on $I,K,p,V,\varrho$ and the fixed
coefficients, such that
\[\label{eq:weighted-transition}
 K^{B,\nu}_{0,H}G(s,B(0)+x)
 \le C(1+x)(1+H)^{-3/2}
 \text{ for }s\in I,\ x\ge0,\ H\ge3,\ \nu\in\{-1,1\}.
\]
The constant is uniform over $B$ and $G$ satisfying
\eqref{eq:wall-modulus} and~\eqref{eq:transition-test}.
\end{lemma}

\begin{proof}
Put $h_1=H/3\ge1$, and fix a starting point
$(s_0,B(0)+x)$ with $s_0\in I$. Let $f\varpi$ be the transverse marginal
of the process surviving to time $h_1$. By
\eqref{eq:compact-density} and Lemma~\ref{lem:survival},
\[\label{eq:first-third-marginal}
 0\le f(s)\le C_I\text{ for a.e. }s\ge0,
 \qquad
 \varpi(f)\le C(1+x)h_1^{-1/2}
 \text{ for }H\ge3,\ x\ge0.
\]
Define $g(s,z)=K^{B,\nu}_{2h_1,3h_1}G(s,z)$ for $z>B(2h_1)$,
extend it by zero below that boundary, and put
$j(s)=\int_\R g(s,z)\,\mathrm dz$. Dropping killing on the last third
and integrating over all axial starting points gives
\[\label{eq:last-third-integral}
 j(s)\le R_0P_{h_1}V(s)\text{ for a.e. }s\ge0,
 \qquad
 \norm j_{L^{p}(\varpi)}\le R_0\norm V_{L^{p}(\varpi)},
 \qquad R_0=\int_0^\infty\varrho(w)\dd w.
\]
Indeed, conditional on the transverse path, an axial translation of the
starting point translates the Gaussian endpoint density by the same
amount, and the integral of that density is one.

Apply~\eqref{eq:killed-duality} on the final third. If
$\widehat B(v)=B(3h_1-v)$ for $0\le v\le h_1$, then for $H\ge3$,
\[\label{eq:last-third-duality}
 \begin{aligned}
 \varpi(j)
 &=\int_0^\infty\int_0^\infty
       G(r,B(3h_1)+w)
       \Prb^{-\nu}_{r,B(3h_1)+w}
                     (\tau_{\widehat B}>h_1)
                  \,\varpi(\mathrm dr)\,\mathrm dw\\
 &\le C h_1^{-1/2}
        \int_0^\infty(1+w)\varrho(w)\dd w
 \le C h_1^{-1/2}.
 \end{aligned}
\]
For the inequality, use Lemma~\ref{lem:survival} with initial transverse
density $V/\varpi(V)$, which belongs to $L^p(\varpi)$, and then multiply by
$\varpi(V)$. Tonelli's theorem justifies the integrals.

On the middle third, let
$E_{\rm bad}=\{\int_{h_1}^{2h_1}q(S_v)\dd v<\kappa h_1\}$,
where $\kappa$ is from Lemma~\ref{lem:clock}. On its complement the
conditional axial density is at most $C h_1^{-1/2}$. Dropping killing
on this third and integrating $g$ over the axial endpoint therefore
bounds its contribution by
$C h_1^{-1/2}\ip{f}{P_{h_1}j}_\varpi$.
Let $p'=p/(p-1)$. Self-adjointness and H\"older's inequality give, for $H\ge3$,
\[\label{eq:middle-third-pairing}
 \begin{aligned}
 \ip{f}{P_{h_1}j}_\varpi
 &=\varpi(f)\varpi(j)
       +\ip{P_{h_1}(f-\varpi(f))}{j}_\varpi\\
 &\le C(1+x)h_1^{-1}
       +\norm{P_{h_1}(f-\varpi(f))}_{L^{p'}(\varpi)}\norm j_{L^{p}(\varpi)}.
 \end{aligned}
\]
The spectral gap bounds the last semigroup factor in $L^2(\varpi)$ by
$Ce^{-gh_1}$. Its $L^\infty$ norm is bounded by
\eqref{eq:first-third-marginal}. Interpolation when $p'>2$, and the
inclusion $L^2(\varpi)\subset L^{p'}(\varpi)$ when $p'\le2$, imply
\[\label{eq:middle-third-mixing}
 \norm{P_{h_1}(f-\varpi(f))}_{L^{p'}(\varpi)}\le Ce^{-\eta_2h_1},
 \qquad
 \ip{f}{P_{h_1}j}_\varpi
 \le C(1+x)h_1^{-1}+Ce^{-\eta_2h_1}
 \text{ for }H\ge3.
\]
The pairing identity first holds for bounded truncations of $j$;
\eqref{eq:last-third-integral} and H\"older's inequality give the
displayed identity on passing to the limit.

For the contribution of $E_{\rm bad}$, drop all killing and use the
boundedness of $\varrho$. If $1<r_0\le p$, then for $s_0\in I$, $x\ge0$ and $H\ge3$,
\[\label{eq:weighted-clock-error}
 \begin{aligned}
 &\E^\nu_{s_0,B(0)+x}
       [G(S_H,Z_H);\,\tau_B>H,\ E_{\rm bad}]\\
 &\quad\le\norm{\varrho}_\infty
        \E_{s_0}[V(S_H)\one_{E_{\rm bad}}]\\
 &\quad\le\norm{\varrho}_\infty
        \bigl(\E_{s_0}V(S_H)^{r_0}\bigr)^{1/r_0}
        \Prb_{s_0}(E_{\rm bad})^{1-1/r_0}
 \le Ce^{-\eta_3h_1}.
 \end{aligned}
\]
Here~\eqref{eq:compact-density} bounds the unconditional transverse
density at time $H$ by $C_I\varpi$, and Lemma~\ref{lem:clock} estimates
the middle clock under the unconditional law at time $h_1$.
Combining~\eqref{eq:middle-third-mixing} and
\eqref{eq:weighted-clock-error} proves~\eqref{eq:weighted-transition}.
\end{proof}

We now extend the estimate to arbitrary transverse starting points for
terminal tests supported in a fixed compact set. Choose a small
$\varepsilon_*>0$ and fix
\[\label{eq:transverse-weight}
 \mathcal W(s)=1+s+e^{\varepsilon_* s}\text{ for }s\ge0,
 \qquad
 \mathcal W\in L^{3/2}(\varpi),
 \qquad
 \sup_{s\ge0}h_\lambda(s)\mathcal W(s)<\infty.
\]
Such a choice follows from~\eqref{eq:eigenfunction-tails}: if
$\psi(s)\le Ce^{-\eta s}$ for large $s$, any sufficiently small
$0<\varepsilon_*<\eta$ has the stated properties. We fix this
$\varepsilon_*$ for the remainder of the proof. The return-time estimate
\eqref{eq:return-time} will be used with $\varepsilon=\varepsilon_*$.

\begin{lemma}\label{lem:compact-terminal}
Fix $R_1,L_1>0$. Let $G$ be a bounded, nonnegative terminal test
whose restriction to $\{(s,z):s\ge0,\ z\ge B(H)\}$ is supported in
$[0,R_1]\times[B(H),B(H)+L_1]$. There is $C>0$, depending on
$R_1,L_1,K,\varepsilon_*$ and the fixed coefficients, such that
\[\label{eq:compact-terminal}
 K^{B,\nu}_{0,H}G(s,B(0)+x)
 \le C\norm G_\infty(1+x)\mathcal W(s)(1+H)^{-3/2}
 \text{ for }s,x,H\ge0,\ \nu\in\{-1,1\}.
\]
The constant is uniform over $B$ satisfying~\eqref{eq:wall-modulus}.
\end{lemma}

\begin{proof}
Bounded horizons are covered by
$K^{B,\nu}_{0,H}G\le\norm G_\infty$, so assume $H\ge6$.
Choose $R>R_1$ large enough for~\eqref{eq:return-time}. For
$0\le s\le R$, Lemma~\ref{lem:weighted-transition} applies with
$V=1$ and $\varrho=\one_{[0,L_1]}$, after factoring out
$\norm G_\infty$.

Fix $s>R$, $x\ge0$ and $H\ge6$, and let $\zeta_R=\inf\{v\ge0:S_v=R\}$. Every path with a
nonzero terminal value of $G$ hits $R$ before time $H$. On
$\{\zeta_R\le H/2,\ \tau_B>\zeta_R\}$, the strong Markov property
and Lemma~\ref{lem:weighted-transition}, applied to the remaining
boundary segment, give
\[\label{eq:compact-terminal-early}
 \begin{aligned}
 &\E^\nu_{s,B(0)+x}
       [G(S_H,Z_H);\,\tau_B>H,\ \zeta_R\le H/2]\\
 &\qquad\le C\norm G_\infty(1+H)^{-3/2}
       \E^\nu_{s,B(0)+x}
       [1+Z_{\zeta_R}-B(\zeta_R);
                         \,\tau_B>\zeta_R,\ \zeta_R\le H/2].
 \end{aligned}
\]
The height at $\zeta_R$ is estimated from~\eqref{eq:sde}. Put
$N_v=Z_v-Z_0-\nu\int_0^v\mathfrak b(S_w)\dd w$ for $v\ge0$.
The martingale isometry and~\eqref{eq:return-time} yield
\[\label{eq:height-at-return}
 \begin{aligned}
 &\E^\nu_{s,B(0)+x}
       [1+Z_{\zeta_R}-B(\zeta_R);
                         \,\tau_B>\zeta_R,\ \zeta_R\le H/2]\\
 &\qquad\le 1+x+\norm{\mathfrak b}_\infty\E_s\zeta_R
                +\E^\nu_{s,B(0)+x}|N_{\zeta_R}|
                +K\E_s\log(1+\zeta_R)\\
 &\qquad\le C(1+x+s),
 \qquad
 \E^\nu_{s,B(0)+x}|N_{\zeta_R}|^2
       \le2D\E_s\zeta_R.
 \end{aligned}
\]
On $\{\zeta_R>H/2\}$, the exponential estimate in~\eqref{eq:return-time} gives
\[\label{eq:compact-terminal-late}
 \begin{aligned}
 &\E^\nu_{s,B(0)+x}
       [G(S_H,Z_H);\,\tau_B>H,\ \zeta_R>H/2]\\
 &\qquad\le C\norm G_\infty e^{\varepsilon_* s-\kappa_{\varepsilon_*}H/2}
 \le C\norm G_\infty e^{\varepsilon_* s}(1+H)^{-3/2}.
 \end{aligned}
\]
Equations~\eqref{eq:compact-terminal-early}--\eqref{eq:compact-terminal-late}
and $1+x+s\le(1+x)(1+s)$ prove~\eqref{eq:compact-terminal}.
\end{proof}

We next obtain lower bounds for the killed kernels. We write
$k^{B,\nu}_{r,t}(s,x;\widetilde s,\widetilde x)$ for the density of
$K^{B,\nu}_{r,t}$ with respect to
$\varpi(\mathrm d\widetilde s)\,\mathrm d\widetilde x$. For the wall
$B\equiv0$, we use the version continuous for axial coordinates above
zero, including at the reflecting transverse boundary. Local uniform
ellipticity and parabolic regularity give this version; reflection across
$s=0$ treats the Neumann boundary. We abbreviate it to $k_H^\nu$, where $H=t-r$,
and put $\tau_0=\inf\{v\ge0:Z_v\le0\}$.

\begin{lemma}\label{lem:meander}
Let $I_0\subset[0,\infty)$ be a compact interval and let $C_0>0$.
There exist $Y>0$, a compact interval $I_1\subset[0,\infty)$,
constants $0<d_1<d_2$, $c_0>0$, and $H_0\ge1$ such that, for $H\ge H_0$, $s\in I_0$, $x\in[Y,Y+C_0]$ and $\nu\in\{-1,1\}$,
\[\label{eq:meander-mass}
 \mathbb P^\nu_{s,x}
 \{\tau_0>H,\ S_H\in I_1,\ d_1\sqrt H\le Z_H\le d_2\sqrt H\}
 \ge\frac{c_0}{\sqrt H}.
\]
\end{lemma}

\begin{proof}
Set $\chi_\nu=\nu\chi$. The process
\[\label{eq:meander-martingale}
 M_v^\nu=Z_v-x+\chi_\nu(S_v)-\chi_\nu(s),\qquad
 \langle M^\nu\rangle_v
 =2\int_0^v(q+D|\chi'|^2)(S_r)\,\mathrm dr
 \quad\text{for }v\ge0
\]
is a continuous martingale. Its quadratic variation has bounded
derivative by Lemma~\ref{lem:diffusion}. Lemma~\ref{lem:survival} gives
\[\label{eq:meander-survival}
 \mathbb P^\nu_{s,x}(\tau_0>H)
 \le C_{I_0}(1+x)(1+H)^{-1/2}
 \quad\text{for }H\ge0,\ s\in I_0,\ x\ge0,\
                   \nu\in\{-1,1\}.
\]
In particular, $\tau_0$ is finite almost surely.

The exponential bounds on the corrector over unit intervals in
Lemma~\ref{lem:diffusion}, followed by a union bound, imply
\[\label{eq:meander-dyadic-corrector}
 \mathbb E_s\sup_{2^j\le v\le2^{j+1}}|\chi(S_v)|^2
 \le C_{I_0}(1+j)^2
 \quad\text{for }j\ge0,\ s\in I_0.
\]
Indeed, the maximum over at most $2^{j+1}$ unit intervals has tail
bounded by the minimum of $1$ and $C2^je^{-\eta R}$ at height $R$, for some $\eta>0$;
integrating this bound gives \eqref{eq:meander-dyadic-corrector}. Cauchy--Schwarz
and~\eqref{eq:meander-survival} now bound the contribution of
$2^j<\tau_0\le2^{j+1}$ to
$\mathbb E^\nu_{s,x}|\chi_\nu(S_{\tau_0})|$ by
$C_{I_0}(1+j)\sqrt{1+x}\,2^{-j/4}$. The first unit interval is
controlled directly by Lemma~\ref{lem:diffusion}. Summation gives
\[\label{eq:meander-hit-corrector}
 \mathbb E^\nu_{s,x}|\chi_\nu(S_{\tau_0})|
 \le C_{I_0}\sqrt{1+x}
 \quad\text{for }s\in I_0,\ x\ge0,\ \nu\in\{-1,1\}.
\]

Optional stopping of~\eqref{eq:meander-martingale} at
$H\wedge\tau_0$ is justified by its bounded quadratic variation
rate. Since $Z_{\tau_0}=0$, it yields
\[\label{eq:meander-first-moment}
 \begin{aligned}
 \mathbb E^\nu_{s,x}[Z_H;\tau_0>H]
 ={}&x+\chi_\nu(s)
 -\mathbb E^\nu_{s,x}[\chi_\nu(S_{\tau_0});\tau_0\le H]\\
 &-\mathbb E^\nu_{s,x}[\chi_\nu(S_H);\tau_0>H]
 \quad\text{for }H\ge0,\ s\in I_0,\ x\ge0.
 \end{aligned}
\]
The transverse second moment of $\chi(S_H)$ is bounded uniformly
for $s\in I_0$ and $H\ge0$, by the compact-start density bound and
the unit-interval estimate in Lemma~\ref{lem:diffusion}. Hence the
last term in~\eqref{eq:meander-first-moment} tends to zero, uniformly
when $s\in I_0$ and $x$ ranges over a fixed bounded interval. The
summable bound preceding~\eqref{eq:meander-hit-corrector} gives the
same uniformity for the stopped term. Consequently the left side
converges to
$V_\nu(s,x)=x+\chi_\nu(s)-\mathbb E^\nu_{s,x}\chi_\nu(S_{\tau_0})$.
Equation~\eqref{eq:meander-hit-corrector} implies
$V_\nu(s,x)\ge x-C_{I_0}\sqrt{1+x}-C_{I_0}$. Choose $Y$ large
enough that this lower bound is positive for $x\ge Y$. There are
then $v_0>0$ and $H_1\ge1$ such that
\[\label{eq:meander-positive-moment}
 \mathbb E^\nu_{s,x}[Z_H;\tau_0>H]\ge v_0
 \quad\text{for }H\ge H_1,\ s\in I_0,\ x\in[Y,Y+C_0],\
                   \nu\in\{-1,1\}.
\]

On this fixed starting rectangle,
\[\label{eq:meander-stopped-second}
 \mathbb E^\nu_{s,x}|M_{H\wedge\tau_0}^\nu|^2
 \le C\int_0^H\mathbb P^\nu_{s,x}(\tau_0>r)\,\mathrm dr
 \le C\sqrt{1+H}
 \quad\text{for }H\ge0.
\]
Using~\eqref{eq:meander-martingale} on $\{\tau_0>H\}$ and the
uniform transverse second moment gives
\[\label{eq:meander-killed-second}
 \mathbb E^\nu_{s,x}[Z_H^2;\tau_0>H]\le C\sqrt{1+H}
 \quad\text{for }H\ge0,\ s\in I_0,\ x\in[Y,Y+C_0],\
                   \nu\in\{-1,1\}.
\]
By~\eqref{eq:meander-survival}, the first moment contributed by
$Z_H<d_1\sqrt H$ is at most $Cd_1$ for $H\ge1$. By
\eqref{eq:meander-killed-second}, the contribution of
$Z_H>d_2\sqrt H$ is at most $C/d_2$. Choose $d_1$ small and
$d_2$ large so that their sum is less than $v_0/4$.

It remains to restrict the transverse endpoint. At time $H/2$, the
transverse marginal of the surviving paths is $f\varpi$, where
$0\le f\le C$ and $\varpi(f)\le C H^{-1/2}$. The first assertion
uses the compact-start density estimate in
Lemma~\ref{lem:diffusion}, and the second uses
\eqref{eq:meander-survival}. Dropping killing after $H/2$ and
applying the transverse spectral gap gives
\[\label{eq:meander-transverse-tail}
 \begin{aligned}
 \mathbb P^\nu_{s,x}(S_H\notin I_1,\tau_0>H)
 &\le\langle f,P_{H/2}\mathbf1_{I_1^c}\rangle_{L^2(\varpi)}\\
 &\le CH^{-1/2}\varpi(I_1^c)+Ce^{-gH/2}
 \quad\text{for }H\ge2.
 \end{aligned}
\]
The constants are uniform on the starting rectangle. Combining
\eqref{eq:meander-transverse-tail} with
\eqref{eq:meander-killed-second} bounds the first moment on
$\{S_H\notin I_1,\tau_0>H\}$ by
$C\sqrt{\varpi(I_1^c)}+o(1)$. Choose the compact interval $I_1$ large
enough to make this less than $v_0/4$ for large $H$. After the three
discarded contributions are removed from
\eqref{eq:meander-positive-moment}, the remaining first moment is
at least $v_0/2$. Since $Z_H\le d_2\sqrt H$ there, division by
$d_2\sqrt H$ proves~\eqref{eq:meander-mass}.
\end{proof}

\begin{lemma}\label{lem:middle-density}
Let $I\subset[0,\infty)$ be a compact interval and let
$0<d_1<d_2$. There are $C>0$ and $H_0\ge1$ such that
\[\label{eq:middle-density-lower}
 k_H^\nu(s,x;r,y)\ge C^{-1}H^{-1/2}
 \quad\text{for }H\ge H_0,\ s,r\in I,\
 x,y\in[d_1\sqrt H,d_2\sqrt H],\ \nu\in\{-1,1\},
\]
where the inequality in the endpoint variables is understood
almost everywhere with respect to $\varpi(\mathrm dr)\,\mathrm dy$.
\end{lemma}

\begin{proof}
Put $\bar q=\varpi(q)=c/(2\lambda)>0$ and
$Q_v=\int_0^vq(S_r)\,\mathrm dr$. If $\mathsf B$ denotes the
Brownian motion driving $S$, It\^o's formula and $A\chi=2\lambda q-c$
give
\[\label{eq:clock-corrector}
 2\lambda(Q_v-\bar qv)
 =\chi(S_v)-\chi(S_0)
       -\sqrt{2D}\int_0^v\chi'(S_r)\,\mathrm d\mathsf B_r
 \quad\text{for }v\ge0.
\]
Under the initial law $\varpi$, the martingale on the right, divided
by $2\lambda$, has quadratic variation whose ratio to time
converges almost surely to
$\sigma_Q^2=D\varpi(|\chi'|^2)/(2\lambda^2)$. This follows from
ergodicity of the transverse diffusion. Monotonicity of quadratic
variation makes the rescaled convergence uniform on $[0,1]$.
The martingale time-change theorem therefore gives convergence to
a Brownian motion of variance $\sigma_Q^2$: after Brownian scaling,
the error caused by the rescaled clock tends to zero in probability
by the Brownian modulus of continuity. The unit-interval
exponential tails in Lemma~\ref{lem:diffusion} also give
\[\label{eq:clock-corrector-negligible}
 \mathbb P_\varpi\left\{\sup_{0\le v\le H}|\chi(S_v)|>\epsilon\sqrt H\right\}
 \le C(1+H)e^{-\eta\epsilon\sqrt H}\to0
 \quad\text{for every }\epsilon>0\text{ as }H\to\infty.
\]
It follows from~\eqref{eq:clock-corrector} that the processes
$u\mapsto(Q_{Hu}-\bar qHu)/\sqrt H$ converge in $C([0,1])$ to
that Brownian motion.

This convergence is uniform for $S_0=s\in I$. To see this, remove
an initial interval of length $\ell_H=H^{1/4}$. Its clock contribution
after division by $\sqrt H$ is at most
$(D+\bar q)H^{-1/4}$. After one unit, the transverse density is
bounded relative to $\varpi$, uniformly for $s\in I$. The spectral
gap therefore makes its distribution at time $\ell_H$ differ from
$\varpi$ in total variation by at most $C_Ie^{-g(\ell_H-1)}$. The Markov
property transfers the stationary convergence to the remaining
path, and the removed time interval has vanishing rescaled length.
Since Brownian motion has positive probability of remaining in
any fixed neighbourhood of zero, for each $\delta>0$ there are
$p_\delta>0$ and $H_\delta\ge1$ such that
\[\label{eq:clock-event}
 \mathbb P_s\left\{
 \sup_{0\le v\le H}|Q_v-\bar qv|<\delta\sqrt H\right\}
 \ge p_\delta
 \quad\text{for }H\ge H_\delta,\ s\in I.
\]
If $\sigma_Q=0$, the limiting process is zero and the same
conclusion holds.

We also need the event in~\eqref{eq:clock-event} with a prescribed
transverse endpoint. Let $\ell_H=H^{1/4}$ and require the clock bound
with threshold $\delta\sqrt H/2$ only up to time $H-\ell_H$; call this
event $E_H$. The preceding convergence gives
$\inf_{s\in I}\mathbb P_s(E_H)\ge p>0$ for large $H$. Because
$q\le D$, the change in $Q_v-\bar qv$ during the last $\ell_H$ units is
at most $(D+\bar q)\ell_H<\delta\sqrt H/2$. Thus $E_H$ implies the
full clock bound in~\eqref{eq:clock-event}.

For the transverse density $p_v(s,r)$ defined in
Lemma~\ref{lem:diffusion}, reversibility and the gap give
\[\label{eq:clock-endpoint-mixing}
 \|p_{\ell_H}(r,\cdot)-1\|_{L^2(\varpi)}
 \le C_Ie^{-g(\ell_H-1)}
 \quad\text{for }\ell_H\ge1,\ r\in I.
\]
The transverse marginal at $H-\ell_H$ weighted by $\mathbf1_{E_H}$
has a density $f_{H,s}$ with $0\le f_{H,s}\le C_I$ and
$\varpi(f_{H,s})=\mathbb P_s(E_H)$. By reversibility, the joint
density of $E_H$ and $S_H=r$ is
$\langle f_{H,s},p_{\ell_H}(r,\cdot)\rangle_{L^2(\varpi)}$.
Equation~\eqref{eq:clock-endpoint-mixing} bounds this below by
$p-C_Ie^{-g(\ell_H-1)}\ge p/2$, uniformly for $s,r\in I$ and large
$H$. The same lower bound holds for the full clock event and
$S_H=r$.

Conditional on the transverse path, the axial process has the
representation
\[\label{eq:middle-axial-representation}
 Z_v=x+\nu(cv-2\lambda Q_v)+\sqrt2\,\mathsf W_{Q_v}
 \quad\text{for }0\le v\le H,\ \nu\in\{-1,1\},
\]
where $\mathsf W$ is independent of $S$. After conditioning also
on $Z_H=y$, its mean is
\[\label{eq:middle-bridge-mean}
 x+(y-x)\frac{Q_v}{Q_H}
       +\nu c\left(v-H\frac{Q_v}{Q_H}\right)
 \quad\text{for }0\le v\le H.
\]
On the full clock event, for large $H$,
\[\label{eq:middle-clock-bounds}
 \frac{\bar qH}{2}\le Q_H\le DH,\qquad
 \left|v-H\frac{Q_v}{Q_H}\right|
       \le\frac{4\delta}{\bar q}\sqrt H
 \quad\text{for }0\le v\le H.
\]
Choose $\delta>0$ so that $4c\delta/\bar q\le d_1/2$.
For the endpoint heights in the statement,
\eqref{eq:middle-bridge-mean} is then at least $d_1\sqrt H/2$.
The centred conditional path has the form
$\sqrt{2Q_H}\,\mathbb B(Q_v/Q_H)$, where $\mathbb B$ is a
standard Brownian bridge. The event
$\sup_{0\le u\le1}|\mathbb B(u)|<d_1/(4\sqrt{2D})$ has a fixed
positive probability and keeps the entire axial path above zero.

On the same clock event, the free endpoint density satisfies, for $x,y\in[d_1\sqrt H,d_2\sqrt H]$ and $\nu\in\{-1,1\}$,
\[\label{eq:middle-free-density}
 \frac{1}{\sqrt{4\pi Q_H}}
 \exp\left[-\frac{(y-x-\nu(cH-2\lambda Q_H))^2}{4Q_H}\right]
 \ge C^{-1}H^{-1/2}.
\]
Here $c=2\lambda\bar q$, so the numerator in the exponent is
bounded by $CH$. Multiplying this density by the conditional
bridge probability and by the lower bound $p/2$ for the joint density of the clock event and $S_H=r$ gives~\eqref{eq:middle-density-lower}.
\end{proof}

\begin{lemma}\label{lem:logarithmic-lower}
Set $\alpha=3/(2\lambda)$, $b(t)=-\alpha\log(1+t)$ for $t\ge0$,
and $B_t(v)=b(t-v)$ for $0\le v\le t$. Let
$I_0\subset[0,\infty)$ be a compact interval of positive length.
There exist $Y>1$, $X>Y$, $C>0$, and $t_0\ge1$ such that
\[\label{eq:logarithmic-kernel-lower}
 k^{B_t,1}_{0,t}(s,b(t)+x;r,y)\ge C^{-1}t^{-3/2}
 \quad\text{for }t\ge t_0,\ s,r\in I_0,\
                 x\in[X,X+1],\ y\in[Y,Y+1],
\]
almost everywhere in $(r,y)$ with respect to
$\varpi(\mathrm dr)\,\mathrm dy$. Consequently, if $G\ge0$ is
bounded, compactly supported, and bounded below by a positive
constant on $I_0\times[Y,Y+1]$, then
\[\label{eq:logarithmic-test-lower}
 K^{B_t,1}_{0,t}G(s,b(t)+x)\ge C_G^{-1}(1+t)^{-3/2}
 \quad\text{for }t\ge t_0,\ s\in I_0,\ x\in[X,X+1].
\]
\end{lemma}

\begin{proof}
Take $C_0=1+\alpha\log2$ and
choose $Y>1,I_1,d_1,d_2$ by Lemma~\ref{lem:meander}.
The adjoint of $\mathcal L^\nu$
relative to $\varpi(\mathrm ds)\,\mathrm dz$ is
$\mathcal L^{-\nu}$. The transverse operator is self-adjoint,
the axial coefficients depend only on $s$, and the Dirichlet
condition at zero is preserved by taking the adjoint. Hence
\[\label{eq:flat-kernel-adjoint}
 k_h^\nu(s,x;r,y)=k_h^{-\nu}(r,y;s,x)
 \quad\text{for }h>0,\ s,r\ge0,\ x,y>0,\
                   \nu\in\{-1,1\}.
\]
For the logarithmic wall, put $h=t/3$ and $X=Y+\alpha\log2$.
On the first third, require the path to remain above the constant
wall $B_t(h)$. This condition implies survival above $B_t$ there.
Since
\[\label{eq:logarithmic-wall-first-gap}
 0\le B_t(h)-B_t(0)
   =\alpha\log\frac{1+3h}{1+2h}
   \le\alpha\log(3/2)
 \quad\text{for }h\ge0,
\]
the initial height relative to this constant wall belongs to
$[Y,Y+C_0]$ whenever $x\in[X,X+1]$. Thus the first endpoint
has mass at least $C^{-1}h^{-1/2}$ in
$I_1\times(B_t(h)+[d_1\sqrt h,d_2\sqrt h])$.

On the final third, require the path to remain above zero.
The actual wall is nonpositive on this interval. Reversing this
constant-wall kernel changes the sign of the axial drift by
\eqref{eq:flat-kernel-adjoint}. A terminal point with
$r\in I_0$ and $y\in[Y,Y+1]$ therefore gives mass at least
$C^{-1}h^{-1/2}$ in $I_1\times[d_1\sqrt h,d_2\sqrt h]$
at the other end of this third.

On the middle third, use the constant wall $B_t(2h)$, which lies
above $B_t$ throughout that interval. The heights of its initial
and final endpoints relative to this constant wall are,
respectively,
\[\label{eq:logarithmic-middle-heights}
 \begin{aligned}
 &[d_1\sqrt h,d_2\sqrt h]
       -\alpha\log\frac{1+2h}{1+h},\\
 &[d_1\sqrt h,d_2\sqrt h]+\alpha\log(1+h)
 \end{aligned}
 \qquad\text{for }h>0.
\]
For sufficiently large $h$, both intervals are contained in
$[d_1\sqrt h/2,2d_2\sqrt h]$. Lemma~\ref{lem:middle-density}
therefore bounds the middle density below by $C^{-1}h^{-1/2}$.
All three imposed constant walls lie above the corresponding
segments of $B_t$. Composing their kernels and multiplying the
two endpoint masses by the middle density gives
\eqref{eq:logarithmic-kernel-lower}. Integrating this inequality
against $G(r,y)\varpi(\mathrm dr)\,\mathrm dy$ over
$I_0\times[Y,Y+1]$ proves~\eqref{eq:logarithmic-test-lower}.
\end{proof}

\subsection{The nonlinear equation and the level positions}\label{sec:nonlinear}

Retain $\alpha=3/(2\lambda)$ and $b(t)=-\alpha\log(1+t)$ from
Lemma~\ref{lem:logarithmic-lower}, and define
\[\label{eq:m-position}
 m(t)=ct+b(t)\text{ for }t\ge0.
\]
For $B_t(v)=b(t-v)$, the elementary inequality
$\log(1+r+h)-\log(1+r)\le\log(1+h)$ for $r,h\ge0$ gives
\[\label{eq:logarithmic-wall-modulus}
 |B_t(v)-B_t(w)|\le\alpha\log(1+|v-w|)
 \text{ for }t\ge0,\ 0\le v,w\le t.
\]
We therefore apply the estimates of Subsection~\ref{sec:absorbing} with
$K=\alpha$.
The critical exponential change of unknowns takes a linear solution
$v(s,y,t)$ to $V(s,z,t)=e^{\lambda z}v(s,ct+z,t)/h_\lambda(s)$.
By~\eqref{eq:reflected-pde} and the eigenvalue equation $\J_\parallel(\lambda)\psi=\Lambda(\lambda)\psi$, it gives
\[\label{eq:critical-transform}
 V_t=\mathcal L^1V
 \text{ for }s>0,\ z\in\R,\ t>0
\]
when $F(s,v)$ is replaced by $a(s)v$. The Neumann condition is preserved
because $h_\lambda'(0)=0$. We use this identity first for the upper bound, and then
for a solution with zero boundary value at $y=m(t)$.

\begin{lemma}\label{lem:nonlinear-upper}
For every compactly supported datum $0\le u_0\le1$, there is $C<\infty$ such
that
\[\label{eq:axis-upper}
 u(0,m(t)+x,t)\le C(1+x)e^{-\lambda x}
 \text{ for }x\ge0,\ t\ge3.
\]
\end{lemma}

\begin{proof}
Choose a smooth, compactly supported function $\bar u_0$, even in $\rho$,
with $0\le\bar u_0\le1$ and $\bar u_0\ge u_0$ almost everywhere.
Let $\bar u$ solve~\eqref{eq:KPP} with this datum. By comparison,
it suffices to prove~\eqref{eq:axis-upper} for $\bar u$.
Its restriction to $s\ge0$ satisfies~\eqref{eq:reflected-pde}.
Put
\[\label{eq:upper-change}
 U(s,z,t)=\frac{e^{\lambda z}}{h_\lambda(s)}\bar u(s,ct+z,t)
 \text{ for }s\ge0,\ z\in\R,\ t\ge0.
\]
By~\eqref{eq:critical-transform} and the KPP upper bound,
$U_t\le\mathcal L^1U$. Moreover, $U$ is bounded by
\eqref{eq:longitudinal-exponential}. Define $W_0(s)=1+\psi(s)^{-1}$ for $s\ge0$.
By~\eqref{eq:eigenfunction-tails},
$W_0\in L^p(\varpi)$ for each $1<p<2$, since
$\int_0^\infty\psi^{2-p}\dd s<\infty$.
The bound on $|\psi'/\psi|$ in~\eqref{eq:eigenfunction-tails} gives $|(\log W_0)'|\le C$.
Since $\bar u(s,y,r)\le C\phiZ(s)/\sqrt{\cosh s}$ by~\eqref{eq:confinement} and $\phiZ$ is bounded, \eqref{eq:upper-change} gives
\[\label{eq:boundary-trace}
 U(s,b(r),r)\le C(1+r)^{-3/2}W_0(s)
 \text{ for }s\ge0,\ r\ge0.
\]

Fix $t\ge3$ and $x>0$, and start the process with generator $\mathcal L^1$
at $(0,b(t)+x)$. For the wall $B_t(v)=b(t-v)$, let
$\tau=\inf\{v\in[0,t]:Z_v\le B_t(v)\}$, with $\inf\varnothing=\infty$.
Bounded stopped parabolic comparison
gives
\[\label{eq:upper-stopped}
 \begin{aligned}
 U(0,b(t)+x,t)
 \le{}&\E^1_{0,b(t)+x}
       [U(S_\tau,b(t-\tau),t-\tau);\ \tau\le t]\\
 &+\E^1_{0,b(t)+x}[G_0(S_t,Z_t);\ \tau>t]
 \quad\text{for }x>0,\ t\ge3,
 \end{aligned}
\]
where $G_0(s,z)=e^{\lambda z}\bar u_0(s,z)/h_\lambda(s)$ is bounded and compactly
supported. Let $\tau_n=\inf\{v\ge0:S_v\ge n\ \text{or}\ |Z_v|\ge n\}$.
Apply It\^o's formula to $U(S_v,Z_v,t-v)$ up to $t\wedge\tau\wedge\tau_n$.
The Neumann condition at $s=0$ eliminates the reflection term.
Since $(S,Z)$ is nonexplosive, $\tau_n\to\infty$ almost surely.
Boundedness of $U$ permits passage to the limit in the stopped expectations
by dominated convergence, giving~\eqref{eq:upper-stopped}.

We estimate the boundary term by intervals of length one. Suppose the
process is at $(s,B_t(j)+w)$ at time $j$, with $w>0$.
For the process restarted from this point at time zero, define
$\tau^{(j)}=\inf\{v\in[0,t-j]:Z_v\le B_t(j+v)\}$.
Over the next unit interval the wall moves by at most $\alpha\log2$.
By~\eqref{eq:sde}, a boundary hit requires an axial martingale displacement
of magnitude at least $w-C$. Its quadratic variation rate is at most $2D$,
so the exponential martingale inequality gives a bound $Ce^{-dw^2}$ for
some $d>0$. On the same interval,
$W_0(S_v)\le W_0(s)\exp(C\sup_{r\le1}|S_r-s|)$.
The exponential moments from~\eqref{eq:unit-increments} and
Cauchy--Schwarz consequently give
\[\label{eq:weighted-unit-hit}
 \E^1_{s,B_t(j)+w}
 [W_0(S_{\tau^{(j)}});\ 0<\tau^{(j)}\le\min\{1,t-j\}]
 \le CW_0(s)e^{-dw^2}
 \text{ for }s\ge0,\ w>0,\ 0\le j<t.
\]
For an integer $j<t$, combine~\eqref{eq:boundary-trace} and
\eqref{eq:weighted-unit-hit}, conditional on survival to $j$.
Since $1+t-j$ and $1+t-\tau$ are comparable on
$j<\tau\le\min\{j+1,t\}$, Lemma~\ref{lem:weighted-transition} gives
\[\label{eq:boundary-contribution}
 \begin{aligned}
 &\E^1_{0,b(t)+x}
 [U(S_\tau,b(t-\tau),t-\tau);\
                 j<\tau\le\min\{j+1,t\}]\\
 &\hspace{1em}\le
 C(1+x)(1+j)^{-3/2}(1+t-j)^{-3/2}
 \text{ for }x>0,\ 3\le j<t,\ j\in\mathbb Z.
 \end{aligned}
\]
For $j=0,1,2$, the same inequality follows from the uniform finite-time
moments of $W_0(S_j)$ starting at zero. These moments are finite because
$W_0$ has at most exponential growth and~\eqref{eq:unit-increments} has
Gaussian tails. The factor $(1+j)^{-3/2}$ is bounded below on these three
values of $j$.
The second term in~\eqref{eq:upper-stopped} is bounded by
$C(1+x)(1+t)^{-3/2}$, again by Lemma~\ref{lem:weighted-transition}.

Splitting the sum at $t/2$ gives
\[\label{eq:discrete-convolution}
 \sum_{j=0}^{\lfloor t\rfloor}
 (1+j)^{-3/2}(1+t-j)^{-3/2}
 \le C(1+t)^{-3/2}
 \text{ for }t\ge0.
\]
By~\eqref{eq:upper-stopped}--\eqref{eq:discrete-convolution},
$U(0,b(t)+x,t)\le C(1+x)(1+t)^{-3/2}$.
Multiplying by $h_\lambda(0)e^{-\lambda(b(t)+x)}$ and using
$e^{-\lambda b(t)}=(1+t)^{3/2}$ proves~\eqref{eq:axis-upper} for $x>0$.
Continuity gives the estimate at $x=0$, and comparison transfers it to $u$.
\end{proof}

\begin{lemma}\label{lem:positive-seed}
For every nonzero compactly supported datum $0\le u_0\le1$, there are
$R_*>0$, $X_*>0$, $\delta>0$ and $T_*>0$ such that
\[\label{eq:fixed-positive-seed}
 u(\rho,m(t)+x,t+1)\ge\delta
 \text{ for }|\rho|\le R_*,\ x\in[X_*,X_*+1],\ t\ge T_*.
\]
\end{lemma}

\begin{proof}
Take $I_0=[0,R_*]$ with $R_*>0$, and choose the constants $X_*,Y$ in
Lemma~\ref{lem:logarithmic-lower}. Choose a smooth nonnegative function
$f_0(s,y)$, compactly supported in $s\ge0$, $y>0$, positive on
$I_0\times[Y,Y+1]$, and satisfying $\partial_sf_0(0,y)=0$.
Choose its even extension in $\rho$ to be smooth.
Strict positivity of $u(\rho,y,1)$ permits a constant $\varepsilon>0$ with
$\varepsilon f_0(|\rho|,y)\le u(\rho,y,1)$ everywhere.
We decrease $\varepsilon$ below as needed, keeping $f_0$ fixed.

Let $v$ solve the linear initial-boundary problem
\[\label{eq:linear-moving-domain}
 \left\{
 \begin{aligned}
 v_t&=D(v_{ss}+\tanh s\,v_s)+q(s)v_{yy}+a(s)v
       &&\text{for }s>0,\ y>m(t),\ t>0,\\
 v_s(0,y,t)&=0
       &&\text{for }y>m(t),\ t>0,\\
 v(s,m(t),t)&=0
       &&\text{for }s\ge0,\ t>0,\\
 v(s,y,0)&=\varepsilon f_0(s,y)
       &&\text{for }s\ge0,\ y>0.
 \end{aligned}
 \right.
\]
Put $G(s,z)=e^{\lambda z}f_0(s,z)/h_\lambda(s)$.
This function is bounded, compactly supported and positive on
$I_0\times[Y,Y+1]$. By~\eqref{eq:critical-transform}, its exact stopped
representation is
\[\label{eq:linear-representation}
 v(s,ct+z,t)=\varepsilon h_\lambda(s)e^{-\lambda z}
                K^{B_t,1}_{0,t}G(s,z)
 \text{ for }s\ge0,\ z>b(t),\ t\ge0.
\]
At $t=0$ the transition operator is the identity.
Define $\omega(x)=(1+x)e^{-\lambda x}$ for $x\ge0$.
Lemma~\ref{lem:compact-terminal}, together with
$\sup_{s\ge0}h_\lambda(s)\mathcal W(s)<\infty$, gives
\[\label{eq:linear-envelope}
 0\le v(s,m(t)+x,t)\le C\varepsilon\omega(x)
 \text{ for }s\ge0,\ x\ge0,\ t\ge0.
\]
The factor $(1+t)^{-3/2}$ in the bound of Lemma~\ref{lem:compact-terminal} cancels
$e^{-\lambda b(t)}=(1+t)^{3/2}$.
Choose $\varepsilon$ so small that the right side of
\eqref{eq:linear-envelope} is at most one.

For $s\in I_0$, $x\in[X_*,X_*+1]$ and $t$ sufficiently large, set
$N_t(s,x)=K^{B_t,1}_{0,t}G(s,b(t)+x)$.
By Lemma~\ref{lem:logarithmic-lower}, there are $c_0>0$ and $T_0\ge1$ such that
\[\label{eq:N-lower}
 N_t(s,x)\ge c_0(1+t)^{-3/2}
 \text{ for }s\in I_0,\ x\in[X_*,X_*+1],\ t\ge T_0.
\]
For the starting points in~\eqref{eq:N-lower}, define a probability measure on paths by
\[\label{eq:path-measure}
 \frac{\mathrm d\mathbb Q_{t,s,x}}
      {\mathrm d\Prb^1_{s,b(t)+x}}
 =\frac{G(S_t,Z_t)\one_{\{\tau_{B_t}>t\}}}{N_t(s,x)}
 \text{ for }s\in I_0,\ x\in[X_*,X_*+1],\ t\ge T_0.
\]
We claim that
\[\label{eq:occupation}
 \E_{\mathbb Q_{t,s,x}}\int_0^t
       \omega(Z_r-B_t(r))\dd r\le C
 \text{ for }s\in I_0,\ x\in[X_*,X_*+1],\ t\ge T_0.
\]

To prove~\eqref{eq:occupation}, apply the Markov property at time $r$.
Lemma~\ref{lem:compact-terminal} on the remaining wall segment gives
\[\label{eq:remaining-segment}
 K^{B_t,1}_{r,t}G(s',B_t(r)+w)
 \le C(1+w)\mathcal W(s')(1+t-r)^{-3/2}
 \text{ for }0\le r\le t,\ s',w\ge0.
\]
Multiplication by $\omega(w)$ leaves the endpoint weight
$\mathcal W(s')(1+w)^2e^{-\lambda w}$ for the first segment.
Its transverse factor is in $L^{3/2}(\varpi)$, and its bounded axial factor has
finite first moment. Applying Lemma~\ref{lem:weighted-transition} with
$p=3/2$ therefore yields
\[\label{eq:occupation-integrand}
 \begin{aligned}
 &\E^1_{s,b(t)+x}
 [\omega(Z_r-B_t(r))G(S_t,Z_t);\ \tau_{B_t}>t]\\
 &\hspace{1em}\le
 C(1+r)^{-3/2}(1+t-r)^{-3/2}
 \text{ for }0\le r\le t,
 \end{aligned}
\]
uniformly on the starting rectangle in~\eqref{eq:N-lower}.
For $r<3$, \eqref{eq:occupation-integrand} follows from
$\sup_{s\in I_0,\,0\le r\le3}\E_s\mathcal W(S_r)<\infty$,
which holds by~\eqref{eq:unit-increments} and the definition of
$\mathcal W$. For the remaining segment,
Lemma~\ref{lem:compact-terminal} applies also when $t-r<3$.
Since
\[\label{eq:continuous-convolution}
 \int_0^t(1+r)^{-3/2}(1+t-r)^{-3/2}\dd r
 \le C(1+t)^{-3/2}
 \text{ for }t\ge0,
\]
integration of~\eqref{eq:occupation-integrand} and division by
\eqref{eq:N-lower} prove~\eqref{eq:occupation}.

Let $C_2$ be the constant in~\eqref{eq:uniform-taylor}.
Let $w$ solve the same initial and boundary conditions as
\eqref{eq:linear-moving-domain}, with equation
\[\label{eq:logistic-moving-domain}
 w_t=D(w_{ss}+\tanh s\,w_s)+q(s)w_{yy}+a(s)w-C_2w^2
 \text{ for }s>0,\ y>m(t),\ t>0.
\]
Linear comparison and~\eqref{eq:linear-envelope} give $0\le w\le v\le1$.
By~\eqref{eq:uniform-taylor}, $w$ is a subsolution of
\eqref{eq:reflected-pde}. Its even extension across $\rho=0$ is therefore
dominated by $u(\rho,y,t+1)$ in $y>m(t)$.

The Feynman--Kac formula for~\eqref{eq:logistic-moving-domain}, followed by
$w\le v$, gives, for $s\ge0$, $z>b(t)$ and $t>0$,
\[\label{eq:axial-nonlinear-FK}
 \begin{aligned}
 w(s,ct+z,t)\ge{}&
 \varepsilon h_\lambda(s)e^{-\lambda z}\,
 \E^1_{s,z}\Bigg[
 G(S_t,Z_t)\\[-2mm]
 &\quad{}\times
 \exp\!\left\{-C_2\int_0^t
        v(S_r,c(t-r)+Z_r,t-r)\dd r\right\};
                \tau_{B_t}>t\Bigg].
 \end{aligned}
\]
These solutions and their stopped representations can be constructed on
bounded domains and passed to the limit by monotone exhaustion.
The coefficients are locally uniformly elliptic, the moving boundary is
smooth, and~\eqref{eq:linear-envelope} bounds the solutions and the
nonlinear killing rate. Thus dominated convergence justifies
\eqref{eq:axial-nonlinear-FK}.

Along every surviving path, its height at physical time $t-r$ is
$Z_r-b(t-r)=Z_r-B_t(r)$. Hence~\eqref{eq:linear-envelope} and
\eqref{eq:occupation} imply
\[\label{eq:nonlinear-occupation}
 \E_{\mathbb Q_{t,s,x}}\int_0^t
     v(S_r,c(t-r)+Z_r,t-r)\dd r\le C\varepsilon
 \text{ for }s\in I_0,\ x\in[X_*,X_*+1],\ t\ge T_0.
\]
Jensen's inequality, applied to the convex exponential in
\eqref{eq:axial-nonlinear-FK}, now yields, for $s\in I_0$, $x\in[X_*,X_*+1]$ and $t\ge T_0$,
\[\label{eq:fixed-amplitude}
 \begin{aligned}
 w(s,m(t)+x,t)
 &\ge \varepsilon h_\lambda(s)e^{-\lambda(b(t)+x)}
                 N_t(s,x)e^{-C\varepsilon}\\
 &\ge\delta>0.
 \end{aligned}
\]
The last inequality follows from~\eqref{eq:N-lower},
$e^{-\lambda b(t)}=(1+t)^{3/2}$, and the positive minimum of
$h_\lambda(s)e^{-\lambda x}$ on the fixed rectangle.
Comparison with the even extension of $w$ proves
\eqref{eq:fixed-positive-seed}.
\end{proof}

\begin{proof}[Proof of Theorem~\ref{thm:longitudinal}(iv)]
By Lemma~\ref{lem:nonlinear-upper}, choose $A_\theta$ so large that
$C(1+x)e^{-\lambda x}<\theta$ for every $x\ge A_\theta$.
Thus $E_\theta^+(t)\subset[0,m(t)+A_\theta]$ for all sufficiently large $t$.

Let $R_*,X_*,\delta,T_*$ be as in Lemma~\ref{lem:positive-seed}.
Put $Q=[-R_*,R_*]\times[X_*,X_*+1]$ in the $(\rho,y)$ coordinates
and $y_*=X_*+1/2$. Let $u_Q$ solve~\eqref{eq:KPP} from $\delta\one_Q$.
Since $p_*(\gamma(y_*))=p(0)>\theta$, estimate~\eqref{eq:liminf-local}
gives $h_0>0$ such that $u_Q(0,y_*,h)>\theta$ for every $h\ge h_0$.
At time $\tau+1$, inequality~\eqref{eq:fixed-positive-seed} bounds $u$
from below by $\delta\one_Q$ translated by $m(\tau)$ along $\G$.
Axial translation preserves the equation. Comparison for the remaining
time $t-1-\tau$ therefore yields
\[\label{eq:recovered-level}
 u(0,m(\tau)+y_*,t)>\theta
 \qquad(T_*\le\tau\le t-1-h_0).
\]
Choose $\tau_0\ge T_*$ so large that $m$ is increasing on
$[\tau_0,\infty)$ and $a=m(\tau_0)+y_*>0$.
At fixed $t\ge1+h_0+\tau_0$, varying $\tau$ from $\tau_0$ to
$t-1-h_0$ in~\eqref{eq:recovered-level} gives
$[a,m(t-1-h_0)+y_*]\subset E_\theta^+(t)$.
On the compact segment $\gamma([0,a])$, the stationary profile equals
$p(0)>\theta$, so~\eqref{eq:liminf-local} gives
$[0,a]\subset E_\theta^+(t)$ for all sufficiently large $t$.
Since $m(t)-m(t-1-h_0)\le c_\parallel(1+h_0)$,
these intervals give the inner inclusion in~\eqref{eq:longitudinal-level-speed}
on the positive half-geodesic.
Applying the same argument to $u(\rho,-y,t)$ proves both inclusions on
the negative half-geodesic. Take the larger position and time constants
and increase the time threshold so that $m(t)\ge C_{\theta,u_0}$ thereafter.
\end{proof}

\begingroup
\fontsize{10}{11}\selectfont
\raggedright

\endgroup
\end{document}